\documentclass[11pt]{ectaart}
\usepackage{amssymb,multirow}
\usepackage{amsfonts}
\usepackage{amsmath,amsthm}
\usepackage{amstext}
\usepackage{graphicx}
\usepackage{fncylab}
\usepackage{natbib}
\usepackage{graphicx}
\usepackage{hypernat}
\usepackage{hyperref}
\usepackage{booktabs}
\usepackage{setspace}

\numberwithin{equation}{section} \theoremstyle{plain}
\newtheorem{assumption}{\textbf{Assumption}}[section]
\newtheorem{theorem}{\textbf{Theorem}}[section]

\newtheorem{lemma}{\textbf{Lemma}}[section]

\newtheorem{remark}{\textbf{Remark}}[section]

\newtheorem{proposition}{\textbf{Proposition}}[section]
\newtheorem{condition}{\textbf{Condition}}[section]

\begin{document}

	\begin{frontmatter}

	\title{Sieve Wald and QLR Inferences on Semi/nonparametric Conditional
		Moment Models}
	
	\runtitle{Sieve Wald and QLR Inference}
\author{Xiaohong Chen and Demian Pouzo \thanksref{r1}}
\thankstext{r1}{
	Earlier versions, some entitled \textquotedblleft On PSMD inference of
	functionals of nonparametric conditional moment
	restrictions\textquotedblright ,\ were presented in April 2009 at the Banff
	conference on seminonparametrics, in June 2009 at the Cemmap conference on
	quantile regression, in July 2009 at the SITE conference on nonparametrics,
	in September 2009 at the Stats in the Chateau/France, in June 2010 at the
	Cemmep workshop on recent developments in nonparametric instrumental
	variable methods, in August 2010 at the Beijing international conference on
	statistics and society, and econometric workshops in numerous universities.
	We thank a co-editor, two referees, Don Andrews, Peter Bickel, Gary
	Chamberlain, Tim Christensen, Michael Jansson, Jim Powell and especially
	Andres Santos for helpful comments. We thank Yinjia Qiu for excellent
	research assistant in simulations using R. Chen acknowledges financial
	support from National Science Foundation grant SES-0838161 and Cowles
	Foundation. Any errors are the responsibility of the authors.}

\address{Chen: Cowles Foundation for Research in Economics, Yale University, Box 208281,	New Haven, CT 06520, USA. Email: xiaohong.chen@yale.edu. Pouzo: Department of Economics, UC Berkeley, 530 Evans Hall 3880, Berkeley, CA	94720, USA. Email: dpouzo@econ.berkeley.edu.}
\runauthor{X. Chen and D. Pouzo}
	
	
	\begin{abstract}
		This paper considers inference on functionals of semi/nonparametric
		conditional moment restrictions with possibly nonsmooth generalized
		residuals, which include all of the (nonlinear) nonparametric instrumental
		variables (IV) as special cases. These models are often ill-posed and hence
		it is difficult to verify whether a (possibly nonlinear) functional is
        root-$n$ estimable or not. We provide computationally simple, unified inference
		procedures that are asymptotically valid regardless of whether a functional
		is root-$n$ estimable or not. We establish the following new useful results:
		(1) the asymptotic normality of a plug-in penalized sieve minimum distance
		(PSMD) estimator of a (possibly nonlinear) functional; (2) the consistency
		of simple sieve variance estimators for the plug-in PSMD estimator, and
		hence the asymptotic chi-square distribution of the sieve Wald statistic;
		(3) the asymptotic chi-square distribution of an optimally weighted sieve
		quasi likelihood ratio (QLR) test under the null hypothesis; (4) the
		asymptotic tight distribution of a non-optimally weighted sieve QLR
		statistic under the null; (5) the consistency of generalized residual
		bootstrap sieve Wald and QLR tests; (6) local power properties of sieve Wald
		and QLR tests and of their bootstrap versions; (7) asymptotic properties of
		sieve Wald and SQLR for functionals of increasing dimension. Simulation
		studies and an empirical illustration of a nonparametric quantile IV
		regression are presented.
	\end{abstract}
	

\end{frontmatter}

\setcounter{page}{0} \thispagestyle{empty} \newpage

\baselineskip=18pt

\section{Introduction}


This paper is about inference on functionals of the unknown true parameters $%
\alpha _{0}\equiv (\theta _{0}^{\prime },h_{0})$ satisfying the
semi/nonparametric conditional moment restrictions
\begin{equation}
E[\rho (Y,X;\theta _{0},h_{0})|X]=0\quad a.s.-X,  \label{semi00}
\end{equation}%
\noindent where $Y$ is a vector of endogenous variables and $X$ is a vector
of conditioning (or instrumental) variables. The conditional distribution of
$Y$ given $X$, $F_{Y|X}$, is not specified beyond that it satisfies (\ref%
{semi00}). $\rho (\cdot ;\theta _{0},h_{0})$ is a $d_{\rho }\times 1-$vector
of generalized residual functions whose functional forms are known up to the
unknown parameters $\alpha _{0}\equiv (\theta _{0}^{\prime },h_{0})\in
\Theta \times \mathcal{H}$, with $\theta _{0}\equiv (\theta _{01},...,\theta
_{0d_{\theta }})^{\prime }\in \Theta $ being a $d_{\theta }\times 1-$vector
of finite dimensional parameters and $h_{0}\equiv (h_{01}(\cdot
),...,h_{0q}(\cdot ))\in \mathcal{H}$ being a $1\times d_{q}-$vector valued
function. The arguments of each unknown function $h_{\ell }(\cdot )$ may
differ across $\ell =1,...,q$, may depend on $\theta ,$ $h_{\ell ^{\prime
}}(\cdot ),$ $\ell ^{\prime }\neq \ell $, $X$ and $Y$. The residual function
$\rho (\cdot ;\alpha )$ could be nonlinear and pointwise non-smooth in the
parameters $\alpha \equiv (\theta ^{\prime },h)\in \Theta \times \mathcal{H}$%
.

The general framework (\ref{semi00}) nests many widely used nonparametric
and semiparametric models in economics and finance. Well known examples
include nonparametric mean instrumental variables regressions (NPIV): $%
E[Y_{1}-h_{0}(Y_{2})|X]=0$ (e.g., \cite{HH_Ann05}, \cite{CFR_bookchp07},
\cite{BCK_Emetrica07}, \cite{DFR_wp10}, \cite{Horowitz_ECMA11});
nonparametric quantile instrumental variables regressions (NPQIV): $%
E[1\{Y_{1}\leq h_{0}(Y_{2})\}-\gamma |X]=0$ (e.g., \cite{CH_Emetrica05},
\cite{CIN_JOE07}, \cite{HL_Emetrica07}, \cite{CP_WP07}, \cite{CGS_WP08});
semi/nonparametric demand models with endogeneity (e.g., \cite%
{BCK_Emetrica07}, \cite{CP_WP07a}, \cite{Souza2012}); semi/nonparametric
random coefficient panel data regressions (e.g., \cite{CHAMBERLAIN_ECMA92},
\cite{GPowell2012}); semi/nonparametric spatial models with endogeneity
(e.g., \cite{Pinkse2002}, \cite{MerloPaula}); semi/nonparametric asset
pricing models (e.g., \cite{Hansen_Richard_ECMA87}, \cite%
{Gallant_Tauchen_ECMA89}, \cite{ChenLudvigson_JAE09}, \cite%
{ChenLudvigson2013}, \cite{Sentana2013}); semi/nonparametric static and
dynamic game models (e.g., \cite{BHN_WP11}); nonparametric optimal
endogenous contract models (e.g., \cite{BMT_WP12}). Additional examples of
the general model (\ref{semi00}) can be found in \cite{CHAMBERLAIN_ECMA92},
\cite{NP_ECMA03}, \cite{AC_Emetrica03}, \cite{CP_WP07}, \cite{CCLN_WP10} and
the references therein. In fact, model (\ref{semi00}) includes all of the
(nonlinear) semi/nonparametric IV regressions when the unknown functions $%
h_{0}$ depend on the endogenous variables $Y$:
\begin{equation}
E[\rho (Y_{1};\theta _{0},h_{0}(Y_{2}))|X]=0\quad a.s.-X,  \label{gnpiv}
\end{equation}%
which could lead to difficult (nonlinear) nonparametric ill-posed inverse
problems with unknown operators.

Let $\left\{ Z_{i}\equiv (Y_{i}^{\prime },X_{i}^{\prime })^{\prime }\right\}
_{i=1}^{n}$ be a random sample from the distribution of $Z\equiv (Y^{\prime
},X^{\prime })^{\prime }$ that satisfies the conditional moment restrictions
(\ref{semi00}) with a unique $\alpha _{0}\equiv (\theta _{0}^{\prime
},h_{0}) $. Let $\phi :\Theta \times \mathcal{H}\rightarrow \mathbb{R}%
^{d_{\phi }}$ be a (possibly nonlinear) functional with a finite $d_{\phi
}\geq 1$. Typical linear functionals include an Euclidean functional $\phi
(\alpha )=\theta $, a point evaluation functional $\phi (\alpha )=h(%
\overline{y}_{2}) $ (for $\overline{y}_{2}\in $ supp$(Y_{2})$), a weighted
derivative functional $\phi (h)=\int w(y_{2})\nabla h(y_{2})dy_{2}$ and many
others. Typical nonlinear functionals include a quadratic functional $\int
w(y_{2})\left\vert h(y_{2})\right\vert ^{2}dy_{2}$, a quadratic derivative
functional $\int w(y_{2})\left\vert \nabla h(y_{2})\right\vert ^{2}dy_{2}$,
a consumer surplus or an average consumer surplus functional of an
endogenous demand function $h$. We are interested in computationally simple,
valid inferences on any $\phi (\alpha _{0})$ of the general model (\ref%
{semi00}) with i.i.d. data.\footnote{%
See our Cowles Foundation Discussion Paper No. 1897 for general theory
allowing for weakly dependent data.}

Although some functionals of the model (\ref{semi00}), such as the (point)
evaluation functional, are known \textit{a priori} to be estimated at slower
than root-$n$ rates, others, such as the weighted derivative functional, are
far less clear without a stare at their semiparametric efficiency bound
expressions. This is because a non-singular efficiency bound is a necessary
condition for $\phi (\alpha _{0})$ to be estimated at a root-$n$ rate.
Unfortunately, as pointed out in \cite{CHAMBERLAIN_ECMA92} and \cite{AC_WP05}%
, there is generally no closed form solution for the efficiency bound of $%
\phi (\alpha _{0})$ (including $\theta _{0}$) of model (\ref{semi00}),
especially so when $\rho (\cdot ;\theta _{0},h_{0})$ contains several
unknown functions and/or when the unknown functions $h_{0}$ of endogenous
variables enter $\rho (\cdot ;\theta _{0},h_{0})$ nonlinearly. It is thus
difficult to verify whether the efficiency bound for $\phi (\alpha _{0})$ is
singular or not. Therefore, it is highly desirable for applied researchers
to be able to conduct simple valid inferences on $\phi (\alpha _{0})$
regardless of whether it is root-$n$ estimable or not. This is the main goal
of our paper.

In this paper, for the general model (\ref{semi00}) that could be
nonlinearly ill-posed and for any $\phi (\alpha _{0})$ that may or may not
be root-$n$ estimable, we first establish the asymptotic normality of the
plug-in penalized sieve minimum distance (PSMD) estimator $\phi (\widehat{%
\alpha }_{n})$ of $\phi (\alpha _{0})$. For the model (\ref{semi00}) with
(pointwise) smooth residuals $\rho (Z;\alpha )$ in $\alpha _{0}$, we propose
two simple consistent sieve variance estimators for possibly slower than
root-$n$ estimator $\phi (\widehat{\alpha }_{n})$, which immediately leads
to the asymptotic chi-square distribution of the sieve Wald statistic.
However, there is no simple variance estimator for $\phi (\widehat{\alpha }%
_{n})$ when $\rho (Z,\alpha )$ is not pointwise smooth in $\alpha _{0}$
(without estimating an extra unknown nuisance function or using numerical
derivatives). We then consider a PSMD criterion based test of the null
hypothesis $\phi (\alpha _{0})=\phi _{0}$. We show that an optimally
weighted sieve quasi likelihood ratio (SQLR) statistic is asymptotically
chi-square distributed under the null hypothesis. This allows us to
construct confidence sets for $\phi (\alpha _{0})$ by inverting the
optimally weighted SQLR statistic, without the need to compute a variance
estimator for $\phi (\widehat{\alpha }_{n})$. Nevertheless, in complicated
real data analysis applied researchers might like to use simple but possibly
non-optimally weighed PSMD procedures for estimation of and inference on $%
\phi (\alpha _{0})$. We show that the non-optimally weighted SQLR statistic
still has a tight limiting distribution under the null regardless of whether
$\phi (\alpha _{0})$ is root-$n$ estimable or not. In addition, we establish
the consistency of the generalized residual bootstrap (possibly
non-optimally weighted) SQLR and sieve Wald tests under virtually the same
conditions as those used to derive the limiting distributions of the
original-sample statistics. The bootstrap SQLR would then lead to
alternative confidence sets construction for $\phi (\alpha _{0})$ without
the need to compute a variance estimator for $\phi (\widehat{\alpha }_{n})$.
To ease notation burden, we present the above listed theoretical results for
a scalar-valued functional in the main text. In Appendix \ref{app:appA} we
present the asymptotic properties of sieve Wald and SQLR for functionals of
increasing dimension (i.e., $d_{\phi }=dim(\phi )$ could grow with sample
size $n$). We also provide the local power properties of sieve Wald and SQLR
tests as well as their bootstrap versions in Appendix \ref{app:appA}.
Regardless of whether a possibly nonlinear functional $\phi (\alpha _{0})$
is root-$n$ estimable or not, we show that the optimally weighted SQLR is
more powerful than the non-optimally weighed SQLR, and that the SQLR and the
sieve Wald using the same weighting matrix have the same local power in
terms of first order asymptotic theory.

To the best of our knowledge, our paper is the first to provide a unified
theory about sieve Wald and SQLR inferences on (possibly nonlinear) $\phi
(\alpha _{0})$ satisfying the general semi/nonparametric model (\ref{semi00}%
) with possibly non-smooth residuals.\footnote{%
We also provide asymptotic properties of sieve score and bootstrap sieve
score statistics in the online Appendix \ref{app:appD}.} Our results allow applied researchers to obtain
limiting distribution of the plug-in PSMD estimator $\phi (\widehat{\alpha }%
_{n})$ and to construct confidence sets for any $\phi (\alpha _{0})$
regardless of whether it is root-$n$ estimable or not. Our paper is also the
first to provide local power properties of sieve Wald and SQLR tests and
their bootstrap versions of general nonlinear hypotheses for the model (\ref%
{semi00}).

Roughly speaking, our results extend the classical theories on Wald and QLR
tests of nonlinear hypothesis based on root-$n$ consistent parametric
minimum distance estimator $\widehat{\alpha }_{n}$ to those based on slower
than root-$n$ consistent nonparametric minimum distance estimator $\widehat{%
\alpha }_{n}\equiv (\widehat{\theta }_{n}^{\prime },\widehat{h}_{n})$ of $%
\alpha _{0}\equiv (\theta _{0}^{\prime },h_{0})$ satisfying the model (\ref%
{semi00}). The implementations of the sieve Wald and SQLR also resemble the
classical Wald and QLR based on parametric extreme estimators and hence are
computationally attractive. For example, our sieve t (Wald) test on a
general nonlinear hypothesis $\phi (h_{0})=\phi _{0}$ of the NPIV model $%
E[Y_{1}-h_{0}(Y_{2})|X]=0$ can be implemented as a standard t (Wald) test
for a parametric linear IV model using two stage least squares (see
Subsection \ref{sec:NPIVex}). The proof techniques are quite different,
however, because one is no longer able to rely on the root-$n$ asymptotic
normality of $\widehat{\alpha }_{n}$ and then a standard \textquotedblleft
delta-method\textquotedblright\ to establish the asymptotic normality of $%
\sqrt{n}\left( \phi (\widehat{\alpha }_{n})-\phi (\alpha _{0})\right) $. In
our framework (\ref{gnpiv}), $\sqrt{n}\left( \phi (\widehat{\alpha }%
_{n})-\phi (\alpha _{0})\right) $ could diverge to infinity under the
combined effects of (i) slower convergence rate of $\widehat{\alpha }_{n}$
to $\alpha _{0}$ due to the ill-posed inverse problem and (ii) nonlinearity
of either the functional $\phi ()$ or the residual function $\rho ()$ in $h$%
. Our proof strategy relies on the convergence rates of the PSMD estimator $%
\widehat{\alpha }_{n}$ to $\alpha _{0}$ in both weak and strong metrics, and
then the local curvatures of the functional $\phi ()$ and the criterion
function under these two metrics. The weak metric is intrinsic to the
variance of the linear approximation to $\phi (\widehat{\alpha }_{n})-\phi
(\alpha _{0})$, while the strong metric controls the nonlinearity (in $%
\alpha $) of the functional $\phi ()$ and of the conditional mean function $%
m(\cdot ,\alpha )=E[\rho (Y,X;\alpha )|X=\cdot ]$. Unfortunately the
convergence rate in the strong metric could be very slow due to the illposed
inverse problem. This explains why it is difficult to establish the
asymptotic normality of $\phi (\widehat{\alpha }_{n})$ for a nonlinear
functional $\phi ()$ even in the NPIV model. Our paper builds upon the
recent results on convergence rates in \cite{CP_WP07} and others. In
particular, under virtually the same conditions as those in \cite{CP_WP07},
we show that our generalized residual bootstrap PSMD estimator of $\alpha
_{0}$ is consistent and achieves the same convergence rates as that of the
original-sample PSMD estimator $\widehat{\alpha }_{n}$. This result is then
used to establish the consistency of the bootstrap sieve Wald and the
bootstrap SQLR statistics under virtually the same conditions as those used
to derive the limiting distributions of the original-sample statistics.%
\footnote{%
The convergence rate of the bootstrap PSMD estimator is also very useful for
the consistency of the bootstrap Wald statistic for semiparametric two-step
GMM estimation of Euclidean parameters when the first-step unknown functions
are estimated via a PSMD procedure. See e.g., \cite{CLvK_Emetrica03}}

There are some published work about estimation of and inference on a
particular linear functional, the Euclidean parameter $\phi (\alpha )=\theta
$, of the general model (\ref{semi00}) when $\theta _{0}$ is assumed to be
root-$n$ estimable; see \cite{AC_Emetrica03}, \cite{CP_WP07a}, \cite%
{OTSU_WP11} and others. None of the existing work allows for $\theta _{0}$
being \textit{irregular} (i.e., slower than root-$n$ estimable),\footnote{%
It is known that $\theta _{0}$ could have singular semiparametric efficiency
bound and could not be root-$n$ estimable; see \cite{CHAMBERLAIN2010}, \cite%
{KhanTamer2010}, \cite{GPowell2012} and the references therein. Following
\cite{KhanTamer2010} and \cite{GPowell2012} we call such a $\theta _{0}$
irregular. Many applied papers on complicated semi/nonparametric models
simply assume that $\theta _{0}$ is root-$n$ estimable.} however. When
specializing our general theory to inference on $\theta _{0}$ of the model (%
\ref{semi00}), we not only recover the results of \cite{AC_Emetrica03} and
\cite{CP_WP07a}, but also provide local power properties of sieve Wald and
SQLR as well as valid bootstrap (possibly non-optimally weighted) SQLR
inference. Moreover, our results remain valid even when $\theta _{0}$ is
irregular.

When specializing our theory to inference on a particular irregular linear
functional, the point evaluation functional $\phi (\alpha )=h(\overline{y}%
_{2})$, of the semi/nonparametric IV model (\ref{gnpiv}), we automatically
obtain the pointwise asymptotic normality of the PSMD estimator of $h_{0}(%
\overline{y}_{2})$ and different ways to construct its confidence set. These
results are directly applicable to the NPIV example with $\rho (Y_{1};\theta
_{0},h_{0}(Y_{2}))=Y_{1}-h_{0}(Y_{2})$ and to the NPQIV example with $\rho
(Y_{1};\theta _{0},h_{0}(Y_{2}))=1\{Y_{1}\leq h_{0}(Y_{2})\}-\gamma $.
Previously, \cite{Horowitz_07} and \cite{CGS_WP08} established the pointwise
asymptotic normality of their kernel based function space Tikhonov
regularization estimators of $h_{0}(\overline{y}_{2})$ for the NPIV and the
NPQIV examples respectively. Immediately after our paper was first presented
in April 2009 Banff/Canada conference on semiparametrics, the authors of
\cite{HL_WP10} informed us that they were concurrently working on confidence
bands for $h_{0}$ using a particular SMD estimator of the NPIV example. To
the best of our knowledge, there is no inference results, in the existing
literature, on any nonlinear functional of $h_{0}$ even for the NPIV and
NPQIV examples. Our paper is the first to provide simple sieve Wald and SQLR
tests for (possibly) nonlinear functionals satisfying the general
semi/nonparametric IV model (\ref{gnpiv}).

The rest of the paper is organized as follows. Section \ref{sec:model}
presents the plug-in PSMD estimator $\phi (\widehat{\alpha }_{n})$ of a
(possibly nonlinear) functional $\phi $ evaluated at $\alpha _{0}\equiv
(\theta _{0}^{\prime },h_{0})$ satisfying the model (\ref{semi00}). It also
provides an overview of the main asymptotic results that will be established
in the subsequent sections, and illustrates the applications through a point
evaluation functional $\phi (\alpha )=h(\overline{y}_{2})$, a weighted
derivative functional $\phi (h)=\int w(y_{2})\nabla h(y_{2})dy_{2}$, and a
quadratic functional $\phi (\alpha )=\int w(y_{2})\left\vert
h(y_{2})\right\vert ^{2}dy_{2}$ of the NPIV and NPQIV examples. Section \ref%
{sec:conditions} states the basic regularity conditions. Section \ref%
{sec:AsymDist} provides the asymptotic properties of sieve t (Wald) and
sieve QLR statistics. Section \ref{sec:bootstrap} establishes the
consistency of the bootstrap sieve t (Wald) and the bootstrap SQLR
statistics. Section \ref{sec-ex} verifies the key regularity conditions for
the asymptotic theories via the three functionals of the NPIV and NPQIV
examples presented in Section \ref{sec:model}. Section \ref%
{sec:sec_simulation} presents simulation studies and an empirical
illustration. Section \ref{sec:conclusion} briefly concludes. Appendix \ref%
{app:appA} consists of several subsections, presenting (1) further results
on sieve Riesz representation of a functional of interest; (2) the
convergence rates of the bootstrap PSMD estimator $\widehat{\alpha }_{n}^{B}$
for model (\ref{semi00}); (3) the local power properties of sieve Wald and
SQLR tests and of their bootstrap versions; (4) asymptotic properties of
sieve Wald and SQLR for functionals of increasing dimension; (5) low level
sufficient conditions with a series least squares (LS) estimated conditional
mean function $m(\cdot ,\alpha )=E[\rho (Y,X;\alpha )|X=\cdot ]$; and (6)
additional useful lemmas with series LS estimated $m(\cdot ,\alpha )$.
Online supplemental materials consist of Appendices \ref{app:appB}, \ref{app:appC} and \ref%
{app:appD}. Appendix \ref{app:appB} contains additional theoretical results
(including other consistent variance estimators and other bootstrap sieve
Wald tests) and proofs of all the results stated in the main text. Appendix %
\ref{app:appC} contains proofs of all the results stated in Appendix \ref%
{app:appA}. The online Appendix \ref{app:appD} provides
computationally attractive sieve score test and sieve score bootstrap.

\textbf{Notation}. We use \textquotedblleft $\equiv $\textquotedblright\ to
implicitly define a term or introduce a notation. For any column vector $A$,
we let $A^{\prime }$ denote its transpose and $||A||_{e}$ its Euclidean norm
(i.e., $||A||_{e}\equiv \sqrt{A^{\prime }A}$, although sometimes we use $%
|A|=||A||_{e}$ for simplicity). Let $||A||_{W}^{2}\equiv A^{\prime }WA$ for
a positive definite weighting matrix $W$. Let $\lambda _{\max }(W)$ and $%
\lambda _{\min }(W)$ denote the maximal and minimal eigenvalues of $W$
respectively. All random variables $Z\equiv (Y^{\prime },X^{\prime
})^{\prime }$, $Z_{i}\equiv (Y_{i}^{\prime },X_{i}^{\prime })^{\prime }$ are
defined on a complete probability space $(\mathcal{Z},\mathcal{B}_{Z},P_{Z})$%
, where $P_{Z}$ is the joint probability distribution of $(Y^{\prime
},X^{\prime })$. We define $(\mathcal{Z}^{\infty },\mathcal{B}_{Z}^{\infty
},P_{Z^{\infty }})$ as the probability space of the sequences $%
(Z_{1},Z_{2},...)$. For simplicity we assume that $Y$ and $X$ are continuous
random variables. Let $f_{X}$ ($F_{X}$) be the marginal density (cdf) of $X$
with support $\mathcal{X}$, and $f_{Y|X}$ ($F_{Y|X}$) be the conditional
density (cdf) of $Y$ given $X$. Let $E_{P}[\cdot ]$ denote the expectation
with respect to a measure $P$. Sometimes we use $P$ for $P_{Z^{\infty }}$
and $E[\cdot ]$ for $E_{P_{Z^{\infty }}}[\cdot ]$. Denote $L^{p}(\Omega
,d\mu )$, $1\leq p<\infty $, as a space of measurable functions with $%
||g||_{L^{p}(\Omega ,d\mu )}\equiv \{\int_{\Omega }|g(t)|^{p}d\mu
(t)\}^{1/p}<\infty $, where $\Omega $ is the support of the sigma-finite
positive measure $d\mu $ (sometimes $L^{p}(d\mu )$ and $||g||_{L^{p}(d\mu )}$
are used). For any (possibly random) positive sequences $\{a_{n}\}_{n=1}^{%
\infty }$ and $\{b_{n}\}_{n=1}^{\infty }$, $a_{n}=O_{P}(b_{n})$ means that $%
\lim_{c\rightarrow \infty }\limsup_{n}\Pr \left( a_{n}/b_{n}>c\right) =0$; $%
a_{n}=o_{P}(b_{n})$ means that for all $\varepsilon >0$, $\lim_{n\rightarrow
\infty }\Pr \left( a_{n}/b_{n}>\varepsilon \right) =0$; and $a_{n}\asymp
b_{n}$ means that there exist two constants $0<c_{1}\leq c_{2}<\infty $ such
that $c_{1}a_{n}\leq b_{n}\leq c_{2}a_{n}$. Also, we use \textquotedblleft
wpa1-$P_{Z^{\infty }}$\textquotedblright\ (or simply wpa1) for an event $%
A_{n}$, to denote that $P_{Z^{\infty }}(A_{n})\rightarrow 1$ as $%
n\rightarrow \infty $. We use $\mathcal{A}_{n}\equiv \mathcal{A}_{k(n)}$ and
$\mathcal{H}_{n}\equiv \mathcal{H}_{k(n)}$ for various sieve spaces. We
assume $\dim (\mathcal{A}_{k(n)})\asymp \dim (\mathcal{H}_{k(n)})\asymp k(n)$
for simplicity, all of which grow to infinity with the sample size $n$. We
use $const.$, $c$ or $C$ to mean a positive finite constant that is
independent of sample size but can take different values at different
places. For sequences, $(a_{n})_{n}$, we sometimes use $a_{n}\nearrow a$ ($%
a_{n}\searrow a$) to denote, that the sequence converges to $a$ and that is
increasing (decreasing) sequence. For any mapping $\digamma :\mathbf{H}%
_{1}\rightarrow \mathbf{H}_{2}$ between two generic Banach spaces, $\frac{%
d\digamma (\alpha _{0})}{d\alpha }[v]\equiv \left. \frac{\partial \digamma
(\alpha _{0}+\tau v)}{\partial \tau }\right\vert _{\tau =0}$ is the pathwise
(or Gateaux) derivative at $\alpha _{0}$ in the direction $v\in \mathbf{H}%
_{1}$. And $\frac{d\digamma (\alpha _{0})}{d\alpha }[\mathbf{v}^{\prime
}]\equiv \left( \frac{d\digamma (\alpha _{0})}{d\alpha }[v_{1}],\cdot \cdot
\cdot ,\frac{d\digamma (\alpha _{0})}{d\alpha }[v_{k}]\right) $ for $\mathbf{%
v}^{\prime }=\left( v_{1},\cdot \cdot \cdot ,v_{k}\right) $ with $v_{j}\in
\mathbf{H}_{1}$ for all $j=1,...,k$.

\section{PSMD Estimation and Inferences: An Overview}

\label{sec:model}

\subsection{The Penalized Sieve Minimum Distance Estimator}

Let $m(X,\alpha )\equiv E\left[ \rho (Y,X;\alpha )|X\right] =\int \rho
(y,X;\alpha )dF_{Y|X}(y)$ be a $d_{\rho }\times 1$ vector valued conditional
mean function, $\Sigma (X)$ be a $d_{\rho }\times d_{\rho }$ positive
definite ($a.s.-X$) weighting matrix, and
\begin{equation*}
Q(\alpha )\equiv E\left[ m(X,\alpha )^{\prime }\Sigma (X)^{-1}m(X,\alpha )%
\right] \equiv E\left[ ||m(X,\alpha )||_{\Sigma ^{-1}}^{2}\right]
\end{equation*}%
be the population minimum distance (MD) criterion function. Then the
semi/nonparametric conditional moment model (\ref{semi00}) can be
equivalently expressed as $m(X,\alpha _{0})=0$ $a.s.-X$, where $\alpha
_{0}\equiv (\theta _{0}^{\prime },h_{0})\in \mathcal{A}\equiv \Theta \times
\mathcal{H}$, or as
\begin{equation*}
\inf_{\alpha \in \mathcal{A}}Q(\alpha )=Q(\alpha _{0})=0.
\end{equation*}%
Let $\Sigma _{0}(X)\equiv Var(\rho (Y,X;\alpha _{0})|X)$ be positive
definite for almost all $X$. In this paper as well as in most applications $%
\Sigma (X)$ is chosen to be either $I_{d_{\rho }}$ (identity) or $\Sigma
_{0}(X)$ for almost all $X$. We call $Q^{0}(\alpha )\equiv E\left[
||m(X,\alpha )||_{\Sigma _{0}^{-1}}^{2}\right] $ the population optimally
weighted MD criterion function.

Let $\phi :\mathcal{A}\rightarrow \mathbb{R}^{d_{\phi }}$ be a functional
with a finite $d_{\phi }\geq 1$. We are interested in inference on $\phi
(\alpha _{0})$. Let
\begin{equation}
\widehat{Q}_{n}(\alpha )\equiv \frac{1}{n}\sum_{i=1}^{n}\widehat{m}%
(X_{i},\alpha )^{\prime }\widehat{\Sigma }(X_{i})^{-1}\widehat{m}%
(X_{i},\alpha )  \label{Qhat}
\end{equation}%
be a sample estimate of $Q(\alpha )$, where $\widehat{m}(X,\alpha )$ and $%
\widehat{\Sigma }(X)$ are any consistent estimators of $m(X,\alpha )$ and $%
\Sigma (X)$ respectively. When $\widehat{\Sigma }(X)=\widehat{\Sigma }%
_{0}(X) $ is a consistent estimator of the optimal weighting matrix $\Sigma
_{0}(X)$, we call the corresponding $\widehat{Q}_{n}(\alpha )$ the\ sample
optimally weighted MD criterion $\widehat{Q}_{n}^{0}(\alpha )$.

We estimate $\phi (\alpha _{0})$ by $\phi (\widehat{\alpha }_{n})$, where $%
\widehat{\alpha }_{n}\equiv (\widehat{\theta }_{n}^{\prime },\widehat{h}%
_{n}) $ is an approximate \textit{penalized sieve minimum distance} (PSMD)
estimator of $\alpha _{0}\equiv (\theta _{0}^{\prime },h_{0})$, defined as
\begin{equation}
\widehat{Q}_{n}(\widehat{\alpha }_{n})+\lambda _{n}Pen(\widehat{h}_{n})\leq
\inf_{\alpha \in \mathcal{A}_{k(n)}}\left\{ \widehat{Q}_{n}(\alpha )+\lambda
_{n}Pen(h)\right\} +o_{P_{Z^{\infty }}}(n^{-1}),  \label{psmd}
\end{equation}%
where $\lambda _{n}Pen(h)\geq 0$ is a penalty term such that $\lambda
_{n}=o(1)$; and $\mathcal{A}_{k(n)}\equiv \Theta \times \mathcal{H}_{k(n)}$
is a finite dimensional sieve for $\mathcal{A}\equiv \Theta \times \mathcal{H%
}$, more precisely, $\mathcal{H}_{k(n)}$ is a finite dimensional \textit{%
linear} sieve for $\mathcal{H}$:
\begin{equation}
\mathcal{H}_{k(n)}=\left\{ h\in \mathcal{H}:h(\cdot )=\sum_{k=1}^{k(n)}\beta
_{k}q_{k}(\cdot )=\beta ^{\prime }q^{k(n)}(\cdot )\right\} ,  \label{sieve}
\end{equation}%
where $\{q_{k}\}_{k=1}^{\infty }$ is a sequence of known basis functions of
a Banach space $(\mathcal{H},\left\Vert \cdot \right\Vert _{\mathbf{H}})$
such as wavelets, splines, Fourier series, Hermite polynomial series, etc.
And $k(n)\rightarrow \infty $ as $n\rightarrow \infty $.

For the purely nonparametric conditional moment models $E\left[ \rho
(Y,X;h_{0})|X\right] =0$, \cite{CP_WP07} proposed more general approximate
PSMD estimators of $h_{0}$ by allowing for possibly infinite dimensional
sieves (i.e., $\dim (\mathcal{H}_{k(n)})=k(n)\leq \infty $). Nevertheless,
both the theoretical properties and Monte Carlo simulations in \cite{CP_WP07}
recommend the use of the PSMD procedures with slowly growing
finite-dimensional linear sieves with a tiny penalty (i.e., $k(n)\rightarrow
\infty ,\frac{k(n)}{n}\rightarrow 0$ as $n$ $\rightarrow \infty $ with a
very small $\lambda _{n}=o(n^{-1})$, and hence the main smoothing parameter
is the sieve dimension $k(n)$). This class of PSMD estimators include the
original SMD estimators of \cite{NP_ECMA03} and \cite{AC_Emetrica03} as
special cases, and has been used in recent empirical estimation of
semiparametric structural models in microeconomics and asset pricing with
endogeneity. See, e.g., \cite{BCK_Emetrica07}, \cite{Horowitz_ECMA11}, \cite%
{CP_WP07a}, \cite{BHN_WP11}, \cite{Souza2012}, \cite{Pinkse2002}, \cite%
{MerloPaula}, \cite{BMT_WP12}, \cite{ChenLudvigson_JAE09}, \cite%
{ChenLudvigson2013}, \cite{Sentana2013} and others.

In this paper we shall develop inferential theory for $\phi (\alpha _{0})$
based on the PSMD procedures with slowly growing finite-dimensional sieves $%
\mathcal{A}_{k(n)}=\Theta \times \mathcal{H}_{k(n)}$. We first establish the
large sample theories under a high level \textquotedblleft local quadratic
approximation\textquotedblright\ (LQA) condition, which allows for any
consistent nonparametric estimator $\widehat{m}(x,\alpha )$ that is linear
in $\rho (Z,\alpha )$:
\begin{equation}
\widehat{m}(x,\alpha )\equiv \sum_{i=1}^{n}\rho (Z_{i},\alpha )A_{n}(X_{i},x)
\label{mhat-linear}
\end{equation}%
where $A_{n}(X_{i},x)$ is a known measurable function of $%
\{X_{j}\}_{j=1}^{n} $ for all $x$, whose expression varies according to
different nonparametric procedures such as kernel, local linear regression,
series and nearest neighbors. In Appendix \ref{app:appA} we provide lower
level sufficient conditions for this LQA assumption when $\widehat{m}%
(x,\alpha )$ is the series least squares (LS) estimator (\ref{mhat}):
\begin{equation}
\widehat{m}(x,\alpha )=\left( \sum_{i=1}^{n}\rho (Z_{i},\alpha
)p^{J_{n}}(X_{i})^{\prime }\right) (P^{\prime }P)^{-}p^{J_{n}}(x),
\label{mhat}
\end{equation}%
which is a linear nonparametric estimator (\ref{mhat-linear}) with $%
A_{n}(X_{i},x)=p^{J_{n}}(X_{i})^{\prime }(P^{\prime }P)^{-}p^{J_{n}}(x)$,
where $\{p_{j}\}_{j=1}^{\infty }$ is a sequence of known basis functions
that can approximate any square integrable functions of $X$ well, $%
p^{J_{n}}(X)=(p_{1}(X),...,p_{J_{n}}(X))^{\prime }$, $%
P=(p^{J_{n}}(X_{1}),...,p^{J_{n}}(X_{n}))^{\prime }$, and $(P^{\prime
}P)^{-} $ is the generalized inverse of the matrix $P^{\prime }P$. Following
\cite{BCK_Emetrica07} and \cite{CP_WP07a}, we let $p^{J_{n}}(X)$ be a
tensor-product linear sieve basis, and $J_{n}$ be the dimension of $%
p^{J_{n}}(X)$ such that $J_{n}\geq d_{\theta }+k(n)\rightarrow \infty $ and $%
\frac{J_{n}}{n}\rightarrow 0$ as $n$ $\rightarrow \infty $.

\subsection{Preview of the Main Results for Inference}

\label{sec:NPIVex}

For simplicity we let $\phi :\mathbb{R}^{d_{\theta }}\times \mathcal{H}%
\rightarrow \mathbb{R}$ be a real-valued functional. Let $\widehat{\phi }%
_{n}\equiv \phi (\widehat{\alpha }_{n})$ be the \textit{plug-in PSMD
estimator} of $\phi (\alpha _{0})$ for $\alpha _{0}=(\theta _{0}^{\prime
},h_{0})\in int(\Theta )\times \mathcal{H}$.

\textbf{Sieve t (or Wald) statistic}. Regardless of whether $\phi (\alpha
_{0})$ is $\sqrt{n}$ estimable or not, Theorem \ref{thm:theta_anorm} shows
that $\frac{\sqrt{n}\left\{ \phi (\widehat{\alpha }_{n})-\phi (\alpha
_{0})\right\} }{||v_{n}^{\ast }||_{sd}}$ is asymptotically standard normal,
and the sieve variance $||v_{n}^{\ast }||_{sd}^{2}$ has a \textit{closed form%
} expression resembling the \textquotedblleft
delta-method\textquotedblright\ variance for a parametric MD problem:
\begin{equation}
||v_{n}^{\ast }||_{sd}^{2}=\left( \frac{d\phi (\alpha _{0})}{d\alpha }[%
\overline{q}^{k(n)}(\cdot )]\right) ^{\prime }D_{n}^{-}\mho
_{n}D_{n}^{-}\left( \frac{d\phi (\alpha _{0})}{d\alpha }[\overline{q}%
^{k(n)}(\cdot )]\right) ,  \label{P-var}
\end{equation}%
where $\overline{q}^{k(n)}(\cdot )\equiv \left( \mathbf{1}_{d_{\theta
}}^{\prime },q^{k(n)}(\cdot )^{\prime }\right) ^{\prime }$ is a $(d_{\theta
}+k(n))\times 1$ vector with $\mathbf{1}_{d_{\theta }}$ a $d_{\theta }\times
1$ vector of $1$'s,
\begin{equation}
\frac{d\phi (\alpha _{0})}{d\alpha }[\overline{q}^{k(n)}(\cdot )]\equiv
\frac{\partial \phi (\theta _{0}+\theta ,h_{0}+\beta ^{\prime
}q^{k(n)}(\cdot ))}{\partial \gamma ^{\prime }} \mid_{\gamma =0}\equiv
\left( \frac{\partial \phi (\alpha _{0})}{\partial \theta ^{\prime }},\frac{%
d\phi (\alpha _{0})}{dh}[q^{k(n)}(\cdot )^{\prime }]\right) ^{\prime }
\label{P_kT_tilde}
\end{equation}%
and $\gamma \equiv (\theta ^{\prime },\beta ^{\prime })^{\prime }$ are $%
(d_{\theta }+k(n))\times 1$ vectors, $\frac{d\phi (\alpha _{0})}{dh}%
[q^{k(n)}(\cdot )^{\prime }]\equiv \frac{\partial \phi (\theta
_{0},h_{0}+\beta ^{\prime }q^{k(n)}(\cdot ))}{\partial \beta } \mid_{\beta =0}$, and%
\begin{equation}
D_{n}=E\left[ \left( \frac{dm(X,\alpha _{0})}{d\alpha }[\overline{q}%
^{k(n)}(\cdot )^{\prime }]\right) ^{\prime }\Sigma (X)^{-1}\left( \frac{%
dm(X,\alpha _{0})}{d\alpha }[\overline{q}^{k(n)}(\cdot )^{\prime }]\right) %
\right] ,  \label{P-R}
\end{equation}%
{\small{\begin{equation}
\mho _{n}=E\left[ \left( \frac{dm(X,\alpha _{0})}{d\alpha }[\overline{q}%
^{k(n)}(\cdot )^{\prime }]\right) ^{\prime }\Sigma (X)^{-1}\rho (Z,\alpha
_{0})\rho (Z,\alpha _{0})^{\prime }\Sigma (X)^{-1}\left( \frac{dm(X,\alpha
_{0})}{d\alpha }[\overline{q}^{k(n)}(\cdot )^{\prime }]\right) \right] ,
\label{P-Omega}
\end{equation}}}
where $\frac{dm(X,\alpha _{0})}{d\alpha }[\overline{q}^{k(n)}(\cdot
)^{\prime }]\equiv \frac{\partial E[\rho (Z,\theta _{0}+\theta ,h_{0}+\beta
^{\prime }q^{k(n)}(\cdot ))|X]}{\partial \gamma } \mid_{\gamma =0}$ is
a $d_{\rho }\times (d_{\theta }+k(n))$ matrix. The closed form expression of
$||v_{n}^{\ast }||_{sd}^{2}$ immediately leads to simple consistent plug-in
sieve variance estimators; one of which is
\begin{equation}
||\widehat{v}_{n}^{\ast }||_{n,sd}^{2}=\widehat{V}_{1}=\left( \frac{d\phi (%
\widehat{\alpha }_{n})}{d\alpha }[\overline{q}^{k(n)}(\cdot )]\right)
^{\prime }\widehat{D}_{n}^{-}\widehat{\mho }_{n}\widehat{D}_{n}^{-}\left(
\frac{d\phi (\widehat{\alpha }_{n})}{d\alpha }[\overline{q}^{k(n)}(\cdot
)]\right) ,  \label{P-var-hat}
\end{equation}%
where $\frac{d\phi (\widehat{\alpha }_{n})}{d\alpha }[\overline{q}%
^{k(n)}(\cdot )]\equiv \frac{\partial \phi (\widehat{\theta }_{n}+\theta ,%
\widehat{h}_{n}+\beta ^{\prime }q^{k(n)}(\cdot ))}{\partial \gamma ^{\prime }%
} \mid_{\gamma =0}$ and
\begin{equation}
\widehat{D}_{n}=\frac{1}{n}\sum_{i=1}^{n}\left[ \left( \frac{d\widehat{m}%
(X_{i},\widehat{\alpha }_{n})}{d\alpha }[\overline{q}^{k(n)}(\cdot )^{\prime
}]\right) ^{\prime }\widehat{\Sigma }(X_{i})^{-1}\left( \frac{d\widehat{m}%
(X_{i},\widehat{\alpha }_{n})}{d\alpha }[\overline{q}^{k(n)}(\cdot )^{\prime
}]\right) \right] ,  \label{P-R-hat}
\end{equation}%
\begin{equation}
\widehat{\mho }_{n}=\frac{1}{n}\sum_{i=1}^{n}\left[ \left( \frac{d\widehat{m}%
(X_{i},\widehat{\alpha }_{n})}{d\alpha }[\overline{q}^{k(n)}(\cdot )^{\prime
}]\right) ^{\prime }\widehat{M}(X_{i})\left( \frac{d\widehat{m}(X_{i},\widehat{\alpha }_{n})}{d\alpha }%
[\overline{q}^{k(n)}(\cdot )^{\prime }]\right) \right] .  \label{P-Omega-hat}
\end{equation}
where $\widehat{M}(X_{i}) \equiv \widehat{\Sigma }(X_{i})^{-1}\rho (Z_{i},\widehat{\alpha
}_{n})\rho (Z_{i},\widehat{\alpha }_{n})^{\prime }\widehat{\Sigma }%
(X_{i})^{-1}$. Theorem \ref{thm:VE} then presents the asymptotic normality of the sieve
(Student's) t statistic:\footnote{%
See Theorems \ref{thm:bootstrap_2} and \ref{thm:waldB_con} for properties of
bootstrap sieve t statistics.}%
\begin{equation*}
\widehat{W}_{n}\equiv \sqrt{n}\frac{\phi (\widehat{\alpha }_{n})-\phi
(\alpha _{0})}{||\widehat{v}_{n}^{\ast }||_{n,sd}}\Rightarrow N(0,1).
\end{equation*}

\textbf{Sieve QLR statistic}. In addition to the sieve t (or sieve Wald)
statistic, we could also use sieve quasi likelihood ratio for constructing
confidence set of $\phi (\alpha _{0})$ and for hypothesis testing of $%
H_{0}:\phi (\alpha _{0})=\phi _{0}$ against $H_{1}:\phi (\alpha _{0})\neq
\phi _{0}$. Denote
\begin{equation}
\widehat{QLR}_{n}(\phi _{0})\equiv n\left( \inf_{\alpha \in \mathcal{A}%
_{k(n)}:\phi (\alpha )=\phi _{0}}\widehat{Q}_{n}(\alpha )-\widehat{Q}_{n}(%
\widehat{\alpha }_{n})\right)  \label{SQLR}
\end{equation}%
as the \textit{sieve quasi likelihood ratio} (SQLR) statistic. It becomes an
\textit{optimally weighted SQLR} statistic, $\widehat{QLR}_{n}^{0}(\phi
_{0}) $, when $\widehat{Q}_{n}(\alpha )$ is the optimally weighted MD
criterion $\widehat{Q}_{n}^{0}(\alpha )$. Regardless of whether $\phi
(\alpha _{0})$ is $\sqrt{n}$ estimable or not, Theorems \ref{thm:chi2}(2)
and \ref{thm:QLR-H1} show that $\widehat{QLR}_{n}^{0}(\phi _{0})$ is
asymptotically chi-square distributed under the null $H_{0}$, and diverges
to infinity under the fixed alternatives $H_{1}$. Theorem \ref%
{thm:chi2_localt} in Appendix \ref{app:appA} states that $\widehat{QLR}%
_{n}^{0}(\phi _{0})$ is asymptotically noncentral chi-square distributed
under local alternatives. One could compute $100(1-\tau )\%$ confidence set
for $\phi (\alpha _{0})$ as
\begin{equation*}
\left\{ r\in \mathbb{R}\colon \text{ }\widehat{QLR}_{n}^{0}(r)\leq c_{\chi
_{1}^{2}}(1-\tau )\right\} ,
\end{equation*}%
where $c_{\chi _{1}^{2}}(1-\tau )$ is the $(1-\tau )$-th quantile of the $%
\chi _{1}^{2}$ distribution.

\textbf{Bootstrap sieve QLR statistic}. Regardless of whether $\phi (\alpha
_{0})$ is $\sqrt{n}$ estimable or not, Theorems \ref{thm:chi2}(1) and \ref%
{thm:QLR-H1} establish that the possibly non-optimally weighted SQLR
statistic $\widehat{QLR}_{n}(\phi _{0})$ is stochastically bounded under the
null $H_{0}$ and diverges to infinity under the fixed alternatives $H_{1}$.
We then consider a bootstrap version of the SQLR statistic. Let $\widehat{QLR%
}_{n}^{B}$ denote a bootstrap SQLR statistic:
\begin{equation}
\widehat{QLR}_{n}^{B}(\widehat{\phi }_{n})\equiv n\left( \inf_{\alpha \in
\mathcal{A}_{k(n)}:\phi (\alpha )=\widehat{\phi }_{n}}\widehat{Q}%
_{n}^{B}(\alpha )-\inf_{\alpha \in \mathcal{A}_{k(n)}}\widehat{Q}%
_{n}^{B}(\alpha )\right) ,  \label{SQLR-B}
\end{equation}%
where $\widehat{\phi }_{n}\equiv \phi (\widehat{\alpha }_{n})$, and $%
\widehat{Q}_{n}^{B}(\alpha )$ is a bootstrap version of $\widehat{Q}%
_{n}(\alpha )$:%
\begin{equation}
\widehat{Q}_{n}^{B}(\alpha )\equiv \frac{1}{n}\sum_{i=1}^{n}\widehat{m}%
^{B}(X_{i},\alpha )^{\prime }\widehat{\Sigma }(X_{i})^{-1}\widehat{m}%
^{B}(X_{i},\alpha ),  \label{QhatB}
\end{equation}%
where $\widehat{m}^{B}(x,\alpha )$ is a bootstrap version of $\widehat{m}%
(x,\alpha )$, which is computed in the same way as that of $\widehat{m}%
(x,\alpha )$ except that we use $\omega _{i,n}\rho (Z_{i},\alpha )$ instead
of $\rho (Z_{i},\alpha )$. Here $\{\omega _{i,n}\geq 0\}_{i=1}^{n}$ is a
sequence of bootstrap weights that has mean 1 and is independent of the
original data $\{Z_{i}\}_{i=1}^{n}$. Typical weights include an i.i.d.
weight $\{\omega _{i}\geq 0\}_{i=1}^{n}$ with $E[\omega _{i}]=1$, $E[|\omega
_{i}-1|^{2}]=1$ and $E[|\omega _{i}-1|^{2+\epsilon }]<\infty $ for some $%
\epsilon >0$, or a multinomial weight (i.e., $(\omega _{1,n},...,\omega
_{n,n})\sim Multinomial(n;n^{-1},...,n^{-1})$). For example, if $\widehat{m}%
(x,\alpha )$ is a series LS estimator (\ref{mhat}) of $m(x,\alpha )$, then $%
\widehat{m}^{B}(x,\alpha )$ is a bootstrap series LS estimator of $%
m(x,\alpha )$, defined as:%
\begin{equation}
\widehat{m}^{B}(x,\alpha )\equiv \left( \sum_{i=1}^{n}\omega _{i,n}\rho
(Z_{i},\alpha )p^{J_{n}}(X_{i})^{\prime }\right) (P^{\prime
}P)^{-}p^{J_{n}}(x).  \label{mhat_B}
\end{equation}%
We sometimes call our bootstrap procedure \textquotedblleft \textit{%
generalized residual bootstrap}\textquotedblright\ since it is based on
randomly perturbing the generalized residual function $\rho (Z,\alpha )$;
see Section \ref{sec:bootstrap} for details. Theorems \ref{thm:bootstrap}
and \ref{thm:BSQLR-loc-alt} establish that under the null $H_{0}$, the fixed
alternatives $H_{1}$ or the local alternatives,\footnote{%
See Section \ref{subsec-A5} for definition of the local alternatives and the
behaviors of $\widehat{QLR}_{n}(\phi _{0})$ and $\widehat{QLR}_{n}^{B}(%
\widehat{\phi }_{n})$ under the local alternatives.} the conditional
distribution of $\widehat{QLR}_{n}^{B}(\widehat{\phi }_{n})$ (given the
data) always converges to the asymptotic null distribution of $\widehat{QLR}%
_{n}(\phi _{0})$. Let $\widehat{c}_{n}(a)$ be the $a-th$ quantile of the
distribution of $\widehat{QLR}_{n}^{B}(\widehat{\phi }_{n})$ (conditional on
the data $\{Z_{i}\}_{i=1}^{n}$). Then for any $\tau \in (0,1)$, we have $%
\lim_{n\rightarrow \infty }\Pr \{\widehat{QLR}_{n}(\phi _{0})>\widehat{c}%
_{n}(1-\tau )\}=\tau $ under the null $H_{0}$, $\lim_{n\rightarrow \infty
}\Pr \{\widehat{QLR}_{n}(\phi _{0})>\widehat{c}_{n}(1-\tau )\}=1$ under the
fixed alternatives $H_{1}$, and $\lim_{n\rightarrow \infty }\Pr \{\widehat{%
QLR}_{n}(\phi _{0})>\widehat{c}_{n}(1-\tau )\}>\tau $ under the local
alternatives. We could also construct a $100(1-\tau )\%$ confidence set
using the bootstrap critical values:
\begin{equation}
\left\{ r\in \mathbb{R}\colon \text{ }\widehat{QLR}_{n}(r)\leq \widehat{c}%
_{n}(1-\tau )\right\} .  \label{boot-cs}
\end{equation}%
The bootstrap consistency holds for possibly non-optimally weighted SQLR
statistic and possibly irregular functionals, without the need to compute
standard errors.

\textbf{Which method to use?} When sieve Wald and SQLR tests are computed
using the same weighting matrix $\widehat{\Sigma }$, there is no local power
difference in terms of first order asymptotic theories; see Appendix \ref%
{app:appA}. As will be demonstrated in simulation Section \ref%
{sec:sec_simulation}, while SQLR and bootstrap SQLR tests are useful for
models (\ref{semi00}) with (pointwise) non-smooth $\rho (Z;\alpha )$, sieve
Wald (or t) statistic is computationally attractive for models with smooth $%
\rho (Z;\alpha )$. Empirical researchers could apply either inference method
depending on whether the residual function $\rho (Z;\alpha )$ in their
specific application is pointwise differentiable with respect to $\alpha $
or not.

\subsubsection{Applications to NPIV and NPQIV models\label{sec:NPIVex1}}

\textbf{An illustration via the NPIV model. }\cite{BCK_Emetrica07} and \cite%
{CR_WP07} established the convergence rate of the identity weighted (i.e., $%
\widehat{\Sigma }=\Sigma =1$) PSMD estimator $\widehat{h}_{n}\in \mathcal{H}%
_{k(n)}$ of the NPIV model:%
\begin{equation}
Y_{1}=h_{0}(Y_{2})+U,\text{\quad }E(U|X)=0.  \label{npiv}
\end{equation}%
By Theorem \ref{thm:theta_anorm} \begin{align*}
	\sqrt{n}\frac{\phi (\widehat{h}_{n})-\phi
		(h_{0})}{||v_{n}^{\ast }||_{sd}}\Rightarrow N(0,1)
\end{align*} with $||v_{n}^{\ast
}||_{sd}^{2}=\frac{d\phi (h_{0})}{dh}[q^{k(n)}(\cdot )]^{\prime
}D_{n}^{-}\mho _{n}D_{n}^{-}\frac{d\phi (h_{0})}{dh}[q^{k(n)}(\cdot )]$,%
\begin{align}
D_{n}&=E\left( E[q^{k(n)}(Y_{2})|X]E[q^{k(n)}(Y_{2})|X]^{\prime }\right)
\text{, }\\
\mho _{n}&=E\left(
E[q^{k(n)}(Y_{2})|X]U^{2}E[q^{k(n)}(Y_{2})|X]^{\prime }\right)
\label{npiv-D}
\end{align}%
and $\frac{d\phi (h_{0})}{dh}[q^{k(n)}(\cdot )]\equiv \frac{\partial \phi
(h_{0}+\beta ^{\prime }q^{k(n)}(\cdot ))}{\partial \beta ^{\prime }} \mid_{\beta =0}$. For example, for a functional $\phi (h)=h(\overline{y}_{2})$,
or $=\int w(y)\nabla h(y)dy$ or $=\int w(y)\left\vert h(y)\right\vert ^{2}dy$%
, we have $\frac{d\phi (h_{0})}{dh}[q^{k(n)}(\cdot )]=q^{k(n)}(\overline{y}%
_{2})$, or $=\int w(y)\nabla q^{k(n)}(y)dy$ or $=2\int
h_{0}(y)w(y)q^{k(n)}(y)dy$.

If $0<\inf_{x}\Sigma _{0}(x)\leq \sup_{x}\Sigma _{0}(x)<\infty $ then \begin{align*}
||v_{n}^{\ast }||_{sd}^{2}\asymp \frac{d\phi (h_{0})}{dh}[q^{k(n)}(\cdot
)]^{\prime }D_{n}^{-}\frac{d\phi (h_{0})}{dh}[q^{k(n)}(\cdot )]
\end{align*}
Without
endogeneity (say $Y_{2}=X$) the model becomes the nonparametric LS
regression
\begin{equation*}
Y_{1}=h_{0}(Y_{2})+U,\text{\quad }E(U|Y_{2})=0,
\end{equation*}%
and the variance satisfies $||v_{n}^{\ast }||_{sd,ex}^{2}\asymp \frac{d\phi
(h_{0})}{dh}[q^{k(n)}(\cdot )]^{\prime }D_{n,ex}^{-}\frac{d\phi (h_{0})}{dh}%
[q^{k(n)}(\cdot )]$, $D_{n,ex}=E[\{q^{k(n)}(Y_{2})\}\{q^{k(n)}(Y_{2})\}^{%
\prime }]$. Since the conditional expectation $E[q^{k(n)}(Y_{2})|X]$ is a
contraction, $D_{n}\leq D_{n,ex}$ and $||v_{n}^{\ast }||_{sd}^{2}\geq
const.||v_{n}^{\ast }||_{sd,ex}^{2}$. Under mild conditions (see, e.g., \cite%
{NP_ECMA03}, \cite{BCK_Emetrica07}, \cite{DFR_wp10}, \cite{Horowitz_ECMA11}%
), the minimal eigenvalue of $D_{n}$, $\lambda _{\min }(D_{n})$, goes to
zero while $\lambda _{\min }(D_{n,ex})$ stays strictly positive as $%
k(n)\rightarrow \infty $. In fact, $D_{n,ex}=I_{k(n)}$ and $\lambda _{\min
}(D_{n,ex})=1$ if $\{q_{j}\}_{j=1}^{\infty }$ is an orthonormal basis of $%
L^{2}(f_{Y_{2}})$, while $\lambda _{\min }(D_{n})\asymp \exp (-k(n))$ if the
conditional density of $Y_{2}$ given $X$ is normal. Therefore, while $%
\lim_{k(n)\rightarrow \infty }||v_{n}^{\ast }||_{sd,ex}^{2}=\infty $ always
implies $\lim_{k(n)\rightarrow \infty }||v_{n}^{\ast }||_{sd}^{2}=\infty $,
it is possible that $\lim_{k(n)\rightarrow \infty }||v_{n}^{\ast
}||_{sd,ex}^{2}<\infty $ but $\lim_{k(n)\rightarrow \infty }||v_{n}^{\ast
}||_{sd}^{2}=\infty $. For example, the point evaluation functional $\phi
(h)=h(\overline{y}_{2})$ is known to be irregular for the nonparametric LS
regression and hence for the NPIV (\ref{npiv}) as well. Under mild
conditions on the weight $w()$ and the smoothness of $h_{0}$, the weighted
derivative functional ($\phi (h)=\int w(y)\nabla h(y)dy$) and the quadratic
functional ($\phi (h)=\int w(y)\left\vert h(y)\right\vert ^{2}dy$) of the
nonparametric LS regression are typically root-$n$ estimable, but they could
be irregular for the NPIV (\ref{npiv}). See Section \ref{sec-ex} for details.

Regardless of whether $\lim_{k(n)\rightarrow \infty }||v_{n}^{\ast
}||_{sd}^{2}$ is finite or infinite, Theorem \ref{thm:VE} shows that the
sieve variance $||v_{n}^{\ast }||_{sd}^{2}$ can be consistently estimated by
a plug-in sieve variance estimator $||\widehat{v}_{n}^{\ast }||_{n,sd}^{2}$,
and that $\sqrt{n}\frac{\phi (\widehat{h}_{n})-\phi (h_{0})}{||\widehat{v}%
_{n}^{\ast }||_{n,sd}}\Rightarrow N(0,1)$.

When the conditional mean function $m(x,h)$ is estimated by the series LS
estimator (\ref{mhat}) as in \cite{NP_ECMA03}, \cite{AC_Emetrica03} and \cite%
{BCK_Emetrica07}, with $\widehat{U}_{i}=Y_{1i}-\widehat{h}_{n}(Y_{2i})$, the
sieve variance estimator $||\widehat{v}_{n}^{\ast }||_{n,sd}^{2}$ given in (%
\ref{P-var-hat}) has a more explicit expression:%
\begin{equation*}
||\widehat{v}_{n}^{\ast }||_{n,sd}^{2}=\widehat{V}_{1}=\left( \frac{d\phi (%
\widehat{h}_{n})}{dh}[q^{k(n)}(\cdot )]\right) ^{\prime }\widehat{D}_{n}^{-}%
\widehat{\mho }_{n}\widehat{D}_{n}^{-}\left( \frac{d\phi (\widehat{h}_{n})}{%
dh}[q^{k(n)}(\cdot )]\right) ,\quad \text{where}
\end{equation*}

$\frac{d\phi (\widehat{h}_{n})}{dh}[q^{k(n)}(\cdot )]\equiv \frac{\partial
\phi (\widehat{h}_{n}+\beta ^{\prime }q^{k(n)}(\cdot ))}{\partial \beta
^{\prime }} \mid_{\beta =0}$ and%
\begin{equation*}
\widehat{D}_{n}=\frac{1}{n}\widehat{C}_{n}(P^{\prime }P)^{-}(\widehat{C}%
_{n})^{\prime },\quad \widehat{C}_{n}\equiv
\sum_{j=1}^{n}q^{k(n)}(Y_{2j})p^{J_{n}}(X_{j})^{\prime },
\end{equation*}%
\begin{equation}
\widehat{\mho }_{n}=\frac{1}{n}\widehat{C}_{n}(P^{\prime }P)^{-}\left(
\sum_{i=1}^{n}p^{J_{n}}(X_{i})\widehat{U}_{i}^{2}p^{J_{n}}(X_{i})^{\prime
}\right) (P^{\prime }P)^{-}(\widehat{C}_{n})^{\prime }.  \label{2sls-vhat}
\end{equation}%
Interestingly, this sieve variance estimator becomes the one computed via
the two stage least squares (2SLS) as if the NPIV model (\ref{npiv}) were a
parametric IV regression:\footnote{This confirms a conjecture of \cite{Newey2013} for the NPIV model (\ref{npiv}%
).} $Y_{1}=q^{k(n)}(Y_{2j})^{\prime }\beta _{0n}+U,$ $%
E[q^{k(n)}(Y_{2})U]\neq 0,$ $E[p^{J_{n}}(X)U]=0$ and $%
E[p^{J_{n}}(X)q^{k(n)}(Y_{2})^{\prime }]$ has a column rank $k(n)\leq J_{n}$%
. See Subsection \ref{sec:sec_simulation1} for simulation studies of finite
sample performances of this sieve variance estimator $\widehat{V}_{1}$ for
both a linear and a nonlinear functional $\phi (h)$.
\medskip

\textbf{An illustration via the NPQIV model. }As an application of their
general theory, \cite{CP_WP07} presented the consistency and the rate of
convergence of the PSMD estimator $\widehat{h}_{n}\in \mathcal{H}_{k(n)}$ of
the NPQIV model:%
\begin{equation}
Y_{1}=h_{0}(Y_{2})+U,\text{\quad }\Pr (U\leq 0|X)=\gamma .  \label{npqiv}
\end{equation}%
In this example we have $\Sigma _{0}(X)=\gamma (1-\gamma )$. So we could use
$\widehat{\Sigma }(X)=\gamma (1-\gamma )$ and $\widehat{Q}_{n}(\alpha )$
given in (\ref{Qhat}) becomes the optimally weighted\ MD criterion.

\noindent By Theorem \ref{thm:theta_anorm}
\begin{align*}
\sqrt{n}\frac{\phi (\widehat{h}%
	_{n})-\phi (h_{0})}{||v_{n}^{\ast }||_{sd}}\Rightarrow N(0,1)
\end{align*}
with $%
||v_{n}^{\ast }||_{sd}^{2}=\left( \frac{d\phi (h_{0})}{dh}[q^{k(n)}(\cdot
)]\right) ^{\prime }D_{n}^{-}\left( \frac{d\phi (h_{0})}{dh}[q^{k(n)}(\cdot
)]\right) $ and%
\begin{equation}
D_{n}=\frac{1}{\gamma (1-\gamma )}E\left(
E[f_{U|Y_{2},X}(0)q^{k(n)}(Y_{2})|X]E[f_{U|Y_{2},X}(0)q^{k(n)}(Y_{2})|X]^{%
\prime }\right) .  \label{npqiv-D}
\end{equation}%
Without endogeneity (say $Y_{2}=X$), the model becomes the nonparametric
quantile regression
\begin{equation*}
Y_{1}=h_{0}(Y_{2})+U,\text{\quad }\Pr (U\leq 0|Y_{2})=\gamma ,
\end{equation*}%
and the sieve variance becomes $||v_{n}^{\ast }||_{sd,ex}^{2}=\left( \frac{%
d\phi (h_{0})}{dh}[q^{k(n)}(\cdot )]\right) ^{\prime }D_{n,ex}^{-}\left(
\frac{d\phi (h_{0})}{dh}[q^{k(n)}(\cdot )]\right) $ with $D_{n,ex}=\frac{1}{%
\gamma (1-\gamma )}E\left[ \{f_{U|Y_{2}}(0)\}^{2}\{q^{k(n)}(Y_{2})\}%
\{q^{k(n)}(Y_{2})\}^{\prime }\right] $. Again $D_{n}\leq D_{n,ex}$ and $%
||v_{n}^{\ast }||_{sd}^{2}\geq ||v_{n}^{\ast }||_{sd,ex}^{2}$. Under mild
conditions (see, e.g., \cite{CP_WP07}, \cite{CCLN_WP10}), $\lambda _{\min
}(D_{n})\rightarrow 0$ while $\lambda _{\min }(D_{n,ex})$ stays strictly
positive as $k(n)\rightarrow \infty $. All of the above discussions for a
functional $\phi (h)$ of the NPIV (\ref{npiv}) now apply to the functional
of the NPQIV (\ref{npqiv}). In particular, a functional $\phi (h)$ could be
root-$n$ estimable for the nonparametric quantile regression ($%
\lim_{k(n)\rightarrow \infty }||v_{n}^{\ast }||_{sd,ex}^{2}<\infty $) but
irregular for the NPQIV (\ref{npqiv}) ($\lim_{k(n)\rightarrow \infty
}||v_{n}^{\ast }||_{sd}^{2}=\infty $). See Section \ref{sec-ex} for details.

By Theorems \ref{thm:chi2}(2) and \ref{thm:QLR-H1}, the optimally weighted
SQLR statistic $\widehat{QLR}_{n}^{0}(\phi _{0})\Rightarrow \chi _{1}^{2}$
under the null of $\phi (h_{0})=\phi _{0}$, and diverges to infinity under
the alternative of $\phi (h_{0})\neq \phi _{0}$. We can compute confidence
set for a functional $\phi (h)$, such as an evaluation or a weighted
derivative functional, as $\left\{ r\in \mathbb{R}\colon \text{ }\widehat{QLR%
}_{n}^{0}(r)\leq c_{\chi _{1}^{2}}(\tau )\right\} $. See Subsection \ref%
{sec:application} for an empirical illustration of this result to the NPQIV
Engel curve regression using the British Family Survey data set that was
first used in \cite{BCK_Emetrica07}. Instead of using the asymptotic
critical values, we could also construct a confidence set using the
bootstrap critical values as in (\ref{boot-cs}).

\section{Basic Regularity Conditions}
\label{sec:conditions}

Before we establish asymptotic properties of sieve t (Wald) and SQLR
statistics, we need to present three sets of basic regularity conditions.
The first set of assumptions allows us to establish the convergence rates of
the PSMD estimator $\widehat{\alpha }_{n}$ to the true parameter value $%
\alpha _{0}$ in both weak and strong metrics, which in turn allows us to
concentrate on some shrinking neighborhood of $\alpha _{0}$ in the
semi/nonparametric model (\ref{semi00}). The second and third regularity
conditions are respectively about the local curvatures of the functional $%
\phi ()$ and of the criterion function under these two metrics. The weak
metric $||\cdot ||$ is closely related to the variance of the linear
approximation to $\phi (\widehat{\alpha }_{n})-\phi (\alpha _{0})$, while
the strong metric $||\cdot ||_{s}$ is used to control the nonlinearity (in $%
\alpha $) of the functional $\phi ()$ and of the conditional mean function $%
m(x,\alpha )$. This section is mostly technical and applied researchers
could skip this and directly go to the subsequent sections on the asymptotic
properties of sieve Wald and SQLR statistics.

\subsection{A brief discussion on the convergence rate of the PSMD estimator
\label{sec:consistency}}

For the purely nonparametric conditional moment model $E\left[ \rho
(Y,X;h_{0}(\cdot ))|X\right] =0$, \cite{CP_WP07} established the consistency
and the convergence rates of their various PSMD estimators of $h_{0}$. Their
results can be trivially extended to establish the corresponding properties
of our PSMD estimator $\widehat{\alpha }_{n}\equiv (\widehat{\theta }%
_{n}^{\prime },\widehat{h}_{n})$ defined in (\ref{psmd}). For the sake of
easy reference and to introduce basic assumptions and notation, we present
some sufficient conditions for consistency and the convergence rate here.
These conditions are also needed to establish the consistency and the
convergence rate of bootstrap PSMD estimators (see Lemma \ref{lem:cons_boot}%
). We first impose three conditions on identification, sieve spaces, penalty
functions and sample criterion function. We equip the parameter space $%
\mathcal{A}\equiv \Theta \times \mathcal{H}$ with a (strong) norm $%
\left\Vert \alpha \right\Vert _{s}\equiv \left\Vert \theta \right\Vert
_{e}+\left\Vert h\right\Vert _{\mathbf{H}}$.

\begin{assumption}[Identification, sieves, criterion]
\label{ass:sieve} (i) $E[\rho (Y,X;\alpha )|X]=0$ if and only if $\alpha \in
(\mathcal{A},\left\Vert \cdot \right\Vert _{s})$ with $\left\Vert \alpha
-\alpha _{0}\right\Vert _{s}=0$; (ii) For all $k\geq 1$, $\mathcal{A}%
_{k}\equiv \Theta \times \mathcal{H}_{k}$, $\Theta $ is a compact subset in $%
\mathbb{R}^{d_{\theta }}$ with a non-empty interior, $\{\mathcal{H}%
_{k}:k\geq 1\}$ is a non-decreasing sequence of non-empty closed linear
subsets of a Banach space $\left( \mathcal{H},\left\Vert \cdot \right\Vert _{%
\mathbf{H}}\right) $ such that $\mathcal{H}=cl\left( \cup _{k}\mathcal{H}%
_{k}\right) $, and there is $\Pi _{n}h_{0}\in \mathcal{H}_{k(n)}$ with $%
||\Pi _{n}h_{0}-h_{0}||_{\mathbf{H}}=o(1)$; (iii) $Q:(\mathcal{A},\left\Vert
\cdot \right\Vert _{s})\rightarrow \lbrack 0,\infty )$ is lower
semicontinuous;\footnote{%
A function $Q$ is lower semicontinuous at a point $\alpha _{o}\in \mathcal{A}
$ iff $\lim_{\left\Vert \alpha -\alpha _{o}\right\Vert _{s}\rightarrow
0}Q(\alpha )\geq Q(\alpha _{o})$; is lower semicontinuous if it is lower
semicontinuous at any point in $\mathcal{A}$.} (iv) $\Sigma (x)$ and $\Sigma
_{0}(x)$ are positive definite, and their smallest and largest eigenvalues
are finite and positive uniformly in $x\in \mathcal{X}$.
\end{assumption}

\begin{assumption}[Penalty]
\label{A_3.6} (i) $\lambda _{n}>0$, $Q(\Pi _{n}\alpha
_{0})+o(n^{-1})=O(\lambda _{n})=o(1)$; (ii) $|Pen(\Pi
_{n}h_{0})-Pen(h_{0})|=O(1)$ with $Pen(h_{0})<\infty $; (iii) $Pen:(\mathcal{%
H},\left\Vert \cdot \right\Vert _{\mathbf{H}})\rightarrow \lbrack 0,\infty )$
is lower semicompact.\footnote{%
A function $Pen$ is lower semicompact iff for all $M$, $\{h\in \mathcal{H}%
\colon Pen(h)\leq M\}$ is a compact subset in $(\mathcal{H},\left\Vert \cdot
\right\Vert _{\mathbf{H}})$.}
\end{assumption}

Let $\Pi _{n}\alpha \equiv (\theta ^{\prime },\Pi _{n}h)\in \mathcal{A}%
_{k(n)}\equiv \Theta \times \mathcal{H}_{k(n)}$. Let $\mathcal{A}%
_{k(n)}^{M_{0}}\equiv \Theta \times \mathcal{H}_{k(n)}^{M_{0}}\equiv
\{\alpha =(\theta ^{\prime },h)\in \mathcal{A}_{k(n)}:\lambda _{n}Pen(h)\leq
\lambda _{n}M_{0}\}$ for a large but finite $M_{0}$ such that $\Pi
_{n}\alpha _{0}\in \mathcal{A}_{k(n)}^{M_{0}}$ and that $\widehat{\alpha }%
_{n}\in \mathcal{A}_{k(n)}^{M_{0}}$ with probability arbitrarily close to
one for all large $n$. Let $\{\bar{\delta}_{m,n}^{2}\}_{n=1}^{\infty }$ be a
sequence of positive real values that decrease to zero as $n\rightarrow
\infty $.

\begin{assumption}[Sample Criterion]
\label{ass:rates} (i) $\widehat{Q}_{n}(\Pi _{n}\alpha _{0})\leq c_{0}Q(\Pi
_{n}\alpha _{0})+o_{P_{Z^{\infty }}}(n^{-1})$ for a finite constant $c_{0}>0$%
; (ii) $\widehat{Q}_{n}(\alpha )\geq cQ(\alpha )-O_{P_{Z^{\infty }}}(\bar{%
\delta}_{m,n}^{2})$ uniformly over $\mathcal{A}_{k(n)}^{M_{0}}$ for some $%
\bar{\delta}_{m,n}^{2}=o(1)$ and a finite constant $c>0$.
\end{assumption}

The following result is a minor modification of Theorem 3.2 of \cite{CP_WP07}%
.

\begin{lemma}
\label{thm:Thm-ill-suff_pencompact3} Let $\widehat{\alpha }_{n}$ be the PSMD
estimator defined in (\ref{psmd}), and Assumptions \ref{ass:sieve}, \ref%
{A_3.6} and \ref{ass:rates} hold. Then: $||\widehat{\alpha }_{n}-\alpha
_{0}||_{s}=o_{P_{Z^{\infty }}}(1)$ and $Pen(\widehat{h}_{n})=O_{P_{Z^{\infty
}}}(1)$.
\end{lemma}

Given the consistency result, the PSMD estimator belongs to any $||\cdot
||_{s}-$neighborhood around $\alpha _{0}$ wpa1. We can restrict our
attention to a convex, $||\cdot ||_{s}-$neighborhood around $\alpha _{0}$,
denoted as $\mathcal{A}_{os}$ such that
\begin{equation*}
\mathcal{A}_{os}\subset \{\alpha \in \mathcal{A}:||\alpha -\alpha
_{0}||_{s}<M_{0},\text{ }\lambda _{n}Pen(h)<\lambda _{n}M_{0}\}
\end{equation*}%
for a positive finite constant $M_{0}$ (the existence of a convex $\mathcal{A%
}_{os}$ is implied by the convexity of $\mathcal{A}$ and quasi-convexity of $%
Pen(\cdot )$). For any $\alpha \in \mathcal{A}_{os}$ we define a pathwise
derivative as
\begin{eqnarray*}
\frac{dm(X,\alpha _{0})}{d\alpha }[\alpha -\alpha _{0}] &\equiv &\left.
\frac{dE[\rho (Z,(1-\tau )\alpha _{0}+\tau \alpha )|X]}{d\tau }\right\vert
_{\tau =0}\quad a.s.~X \\
&=&\frac{dE[\rho (Z,\alpha _{0})|X]}{d\theta ^{\prime }}(\theta -\theta
_{0})\\
&& +\frac{dE[\rho (Z,\alpha _{0})|X]}{dh}[h-h_{0}]\quad a.s.~X.
\end{eqnarray*}%
Following \cite{AC_Emetrica03} and \cite{CP_WP07a}, we introduce two
pseudo-metrics $||\cdot ||$ and $||\cdot ||_{0}$ on $\mathcal{A}_{os}$ as:
for any $\alpha _{1},$ $\alpha _{2}\in \mathcal{A}_{os}$,
\begin{equation}
||\alpha _{1}-\alpha _{2}||^{2}\equiv E\left[ \left( \frac{%
dm(X,\alpha _{0})}{d\alpha }[\alpha _{1}-\alpha _{2}]\right) ^{\prime }\Sigma (X)^{-1}\left( \frac{dm(X,\alpha _{0})}{d\alpha }%
[\alpha _{1}-\alpha _{2}]\right) \right] ;  \label{fmetric}
\end{equation}%
\begin{equation}
||\alpha _{1}-\alpha _{2}||_{0}^{2}\equiv E\left[ \left(
\frac{dm(X,\alpha _{0})}{d\alpha }[\alpha _{1}-\alpha _{2}]\right) ^{\prime }\Sigma _{0}(X)^{-1}\left( \frac{dm(X,\alpha _{0})}{%
d\alpha }[\alpha _{1}-\alpha _{2}]\right) \right] .
\label{fmetric0}
\end{equation}%
It is clear that, under Assumption \ref{ass:sieve}(iv), these two
pseudo-metrics are equivalent, i.e., $||\cdot ||\asymp ||\cdot ||_{0}$ on $%
\mathcal{A}_{os}$. This is why Assumption \ref{ass:sieve}(iv) is imposed
throughout the paper.

Let $\mathcal{A}_{osn}=\mathcal{A}_{os}\cap \mathcal{A}_{k(n)}$. Let $%
\{\delta _{n}\}_{n=1}^{\infty }$ be a sequence of positive real values such
that $\delta _{n}=o(1)$ and $\delta _{n}\leq \bar{\delta}_{m,n}$.

\begin{assumption}
\label{ass:weak_equiv} (i) There exists a convex $||\cdot ||_{s}-$%
neighborhood of $\alpha _{0}$, $\mathcal{A}_{os}$, such that $m(\cdot
,\alpha )$ is continuously pathwise differentiable with respect to $\alpha
\in \mathcal{A}_{os}$, and there is a finite constant $C>0$ such that $%
||\alpha -\alpha _{0}||\leq C||\alpha -\alpha _{0}||_{s}$ for all $\alpha
\in \mathcal{A}_{os}$; (ii) $Q(\alpha )\asymp ||\alpha -\alpha _{0}||^{2}$
for all $\alpha \in \mathcal{A}_{os}$; (iii) $\widehat{Q}_{n}(\alpha )\geq
cQ(\alpha )-O_{P_{Z^{\infty }}}(\delta _{n}^{2})$ uniformly over $\mathcal{A}%
_{osn}$, and $\max \{\delta _{n}^{2},Q(\Pi _{n}\alpha _{0}),\lambda
_{n},o(n^{-1})\}=\delta _{n}^{2}$; (iv) $\lambda _{n}\times \sup_{\alpha
,\alpha ^{\prime }\in \mathcal{A}_{os}}\left\vert Pen(h)-Pen(h^{\prime
})\right\vert =o(n^{-1})$ or $\lambda _{n}=o(n^{-1})$.
\end{assumption}

Assumption \ref{ass:weak_equiv}(ii) is about the local curvature of the
population criterion $Q(\alpha )$ at $\alpha _{0}$. It can be weakened to Assumption 4.1(ii) in \cite{CP_WP07}. When $\widehat{Q}%
_{n}(\alpha )$ is computed using the series LS estimator (\ref{mhat}), Lemma
C.2 of \cite{CP_WP07} shows that $\widehat{Q}_{n}(\alpha )\asymp Q(\alpha
)-O_{P_{Z^{\infty }}}(\delta _{n}^{2})$ uniformly over $\mathcal{A}_{osn}$
and hence Assumption \ref{ass:weak_equiv}(iii) is satisfied.

Recall the definition of the \textit{sieve measure of local ill-posedness}
\begin{equation}
\tau _{n}\equiv \sup_{\alpha \in \mathcal{A}_{osn}:||\alpha -\Pi _{n}\alpha
_{0}||\neq 0}\frac{||\alpha -\Pi _{n}\alpha _{0}||_{s}}{||\alpha -\Pi
_{n}\alpha _{0}||}.  \label{tau-n}
\end{equation}%
The problem of estimating $\alpha _{0}$ under $||\cdot ||_{s}$ is \textit{%
locally} \textit{ill-posed in rate} if and only if $\limsup_{n\rightarrow
\infty }\tau _{n}=\infty $. We say the problem is \textit{mildly ill-posed}
if $\tau _{n}=O([k(n)]^{a})$, and \textit{severely ill-posed} if $\tau
_{n}=O(\exp \{\frac{a}{2}k(n)\})$ for some finite $a>0$. The following
general rate result is a minor modification of Theorem 4.1 and Remark 4.1(i)
of \cite{CP_WP07}, and hence we omit its proof.

\begin{lemma}
\label{thm:ThmCONVRATEGRAL} Let $\widehat{\alpha }_{n}$ be the PSMD
estimator defined in (\ref{psmd}), and Assumptions \ref{ass:sieve}, \ref%
{A_3.6}(ii)(iii), \ref{ass:rates} and \ref{ass:weak_equiv}(i)(ii)(iii) hold.
Then:%
\begin{equation*}
||\widehat{\alpha }_{n}-\alpha _{0}||=O_{P_{Z^{\infty }}}\left( \delta
_{n}\right) \quad \text{and}\quad ||\widehat{\alpha }_{n}-\alpha
_{0}||_{s}=O_{P_{Z^{\infty }}}\left( ||\alpha _{0}-\Pi _{n}\alpha
_{0}||_{s}+\tau _{n}\delta _{n}\right) .
\end{equation*}
\end{lemma}

The above convergence rate result is applicable to any nonparametric
estimator $\widehat{m}(X,\alpha )$ of $m(X,\alpha )$ as soon as one could
compute $\delta _{n}^{2}$, the rate at which $\widehat{Q}_{n}(\alpha )$ goes
to $Q(\alpha )$. See \cite{CP_WP07} and \cite{CP_WP07a} for low level
sufficient conditions in terms of the series LS estimator (\ref{mhat}) of $%
m(X,\alpha )$.

Let $\left\{ \delta _{s,n}:n\geq 1\right\} $ be a sequence of real positive
numbers such that $\delta _{s,n}=||h_{0}-\Pi _{n}h_{0}||_{s}+\tau _{n}\delta
_{n}=o(1)$. Lemma \ref{thm:ThmCONVRATEGRAL} implies that $\widehat{\alpha }%
_{n}\in \mathcal{N}_{osn}\subseteq \mathcal{N}_{os}$ wpa1-$P_{Z^{\infty }}$,
where
\begin{eqnarray*}
\mathcal{N}_{os} &\equiv &\left\{ \alpha \in \mathcal{A}\colon \text{ }%
||\alpha -\alpha _{0}||\leq M_{n}\delta _{n},\text{ }||\alpha -\alpha
_{0}||_{s}\leq M_{n}\delta _{s,n},\text{ }\lambda _{n}Pen(h)\leq \lambda
_{n}M_{0}\right\} , \\
\mathcal{N}_{osn} &\equiv &\mathcal{N}_{os}\cap \mathcal{A}_{k(n)},\quad
\text{with }M_{n}\equiv \min \left\{ \log (\log (n+1)),\log ((\delta
_{s,n}^{-1}+1))\right\} .
\end{eqnarray*}%
We can regard $\mathcal{N}_{os}$ as the effective parameter space and $%
\mathcal{N}_{osn}$ as its sieve space in the rest of the paper. Assumption %
\ref{ass:weak_equiv}(iv) is not needed for establishing a convergence rate
in Lemma \ref{thm:ThmCONVRATEGRAL}. but, it will be imposed in the rest of
the paper so that we can ignore penalty effect in the first order local
asymptotic analysis.

\subsection{(Sieve) Riesz representation and (sieve) variance}

\label{sec:phi}

We first introduce a representation of the functional of interest $\phi ()$
at $\alpha _{0}$ that is crucial for all the subsequent local asymptotic
theories. Let $\phi :\mathbb{R}^{d_{\theta }}\times \mathcal{H}\rightarrow
\mathbb{R}$ be continuous in $||\cdot ||_{s}$. We assume that $\frac{d\phi
(\alpha _{0})}{d\alpha }[\cdot ]:\left( \mathbb{R}^{d_{\theta }}\times
\mathcal{H},||\cdot ||_{s}\right) \rightarrow \mathbb{R}$ is a $||\cdot
||_{s}-$bounded linear functional (i.e., $\left\vert \frac{d\phi (\alpha
_{0})}{d\alpha }[v]\right\vert \leq c||v||_{s}$ uniformly over $v\in \mathbb{%
R}^{d_{\theta }}\times \mathcal{H}$ for a finite positive constant $c$),
which could be computed as a pathwise\ (directional) derivative of the
functional $\phi \left( \cdot \right) $ at $\alpha _{0}$ in the direction of
$v=\alpha -\alpha _{0}\in \mathbb{R}^{d_{\theta }}\times \mathcal{H}:$

\begin{equation*}
\frac{d\phi (\alpha _{0})}{d\alpha }\left[ v\right] =\left. \frac{\partial
\phi (\alpha _{0}+\tau v)}{\partial \tau }\right\vert _{\tau =0}.
\end{equation*}

Let $\mathbf{V}$ be a linear span of $\mathcal{A}_{os}-\{\alpha _{0}\}$,
which is endowed with both $||\cdot ||_{s}$ and $||\cdot ||$ (in equation (%
\ref{fmetric})) norms, and $||v||\leq C||v||_{s}$ for all $v\in \mathbf{V}$
(under Assumption \ref{ass:weak_equiv}(i)). Let $\overline{\mathbf{V}}\equiv
clsp(\mathcal{A}_{os}-\{\alpha _{0}\})$, where $clsp(\cdot )$ is the closure
of the linear span under $||\cdot ||$. For any $v_{1},v_{2}\in \overline{%
\mathbf{V}}$, we define an inner product induced by the metric $||\cdot ||$:
\begin{equation*}
\left\langle v_{1},v_{2}\right\rangle =E\left[ \left( \frac{dm(X,\alpha _{0})%
}{d\alpha }[v_{1}]\right) ^{\prime }\Sigma (X)^{-1}\left( \frac{dm(X,\alpha
_{0})}{d\alpha }[v_{2}]\right) \right] ,
\end{equation*}%
and for any $v\in \overline{\mathbf{V}}$ we call $v=0$ if and only if $%
||v||=0$ (i.e., functions in $\overline{\mathbf{V}}$ are defined in an
equivalent class sense according to the metric $||\cdot ||$). It is clear
that $(\overline{\mathbf{V}},||\cdot ||)$ is an infinite dimensional Hilbert
space (under Assumptions \ref{ass:sieve}(i)(iii)(iv) and \ref{ass:weak_equiv}%
(i)(ii)).

If the linear functional $\frac{d\phi (\alpha _{0})}{d\alpha }[\cdot ]$ is
\textit{bounded} on $(\mathbf{V},||\cdot ||)$, i.e.
\begin{equation*}
\sup_{v\in \mathbf{V},v\neq 0}\frac{\left\vert \frac{d\phi (\alpha _{0})}{%
d\alpha }\left[ v\right] \right\vert }{\left\Vert v\right\Vert }<\infty ,
\end{equation*}%
then there is a unique extension of\ $\frac{d\phi (\alpha _{0})}{d\alpha }%
[\cdot ]$ from $(\mathbf{V},||\cdot ||)$ to $(\overline{\mathbf{V}},||\cdot
||)$, and a unique Riesz representer $v^{\ast }\in \overline{\mathbf{V}}$ of
$\frac{d\phi (\alpha _{0})}{d\alpha }[\cdot ]$ on $(\overline{\mathbf{V}}%
,||\cdot ||)$ such that\footnote{%
See, e.g., page 206-207 and theorem 3.10.1 in \cite{Debnath-Hilbert}.}%
\begin{align}\label{RRT-0}
&\frac{d\phi (\alpha _{0})}{d\alpha }\left[ v\right] =\left\langle v^{\ast
},v\right\rangle \text{ for all }v\in \overline{\mathbf{V}}\text{\quad
and\quad }\\ \notag
&\left\Vert v^{\ast }\right\Vert \equiv \sup_{v\in \overline{%
\mathbf{V}},v\neq 0}\frac{\left\vert \frac{d\phi (\alpha _{0})}{d\alpha }%
\left[ v\right] \right\vert }{\left\Vert v\right\Vert }=\sup_{v\in \mathbf{V}%
,v\neq 0}\frac{\left\vert \frac{d\phi (\alpha _{0})}{d\alpha }\left[ v\right]
\right\vert }{\left\Vert v\right\Vert }<\infty .
\end{align}%
If $\frac{d\phi (\alpha _{0})}{d\alpha }[\cdot ]$ is \textit{unbounded} on $(%
\mathbf{V},||\cdot ||)$, i.e.
\begin{equation*}
\sup_{v\in \mathbf{V},v\neq 0}\frac{\left\vert \frac{d\phi (\alpha _{0})}{%
d\alpha }\left[ v\right] \right\vert }{\left\Vert v\right\Vert }=\infty ,
\end{equation*}%
then there is no unique extension of the mapping\ $\frac{d\phi (\alpha _{0})%
}{d\alpha }[\cdot ]$ from $(\mathbf{V},||\cdot ||)$ to $(\overline{\mathbf{V}%
},||\cdot ||)$, and nor existing any Riesz representer of\ $\frac{d\phi
(\alpha _{0})}{d\alpha }[\cdot ]$ on $(\overline{\mathbf{V}},||\cdot ||)$.

Since $||v||\leq C||v||_{s}$ for all $v\in \mathbf{V}$, it is clear that a $%
||\cdot ||_{s}-$bounded linear functional $\frac{d\phi (\alpha _{0})}{%
d\alpha }[\cdot ]$ could be either bounded or unbounded on $(\mathbf{V}%
,||\cdot ||)$.

\textbf{Sieve Riesz representation}. Let $\alpha _{0,n}\in \mathbb{R}%
^{d_{\theta }}\times \mathcal{H}_{k(n)}$ be such that%
\begin{equation}
||\alpha _{0,n}-\alpha _{0}||\equiv \min_{\alpha \in \mathbb{R}^{d_{\theta
}}\times \mathcal{H}_{k(n)}}||\alpha -\alpha _{0}||.  \label{SP-1}
\end{equation}%
Let $\overline{\mathbf{V}}_{k(n)}\equiv clsp\left( \mathcal{A}%
_{osn}-\{\alpha _{0,n}\}\right) $, where $clsp\left( .\right) $ denotes the
closed linear span under $\left\Vert \cdot \right\Vert $. Then $\overline{%
\mathbf{V}}_{k(n)}$ is a finite dimensional Hilbert space under $\left\Vert
\cdot \right\Vert $. Moreover, $\overline{\mathbf{V}}_{k(n)}$ is dense in $%
\overline{\mathbf{V}}$ under $\left\Vert \cdot \right\Vert $. To simplify
the presentation, we assume that $\dim (\overline{\mathbf{V}}_{k(n)})=\dim (%
\mathcal{A}_{k(n)})\asymp k(n)$, all of which grow to infinity with $n$. By
definition we have $\left\langle v_{n},\alpha _{0,n}-\alpha
_{0}\right\rangle =0$ for all $v_{n}\in \overline{\mathbf{V}}_{k(n)}$.

Note that $\overline{\mathbf{V}}_{k(n)}$ is a finite dimensional Hilbert
space. As any linear functional on a finite dimensional Hilbert space is
bounded, we can invoke the Riesz representation theorem to deduce that there
is a $v_{n}^{\ast }\in \overline{\mathbf{V}}_{k(n)}$ such that
\begin{equation}
\frac{d\phi (\alpha _{0})}{d\alpha }[v]=\left\langle v_{n}^{\ast
},v\right\rangle,~\forall v\in \overline{\mathbf{V}}_{k(n)},~and~\left\Vert v_{n}^{\ast }\right\Vert \equiv \sup_{v\in
\overline{\mathbf{V}}_{k(n)}:\left\Vert v\right\Vert \neq 0}\frac{\left\vert
\frac{d\phi (\alpha _{0})}{d\alpha }[v]\right\vert }{\left\Vert v\right\Vert
}<\infty .  \label{RRT-1}
\end{equation}%
We call $v_{n}^{\ast }$ the \textit{sieve Riesz representer} of the
functional\ $\frac{d\phi (\alpha _{0})}{d\alpha }[\cdot ]$ on $\overline{%
\mathbf{V}}_{k(n)}$. By definition, for any non-zero linear functional $%
\frac{d\phi (\alpha _{0})}{d\alpha }[\cdot ]$, we have:%
\begin{equation*}
0<\left\Vert v_{n}^{\ast }\right\Vert ^{2}=E\left[ \left( \frac{dm(X,\alpha
_{0})}{d\alpha }[v_{n}^{\ast }]\right) ^{\prime }\Sigma (X)^{-1}\left( \frac{%
dm(X,\alpha _{0})}{d\alpha }[v_{n}^{\ast }]\right) \right]
\end{equation*}
is non-decreasing in $k(n)$.

We emphasize that the sieve Riesz representer $v_{n}^{\ast }$ of a linear
functional\ $\frac{d\phi (\alpha _{0})}{d\alpha }[\cdot ]$ on $\overline{%
\mathbf{V}}_{k(n)}$ always exists regardless of whether $\frac{d\phi (\alpha
_{0})}{d\alpha }[\cdot ]$ is bounded on the infinite dimensional space $(%
\mathbf{V},||\cdot ||)$ or not. Moreover, $v_{n}^{\ast }\in \overline{%
\mathbf{V}}_{k(n)}$ and its norm $\left\Vert v_{n}^{\ast }\right\Vert $ can
be computed in closed form (see Subsection \ref{closed_form_sieve_variance}%
). The next Lemma allows us to verify whether or not $\frac{d\phi (\alpha
_{0})}{d\alpha }[\cdot ]$ is bounded on $(\mathbf{V},||\cdot ||)$ by
checking whether or not $\lim_{k(n)\rightarrow \infty }\left\Vert
v_{n}^{\ast }\right\Vert <\infty $.

\begin{lemma}
\label{lem:sieve-Riesz-property} Let $\{\overline{\mathbf{V}}%
_{k}\}_{k=1}^{\infty }$ be an increasing sequence of finite dimensional
Hilbert spaces that is dense in $(\overline{\mathbf{V}},\left\Vert \cdot
\right\Vert )$, and $v_{n}^{\ast }\in \overline{\mathbf{V}}_{k(n)}$ be
defined in (\ref{RRT-1}). (1) If $\frac{d\phi (\alpha _{0})}{d\alpha }[\cdot
]$ is bounded on $(\mathbf{V},||\cdot ||)$, then (\ref{RRT-0}) holds, $%
v_{n}^{\ast }=\arg \min_{v\in \overline{\mathbf{V}}_{k(n)}}\left\Vert
v^{\ast }-v\right\Vert $ and $\left\Vert v^{\ast }-v_{n}^{\ast }\right\Vert
\rightarrow 0$, $\lim_{k(n)\rightarrow \infty }\left\Vert v_{n}^{\ast
}\right\Vert =\left\Vert v^{\ast }\right\Vert <\infty $; (2) Let $\frac{%
d\phi (\alpha _{0})}{d\alpha }[\cdot ]$ be bounded on $(\mathbf{V},||\cdot
||_{s})$ and $\{\overline{\mathbf{V}}_{k}\}_{k=1}^{\infty }$ be dense in $(%
\mathbf{V},\left\Vert \cdot \right\Vert _{s})$. If $\frac{d\phi (\alpha _{0})%
}{d\alpha }[\cdot ]$ is unbounded on $(\mathbf{V},||\cdot ||)$ then $%
\lim_{k(n)\rightarrow \infty }\left\Vert v_{n}^{\ast }\right\Vert =\infty $.
\end{lemma}

\textbf{Sieve score and sieve variance}. For each sieve dimension $k(n)$, we
call
\begin{equation}
S_{n,i}^{\ast }\equiv \left( \frac{dm(X_{i},\alpha _{0})}{d\alpha }%
[v_{n}^{\ast }]\right) ^{\prime }\Sigma (X_{i})^{-1}\rho (Z_{i},\alpha _{0})
\label{score}
\end{equation}%
the \textit{sieve score} associated with the $i$-th observation, and $%
\left\Vert v_{n}^{\ast }\right\Vert _{sd}^{2}\equiv Var\left( S_{n,i}^{\ast
}\right) $ as the \textit{sieve variance}. Recall that $\Sigma _{0}(X)\equiv
Var(\rho (Z;\alpha _{0})|X)$ a.s.-$X$. Then%
\begin{align} \label{svar}
\left\Vert v_{n}^{\ast }\right\Vert _{sd}^{2} &=  E[S_{n,i}^{\ast }S_{n,i}^{\ast
\prime }]\\ \notag
& =  E\left[ \left( \frac{dm(X,\alpha _{0})}{d\alpha }[v_{n}^{\ast
}]\right) ^{\prime }\Sigma (X)^{-1}\Sigma _{0}(X)\Sigma (X)^{-1}\left( \frac{%
dm(X,\alpha _{0})}{d\alpha }[v_{n}^{\ast }]\right) \right] .
\end{align}%
(See Subsection \ref{closed_form_sieve_variance} for closed form expressions
of $\left\Vert v_{n}^{\ast }\right\Vert _{sd}^{2}$.) Under Assumption \ref%
{ass:sieve}(iv), we have $\left\Vert v_{n}^{\ast }\right\Vert
_{sd}^{2}\asymp \left\Vert v_{n}^{\ast }\right\Vert ^{2}$, and hence $%
\lim_{k(n)\rightarrow \infty }\left\Vert v_{n}^{\ast }\right\Vert
_{sd}<\infty $ (or $=\infty $) iff $\lim_{k(n)\rightarrow \infty }\left\Vert
v_{n}^{\ast }\right\Vert <\infty $ (or $=\infty $). Therefore, in this paper
we call $\phi ()$ \textit{regular} (or \textit{irregular}) at $\alpha _{0}$
whenever $\lim_{k(n)\rightarrow \infty }\left\Vert v_{n}^{\ast }\right\Vert
<\infty $ (or $=\infty $), which, by Lemma \ref{lem:sieve-Riesz-property},
is also whenever $\frac{d\phi (\alpha _{0})}{d\alpha }[\cdot ]$ is bounded
(or unbounded) on $(\mathbf{V},||\cdot ||)$. It is clear that our notion of
a regular $\phi \left( \cdot \right) $ at $\alpha _{0}$ is only necessary
but not sufficient for the existence of root-$n$ asymptotically normal
regular estimators of $\phi \left( \alpha _{0}\right) $. Moreover, if $\phi
\left( \cdot \right) $ is regular at $\alpha _{0}$ then we can define
\begin{equation*}
S_{i}^{\ast }\equiv \left( \frac{dm(X_{i},\alpha _{0})}{d\alpha }[v^{\ast
}]\right) ^{\prime }\Sigma (X_{i})^{-1}\rho (Z_{i},\alpha _{0})
\end{equation*}%
as the \textit{score} associated with the $i$-th observation, and $%
\left\Vert v^{\ast }\right\Vert _{sd}^{2}\equiv Var\left( S_{i}^{\ast
}\right) $ as the \textit{asymptotic variance}. By Lemma \ref%
{lem:sieve-Riesz-property}(1) for a regular functional we have: $\left\Vert
v^{\ast }\right\Vert _{sd}^{2}\asymp \left\Vert v^{\ast }\right\Vert <\infty
$ and $Var\left( S_{i}^{\ast }-S_{n,i}^{\ast }\right) \asymp \left\Vert
v^{\ast }-v_{n}^{\ast }\right\Vert ^{2}\rightarrow 0$ as $k(n)\rightarrow
\infty $. See Appendix \ref{app:appA} for further discussions.

\subsection{Two key local conditions}

\label{sec:LAQ}

For all $k(n)$, let
\begin{equation}
u_{n}^{\ast }\equiv \frac{v_{n}^{\ast }}{\left\Vert v_{n}^{\ast }\right\Vert
_{sd}}  \label{u*}
\end{equation}%
be the \textquotedblleft scaled sieve Riesz representer\textquotedblright .
Since $\left\Vert v_{n}^{\ast }\right\Vert _{sd}^{2}\asymp \left\Vert
v_{n}^{\ast }\right\Vert ^{2}$ (under Assumption \ref{ass:sieve}(iv)), we
have: $\left\Vert u_{n}^{\ast }\right\Vert \asymp 1$ and $\left\Vert
u_{n}^{\ast }\right\Vert _{s}\leq c\tau _{n}$ for $\tau _{n}$ defined in (%
\ref{tau-n}) and a finite constant $c>0$.

Let $\mathcal{T}_{n}\equiv \{t\in \mathbb{R}\colon |t|\leq 4M_{n}^{2}\delta
_{n}\}$ with $M_{n}$ and $\delta _{n}$ given in the definition of $\mathcal{N%
}_{osn}$.

\begin{assumption}[Local behavior of $\protect\phi $]
\label{ass:phi} (i) $v\mapsto \frac{d\phi (\alpha _{0})}{d\alpha }[v]$ is a
non-zero linear functional mapping from $\mathbf{V}$ to $\mathbb{R}$; $\{%
\overline{\mathbf{V}}_{k}\}_{k=1}^{\infty }$ is an increasing sequence of
finite dimensional Hilbert spaces that is dense in $(\overline{\mathbf{V}}%
,\left\Vert \cdot \right\Vert )$; and $\frac{\left\Vert v_{n}^{\ast
}\right\Vert }{\sqrt{n}}=o(1)$;%
\begin{equation*}
\text{(ii)\quad }\sup_{(\alpha ,t)\in \mathcal{N}_{osn}\times \mathcal{T}%
_{n}}\frac{\sqrt{n}\left\vert \phi \left( \alpha +tu_{n}^{\ast }\right)
-\phi (\alpha _{0})-\frac{d\phi (\alpha _{0})}{d\alpha }[\alpha
+tu_{n}^{\ast }-\alpha _{0}]\right\vert }{\left\Vert v_{n}^{\ast
}\right\Vert }=o\left( 1\right) ;
\end{equation*}%
(iii) $\frac{\sqrt{n}\left\vert \frac{d\phi (\alpha _{0})}{d\alpha }[\alpha
_{0,n}-\alpha _{0}]\right\vert }{\left\Vert v_{n}^{\ast }\right\Vert }%
=o\left( 1\right) .$
\end{assumption}

Since $\left\Vert v_{n}^{\ast }\right\Vert _{sd}^{2}\asymp \left\Vert
v_{n}^{\ast }\right\Vert ^{2}$ (under Assumption \ref{ass:sieve}(iv)), we
could rewrite Assumption \ref{ass:phi} using $\left\Vert v_{n}^{\ast
}\right\Vert _{sd}$ instead $\left\Vert v_{n}^{\ast }\right\Vert $. As it
will become clear in Theorem \ref{thm:theta_anorm} that $\frac{\left\Vert
v_{n}^{\ast }\right\Vert _{sd}^{2}}{n}$ is the variance of $\phi (\widehat{%
\alpha }_{n})-\phi (\alpha _{0})$, Assumption \ref{ass:phi}(i) puts a
restriction on how fast the sieve dimension $k(n)$ could grow with the
sample size $n$.

Assumption \ref{ass:phi}(ii) controls the nonlinearity bias of $\phi \left(
\cdot \right) $ (i.e., the linear approximation error of a possibly
nonlinear functional $\phi \left( \cdot \right) $). It is automatically
satisfied when $\phi \left( \cdot \right) $ is a linear functional. For a
nonlinear functional $\phi \left( \cdot \right) $ (such as the quadratic
functional), it can be verified using the smoothness of $\phi \left( \cdot
\right) $ and the convergence rates in both $||\cdot ||$ and $||\cdot ||_{s}$
metrics (the definition of $\mathcal{N}_{osn}$). See Section \ref{sec-ex}
for verification.

Assumption \ref{ass:phi}(iii) controls the linear bias part due to the
finite dimensional sieve approximation of $\alpha _{0,n}$ to $\alpha _{0}$.
It is a condition imposed on the growth rate of the sieve dimension $k(n)$.
When $\phi \left( \cdot \right) $ is an irregular functional, we have $%
\left\Vert v_{n}^{\ast }\right\Vert \nearrow \infty $. Assumption \ref%
{ass:phi}(iii) requires that the sieve bias term, $\left\vert \frac{d\phi
(\alpha _{0})}{d\alpha }[\alpha _{0,n}-\alpha _{0}]\right\vert $, is of a
smaller order than that of the sieve standard deviation term, $%
n^{-1/2}\left\Vert v_{n}^{\ast }\right\Vert _{sd}$. This is a standard
condition imposed for the asymptotic normality of any plug-in nonparametric
estimator of an irregular functional (such as a point evaluation functional
of a nonparametric mean regression).

\begin{remark}
\label{remark-bias} When $\phi \left( \cdot \right) $ is regular at $\alpha
_{0}$ (i.e., $\left\Vert v_{n}^{\ast }\right\Vert \nearrow \left\Vert
v^{\ast }\right\Vert <\infty $), since $\left\langle v_{n}^{\ast },\alpha
_{0,n}-\alpha _{0}\right\rangle =0$ (by definition of $\alpha _{0,n}$) we
have $\left\vert \frac{d\phi (\alpha _{0})}{d\alpha }[\alpha _{0,n}-\alpha
_{0}]\right\vert \leq \left\Vert v^{\ast }-v_{n}^{\ast }\right\Vert \times
\left\Vert \alpha _{0,n}-\alpha _{0}\right\Vert $. And Assumption \ref%
{ass:phi}(iii) is satisfied if
\begin{equation}
||v^{\ast }-v_{n}^{\ast }||\times ||\alpha _{0,n}-\alpha _{0}||=o(n^{-1/2}).
\label{regular-orth}
\end{equation}%
This is similar to assumption 4.2 in \cite{AC_Emetrica03} and assumption
3.2(iii) in \cite{CP_WP07a} for the root-$n$ estimable Euclidean parameter $%
\theta _{0}$ of the model (\ref{semi00}). As pointed out by \cite{CP_WP07a},
Condition (\ref{regular-orth}) could be satisfied when $\dim (\mathcal{A}%
_{k(n)})\asymp k(n)$ is chosen to obtain optimal nonparametric convergence
rate in $||\cdot ||_{s}$ norm. But this nice feature only applies to regular
functionals.
\end{remark}

The next assumption is about the local quadratic approximation (LQA) to the
sample criterion difference along the scaled sieve Riesz representer
direction $u_{n}^{\ast }=v_{n}^{\ast }/\left\Vert v_{n}^{\ast }\right\Vert
_{sd}$.

For any $(\alpha ,t)\in \mathcal{N}_{osn}\times \mathcal{T}_{n}$, we let $%
\widehat{\Lambda }_{n}(\alpha (t),\alpha )\equiv 0.5\{\widehat{Q}_{n}(\alpha
(t))-\widehat{Q}_{n}(\alpha )\}$ with $\alpha (t)\equiv \alpha +tu_{n}^{\ast
}$. Denote
\begin{equation}
\mathbb{Z}_{n}\equiv n^{-1}\sum_{i=1}^{n}\left( \frac{dm(X_{i},\alpha _{0})}{%
d\alpha }[u_{n}^{\ast }]\right) ^{\prime }\Sigma (X_{i})^{-1}\rho
(Z_{i},\alpha _{0})=n^{-1}\sum_{i=1}^{n}\frac{S_{n,i}^{\ast }}{\left\Vert
v_{n}^{\ast }\right\Vert _{sd}}.  \label{effscore}
\end{equation}

\begin{assumption}[LQA]
\label{ass:LAQ} (i) $\alpha (t)\in \mathcal{A}_{k(n)}$ for any $(\alpha
,t)\in \mathcal{N}_{osn}\times \mathcal{T}_{n}$; and with $%
r_{n}(t_{n})=\left( \max \{t_{n}^{2},t_{n}n^{-1/2},o(n^{-1})\}\right) ^{-1}$%
,
\begin{equation*}
\sup_{(\alpha ,t_{n})\in \mathcal{N}_{osn}\times \mathcal{T}%
_{n}}r_{n}(t_{n})\left\vert \widehat{\Lambda }_{n}(\alpha (t_{n}),\alpha
)-t_{n}\left\{ \mathbb{Z}_{n}+\langle u_{n}^{\ast },\alpha -\alpha
_{0}\rangle \right\} -\frac{B_{n}}{2}t_{n}^{2}\right\vert =o_{P_{Z^{\infty
}}}(1),
\end{equation*}%
where, for each $n$, $B_{n}$ is a $Z^{n}$ measurable positive random
variable, and $B_{n}=O_{P_{Z^{\infty }}}(1)$;

(ii) $\sqrt{n}\mathbb{Z}_{n}\Rightarrow N(0,1)$.
\end{assumption}

Assumption \ref{ass:LAQ}(ii) is a standard one, and is implied by the
following Lindeberg condition: For all $\epsilon >0$,
\begin{equation}
\limsup_{n\rightarrow \infty }E\left[ \left( \frac{S_{n,i}^{\ast }}{%
\left\Vert v_{n}^{\ast }\right\Vert _{sd}}\right) ^{2}1\left\{ \left\vert
\frac{S_{n,i}^{\ast }}{\epsilon \sqrt{n}\left\Vert v_{n}^{\ast }\right\Vert
_{sd}}\right\vert >1\right\} \right] =0,  \label{LF}
\end{equation}%
which, under Lemma \ref{lem:sieve-Riesz-property}(1) and Assumption \ref%
{ass:sieve}(iv), is satisfied when the functional $\phi (\cdot )$ is regular
($\left\Vert v_{n}^{\ast }\right\Vert _{sd}\asymp \left\Vert v_{n}^{\ast
}\right\Vert \rightarrow \left\Vert v^{\ast }\right\Vert <\infty $). This is
why Assumption \ref{ass:LAQ}(ii) is not imposed in \cite{AC_Emetrica03} and
\cite{CP_WP07a} in their root-$n$ asymptotically normal estimation of the
regular functional $\phi (\alpha )=\lambda ^{\prime }\theta $.

Assumption \ref{ass:LAQ}(i) implicitly imposes restrictions on the
nonparametric estimator $\widehat{m}(x,\alpha )$ of $m(x,\alpha )=E[\rho
(Z,\alpha )|X=x]$ in a shrinking neighborhood of $\alpha _{0}$, so that the
criterion difference could be well approximated by a quadratic form. It is
trivially satisfied when $\widehat{m}(x,\alpha )$ is linear in $\alpha $,
such as the series LS estimator (\ref{mhat}) when $\rho (Z,\alpha )$ is
linear in $\alpha $. There are two potential difficulties in verifying this
assumption for nonlinear conditional moment models with nonparametric
endogeneity (such as the NPQIV\ model). First, due to the non-smooth
residual function $\rho (Z,\alpha )$, the estimator $\widehat{m}(x,\alpha )$
(and hence the sample criterion $\widehat{Q}_{n}(\alpha )$) could be
pointwise non-smooth with respect to $\alpha $. Second, due to the slow
convergence rates in the strong norm $||\cdot ||_{s}$ present in nonlinear
nonparametric ill-posed inverse problems, it could be challenging to control
the remainder of a quadratic approximation. When $\widehat{m}(x,\alpha )$ is
the series LS estimator (\ref{mhat}), Lemma \ref{lem:Qdiff_B} in Section \ref%
{sec:bootstrap} shows that Assumption \ref{ass:LAQ}(i) is satisfied by a set
of relatively low level sufficient conditions (Assumptions \ref{ass:m_ls} - %
\ref{ass:cont_diffm} in Appendix \ref{app:appA}). See Section \ref{sec-ex}
for verification of these sufficient conditions for functionals of the NPQIV
model.

\section{Asymptotic Properties of Sieve Wald and SQLR Statistics}

\label{sec:AsymDist}

In this section, we first establish the asymptotic normality of the plug-in
PSMD estimator $\phi (\widehat{\alpha }_{n})$ of $\phi (\alpha _{0})$ for
the model (\ref{semi00}), regardless of whether it is root-$n$ estimable or
not. We then provide a simple consistent variance estimator and hence the
asymptotic standard normality of the corresponding sieve t statistic for a
real-valued functional $\phi :\mathbb{R}^{d_{\theta }}\times \mathcal{H}%
\rightarrow \mathbb{R}$. We finally derive the asymptotic properties of SQLR
tests for the hypothesis $\phi (\alpha _{0})=\phi _{0}$. See Appendix \ref%
{app:appA} for the case of a vector-valued functional $\phi :\mathbb{R}%
^{d_{\theta }}\times \mathcal{H}\rightarrow \mathbb{R}^{d_{\phi }}$ (where $%
d_{\phi }$ could grow slowly with $n$).

\subsection{Asymptotic normality of the plug-in PSMD estimator}

\label{sec:anormality}

\renewcommand{\theassumption}{\thesection.\arabic{assumption}}

The next result allows for a (possibly) nonlinear irregular functional $\phi
()$ of the general model (\ref{semi00}).

\begin{theorem}
\label{thm:theta_anorm} Let $\widehat{\alpha }_{n}$ be the PSMD estimator (%
\ref{psmd}) and Assumptions \ref{ass:sieve} - \ref{ass:weak_equiv} hold. If
Assumptions \ref{ass:phi} and \ref{ass:LAQ} hold, then:
\begin{equation*}
\sqrt{n}\frac{\phi (\widehat{\alpha }_{n})-\phi (\alpha _{0})}{||v_{n}^{\ast
}||_{sd}}=-\sqrt{n}\mathbb{Z}_{n}+o_{P_{Z^{\infty }}}(1)\Rightarrow N(0,1).
\end{equation*}
\end{theorem}

When the functional $\phi (\cdot )$ is regular at $\alpha =\alpha _{0}$, we
have $\left\Vert v_{n}^{\ast }\right\Vert _{sd}\asymp \left\Vert v_{n}^{\ast
}\right\Vert =O(1)$ and $\phi (\widehat{\alpha }_{n})$ converges to $\phi
(\alpha _{0})$ at the parametric rate of $1/\sqrt{n}$. When the functional $%
\phi (\cdot )$ is irregular at $\alpha =\alpha _{0}$, we have $\left\Vert
v_{n}^{\ast }\right\Vert _{sd}\asymp \left\Vert v_{n}^{\ast }\right\Vert
\rightarrow \infty $; so the convergence rate of $\phi (\widehat{\alpha }%
_{n})$ becomes slower than $1/\sqrt{n}$.

For any regular functional of the semi/nonparametric model (\ref{semi00}),
Theorem \ref{thm:theta_anorm} implies that
\begin{equation*}
\sqrt{n}\left( \phi (\widehat{\alpha }_{n})-\phi (\alpha _{0})\right)
=-n^{-1/2}\sum_{i=1}^{n}S_{n,i}^{\ast }+o_{P_{Z^{\infty }}}(1)\Rightarrow
N(0,\sigma _{v^{\ast }}^{2}),~~\text{with}
\end{equation*}%
\begin{eqnarray*}
\sigma _{v^{\ast }}^{2} &= &\lim_{n\rightarrow \infty }\left\Vert v_{n}^{\ast
}\right\Vert _{sd}^{2}=\left\Vert v^{\ast }\right\Vert _{sd}^{2}\\
& = & E\left[
\left( \frac{dm(X,\alpha _{0})}{d\alpha }[v^{\ast }]\right) ^{\prime }\Sigma
(X)^{-1}\Sigma _{0}(X)\Sigma (X)^{-1}\left( \frac{dm(X,\alpha _{0})}{d\alpha
}[v^{\ast }]\right) \right] .
\end{eqnarray*}%
Thus, Theorem \ref{thm:theta_anorm} is a natural extension of the asymptotic
normality results of \cite{AC_Emetrica03} and \cite{CP_WP07a} for the
specific regular functional $\phi (\alpha _{0})=\lambda ^{\prime }\theta
_{0} $ of the model (\ref{semi00}). See Remark \ref{remark-norm} in Appendix %
\ref{app:appA} for further discussions.

\subsubsection{Closed form expressions of sieve Riesz representer and sieve
variance\label{closed_form_sieve_variance}}

To apply Theorem \ref{thm:theta_anorm}, one needs to know the sieve Riesz
representer $v_{n}^{\ast }$ defined in (\ref{RRT-1}) and the sieve variance $%
\left\Vert v_{n}^{\ast }\right\Vert _{sd}^{2}$ given in (\ref{svar}). It
turns out that both can be computed in closed form.

\begin{lemma}
\label{lem:P-sieve} Let $\overline{\mathbf{V}}_{k(n)}=\mathbb{R}^{d_{\theta
}}\times \{v_{h}(\cdot )=\psi ^{k(n)}(\cdot )^{\prime }\beta :\beta \in
\mathbb{R}^{k(n)}\}=\{v(\cdot )=\overline{\psi }^{k(n)}(\cdot )^{\prime
}\gamma :\gamma \in \mathbb{R}^{d_{\theta }+k(n)}\}$ be dense in the
infinite dimensional Hilbert space $(\overline{\mathbf{V}},\left\Vert \cdot
\right\Vert )$ with the norm $\left\Vert \cdot \right\Vert $ defined in (\ref%
{fmetric}). Then: the sieve Riesz representer $v_{n}^{\ast }=(v_{\theta
,n}^{\ast \prime },v_{h,n}^{\ast }\left( \cdot \right) )^{\prime }\in
\overline{\mathbf{V}}_{k(n)}$ of $\frac{d\phi (\alpha _{0})}{d\alpha }[\cdot
]$ has a closed form expression:
\begin{equation}
v_{n}^{\ast }=(v_{\theta ,n}^{\ast \prime },\psi ^{k(n)}(\cdot )^{\prime
}\beta _{n}^{\ast })^{\prime }=\overline{\psi }^{k(n)}(\cdot )^{\prime
}\gamma _{n}^{\ast }\text{, and }\gamma _{n}^{\ast }=D_{n}^{-}\digamma _{n}
\label{finite-riesz-s}
\end{equation}%
with $D_{n}=E\left[ \left( \frac{dm(X,\alpha _{0})}{d\alpha }[\overline{\psi
}^{k(n)}(\cdot )^{\prime }]\right) ^{\prime }\Sigma (X)^{-1}\left( \frac{%
dm(X,\alpha _{0})}{d\alpha }[\overline{\psi }^{k(n)}(\cdot )^{\prime
}]\right) \right] $ and $\digamma _{n}=\frac{d\phi (\alpha _{0})}{d\alpha }[%
\overline{\psi }^{k(n)}(\cdot )]$. Thus
\begin{equation}
\left\Vert v_{n}^{\ast }\right\Vert ^{2}=\gamma _{n}^{\ast \prime
}D_{n}\gamma _{n}^{\ast }=\digamma _{n}^{\prime }D_{n}^{-}\digamma _{n}\text{%
.}  \label{finite-riesz-s1}
\end{equation}%
The sieve variance (\ref{svar}) also has a closed form expression:%
\begin{equation}
||v_{n}^{\ast }||_{sd}^{2}=\digamma _{n}^{\prime }D_{n}^{-}\mho
_{n}D_{n}^{-}\digamma _{n},  \label{P-svar}
\end{equation}%
{\small{\begin{eqnarray*}
& &\mho _{n}\equiv\\
 & & E\left[ \left( \frac{dm(X,\alpha _{0})}{d\alpha }[\overline{%
\psi }^{k(n)}(\cdot )^{\prime }]\right) ^{\prime }\Sigma (X)^{-1}\rho
(Z,\alpha _{0})\rho (Z,\alpha _{0})^{\prime }\Sigma (X)^{-1}\left( \frac{%
dm(X,\alpha _{0})}{d\alpha }[\overline{\psi }^{k(n)}(\cdot )^{\prime
}]\right) \right] .
\end{eqnarray*}}}
\end{lemma}

Let $\mathcal{A}_{k(n)}=\Theta \times \mathcal{H}_{k(n)}$ with $\mathcal{H}%
_{k(n)}$ given in (\ref{sieve}). Then $\overline{\mathbf{V}}%
_{k(n)}=clsp\left( \mathcal{A}_{k(n)}-\{\alpha _{0,n}\}\right) $ and one
could let ${\psi }^{k(n)}(\cdot )= {q}^{k(n)}(\cdot )$ in
Lemma \ref{lem:P-sieve}, and (\ref{P-svar}) becomes the sieve variance
expression given in (\ref{P-var}).

\noindent Lemmas \ref{lem:sieve-Riesz-property} and \ref{lem:P-sieve} imply
that $\phi \left( \cdot \right) $ is \textit{regular (or irregular) at} $%
\alpha =\alpha _{0}$ \textit{iff} $\lim_{k(n)\rightarrow \infty }\left(
\digamma _{n}^{\prime }D_{n}^{-}\digamma _{n}\right) <\infty $ \textit{(or }$%
=\infty $\textit{)}.

According to Lemma \ref{lem:P-sieve} we could use different finite
dimensional linear sieve basis $\psi ^{k(n)}$ to compute sieve Riesz
representer $v_{n}^{\ast }=(v_{\theta ,n}^{\ast \prime },v_{h,n}^{\ast
}\left( \cdot \right) )^{\prime }\in \overline{\mathbf{V}}_{k(n)}$, $%
\left\Vert v_{n}^{\ast }\right\Vert ^{2}$ and $||v_{n}^{\ast }||_{sd}^{2}$.
Most typical choices include orthonormal bases and the original sieve basis $%
q^{k(n)}$ (used to approximate unknown function $h_{0}$). It is typically
easier to characterize the speed of $\left\Vert v_{n}^{\ast }\right\Vert
^{2}=\digamma _{n}^{\prime }D_{n}^{-}\digamma _{n}$ as a function of $k(n)$
when an orthonormal basis is used, while there is a nice interpretation in
terms of sieve variance estimation when the original sieve basis $q^{k(n)}$
is used. See Sections \ref{sec:NPIVex}, \ref{sec:est_avar} and \ref{sec-ex}
for related discussions.

\subsection{Consistent estimator of sieve variance of \protect$\phi (\widehat{\alpha }_{n})$}

\label{sec:est_avar}

In order to apply the asymptotic normality Theorem \ref{thm:theta_anorm}, we
need an estimator of the sieve variance $\left\Vert v_{n}^{\ast }\right\Vert
_{sd}^{2}$ defined in (\ref{svar}). We now provide one simple consistent
estimator of the sieve variance when the residual function $\rho ()$ is
pointwise smooth with respect to $\alpha _{0}$. See Appendix \ref{app:appB}
for additional consistent variance estimators.

The theoretical sieve Riesz representer $v_{n}^{\ast }$ is unknown but can
be estimated easily. Let $\left\Vert \cdot \right\Vert _{n,M}$ denote the
empirical norm induced by the following empirical inner product%
\begin{equation}
\langle v_{1},v_{2}\rangle _{n,M}\equiv \frac{1}{n}\sum_{i=1}^{n}\left(
\frac{d\widehat{m}(X_{i},\widehat{\alpha }_{n})}{d\alpha }[v_{1}]\right)
^{\prime }M_{n,i}\left( \frac{d\widehat{m}(X_{i},\widehat{\alpha }_{n})}{%
d\alpha }[v_{2}]\right) ,  \label{RRT-3}
\end{equation}%
for any $v_{1},v_{2}\in \overline{\mathbf{V}}_{k(n)}$, where $M_{n,i}$ is
some (almost surely) positive definite weighting matrix.

We define an \textit{empirical sieve Riesz representer} $\widehat{v}%
_{n}^{\ast }$ of the functional $\frac{d\phi (\widehat{\alpha }_{n})}{%
d\alpha }[\cdot ]$ with respect to the empirical norm $||\cdot ||_{n,%
\widehat{\Sigma }^{-1}}$ as
\begin{equation}
\frac{d\phi (\widehat{\alpha }_{n})}{d\alpha }[\widehat{v}_{n}^{\ast
}]=\sup_{v\in \overline{\mathbf{V}}_{k(n)},v\neq 0}\frac{|\frac{d\phi (%
\widehat{\alpha }_{n})}{d\alpha }[v]|^{2}}{||v||_{n,\widehat{\Sigma }%
^{-1}}^{2}}<\infty  \label{RRT-4}
\end{equation}%
and
\begin{equation}
\frac{d\phi (\widehat{\alpha }_{n})}{d\alpha }[v]=\langle \widehat{v}%
_{n}^{\ast },v\rangle _{n,\widehat{\Sigma }^{-1}}\text{\quad for any }v\in
\overline{\mathbf{V}}_{k(n)}.  \label{RRT-5}
\end{equation}

For $\left\Vert v_{n}^{\ast }\right\Vert _{sd}^{2}=E\left( S_{n,i}^{\ast
}S_{n,i}^{\ast \prime }\right) $ given in (\ref{svar}) we can define a
simple plug-in sieve variance estimator:%
\begin{align}\notag
||\widehat{v}_{n}^{\ast }||_{n,sd}^{2}= &\frac{1}{n}\sum_{i=1}^{n}\widehat{S}%
_{n,i}^{\ast }\widehat{S}_{n,i}^{\ast \prime }\\
= & \frac{1}{n}%
\sum_{i=1}^{n}\left( \frac{d\widehat{m}(X_{i},\widehat{\alpha }_{n})}{%
d\alpha }[\widehat{v}_{n}^{\ast }]\right) ^{\prime }\widehat{\Sigma }%
_{i}^{-1}\left( \widehat{\rho }_{i}\widehat{\rho }_{i}^{\prime }\right)
\widehat{\Sigma }_{i}^{-1}\left( \frac{d\widehat{m}(X_{i},\widehat{\alpha }%
_{n})}{d\alpha }[\widehat{v}_{n}^{\ast }]\right)  \label{svar-hat1}
\end{align}%
with $\widehat{\rho }_{i}=\rho (Z_{i},\widehat{\alpha }_{n})$ and $\widehat{%
\Sigma }_{i}=\widehat{\Sigma }(X_{i})$.

Under the condition stated in Lemma \ref{lem:P-sieve}, $\widehat{v}_{n}^{\ast }$
defined in (\ref{RRT-4}-\ref{RRT-5}) also has a closed form solution:%
\begin{equation}
\widehat{v}_{n}^{\ast }=\overline{\psi }^{k(n)}(\cdot )^{\prime }\widehat{%
\gamma }_{n}^{\ast }\text{,\quad and\quad }\widehat{\gamma }_{n}^{\ast }=%
\widehat{D}_{n}^{-}\widehat{\digamma }_{n},  \label{RRT-6}
\end{equation}%
with $\widehat{D}_{n}=\frac{1}{n}\sum_{i=1}^{n}\left( \frac{d\widehat{m}%
(X_{i},\widehat{\alpha }_{n})}{d\alpha }[\overline{\psi }^{k(n)}(\cdot
)^{\prime }]\right) ^{\prime }\widehat{\Sigma }_{i}^{-1}\left( \frac{d%
\widehat{m}(X_{i},\widehat{\alpha }_{n})}{d\alpha }[\overline{\psi }%
^{k(n)}(\cdot )^{\prime }]\right) $ and $\widehat{\digamma }_{n}=\frac{d\phi
(\widehat{\alpha }_{n})}{d\alpha }[\overline{\psi }^{k(n)}(\cdot )]$. Hence
the sieve variance estimator given in (\ref{svar-hat1}) now becomes%
\begin{equation}
||\widehat{v}_{n}^{\ast }||_{n,sd}^{2}=\widehat{V}_{1}\equiv \widehat{%
\digamma }_{n}^{\prime }\widehat{D}_{n}^{-}\widehat{\mho }_{n}\widehat{D}%
_{n}^{-}\widehat{\digamma }_{n}\text{\quad with}  \label{P-svar-hat1}
\end{equation}%
\begin{equation*}
\widehat{\mho }_{n}=\frac{1}{n}\sum_{i=1}^{n}\left( \frac{d\widehat{m}(X_{i},%
\widehat{\alpha }_{n})}{d\alpha }[\overline{\psi }^{k(n)}(\cdot )^{\prime
}]\right) ^{\prime }\widehat{\Sigma }_{i}^{-1}\left( \widehat{\rho }_{i}%
\widehat{\rho }_{i}^{\prime }\right) \widehat{\Sigma }_{i}^{-1}\left( \frac{d%
\widehat{m}(X_{i},\widehat{\alpha }_{n})}{d\alpha }[\overline{\psi }%
^{k(n)}(\cdot )^{\prime }]\right) .
\end{equation*}%
In particular, with $\psi ^{k(n)}=q^{k(n)}$ the sieve variance estimator $||%
\widehat{v}_{n}^{\ast }||_{n,sd}^{2}$ given in (\ref{P-svar-hat1}) becomes
the one given in (\ref{P-var-hat}) in Subsection \ref{sec:NPIVex}.

Let $\langle v_{1},v_{2}\rangle _{M}\equiv E\left[ \left( \frac{dm(X,\alpha
_{0})}{d\alpha }[v_{1}]\right) ^{\prime }M\left( \frac{dm(X,\alpha _{0})}{%
d\alpha }[v_{2}]\right) \right] $. Then $\langle v_{1},v_{2}\rangle _{\Sigma
^{-1}}\equiv \langle v_{1},v_{2}\rangle $ for all $v_{1},v_{2}\in \overline{%
\mathbf{V}}_{k(n)}$. Denote $\overline{\mathbf{V}}_{k(n)}^{1}\equiv \{v\in
\overline{\mathbf{V}}_{k(n)}\colon ||v||=1\}$.

\begin{assumption}
\label{ass:VE} (i) $\sup_{\alpha \in \mathcal{N}_{osn}}\sup_{v\in \overline{%
\mathbf{V}}_{k(n)}^{1}}\left\vert \frac{d\phi (\alpha )}{d\alpha }[v]-\frac{%
d\phi (\alpha _{0})}{d\alpha }[v]\right\vert =o(1)$;

\noindent (ii) for each $k(n)$ and any $\alpha \in \mathcal{N}_{osn}$, $v\in
\overline{\mathbf{V}}_{k(n)}\mapsto \frac{d\widehat{m}(\cdot ,\alpha )}{%
d\alpha }[v]\in L^{2}(f_{X})$ is a linear functional measurable with respect
to $Z^{n}$; and\\
 $\sup_{v_{1},v_{2}\in \overline{\mathbf{V}}%
_{k(n)}^{1}}\left\vert \langle v_{1},v_{2}\rangle _{n,\Sigma ^{-1}}-\langle
v_{1},v_{2}\rangle _{\Sigma ^{-1}}\right\vert =o_{P_{Z^{\infty }}}(1)$;

\noindent(iii) $\sup_{x\in \mathcal{X}}||\widehat{\Sigma }(x)-\Sigma
(x)||_{e}=o_{P_{Z^{\infty }}}(1)$;

\noindent (iv) $\sup_{x\in \mathcal{X}}E\left[ \sup_{\alpha \in \mathcal{N}%
_{osn}}||\rho (Z,\alpha )\rho (Z,\alpha )^{\prime }-\rho (Z,\alpha _{0})\rho
(Z,\alpha _{0})^{\prime }||_{e}|X=x\right] =o(1)$.

\noindent(v) $\sup_{v\in \overline{\mathbf{V}}_{k(n)}^{1}}\left\vert \langle
v,v\rangle _{n,M}-\langle v,v\rangle _{M}\right\vert =o_{P_{Z^{\infty }}}(1)$
with $M=\Sigma ^{-1}\rho (Z,\alpha _{0})\rho (Z,\alpha _{0})^{\prime }\Sigma
^{-1}$.
\end{assumption}

Assumption \ref{ass:VE}(i) becomes vacuous if $\phi $ is linear; otherwise
it requires smoothness of the family $\{\frac{d\phi (\alpha )}{d\alpha }%
[v]:\alpha \in \mathcal{N}_{osn}\}$ uniformly in $v\in \overline{\mathbf{V}}%
_{k(n)}^{1}$. Assumption \ref{ass:VE}(ii) implicitly assumes that the
residual function $\rho (z,\cdot )$ is \textquotedblleft
smooth\textquotedblright\ in $\alpha \in \mathcal{N}_{osn}$ (see, e.g., \cite%
{AC_Emetrica03}) or that $\frac{d\widehat{m}(X,\widehat{\alpha }_{n})}{%
d\alpha }[v]$ can be well approximated by numerical derivatives (see, e.g.,
\cite{HMN_WP10}). Assumption \ref{ass:VE}(iii) assumes the existence of
consistent estimators for $\Sigma $. In most applications, $\Sigma (\cdot )$
is either completely known (such as the identity matrix) or $\Sigma _{0}$;
while $\Sigma _{0}(x)$ could be consistently estimated via kernel, series
LS, local linear regression and other nonparametric procedures (see, e.g.,
\cite{AC_Emetrica03} and \cite{CP_WP07a})

\begin{theorem}
\label{thm:VE} Let Assumptions \ref{ass:sieve} - \ref{ass:weak_equiv} hold.
If Assumption \ref{ass:VE} is satisfied, then:

(1) $\left\vert \frac{||\widehat{v}_{n}^{\ast }||_{n,sd}}{||v_{n}^{\ast
}||_{sd}}-1\right\vert =o_{P_{Z^{\infty }}}(1)$ for $||\widehat{v}_{n}^{\ast
}||_{n,sd}$ given in (\ref{svar-hat1}).

(2) If, in addition, Assumptions \ref{ass:phi} and \ref{ass:LAQ} hold, then:
\begin{equation*}
\widehat{W}_{n}\equiv \sqrt{n}\frac{\phi (\widehat{\alpha }_{n})-\phi
(\alpha _{0})}{||\widehat{v}_{n}^{\ast }||_{n,sd}}=-\sqrt{n}\mathbb{Z}%
_{n}+o_{P_{Z^{\infty }}}(1)\Rightarrow N(0,1).
\end{equation*}
\end{theorem}

Theorem \ref{thm:VE}(2) allows us to construct confidence sets for $\phi
(\alpha _{0})$ based on a possibly non-optimally weighted plug-in PSMD
estimator $\phi (\widehat{\alpha }_{n})$. A potential drawback, is that it
requires a consistent estimator for $v\mapsto \frac{dm(\cdot ,\alpha _{0})}{%
d\alpha }[v]$, which may be hard to compute in practice when the residual
function $\rho (Z,\alpha )$ is not pointwise smooth in $\alpha \in \mathcal{N%
}_{osn}$ such as in the NPQIV (\ref{npqiv}) example.

\begin{remark}
\label{remark-Wald} Let $\mathcal{W}_{n}\equiv \left( \sqrt{n}\frac{\phi (%
\widehat{\alpha }_{n})-\phi _{0}}{||\widehat{v}_{n}^{\ast }||_{n,sd}}\right)
^{2}=\left( \widehat{W}_{n}+\sqrt{n}\frac{\phi (\alpha _{0})-\phi _{0}}{||%
\widehat{v}_{n}^{\ast }||_{n,sd}}\right) ^{2}$ be the Wald test statistic.
Then Theorem \ref{thm:VE} (with $\frac{||v_{n}^{\ast }||_{sd}}{\sqrt{n}}%
\asymp \frac{||v_{n}^{\ast }||}{\sqrt{n}}=o(1)$) immediately implies the
following results:

\noindent Under $H_{0}:$ $\phi (\alpha _{0})=\phi _{0}$, $\mathcal{W}%
_{n}=\left( \widehat{W}_{n}\right) ^{2}\Rightarrow \chi _{1}^{2}$.

\noindent Under $H_{1}:$ $\phi (\alpha _{0})\neq \phi _{0}$, $\mathcal{W}%
_{n}=\left( O_{P}(1)+\sqrt{n}||v_{n}^{\ast }||_{sd}^{-1}[\phi (\alpha
_{0})-\phi _{0}]\left( 1+o_{P}(1)\right) \right) ^{2}\rightarrow \infty $ in
probability.

\noindent See Theorem \ref{thm:t-localt} in Appendix \ref{app:appA} for
asymptotic properties of $\mathcal{W}_{n}$ under local alternatives.
\end{remark}

\subsection{Sieve QLR statistics\label{sec:asym_SQLR}}

We now characterize the asymptotic behaviors of the possibly \emph{%
non-optimally weighted} SQLR statistic $\widehat{QLR}_{n}(\phi _{0})$
defined in (\ref{SQLR}).

Let $\mathcal{A}_{k(n)}^{R}\equiv \{\alpha \in \mathcal{A}_{k(n)}\colon \phi
(\alpha )=\phi _{0}\}$ be the restricted sieve space, and $\widehat{\alpha }%
_{n}^{R}\in \mathcal{A}_{k(n)}^{R}$ be a restricted approximate PSMD
estimator, defined as
\begin{equation}
\widehat{Q}_{n}(\widehat{\alpha }_{n}^{R})+\lambda _{n}Pen(\widehat{h}%
_{n}^{R})\leq \inf_{\alpha \in \mathcal{A}_{k(n)}^{R}}\left\{ \widehat{Q}%
_{n}(\alpha )+\lambda _{n}Pen(h)\right\} +o_{P_{Z^{\infty }}}(n^{-1}).
\label{Rpsmd}
\end{equation}%
Then:
\begin{align*}
\widehat{QLR}_{n}(\phi _{0})= &n\left( \widehat{Q}_{n}(\widehat{\alpha }%
_{n}^{R})-\widehat{Q}_{n}(\widehat{\alpha }_{n})\right)\\
 = & n\left(
\inf_{\alpha \in \mathcal{A}_{k(n)}^{R}}\widehat{Q}_{n}(\alpha
)-\inf_{\alpha \in \mathcal{A}_{k(n)}}\widehat{Q}_{n}(\alpha )\right) +o_{P_{Z^{\infty }}}(1).
\end{align*}

Recall that $u_{n}^{\ast }\equiv v_{n}^{\ast }/\left\Vert v_{n}^{\ast
}\right\Vert _{sd}$, and that $\widehat{QLR}_{n}^{0}(\phi _{0})$ denotes the
optimally weighted (i.e., $\Sigma =\Sigma _{0}$) SQLR statistic in
Subsection \ref{sec:NPIVex}. We note that $||u_{n}^{\ast }||=1$ for the
optimally weighted case.

\begin{theorem}
\label{thm:chi2} Let Assumptions \ref{ass:sieve} - \ref{ass:LAQ} hold with $%
\left\vert B_{n}-||u_{n}^{\ast }||^{2}\right\vert =o_{P_{Z^{\infty }}}(1)$.
If $\widehat{\alpha }_{n}^{R}\in \mathcal{N}_{osn}$ wpa1-$P_{Z^{\infty }}$,
then: (1) under the null $H_{0}:\phi (\alpha _{0})=\phi _{0}$,
\begin{equation*}
||u_{n}^{\ast }||^{2}\times \widehat{QLR}_{n}(\phi _{0})=\left( \sqrt{n}%
\mathbb{Z}_{n}\right) ^{2}+o_{P_{Z^{\infty }}}(1)\Rightarrow \chi _{1}^{2}.
\end{equation*}%
(2) Further, let $\widehat{\alpha }_{n}$ be the optimally weighted PSMD
estimator (\ref{psmd}) with $\Sigma =\Sigma _{0}$. Then: under $H_{0}:$ $%
\phi (\alpha _{0})=\phi _{0}$,
\begin{equation*}
\widehat{QLR}_{n}^{0}(\phi _{0})=\left( \sqrt{n}\mathbb{Z}_{n}\right)
^{2}+o_{P_{Z^{\infty }}}(1)\Rightarrow \chi _{1}^{2}.
\end{equation*}

\noindent \textit{See Theorem \ref{thm:chi2_localt} in Appendix \ref%
{app:appA} for the asymptotic behavior under local alternatives}.
\end{theorem}

Compared to Theorem \ref{thm:theta_anorm} on the asymptotic normality of $%
\phi (\widehat{\alpha }_{n})$, Theorem \ref{thm:chi2} on the asymptotic null
distribution of the SQLR statistic requires two extra conditions: $%
\left\vert B_{n}-||u_{n}^{\ast }||^{2}\right\vert =o_{P_{Z^{\infty }}}(1)$
and $\widehat{\alpha }_{n}^{R}\in \mathcal{N}_{osn}$ wpa1-$P_{Z^{\infty }}$.
Both conditions are also needed even for QLR statistics in parametric
extremum estimation and testing problems. Lemma \ref{lem:Qdiff_B} in Section %
\ref{sec:bootstrap} provides a simple sufficient condition (Assumption \ref%
{ass:LLN_triangular}) for $\left\vert B_{n}-||u_{n}^{\ast }||^{2}\right\vert
=o_{P_{Z^{\infty }}}(1)$. Proposition \ref{pro:conv-rate-RPSMDE} in Appendix %
\ref{app:appB} establishes $\widehat{\alpha }_{n}^{R}\in \mathcal{N}_{osn}$
wpa1-$P_{Z^{\infty }}$ under the null $H_{0}:\phi (\alpha _{0})=\phi _{0}$
and other conditions virtually the same as those for Lemma \ref%
{thm:ThmCONVRATEGRAL} (i.e., $\widehat{\alpha }_{n}\in \mathcal{N}_{osn}$
wpa1-$P_{Z^{\infty }}$).

Theorem \ref{thm:chi2}(2) recommends to construct an asymptotic $100(1-\tau
)\%$ confidence set for $\phi (\alpha )$ by inverting the optimally weighted
SQLR statistic: $\left\{ r\in \mathbb{R}\colon \text{ }\widehat{QLR}%
_{n}^{0}(r)\leq c_{\chi _{1}^{2}}(1-\tau )\right\} $. This result extends
that of \cite{CP_WP07a} for a regular Euclidean functional $\phi (\alpha
)=\lambda ^{\prime }\theta $ to possibly irregular nonlinear functionals.

Next, we consider the asymptotic behavior of $\widehat{QLR}_{n}(\phi _{0})$
under the fixed alternatives $H_{1}:$ $\phi (\alpha _{0})\neq \phi _{0}$.

\begin{theorem}
\label{thm:QLR-H1} Let Assumptions \ref{ass:sieve}, \ref{A_3.6} and \ref%
{ass:rates} hold. Suppose that $\sup_{h\in \mathcal{H}}Pen(h)<\infty $ and $%
\phi $ is continuous in $||\cdot ||_{s}$. Then: under $H_{1}:$ $\phi (\alpha
_{0})\neq \phi _{0}$, there is a constant $C>0$ such that
\begin{equation*}
\frac{\widehat{QLR}_{n}(\phi _{0})}{n}\geq C>0\quad \text{wpa1.}
\end{equation*}
\end{theorem}

\section{Inference Based on Generalized Residual Bootstrap}

\label{sec:bootstrap}

\renewcommand{\theassumption}{Boot.\arabic{assumption}} %
\setcounter{assumption}{0}

The inference procedures described in Subsections \ref{sec:est_avar} and \ref%
{sec:asym_SQLR} are based on the asymptotic critical values. For many
parametric models it is known that bootstrap based procedures could
approximate finite sample distributions more accurately. In this section we
establish the consistency of the bootstrap sieve Wald and SQLR statistics
under virtually the same conditions as those imposed for the original-sample
sieve Wald and SQLR statistics.

A bootstrap procedure is described by an array of \textquotedblleft
weights\textquotedblright\ $\left\{ \omega _{i,n}\right\} _{i=1}^{n}$ for
each $n$, where each bootstrap sample is drawn independently of the original
data $\left\{ Z_{i}\right\} _{i=1}^{n}$. Different bootstrap procedures
correspond to different choices of the weights $\left\{ \omega
_{i,n}\right\} _{i=1}^{n}$ but all satisfy $\omega _{i,n}\geq 0$ and $%
E[\omega _{i,n}]=1$. For the time being we assume that $\lim_{n\rightarrow
\infty }Var(\omega _{i,n})=\sigma _{\omega }^{2}\in (0,\infty )$ for all $i$.

In this paper we focus on two types of bootstrap weights:

\begin{assumption}[I.i.d Weights]
\label{ass:Wboot} Let $(\omega _{i})_{i=1}^{n}$ be a sequence such that $%
\omega _{i}\in \mathbb{R}_{+}$, $\omega _{i}\sim iidP_{\omega }$, $E[\omega
]=1$, $Var(\omega )=\sigma _{\omega }^{2}$, and $\int_{0}^{\infty }\sqrt{%
P(|\omega -1|\geq t)}dt<\infty $.
\end{assumption}

The condition $\int_{0}^{\infty }\sqrt{P(|\omega -1|\geq t)}dt<\infty $ is
implied by $E[|\omega -1|^{2+\epsilon }]<\infty $ for some $\epsilon >0$.

\begin{assumption}[Multinomial Weights]
\label{ass:Wboot_e} Let $(\omega _{i,n})_{i=1}^{n}$ be a triangular array of
random variables such that $(\omega _{1,n},...,\omega _{n,n})\sim
Multinomial(n;n^{-1},...,n^{-1})$.
\end{assumption}

We sometimes omit the $n$ subscript from the weight series. Note that under
Assumption \ref{ass:Wboot_e}, $E[\omega _{1}]=1$, $Var(\omega
_{1})=(1-1/n)\rightarrow 1\equiv \sigma _{\omega }^{2}$ and $Cov(\omega
_{i},\omega _{j})=-n^{-1}$ (for $i\neq j$). Finally, $n^{-1}\max_{1\leq
i\leq n}(\omega _{i}-1)^{2}=o_{P_{\omega }}(1)$. We use these facts in the
proofs.

Let $V_{i}\equiv (Z_{i},\omega _{i,n})$ and
\begin{equation*}
\rho ^{B}(V_{i},\alpha )\equiv \omega _{i,n}\rho (Z_{i},\alpha ),
\end{equation*}%
be the bootstrap residual function. Let $\widehat{m}^{B}(x,\alpha )$ be a
bootstrap version of $\widehat{m}(x,\alpha )$, that is, $\widehat{m}%
^{B}(x,\alpha )$ is computed in the same way as that of $\widehat{m}%
(x,\alpha )$ except that we use $\rho ^{B}(V_{i},\alpha )$ instead of $\rho
(Z_{i},\alpha )$. In particular, $\widehat{m}^{B}(x,\alpha
)=\sum_{i=1}^{n}\omega _{i,n}\rho (Z_{i},\alpha )A_{n}(X_{i},x)$ for any
linear estimator $\widehat{m}(x,\alpha )$ (\ref{mhat-linear}) of $m(x,\alpha
)$. For example, if $\widehat{m}(x,\alpha )$ is a series LS estimator (\ref%
{mhat}), then $\widehat{m}^{B}(x,\alpha )$ is the bootstrap series LS
estimator (\ref{mhat_B}) defined in Subsection \ref{sec:NPIVex}.

Let $\widehat{Q}_{n}^{B}(\alpha )\equiv \frac{1}{n}\sum_{i=1}^{n}\widehat{m}%
^{B}(X_{i},\alpha )^{\prime }\widehat{\Sigma }(X_{i})^{-1}\widehat{m}%
^{B}(X_{i},\alpha )$ be a bootstrap version of $\widehat{Q}_{n}(\alpha )$,
and $\widehat{\alpha }_{n}^{B}$ be the bootstrap PSMD estimator, i.e., $%
\widehat{\alpha }_{n}^{B}$ is an approximate minimizer of $\left\{ \widehat{Q%
}_{n}^{B}(\alpha )+\lambda _{n}Pen(h)\right\} $ on $\mathcal{A}_{k(n)}$.
Denote $\widehat{\phi }_{n}\equiv \phi (\widehat{\alpha }_{n})$. Then
\begin{equation*}
\widehat{QLR}_{n}^{B}(\widehat{\phi }_{n})=n\left( \inf_{\{\mathcal{A}%
_{k(n)}\colon \phi (\alpha )=\widehat{\phi }_{n}\}}\widehat{Q}%
_{n}^{B}(\alpha )-\widehat{Q}_{n}^{B}(\widehat{\alpha }_{n}^{B})\right)
\end{equation*}%
is the (generalized residual) bootstrap SQLR test statistic. And $\mathcal{W}%
_{1,n}^{B}\equiv \left( \sqrt{n}\frac{\phi (\widehat{\alpha }_{n}^{B})-%
\widehat{\phi }_{n}}{\sigma _{\omega }||\widehat{v}_{n}^{\ast }||_{n,sd}}%
\right) ^{2}$ is one simple bootstrap Wald test statistic (see Subsection %
\ref{sub-boots-t} for another simple bootstrap Wald statistic).

\textbf{Additional notation}. To be more precise, we introduce some
definitions associated with the new random variables $V_{i}\equiv
(Z_{i},\omega _{i,n})$ and the enlarged probability spaces. Let $\Omega
=\{\omega _{i,n}\colon i=1,...,n;~n=1,...\}$ be the space of weights,
defined as a triangle array with elements in $\mathbb{R}$, the corresponding
$\sigma $-algebra and probability are $(\mathcal{B}_{\Omega },P_{\Omega })$.
Let $\mathcal{V}^{\infty }\equiv \mathcal{Z}^{\infty }\times \Omega $, $%
\mathcal{B}^{\infty }\equiv \mathcal{B}_{Z}^{\infty }\times \mathcal{B}%
_{\Omega }$ be the $\sigma $-algebra, and $P_{V^{\infty }}$ be the joint
probability over $\mathcal{V}^{\infty }$. Finally, for each $n$, let $%
\mathcal{B}^{n}$ be the $\sigma $-algebra generated by $V^{n}\equiv
Z^{n}\times (\omega _{1,n},...,\omega _{n,n})$, where each $\omega _{i,n}$
acts as a \textquotedblleft weight\textquotedblright\ of $Z_{i}$. Let $A_{n}$
be a random variable that is measurable with respect to $\mathcal{B}^{n}$,
and $\mathcal{L}_{V^{\infty }|Z^{\infty }}(A_{n}|Z^{n})$ (or $P_{V^{\infty
}|Z^{\infty }}\left( A_{n}\leq \cdot \mid Z^{n}\right) $) be the conditional
law (or conditional distribution) of $A_{n}$ given $Z^{n}$. Let $B_{n}$ be a
random variable measurable with respect to $\mathcal{B}_{Z}^{\infty }$, and $%
\mathcal{L}(B_{n})$ (or $P_{Z^{\infty }}\left( B_{n}\leq \cdot \right) $) be
the law (or distribution) of $B_{n}$. For two real valued random variables, $%
A_{n}$ (measurable with respect to $\mathcal{B}^{n}$) and $B$ (measurable
with respect to some $\sigma $-algebra $\mathcal{B}_{B}$), we say $%
\left\vert \mathcal{L}_{V^{\infty }|Z^{\infty }}(A_{n}|Z^{n})-\mathcal{L}%
(B)\right\vert =o_{P_{Z^{\infty }}}(1)$ if for any $\delta >0$, there exists
a $N(\delta )$ such that
\begin{equation*}
P_{Z^{\infty }}\left( \sup_{f\in BL_{1}}\left\vert
E[f(A_{n})|Z^{n}]-E[f(B)]\right\vert \leq \delta \right) \geq 1-\delta \text{%
\quad for all }n\geq N(\delta ),
\end{equation*}%
(i.e., $\sup_{f\in BL_{1}}\left\vert E[f(A_{n})|Z^{n}]-E[f(B)]\right\vert
=o_{P_{Z^{\infty }}}(1)$), where $BL_{1}$ denotes the class of uniformly
bounded Lipschitz functions $f:\mathbb{R}\rightarrow \mathbb{R}$ such that $%
||f||_{L^{\infty }}\leq 1$ and $|f(z)-f(z^{\prime })|\leq |z-z^{\prime }|$.
See chapter 1.12 of \cite{VdV-W_book96} (henceforth, VdV-W) for more details.

We say $\Delta _{n}$ is of order $o_{P_{V^{\infty }|Z^{\infty }}}(1)$ in $%
P_{Z^{\infty }}$ probability, and denote it as $\Delta _{n}=o_{P_{V^{\infty
}|Z^{\infty }}}(1)~wpa1(P_{Z^{\infty }})$, if for any $\epsilon >0$, $%
P_{Z^{\infty }}\left( P_{V^{\infty }|Z^{\infty }}\left( |\Delta
_{n}|>\epsilon \mid Z^{n}\right) >\epsilon \right) \rightarrow 0$ as $%
n\rightarrow \infty $.

We say $\Delta _{n}$ is of order $O_{P_{V^{\infty }|Z^{\infty }}}(1)$ in $%
P_{Z^{\infty }}$ probability, and denote it as $\Delta _{n}=O_{P_{V^{\infty
}|Z^{\infty }}}(1)~wpa1(P_{Z^{\infty }})$, if for any $\epsilon >0$ there
exists a $M\in (0,\infty )$, such that $P_{Z^{\infty }}\left( P_{V^{\infty
}|Z^{\infty }}\left( |\Delta _{n}|>M\mid Z^{n}\right) >\epsilon \right)
\rightarrow 0$ as $n\rightarrow \infty $.

\subsection{Bootstrap local quadratic approximation (LQA$^{B}$)}

Lemma \ref{lem:cons_boot} in Appendix \ref{app:appA} shows that the
bootstrap PSMD estimator $\widehat{\alpha }_{n}^{B}\in \mathcal{N}_{osn}$
wpa1 under Assumptions \ref{ass:rates_B} and \ref{ass:sieve} - \ref%
{ass:weak_equiv}. This allows us to introduce a condition that is a
bootstrap version of the LQA Assumption \ref{ass:LAQ}. For any $\alpha \in
\mathcal{N}_{osn}$, we let $\widehat{\Lambda }_{n}^{B}(\alpha (t_{n}),\alpha
)\equiv 0.5\{\widehat{Q}_{n}^{B}(\alpha (t_{n}))-\widehat{Q}_{n}^{B}(\alpha
)\}$ with $\alpha (t_{n})\equiv \alpha +t_{n}u_{n}^{\ast }$ for $t_{n}\in
\mathcal{T}_{n}$. For any sequence of non-negative weights $(b_{i})_{i}$,
let
\begin{equation*}
\mathbb{Z}_{n}^{b}\equiv n^{-1}\sum_{i=1}^{n}b_{i}\left( \frac{%
dm(X_{i},\alpha _{0})}{d\alpha }[u_{n}^{\ast }]\right) ^{\prime }\Sigma
(X_{i})^{-1}\rho (Z_{i},\alpha _{0})=n^{-1}\sum_{i=1}^{n}b_{i}\frac{%
S_{n,i}^{\ast }}{\left\Vert v_{n}^{\ast }\right\Vert _{sd}}.
\end{equation*}

\begin{assumption}[LQA$^{B}$]
\label{ass:LAQ_B} (i) $\alpha (t)\in \mathcal{A}_{k(n)}$ for any $(\alpha
,t)\in \mathcal{N}_{osn}\times \mathcal{T}_{n}$, and with $%
r_{n}(t_{n})=\left( \max \{t_{n}^{2},t_{n}n^{-1/2},o(n^{-1})\}\right) ^{-1}$%
,
\begin{eqnarray*}
& & \sup_{(\alpha ,t_{n})\in \mathcal{N}_{osn}\times \mathcal{T}%
_{n}}r_{n}(t_{n})\left\vert \widehat{\Lambda }_{n}^{B}(\alpha (t_{n}),\alpha
)-t_{n}\left\{ \mathbb{Z}_{n}^{\omega }+\langle u_{n}^{\ast },\alpha -\alpha
_{0}\rangle \right\} -\frac{B_{n}^{\omega }}{2}t_{n}^{2}\right\vert\\
& & =  o_{P_{V^{\infty }|Z^{\infty }}}(1)~wpa1(P_{Z^{\infty }})
\end{eqnarray*}%
where $B_{n}^{\omega }$ is a $V^{n}$ measurable positive random variable
such that $B_{n}^{\omega }=O_{P_{V^{\infty }|Z^{\infty
}}}(1)~wpa1(P_{Z^{\infty }})$;
\begin{equation*}
\text{(ii)\quad }\left\vert \mathcal{L}_{V^{\infty }|Z^{\infty }}\left(
\sqrt{n}\frac{\mathbb{Z}_{n}^{\omega -1}}{\sigma _{\omega }}\mid
Z^{n}\right) -\mathcal{L}\left( \mathbb{Z}\right) \right\vert
=o_{P_{Z^{\infty }}}(1),
\end{equation*}%
where $\mathbb{Z}$ is a standard normal random variable.
\end{assumption}

Assumption \ref{ass:LAQ_B}(i) implicitly imposes restrictions on the
bootstrap estimator $\widehat{m}^{B}(x,\alpha )$ of the conditional mean
function $m(x,\alpha )$. Below we provide low level sufficient conditions
for Assumption \ref{ass:LAQ_B}(i) when $\widehat{m}^{B}(x,\alpha )$ is a
bootstrap series LS estimator.

Let $g(X,u_{n}^{\ast })\equiv \{\frac{dm(X,\alpha _{0})}{d\alpha }%
[u_{n}^{\ast }]\}^{\prime }\Sigma (X)^{-1}$. Then $E\left[
g(X_{i},u_{n}^{\ast })\Sigma (X_{i})g(X_{i},u_{n}^{\ast })^{\prime }\right]
=||u_{n}^{\ast }||^{2}$.

\renewcommand{\theassumption}{B} \setcounter{assumption}{0}

\begin{assumption} \label{ass:LLN_triangular} For $\Gamma (\cdot )\in \{\Sigma (\cdot ),\Sigma
_{0}(\cdot )\}$,
\begin{equation*}
\left\vert n^{-1}\sum_{i=1}^{n}g(X_{i},u_{n}^{\ast })\Gamma
(X_{i})g(X_{i},u_{n}^{\ast })^{\prime }-E\left[ g(X_{i},u_{n}^{\ast })\Gamma
(X_{i})g(X_{i},u_{n}^{\ast })^{\prime }\right] \right\vert =o_{P_{Z^{\infty
}}}(1).
\end{equation*}
\end{assumption}

\begin{lemma}
\label{lem:Qdiff_B} Let Assumptions \ref{ass:sieve} - \ref{ass:weak_equiv}
and \ref{ass:m_ls} - \ref{ass:cont_diffm} hold.

(1) Let $\widehat{m}$ be the series LS estimator (\ref{mhat}). Then
Assumption \ref{ass:LAQ}(i) is satisfied. Further, if Assumption \ref%
{ass:LLN_triangular} holds then $\left\vert B_{n}-||u_{n}^{\ast
}||^{2}\right\vert =o_{P_{Z^{\infty }}}(1)$.

(2) Let $\widehat{m}^{B}(\cdot ,\alpha )$ be the bootstrap series LS
estimator (\ref{mhat_B}), Assumption \ref{ass:rates_B}, and either
Assumption \ref{ass:Wboot} or \ref{ass:Wboot_e} hold. Then Assumption \ref%
{ass:LAQ_B}(i) holds with $B_{n}^{\omega }=B_{n}$. Further, if Assumption %
\ref{ass:LLN_triangular} holds then $\left\vert B_{n}^{\omega
}-||u_{n}^{\ast }||^{2}\right\vert =o_{P_{V^{\infty }|Z^{\infty
}}}(1)~wpa1(P_{Z^{\infty }})$.
\end{lemma}

Lemma \ref{lem:Qdiff_B} indicates that the low level Assumptions \ref%
{ass:m_ls} - \ref{ass:cont_diffm} are sufficient for both the
original-sample LQA Assumption \ref{ass:LAQ}(i) and the bootstrap LQA
Assumption \ref{ass:LAQ_B}(i).

Assumption \ref{ass:LAQ_B}(ii) can be easily verified by applying some
central limit theorems. For example, if the weights are independent
(Assumption \ref{ass:Wboot}), we can use Lindeberg-Feller CLT; if the
weights are multinomial (Assumption \ref{ass:Wboot_e}) we can apply Hayek
CLT (see \cite{VdV-W_book96} p. 458 ). The next lemma provides some simple
sufficient conditions for Assumption \ref{ass:LAQ_B}(ii).

\begin{lemma}
\label{lem:clt_B} Let either Assumption \ref{ass:Wboot} or Assumption \ref%
{ass:Wboot_e} hold. If there is a positive real sequence $(b_{n})_{n}$ such
that $b_{n}=o\left( \sqrt{n}\right) $ and
\begin{equation}
\limsup_{n\rightarrow \infty }E\left[ \left( g(X,u_{n}^{\ast })\rho
(Z,\alpha _{0})\right) ^{2}1\left\{ \frac{(g(X,u_{n}^{\ast })\rho (Z,\alpha
_{0}))^{2}}{b_{n}}>1\right\} \right] =0,  \label{CLT_triangular}
\end{equation}%
then Assumptions \ref{ass:LAQ_B}(ii) and \ref{ass:LAQ}(ii) hold.
\end{lemma}

\subsection{Bootstrap sieve Student t statistic\label{sub-boots-t}}

\renewcommand{\theassumption}{Boot.\arabic{assumption}} %
\setcounter{assumption}{3}

Lemma \ref{lem:cons_boot} shows that $\widehat{\alpha }_{n}^{B}\in \mathcal{N%
}_{osn}$ wpa1 under virtually the same conditions as those for the
original-sample estimator $\widehat{\alpha }_{n}\in \mathcal{N}_{osn}$ wpa1.
This would easily lead to the consistency of the simplest bootstrap sieve t
statistic $\widehat{W}_{1,n}^{B}\equiv \sqrt{n}\frac{\phi (\widehat{\alpha }%
_{n}^{B})-\phi (\widehat{\alpha }_{n})}{\sigma _{\omega }||\widehat{v}%
_{n}^{\ast }||_{n,sd}}$.

We now establish the consistency of another bootstrap sieve t statistic $%
\widehat{W}_{2,n}^{B}\equiv \sqrt{n}\frac{\phi (\widehat{\alpha }%
_{n}^{B})-\phi (\widehat{\alpha }_{n})}{||\widehat{v}_{n}^{\ast }||_{B,sd}}$%
, where $||\widehat{v}_{n}^{\ast }||_{B,sd}^{2}$ is a bootstrap sieve
variance estimator:
\begin{equation}
||\widehat{v}_{n}^{\ast }||_{B,sd}^{2}\equiv \frac{1}{n}\sum_{i=1}^{n}\left(
\frac{d\widehat{m}(X_{i},\widehat{\alpha }_{n})}{d\alpha }[\widehat{v}%
_{n}^{\ast }]\right) ^{\prime }\widehat{\Sigma }_{i}^{-1}\varrho (V_{i},%
\widehat{\alpha }_{n})\varrho (V_{i},\widehat{\alpha }_{n})^{\prime }%
\widehat{\Sigma }_{i}^{-1}\left( \frac{d\widehat{m}(X_{i},\widehat{\alpha }%
_{n})}{d\alpha }[\widehat{v}_{n}^{\ast }]\right)  \label{svar-hat1-boot}
\end{equation}%
with $\varrho (V_{i},\alpha )\equiv (\omega _{i,n}-1)\rho (Z_{i},\alpha
)\equiv \rho ^{B}(V_{i},\alpha )-\rho (Z_{i},\alpha )$ for any $\alpha $.

We note that $||\widehat{v}_{n}^{\ast }||_{B,sd}^{2}$ is an analog to $||%
\widehat{v}_{n}^{\ast }||_{n,sd}^{2}$ defined in (\ref{svar-hat1}) but using
the bootstrapped generalized residual $\varrho (V_{i},\widehat{\alpha }_{n})$
instead of the original sample fitted residual $\rho (Z_{i},\widehat{\alpha }%
_{n})$. It also has a closed form expression: $||\widehat{v}_{n}^{\ast
}||_{B,sd}^{2}=\widehat{\digamma }_{n}^{\prime }\widehat{D}_{n}^{-}\widehat{%
\mho }_{n}^{B}\widehat{D}_{n}^{-}\widehat{\digamma }_{n}$ with
\begin{equation*}
\widehat{\mho }_{n}^{B}=\frac{1}{n}\sum_{i=1}^{n}\left( \frac{d\widehat{m}%
(X_{i},\widehat{\alpha }_{n})}{d\alpha }[\overline{\psi }^{k(n)}(\cdot
)^{\prime }]\right) ^{\prime }(\omega
_{i,n}-1)^{2} \widehat{M}_i \left( \frac{d\widehat{m}(X_{i},%
\widehat{\alpha }_{n})}{d\alpha }[\overline{\psi }^{k(n)}(\cdot )^{\prime
}]\right)
\end{equation*}%
where $\widehat{M}_i \equiv \widehat{\Sigma }_{i}^{-1}\rho (Z_{i},\widehat{\alpha }_{n})\rho (Z_{i},\widehat{\alpha }%
_{n})^{\prime }\widehat{\Sigma }_{i}^{-1}$. That is, $||\widehat{v}_{n}^{\ast }||_{B,sd}^{2}$ is computed in the same
way as $||\widehat{v}_{n}^{\ast }||_{n,sd}^{2}=\widehat{\digamma }%
_{n}^{\prime }\widehat{D}_{n}^{-}\widehat{\mho }_{n}\widehat{D}_{n}^{-}%
\widehat{\digamma }_{n}$ given in (\ref{P-svar-hat1}) except using $\widehat{%
\mho }_{n}^{B}$ instead of $\widehat{\mho }_{n}$.

\begin{assumption}
\label{ass:VE_boot2} $\sup_{v\in \overline{\mathbf{V}}_{k(n)}^{1}}|\langle
v,v\rangle _{n,\widehat{M}^{B}}-\sigma _{\omega }^{2}\langle v,v\rangle _{n,\widehat{%
M}}|=o_{P_{V^{\infty }|Z^{\infty }}}(1)~wpa1(P_{Z^{\infty }})$ with $\widehat{M}%
_{i}^{B}=(\omega _{i,n}-1)^{2}\widehat{M}_{i}$.
\end{assumption}

This assumption can be verified given Assumptions \ref{ass:Wboot} or \ref%
{ass:Wboot_e}. The following result is a bootstrap version of Theorem \ref%
{thm:VE}(1).

\begin{theorem}
\label{thm:VE-boot2} Let Assumptions \ref{ass:sieve} - \ref{ass:weak_equiv}, %
\ref{ass:VE} and \ref{ass:VE_boot2} hold. Then:%
\begin{equation*}
\left\vert \frac{||\widehat{v}_{n}^{\ast }||_{B,sd}}{\sigma _{\omega
}||v_{n}^{\ast }||_{sd}}-1\right\vert =o_{P_{V^{\infty }|Z^{\infty
}}}(1)~wpa1(P_{Z^{\infty }}).
\end{equation*}
\end{theorem}

Recall that $\widehat{W}_{n}\equiv \sqrt{n}\frac{\phi (\widehat{\alpha }%
_{n})-\phi (\alpha _{0})}{||\widehat{v}_{n}^{\ast }||_{n,sd}}$, whose
probability distribution $P_{Z^{\infty }}\left( \widehat{W}_{n}\leq \cdot
\right) $ converges to the standard normal cdf $\Phi (\cdot )$. The next
result is about the consistency of the bootstrap sieve t statistic $\widehat{%
W}_{2,n}^{B}$.

\begin{theorem}
\label{thm:bootstrap_2} Let $\widehat{\alpha }_{n}$ be the PSMD estimator (%
\ref{psmd}) and $\widehat{\alpha }_{n}^{B}$ the bootstrap PSMD estimator.
Let Assumptions \ref{ass:sieve} - \ref{ass:weak_equiv} and \ref{ass:rates_B}
hold. Let Assumptions \ref{ass:phi}, \ref{ass:LAQ} and \ref{ass:LAQ_B} hold.

(1) Let Assumptions \ref{ass:VE} and \ref{ass:VE_boot2} hold. Then:%
\begin{equation*}
\sup_{t\in \mathbb{R}}\left\vert P_{V^{\infty }|Z^{\infty }}\left( \widehat{W%
}_{2,n}^{B}\leq t\mid Z^{n}\right) -P_{Z^{\infty }}\left( \widehat{W}%
_{n}\leq t\right) \right\vert =o_{P_{V^{\infty }|Z^{\infty
}}}(1)~wpa1(P_{Z^{\infty }}).
\end{equation*}%
(2) If $\phi ()$ is regular at $\alpha _{0}$, without imposing Assumptions %
\ref{ass:VE} and \ref{ass:VE_boot2}, we have:%
\begin{eqnarray*}
& & \sup_{t\in \mathbb{R}}\left\vert P_{V^{\infty }|Z^{\infty }}\left( \sqrt{n}%
\frac{\phi (\widehat{\alpha }_{n}^{B})-\phi (\widehat{\alpha }_{n})}{\sigma
_{\omega }}\leq t\mid Z^{n}\right) -P_{Z^{\infty }}\left( \sqrt{n}\left(
\phi (\widehat{\alpha }_{n})-\phi (\alpha _{0})\right) \leq t\right)
\right\vert \\
&& =  o_{P_{V^{\infty }|Z^{\infty }}}(1)~wpa1(P_{Z^{\infty }}).
\end{eqnarray*}
\end{theorem}

For a regular functional, Theorem \ref{thm:bootstrap_2}(2) provides one way
to construct its confidence sets without the need to compute any variance
estimator. This extends the result in \cite{CP_WP07a} for a regular
Euclidean parameter $\lambda ^{\prime }\theta $ to a general regular
functional $\phi (\alpha )$. Unfortunately for an irregular functional, we
need to compute a consistent bootstrap sieve variance estimator $||\widehat{v%
}_{n}^{\ast }||_{B,sd}^{2}$ to apply Theorem \ref{thm:bootstrap_2}(1).
Luckily $||\widehat{v}_{n}^{\ast }||_{B,sd}^{2}$ is easy to compute when the
residual function $\rho (Z_{i},\alpha )$ is pointwise smooth in $\alpha _{0}$%
. Moreover, since $E\left( ||\widehat{v}_{n}^{\ast }||_{B,sd}^{2}\mid
Z^{n}\right) =\sigma _{\omega }^{2}||\widehat{v}_{n}^{\ast }||_{n,sd}^{2}$
we suspect that the bootstrap sieve t statistic $\widehat{W}_{2,n}^{B}$
might have second order refinement property by choices of bootstrap weights $%
\{\omega _{i,n}\}$. This will be a subject of future research.

The bootstrap sieve t statistic $\widehat{W}_{2,n}^{B}$ requires to compute
the original sample PSMD estimator $\widehat{\alpha }_{n}$ and the bootstrap
PSMD estimator $\widehat{\alpha }_{n}^{B}$. In the online Appendix \ref{app:appD} we present a sieve score test and its bootstrap version,
which only use the original sample restricted PSMD estimator $\widehat{%
\alpha }_{n}^{R}$ and do not use $\widehat{\alpha }_{n}^{B}$, and hence are
computationally simple.

\begin{remark}
\label{remark-Wald-B} Theorems \ref{thm:VE}(2) and \ref{thm:bootstrap_2}(1)
imply that the bootstrap Wald test statistic $\mathcal{W}_{2,n}^{B}\equiv
\left( \widehat{W}_{2,n}^{B}\right) ^{2}$ always has the same limiting
distribution $\chi _{1}^{2}$ (conditional on the data) under the null and
the alternatives. Let $\widehat{c}_{2,n}(a)$ be the $a-th$ quantile of the
distribution of $\mathcal{W}_{2,n}^{B}$ (conditional on the data $%
\{Z_{i}\}_{i=1}^{n}$). Let $\mathcal{W}_{n}\equiv \left( \sqrt{n}\frac{\phi (%
\widehat{\alpha }_{n})-\phi _{0}}{||\widehat{v}_{n}^{\ast }||_{n,sd}}\right)
^{2}$ be the original sample Wald test statistic. Then Remark \ref%
{remark-Wald} and Theorem \ref{thm:bootstrap_2}(1) immediately imply that
for any $\tau \in (0,1)$,

under $H_{0}:$ $\phi (\alpha _{0})=\phi _{0}$, $\lim_{n\rightarrow \infty
}\Pr \left( \mathcal{W}_{n}\geq \widehat{c}_{2,n}(1-\tau )\right) =\tau $;

under $H_{1}:$ $\phi (\alpha _{0})\neq \phi _{0}$, $\lim_{n\rightarrow
\infty }\Pr \left( \mathcal{W}_{n}\geq \widehat{c}_{2,n}(1-\tau )\right) =1.$

\noindent See Theorem \ref{thm:waldB_con} in Appendix \ref{app:appA} for
properties under local alternatives.
\end{remark}

See online supplemental Appendix \ref{app:appB} for consistency of $\mathcal{%
W}_{1,n}^{B}\equiv \left( \sqrt{n}\frac{\phi (\widehat{\alpha }_{n}^{B})-%
\widehat{\phi }_{n}}{\sigma _{\omega }||\widehat{v}_{n}^{\ast }||_{n,sd}}%
\right) ^{2}$ and other bootstrap sieve Wald (t) statistics based on
different sieve variance estimators.

\subsection{Bootstrap SQLR statistic}

\renewcommand{\theassumption}{\thesection.\arabic{assumption}}

If $\Sigma \neq \Sigma _{0}$, the SQLR statistic $\widehat{QLR}_{n}(\phi
_{0})=n\left( \widehat{Q}_{n}(\widehat{\alpha }_{n}^{R})-\widehat{Q}_{n}(%
\widehat{\alpha }_{n})\right) $ is no longer asymptotically chi-square even
under the null; Theorem \ref{thm:chi2}(1), however, implies that the SQLR
statistic converges weakly to a tight limit under the null. In this
subsection we show that the asymptotic null distribution of the SQLR can be
consistently approximated by that of the (generalized residual) bootstrap
SQLR statistic $\widehat{QLR}_{n}^{B}(\widehat{\phi }_{n})$. Recall that%
\begin{equation*}
\widehat{QLR}_{n}^{B}(\widehat{\phi }_{n})=n\left( \widehat{Q}_{n}^{B}(%
\widehat{\alpha }_{n}^{R,B})-\widehat{Q}_{n}^{B}(\widehat{\alpha }%
_{n}^{B})\right) +o_{P_{V^{\infty }|Z^{\infty }}}(1)~wpa1(P_{Z^{\infty }})
\end{equation*}%
where $\widehat{\phi }_{n}\equiv \phi (\widehat{\alpha }_{n})$, and $%
\widehat{\alpha }_{n}^{R,B}$ is the \emph{restricted} bootstrap PSMD
estimator, defined as
\begin{eqnarray}
\widehat{Q}_{n}^{B}(\widehat{\alpha }_{n}^{R,B})+\lambda _{n}Pen(\widehat{h}%
_{n}^{R,B}) &\leq& \inf_{\alpha \in \mathcal{A}_{k(n)}:\phi (\alpha )=\widehat{%
\phi }_{n}}\left\{ \widehat{Q}_{n}^{B}(\alpha )+\lambda _{n}Pen(h)\right\}\\
&& +o_{P_{V^{\infty }|Z^{\infty }}}(\frac{1}{n})~wpa1(P_{Z^{\infty }}).
\label{R-B-PSMD}
\end{eqnarray}

Lemma \ref{lem:cons_boot} in Appendix \ref{app:appA} implies that $\widehat{%
\alpha }_{n}^{R,B},$ $\widehat{\alpha }_{n}^{B}\in \mathcal{N}_{osn}$ wpa1
under both the null $H_{0}:$ $\phi (\alpha _{0})=\phi _{0}$ and the
alternatives $H_{1}:$ $\phi (\alpha _{0})\neq \phi _{0}$. This indicates
that the bootstrap SQLR statistic $\widehat{QLR}_{n}^{B}(\widehat{\phi }%
_{n}) $ is always properly centered and should be stochastically bounded
under both the null and the alternatives, as shown in the next theorem. Let $%
P_{Z^{\infty }}\left( \widehat{QLR}_{n}(\phi _{0})\leq \cdot \mid
H_{0}\right) $ denote the probability distribution of $\widehat{QLR}%
_{n}(\phi _{0})$ under the null $H_{0}:$ $\phi (\alpha _{0})=\phi _{0}$,
which would converge to the cdf of $\chi _{1}^{2}$ when $\widehat{QLR}%
_{n}(\phi _{0})=\widehat{QLR}_{n}^{0}(\phi _{0})$ (the optimally weighted
SQLR).

\begin{theorem}
\label{thm:bootstrap} Let Assumptions \ref{ass:sieve} - \ref{ass:weak_equiv}
and \ref{ass:rates_B} hold. Let Assumptions \ref{ass:phi}, \ref{ass:LAQ} and %
\ref{ass:LAQ_B} hold with $\left\vert B_{n}^{\omega }-||u_{n}^{\ast
}||^{2}\right\vert =o_{P_{V^{\infty }|Z^{\infty }}}(1)~wpa1(P_{Z^{\infty }})$%
. Then:
\begin{equation*}
\text{(1) }\frac{\widehat{QLR}_{n}^{B}(\widehat{\phi }_{n})}{\sigma _{\omega
}^{2}}=\left( \sqrt{n}\frac{\mathbb{Z}_{n}^{\omega -1}}{\sigma _{\omega
}||u_{n}^{\ast }||}\right) ^{2}+o_{P_{V^{\infty }|Z^{\infty
}}}(1)=O_{P_{V^{\infty }|Z^{\infty }}}(1)~wpa1(P_{Z^{\infty }});
\end{equation*}
and
\begin{eqnarray*}
\text{(2) } && \sup_{t\in \mathbb{R}}\left\vert P_{V^{\infty }|Z^{\infty
}}\left( \frac{\widehat{QLR}_{n}^{B}(\widehat{\phi }_{n})}{\sigma _{\omega
}^{2}}\leq t\mid Z^{n}\right) -P_{Z^{\infty }}\left( \widehat{QLR}_{n}(\phi
_{0})\leq t\mid H_{0}\right) \right\vert \\
&& =o_{P_{V^{\infty }|Z^{\infty
}}}(1)~wpa1(P_{Z^{\infty }}).
\end{eqnarray*}
\end{theorem}

Theorem \ref{thm:bootstrap} allows us to construct valid confidence sets
(CS) for $\phi (\alpha _{0})$ based on inverting possibly \emph{non}%
-optimally weighted SQLR statistic without the need to compute a variance
estimator. We recommend this procedure when it is difficult to compute any
consistent variance estimator for $\phi (\widehat{\alpha })$, such as in the
cases when the residual function $\rho (Z;\alpha )$ is pointwise non-smooth
in $\alpha _{0}$. See, e.g., \cite{AB_Emetrica00} for a thorough discussion
about how to construct CS via bootstrap.

\begin{remark}
\label{remark-QLR-B} Let $\widehat{c}_{n}(a)$ be the $a-th$ quantile of the
distribution of $\frac{\widehat{QLR}_{n}^{B}(\widehat{\phi }_{n})}{\sigma
_{\omega }^{2}}$ (conditional on the data $\{Z_{i}\}_{i=1}^{n}$). Then
Theorems \ref{thm:chi2}, \ref{thm:QLR-H1} and \ref{thm:bootstrap}
immediately imply that for any $\tau \in (0,1)$,

under $H_{0}:$ $\phi (\alpha _{0})=\phi _{0}$, $\lim_{n\rightarrow \infty
}\Pr \left( \widehat{QLR}_{n}(\phi _{0})\geq \widehat{c}_{n}(1-\tau )\right)
=\tau $;

under $H_{1}:$ $\phi (\alpha _{0})\neq \phi _{0}$, $\lim_{n\rightarrow
\infty }\Pr \left( \widehat{QLR}_{n}(\phi _{0})\geq \widehat{c}_{n}(1-\tau
)\right) =1.$

\noindent See Theorem \ref{thm:BSQLR-loc-alt} in Appendix \ref{app:appA} for
properties under local alternatives.
\end{remark}

\section{Verification of Assumptions 3.5 and 3.6}
\label{sec-ex}

In this section, we illustrate the verification of the two key regularity
conditions, Assumption \ref{ass:phi} and Assumption \ref{ass:LAQ}(i), via
some functionals $\phi (h)$ of the (nonlinear) nonparametric IV regressions:
\begin{equation}
E[\rho (Y_{1};h_{0}(Y_{2}))|X]=0\quad a.s.-X,  \label{gnpiv1}
\end{equation}%
where the scalar valued residual function $\rho ()$ could be nonlinear and
pointwise non-smooth in $h$. This model includes the NPIV and NPQIV as
special cases. To be concrete, we consider a PSMD estimator $\widehat{h}\in
\mathcal{H}_{k(n)}$ of $h_{0}$ with $\widehat{\Sigma }=\Sigma =1$, and $%
\widehat{m}(\cdot ,h)$ being the series LS estimator (\ref{mhat}) of $%
m(\cdot ,h)=E[\rho (Y_{1};h(Y_{2}))|X=\cdot ]$ with $J_{n}=ck(n)$ for a
finite constant $c\geq 1$. We assume that $h_{0}\in \mathcal{H}=\Lambda
_{c}^{\varsigma }\left( [-1,1]\right) $ with smoothness $\varsigma >1/2$ (a H%
\"{o}lder ball with support $[-1,1]$, see, e.g., \cite{CLvK_Emetrica03}).%
\footnote{%
This H\"{o}lder ball condition and several other conditions assumed in this
subsection are for illustration only, and can be replaced by weaker
sufficient conditions.} By definition, $\mathcal{H}\subset L^{2}(f_{Y_{2}})$
and we let $||\cdot ||_{s}=||\cdot ||_{L^{2}(f_{Y_{2}})}.$ We assume that $%
\mathcal{H}_{k(n)}=clsp\{q_{1},...,q_{k(n)}\}$ with $\{q_{k}\}_{k=1}^{\infty
}$ being a Riesz basis of $(\mathcal{H},||\cdot ||_{s})$. The convergence
rates of $\widehat{h}$ to $h_{0}$ in both $||\cdot ||$ and $||\cdot
||_{s}=||\cdot ||_{L^{2}(f_{Y_{2}})}$ metrics have already been established
in \cite{CP_WP07}, and hence will not be repeated here.

We use $\mathcal{H}_{os}$ and $\mathcal{H}_{osn}$ for $\mathcal{A}_{os}$ and
$\mathcal{A}_{osn}$ defined in Subsection \ref{sec:consistency} (since there
is no $\theta $ here). Denote $T\equiv \frac{dm(\cdot ,h_{0})}{dh}:\mathcal{H%
}_{os}\subset L^{2}(f_{Y_{2}})\rightarrow L^{2}(f_{X})$\textbf{, }i.e., for
any $h\in \mathcal{H}_{os}\subset L^{2}(f_{Y_{2}})$,
\begin{equation*}
Th\equiv \left. \frac{dE[\rho (Y_{1};h_{0}(Y_{2})+\tau h(Y_{2}))|X=\cdot ]}{%
d\tau }\right\vert _{\tau =0}.
\end{equation*}%
Let $T^{\ast }$ be the adjoint of $T$. Then for all $h\in \mathcal{H}_{os}$,
we have $||h||^{2}\equiv ||Th||_{L^{2}(f_{X})}^{2}=||(T^{\ast
}T)^{1/2}h||_{L^{2}(f_{Y_{2}})}^{2}$. Under mild conditions as stated in
\cite{CP_WP07}, $T$ and $T^{\ast }$ are compact. Then $T$ has a singular
value decomposition $\{\mu _{k};\psi _{k},\phi _{0k}\}_{k=1}^{\infty }$,
where $\{\mu _{k}>0\}_{k=1}^{\infty }$ is the sequence of singular values in
non-increasing order ($\mu _{k}\geq \mu _{k+1}\geq ...$) with $%
\liminf_{k\rightarrow \infty }\mu _{k}=0$, $\{\psi _{k}\in
L^{2}(f_{Y_{2}})\}_{k=1}^{\infty }$ and $\{\phi _{0k}\in
L^{2}(f_{X})\}_{k=1}^{\infty }$ are sequences of eigenfunctions of the
operators $(T^{\ast }T)^{1/2}$ and $(TT^{\ast })^{1/2}$:
\begin{equation*}
T\psi _{k}=\mu _{k}\phi _{0k},\text{\quad }(T^{\ast }T)^{1/2}\psi _{k}=\mu
_{k}\psi _{k}\quad \text{and\quad }(TT^{\ast })^{1/2}\phi _{0k}=\mu _{k}\phi
_{0k}\quad \text{for all }k.
\end{equation*}%
Since $\{q_{k}\}_{k=1}^{\infty }$ is a Riesz basis of $(\mathcal{H},||\cdot
||_{s})$ we could also have $\mathcal{H}_{k(n)}=clsp\{\psi _{1},...,\psi
_{k(n)}\}$. The sieve measure of local ill-posedness now becomes $\tau
_{n}=\mu _{k(n)}^{-1}$ (see, e.g., \cite{BCK_Emetrica07} and \cite{CP_WP07}%
), and hence $\left\Vert u_{n}^{\ast }\right\Vert _{s}\leq c\mu _{k(n)}^{-1}$
for a finite constant $c>0$. Also, $\Pi _{n}h_{0}\equiv \arg \min_{h\in
\mathcal{H}_{k(n)}}||h-h_{0}||_{s}=\sum_{k=1}^{k(n)}\langle h_{0},\psi
_{k}\rangle _{s}\psi _{k}$ is the LS projection of $h_{0}$ onto the sieve
space $\mathcal{H}_{n}$ under the strong norm $||\cdot ||_{s}=||\cdot
||_{L^{2}(f_{Y_{2}})}$. Recall that $h_{0,n}\equiv \arg \min_{h\in \mathcal{H%
}_{k(n)}}||h-h_{0}||^{2}\equiv \arg \min_{h\in \mathcal{H}%
_{k(n)}}||T[h-h_{0}]||_{L^{2}(f_{X})}^{2}$. We have:%
\begin{eqnarray} \notag
h_{0,n} &=&\arg \min_{\{a_{k}\}}\left[ \sum_{k=1}^{k(n)}\left( \langle
h_{0},\psi _{k}\rangle _{s}-a_{k}\right) ^{2}\mu
_{k}^{2}+\sum_{k=k(n)+1}^{\infty }\langle h_{0},\psi _{k}\rangle _{s}^{2}\mu
_{k}^{2}\right]\\
 &= &\sum_{k=1}^{k(n)}\langle h_{0},\psi _{k}\rangle _{s}\psi
_{k}=\Pi _{n}h_{0}.  \label{h0n}
\end{eqnarray}

The next remark specializes Theorem \ref{thm:theta_anorm} to a general
functional $\phi (h)$ of the model (\ref{gnpiv1}).

\begin{remark}
\label{remark-normality-gnpiv1} Let $\widehat{m}$ be the series LS estimator
(\ref{mhat}) for the model (\ref{gnpiv1}) with $\widehat{\Sigma }=\Sigma =1$%
, and Assumptions \ref{ass:sieve}(i)(ii), \ref{A_3.6}(ii)(iii), and \ref%
{ass:weak_equiv} hold with $\delta _{n}=O\left( \sqrt{\frac{k(n)}{n}}\right)
=o(n^{-1/4})$ and $\delta _{s,n}=O\left( \{k(n)\}^{-\varsigma }+\mu
_{k(n)}^{-1}\sqrt{\frac{k(n)}{n}}\right) =o(1)$. Let Assumption \ref{ass:phi}%
, equation (\ref{LF}) and Assumptions \ref{ass:m_ls} - \ref{ass:cont_diffm}
hold. Then:
\begin{equation}
\sqrt{n}\frac{\phi (\widehat{h}_{n})-\phi (h_{0})}{||v_{n}^{\ast }||_{sd}}%
\Rightarrow N(0,1),  \label{P-var-gnpiv1}
\end{equation}
with $||v_{n}^{\ast }||_{sd}^{2}=(\frac{%
	d\phi (h_{0})}{dh} [q^{k(n)}(\cdot )])^{\prime }D_{n}^{-1}{\mho }_{n}D_{n}^{-1}(\frac{d\phi (h_{0})}{dh}%
 [q^{k(n)} (\cdot )])$, and $D_{n}=E\left[ \left( T[q^{k(n)}(\cdot )^{\prime }]\right) ^{\prime }\left(
T[q^{k(n)}(\cdot )^{\prime }]\right) \right] $ and $\mho _{n}=E\left[ \left(
T[q^{k(n)}(\cdot )^{\prime }]\right) ^{\prime }\rho (Z,h_{0})^{2}\left(
T[q^{k(n)}(\cdot )^{\prime }]\right) \right] .$
\end{remark}

Remark \ref{remark-normality-gnpiv1} includes the NPIV\ and NPQIV examples
in Subsection \ref{sec:NPIVex} as special cases. In particular, the sieve
variance expression (\ref{P-var-gnpiv1}) reproduces the one for the NPIV
model (\ref{npiv}) with $T[q^{k(n)}(\cdot )^{\prime
}]=E[q^{k(n)}(Y_{2})^{\prime }|X]$, and the one for the NPQIV model (\ref%
{npqiv}) with $T[q^{k(n)}(\cdot )^{\prime
}]=E[f_{U|Y_{2},X}(0)q^{k(n)}(Y_{2})^{\prime }|X]$.

By the result in \cite{CP_WP07}, the sieve dimension $k_{n}^{\ast }$
satisfying $\{k_{n}^{\ast }\}^{-\varsigma }\asymp \mu _{k_{n}^{\ast
}}^{-1}\times \sqrt{\frac{k_{n}^{\ast }}{n}}$ leads to the nonparametric
optimal convergence rate of $||\widehat{h}-h_{0}||_{s}=O_{P_{Z^{\infty
}}}(\delta _{s,n}^{\ast })=o(1)$ in strong norm, where $\delta _{s,n}^{\ast
}\asymp \{k_{n}^{\ast }\}^{-\varsigma }$. In particular, $k_{n}^{\ast
}\asymp n^{\frac{1}{2(\varsigma +a)+1}}$ and $\delta _{s,n}^{\ast }=n^{-%
\frac{\varsigma }{2(\varsigma +a)+1}}$ for the \textit{mildly ill-posed case}
$\mu _{k}\asymp k^{-a}$ for a finite $a>0$; and $\delta _{s,n}^{\ast }=\{\ln
n\}^{-\varsigma }$ for the \textit{severely ill-posed case} $\mu _{k}\asymp
\exp \{-0.5ak\}$ for a finite $a>0$. However this paper aims at simple valid
inferences on functional $\phi (h_{0})$. As will be illustrated in the next
subsection, although the nonparametric optimal choice $k_{n}^{\ast }$ is
compatible with the sufficient conditions for the asymptotic normality of $%
\sqrt{n}(\phi (\widehat{h})-\phi (h_{0}))$ for a regular linear functional $%
\phi (h_{0})$ (see Remark \ref{remark-bias}), it is typically ruled out by
Assumption \ref{ass:phi}(iii) for irregular functionals.

\subsection{Verification of Assumption \protect\ref{ass:phi}}

Let $b_{j}\equiv \frac{d\phi (h_{0})}{dh}[\psi _{j}(\cdot )]$ for all $j$.
By Lemma \ref{lem:P-sieve} $D_{n}=E\left[ \left( T[q^{k(n)}(\cdot )^{\prime
}]\right) ^{\prime }\left( T[q^{k(n)}(\cdot )^{\prime }]\right) \right]
=Diag\left\{ \mu _{1}^{2},...,\mu _{k(n)}^{2}\right\} $ and
\begin{equation}
||v_{n}^{\ast }||^{2}=\left( \frac{d\phi (h_{0})}{dh}[q^{k(n)}(\cdot
)]\right) ^{\prime }D_{n}^{-1}\left( \frac{d\phi (h_{0})}{dh}[q^{k(n)}(\cdot
)]\right) =\sum_{j=1}^{k(n)}\mu _{j}^{-2}b_{j}^{2}.
\label{sieve-riesz-gnpiv1}
\end{equation}%
By Lemma \ref{lem:sieve-Riesz-property}, $\phi (h)$ of the model (\ref%
{gnpiv1}) is regular (at $h=h_{0}$) iff $\sum_{j=1}^{\infty }\mu
_{j}^{-2}b_{j}^{2}<\infty $, and is irregular (at $h=h_{0}$) iff $%
\sum_{j=1}^{\infty }\mu _{j}^{-2}b_{j}^{2}=\infty $.

For the same functional $\phi (h)$ of a model (\ref{ex-gnpiv}) without
endogeneity:%
\begin{equation}
E[\rho (Y_{1};h_{0}(Y_{2}))|Y_{2}]=0\quad a.s.-Y_{2},  \label{ex-gnpiv}
\end{equation}%
we have $D_{n}\asymp I_{k(n)}$ and $||v_{n}^{\ast }||^{2}\asymp
\sum_{j=1}^{k(n)}b_{j}^{2}$. Thus, $\phi (h)$ of the model (\ref{ex-gnpiv})
is regular (or irregular) iff $\sum_{j=1}^{\infty }b_{j}^{2}<\infty $ (or $%
=\infty $).

Since $\mu _{k(n)}\rightarrow 0$ as $k(n)\rightarrow \infty $, if a
functional $\phi (h)$ is irregular for the model (\ref{ex-gnpiv}) without
endogeneity, then it is irregular for the model (\ref{gnpiv1}). But, even if
a functional $\phi (h)$ is regular for the model (\ref{ex-gnpiv}) without
endogeneity, it could still be irregular for the model (\ref{gnpiv1}) with
endogeneity.

\subsubsection{Linear functionals of the model (\protect\ref{gnpiv1})}

For a linear functional $\phi (h)$ of the model (\ref{gnpiv1}), given
relation (\ref{h0n}), Assumption \ref{ass:phi} is satisfied provided that
the sieve dimension $k(n)$ satisfies (\ref{a3.1-linear}):%
\begin{equation}
\frac{||v_{n}^{\ast }||}{\sqrt{n}}=o(1)\text{\quad and\quad }\sqrt{n}\frac{%
\left\vert \frac{d\phi (h_{0})}{dh}[\Pi _{n}h_{0}-h_{0}]\right\vert }{%
||v_{n}^{\ast }||}=o(1).  \label{a3.1-linear}
\end{equation}%
When $\phi (h)$ of the model (\ref{gnpiv1}) is regular, Remark \ref%
{remark-bias} implies that (\ref{a3.1-linear}) is satisfied provided
\begin{equation}
\sum_{j=1}^{\infty }\mu _{j}^{-2}b_{j}^{2}<\infty \text{\quad and\quad }%
n\times \sum_{j=k(n)+1}^{\infty }\mu _{j}^{-2}b_{j}^{2}\times ||\Pi
_{n}h_{0}-h_{0}||^{2}=o(1).  \label{a3.1-linear-r}
\end{equation}%
We shall illustrate below that both these sufficient conditions allow for
severely ill-posed problems.

\smallskip

\textbf{Example 1 (evaluation functional).} For $\phi (h)=h(\overline{y}%
_{2}) $, we have: $||v_{n}^{\ast }||^{2}=\sum_{j=1}^{k(n)}\mu _{j}^{-2}[\psi
_{j}(\overline{y}_{2})]^{2}$. Let $\mathcal{H}_{k(n)}$ be the spline or the CDV wavelet sieve as described in \cite{CC_WP13}, say. Then
\begin{equation*}
\left\vert \frac{d\phi (h_{0})}{dh}[\Pi _{n}h_{0}-h_{0}]\right\vert =|(\Pi
_{n}h_{0})(\overline{y}_{2})-h_{0}(\overline{y}_{2})|\leq ||\Pi
_{n}h_{0}-h_{0}||_{\infty }\leq const.\{k(n)\}^{-\varsigma }.
\end{equation*}%
To provide concrete sufficient condition for (\ref{a3.1-linear}), we assume $%
||v_{n}^{\ast }||^{2}\asymp E\left( \sum_{j=1}^{k(n)}\mu _{j}^{-2}[\psi
_{j}(Y_{2})]^{2}\right) =\sum_{k=1}^{k(n)}\mu _{k}^{-2}$. Since $%
\lim_{k(n)\rightarrow \infty }||v_{n}^{\ast }||^{2}=\infty $, the evaluation
functional is irregular. Condition (\ref{a3.1-linear}) is satisfied provided
that%
\begin{equation}
\frac{||v_{n}^{\ast }||^{2}}{n}=\frac{\sum_{k=1}^{k(n)}\mu _{k}^{-2}}{n}=o(1)%
\text{\quad and\quad }\frac{\{k(n)\}^{-2\varsigma }}{\frac{1}{n}%
||v_{n}^{\ast }||^{2}}=\frac{\{k(n)\}^{-2\varsigma }}{\frac{1}{n}%
\sum_{k=1}^{k(n)}\mu _{k}^{-2}}=o(1).  \label{Ex1-both}
\end{equation}%
Condition (\ref{Ex1-both}) allows for both mildly and severely ill-posed
cases.

(a) \textit{Mildly ill-posed}: $\mu _{k}\asymp k^{-a}$ for a finite $a>0$.
Then $||v_{n}^{\ast }||^{2}\asymp \{k(n)\}^{2a+1}$. Condition (\ref{Ex1-both}%
) is satisfied by a wide range of sieve dimensions, such as $k(n)\asymp n^{%
\frac{1}{2(\varsigma +a)+1}}(\ln \ln n)^{\varpi }$ or $n^{\frac{1}{%
2(\varsigma +a)+1}}(\ln n)^{\varpi }$ for any finite $\varpi >0$, or $%
k(n)\asymp n^{\epsilon }$ for any $\epsilon \in (\frac{1}{2(\varsigma +a)+1},%
\frac{1}{2a+1})$. Note that any $k(n)$ satisfying Condition (\ref{Ex1-both})
also ensures $\delta _{s,n}=o(1)$. However, it does require $%
k(n)/k_{n}^{\ast }\rightarrow \infty $, where $k_{n}^{\ast }\asymp n^{\frac{1%
}{2(\varsigma +a)+1}}$ is the choice for the nonparametric optimal
convergence rate in strong norm.

(b) \textit{Severely ill-posed}: $\mu _{k}\asymp \exp \{-0.5ak\}$ for a
finite $a>0$. Then $||v_{n}^{\ast }||^{2}\asymp \exp \{ak(n)\}$. Condition (%
\ref{Ex1-both}) is satisfied with $k(n)\asymp a^{-1}\left[ \ln n-\varpi \ln
(\ln n)\right] $ for $0<\varpi <2\varsigma $. In addition we need $\varpi >1$
(and hence $\varsigma >1/2$) to ensure $\delta _{s,n}=O\left(
\{k(n)\}^{-\varsigma }+\mu _{k(n)}^{-1}\sqrt{\frac{k(n)}{n}}\right) =o(1)$.

\smallskip

\textbf{Example 2 (weighted derivative functional).} For $\phi (h)=\int
w(y)\nabla h(y)dy$, where $w(y)$ is a weight satisfying the integration by
part formula: $\phi (h)=\int w(y)\nabla h(y)dy=-\int h(y)\nabla w(y)dy$, we
have: $||v_{n}^{\ast }||^{2}=\sum_{j=1}^{k(n)}\mu _{j}^{-2}b_{j}^{2}$ with $%
b_{j}=\int \psi _{j}(y)\nabla w(y)dy$ for all $j$, and%
\begin{eqnarray*}
\left\vert \frac{d\phi (h_{0})}{dh}[\Pi _{n}h_{0}-h_{0}]\right\vert
&= &\left\vert \int [\Pi _{n}h_{0}(y)-h_{0}(y)]\nabla w(y)dy\right\vert \\
&\leq &
C\times ||\Pi _{n}h_{0}-h_{0}||_{L^{2}(f_{Y_{2}})}\leq
const.\{k(n)\}^{-\varsigma }
\end{eqnarray*}%
provided that $E\left( \left[ \frac{\nabla w(Y_{2})}{f_{Y_{2}}(Y_{2})}\right]
^{2}\right) =\sum_{j=1}^{\infty }b_{j}^{2}=C<\infty $. That is, the weighted
derivative is assumed to be regular for the model (\ref{ex-gnpiv}) without
endogeneity.

\textbf{(i)} When the weighted derivative is regular (i.e., $%
\sum_{j=1}^{\infty }\mu _{j}^{-2}b_{j}^{2}<\infty $) for the model (\ref%
{gnpiv1}), Condition (\ref{a3.1-linear-r}) is satisfied provided that $%
n\times \sum_{j=k(n)+1}^{\infty }\mu _{j}^{-2}b_{j}^{2}\times \delta
_{n}^{2}=o(1)$, which is the condition imposed in \cite{AC_JOE07} for their
root-$n$ estimation of an average derivative of NPIV example, and is shown
to allow for severely ill-posed inverse case in \cite{AC_JOE07}.

\textbf{(ii)} When the weighted derivative is irregular (i.e., $%
\sum_{j=1}^{\infty }\mu _{j}^{-2}b_{j}^{2}=\infty $) for the model (\ref%
{gnpiv1}), Condition (\ref{a3.1-linear}) is satisfied provided that%
\begin{equation}
\frac{||v_{n}^{\ast }||^{2}}{n}=\frac{\sum_{j=1}^{k(n)}\mu _{j}^{-2}b_{j}^{2}%
}{n}=o(1)\text{\quad and\quad }\frac{\{k(n)\}^{-2\varsigma }}{\frac{1}{n}%
||v_{n}^{\ast }||^{2}}=\frac{\{k(n)\}^{-2\varsigma }}{\frac{1}{n}%
\sum_{j=1}^{k(n)}\mu _{j}^{-2}b_{j}^{2}}=o(1).  \label{Ex2-severe}
\end{equation}%
Condition (\ref{Ex2-severe}) allows for both mildly and severely ill-posed
cases. To provide concrete sufficient conditions for (\ref{Ex2-severe}) we
assume $b_{j}^{2}\asymp \left( j\ln (j)\right) ^{-1}$ in the following
calculations.

(a) \textit{Mildly ill-posed}: $\mu _{k}\asymp k^{-a}$ for a finite $a>0$.
Then $||v_{n}^{\ast }||^{2}\in \lbrack c\frac{k(n)^{2a}}{\ln (k(n))}%
,c^{\prime }k(n)^{2a}]$ for some $0<c\leq c^{\prime }<\infty $. Condition (%
\ref{Ex2-severe}) and $\delta _{s,n}=o(1)$ are jointly satisfied by a wide
range of sieve dimensions, such as $k(n)\asymp n^{\frac{1}{2(\varsigma +a)}%
}(\ln n)^{\varpi }$ for any finite $\varpi >\frac{1}{2(\varsigma +a)}$, or $%
k(n)\asymp n^{\epsilon }$ for any $\epsilon \in (\frac{1}{2(\varsigma +a)},%
\frac{1}{2a+1})$ and $\varsigma >1/2$.

(b) \textit{Severely ill-posed}: $\mu _{k}\asymp \exp \{-0.5ak\}$ for $a>0$.
Then $||v_{n}^{\ast }||^{2}\in \lbrack c\frac{\exp \{ak(n)\}}{k(n)\ln (k(n))}%
,c^{\prime }\frac{\exp \{ak(n)\}}{\ln (k(n))}]$ for some $0<c\leq c^{\prime
}<\infty $. Condition (\ref{Ex2-severe}) and $\delta _{s,n}=o(1)$ are
jointly satisfied by $k(n)\asymp a^{-1}\left[ \ln (n)-\varpi \ln (\ln (n))%
\right] $ for $\varpi \in (1,2\varsigma -1)$ and $\varsigma >1$.

\subsubsection{Nonlinear functionals}

For a nonlinear functional $\phi (h)$ of the model (\ref{gnpiv1}),
Assumption \ref{ass:phi} is satisfied provided that the sieve dimension $%
k(n) $ satisfies (\ref{a3.1-linear}) (or (\ref{a3.1-linear-r}) if $\phi (h)$
is regular) and Assumption \ref{ass:phi}(ii), which is implied by the
following condition:

\noindent \textbf{Assumption \ref{ass:phi}(ii)'}: \textit{there are finite
non-negative constants }$C\geq 0,\omega _{1},\omega _{2}\geq 0$\textit{\
such that for all }$(\alpha ,t)\in \mathcal{N}_{osn}\times \mathcal{T}_{n}$,
\begin{eqnarray*}
&& \left\vert \phi (\alpha +tu_{n}^{\ast })-\phi (\alpha _{0})-\frac{d\phi
(\alpha _{0})}{d\alpha }[\alpha +tu_{n}^{\ast }-\alpha _{0}]\right\vert \\
& &\leq
C\times (||\alpha -\alpha _{0}+tu_{n}^{\ast }||^{\omega _{1}}\times ||\alpha
-\alpha _{0}+tu_{n}^{\ast }||_{s}^{\omega _{2}}),
\end{eqnarray*}
and
\begin{equation*}
C\times \frac{\sqrt{n}\times (\delta _{n}(1+M_{n}^{2}))^{\omega _{1}}\times
(\delta _{s,n}+M_{n}^{2}\delta _{n}||u_{n}^{\ast }||_{s})^{\omega _{2}}}{%
||v_{n}^{\ast }||}=o\left( 1\right) .
\end{equation*}

Assumption \ref{ass:phi}(ii) or (ii)' controls the nonlinearity bias of $%
\phi \left( \cdot \right) $ (i.e., the linear approximation error of a
nonlinear functional $\phi \left( \cdot \right) $). It typically rules out
nonlinear regular functionals of severely illposed inverse problems, but
allows for nonlinear irregular functionals of severely illposed inverse
problems.

\textbf{Example 3 (weighted quadratic functional).} For $\phi (h)=\frac{1}{2}%
\int w(y)\left\vert h(y)\right\vert ^{2}dy$, we have $||v_{n}^{\ast
}||^{2}=\sum_{j=1}^{k(n)}\mu _{j}^{-2}b_{j}^{2}$ with $b_{j}=\int
h_{0}(y)w(y)\psi _{j}(y)dy$ for all $j$, and%
\begin{equation*}
\left\vert \frac{d\phi (h_{0})}{dh}[\Pi _{n}h_{0}-h_{0}]\right\vert
=\left\vert \int w(y)h_{0}(y)[\Pi _{n}h_{0}(y)-h_{0}(y)]dy\right\vert \leq
const.\times ||\Pi _{n}h_{0}-h_{0}||_{L^{2}(f_{Y_{2}})}
\end{equation*}%
provided that $\sup_{y}\frac{w(y)}{f_{Y_{2}}(y)}<\infty $. This and $E\left( %
\left[ h_{0}(Y_{2})\right] ^{2}\right) <\infty $ imply that $%
\sum_{j=1}^{\infty }b_{j}^{2}<\infty $. That is, the weighted quadratic
functional is regular for the model (\ref{ex-gnpiv}) without endogeneity.
Also,%
\begin{equation*}
\left\vert \phi (h)-\phi (h_{0})-\frac{d\phi (h_{0})}{dh}[h-h_{0}]\right%
\vert =\frac{1}{2}\int w(y)\left\vert h(y)-h_{0}(y)\right\vert ^{2}dy\leq
const.\times ||h-h_{0}||_{L^{2}(f_{Y_{2}})}^{2}.
\end{equation*}

\textbf{(i)} When the weighted quadratic functional is regular (i.e., $%
\sum_{j=1}^{\infty }\mu _{j}^{-2}b_{j}^{2}<\infty $) for the model (\ref%
{gnpiv1}), Condition (\ref{a3.1-linear-r}) is satisfied provided that $%
n\times \sum_{j=k(n)+1}^{\infty }\mu _{j}^{-2}b_{j}^{2}\times \delta
_{n}^{2}=o(1)$, which allows for severely ill-posed cases. But Assumption %
\ref{ass:phi}(ii)' requires that $\sqrt{n}\times \delta _{s,n}^{2}=\sqrt{n}%
\times \left( \{k(n)\}^{-\varsigma }+\mu _{k(n)}^{-1}\sqrt{\frac{k(n)}{n}}%
\right) ^{2}=o(1)$, which clearly rules out severely ill-posed inverse case
where $\mu _{k}\asymp \exp \{-0.5ak\}$ for some finite $a>0$.

\textbf{(ii)} When the weighted quadratic functional is irregular (i.e., $%
\sum_{j=1}^{\infty }\mu _{j}^{-2}b_{j}^{2}=\infty $) for the model (\ref%
{gnpiv1}), Condition (\ref{a3.1-linear}) is satisfied provided that
Condition (\ref{Ex2-severe}) holds with $b_{j}=\int h_{0}(y)w(y)\psi
_{j}(y)dy$ for Example 3. Assumption \ref{ass:phi}(ii)' is satisfied
provided that%
\begin{equation}
\sqrt{n}\frac{\delta _{s,n}^{2}}{||v_{n}^{\ast }||}=\frac{\sqrt{n}\times
\left( \{k(n)\}^{-\varsigma }+\mu _{k(n)}^{-1}\sqrt{\frac{k(n)}{n}}\right)
^{2}}{||v_{n}^{\ast }||}\leq n^{-1/2}\frac{\mu _{k(n)}^{-2}k(n)}{\sqrt{%
\sum_{j=1}^{k(n)}\mu _{j}^{-2}b_{j}^{2}}}=o(1).  \label{Ex3-severe2}
\end{equation}%
Any $k(n)$ satisfying Conditions (\ref{Ex2-severe}) and (\ref{Ex3-severe2})
automatically satisfies $\delta _{s,n}=o(1)$. In addition, both conditions
allow for mildly and severely ill-posed cases. To provide concrete
sufficient conditions we assume $b_{j}^{2}\asymp \left( j\ln (j)\right)
^{-1} $ in the following calculations.

(a) \textit{Mildly ill-posed}: $\mu _{k}\asymp k^{-a}$ for a finite $a>0$.
Then $||v_{n}^{\ast }||^{2}\in \lbrack c\frac{k(n)^{2a}}{\ln (k(n))}%
,c^{\prime }k(n)^{2a}]$ for some $0<c\leq c^{\prime }<\infty $. Conditions (%
\ref{Ex2-severe}) and (\ref{Ex3-severe2}) are satisfied by a wide range of
sieve dimensions, such as $k(n)\asymp n^{\frac{1}{2(\varsigma +a)}}(\ln
n)^{\varpi }$ for any finite $\varpi >\frac{1}{2(\varsigma +a)}$, or $%
k(n)\asymp n^{\epsilon }$ for any $\epsilon \in (\frac{1}{2(\varsigma +a)},%
\frac{1}{2a+2})$ and $\varsigma >1$.

(b) \textit{Severely ill-posed}: $\mu _{k}\asymp \exp \{-0.5ak\}$ for $a>0$.
Then $||v_{n}^{\ast }||^{2}\in \lbrack c\frac{\exp \{ak(n)\}}{k(n)\ln (k(n))}%
,c^{\prime }\frac{\exp \{ak(n)\}}{\ln (k(n))}]$ for some $0<c\leq c^{\prime
}<\infty $. Conditions (\ref{Ex2-severe}) and (\ref{Ex3-severe2}) are
satisfied with $k(n)\asymp a^{-1}\left[ \ln (n)-\varpi \ln (\ln (n))\right] $
and $\varpi \in (3,2\varsigma -1)$ for $\varsigma >2$.

\subsection{Verification of Assumption \protect\ref{ass:LAQ}(i)}

By Lemma \ref{lem:Qdiff_B}(1), to verify Assumption \ref{ass:LAQ}(i), it
suffices to verify Assumptions \ref{ass:m_ls} - \ref{ass:cont_diffm} in
Appendix \ref{app:appA}. Note that Assumptions \ref{ass:m_ls} and \ref%
{ass:rho_Donsker} do not depend on sieve Riesz representer at all, and have
already been verified in \cite{CP_WP07a}, \cite{AC_JOE07} and others for
(penalized) SMD estimators for the model (\ref{gnpiv1}). Assumptions \ref%
{ass:anor-mtilde} and \ref{ass:cont_diffm} do depend on the scaled sieve
Riesz representer $u_{n}^{\ast }\equiv v_{n}^{\ast }/||v_{n}^{\ast }||_{sd}$%
. Both these assumptions are also verified in \cite{AC_Emetrica03}, \cite%
{CP_WP07a}, \cite{AC_JOE07} for examples of regular functionals of the model
(\ref{gnpiv1}). Here, we present simple (albeit somewhat strong) sufficient conditions for Assumptions \ref%
{ass:anor-mtilde} and \ref{ass:cont_diffm} for irregular functionals of the
NPIV and NPQIV examples.

\begin{condition}
\label{eqn:NPIV-holder} (i) $\{E[h(Y_{2})|\cdot ]:h\in \mathcal{H}%
\}\subseteq \Lambda _{c}^{\gamma }(\mathcal{X})$, with $\gamma >0.5$; (ii) $%
\sup_{x,y_{2}}\frac{f_{Y_{2}X}(y_{2},x)}{f_{Y_{2}}(y_{2})f_{X}(x)}\leq
Const.<\infty $.
\end{condition}

\begin{proposition}
\label{pro:NPIV-suff} Let all conditions for Remark \ref%
{remark-normality-gnpiv1} hold. Under Condition \ref{eqn:NPIV-holder},
Assumptions \ref{ass:anor-mtilde} and \ref{ass:cont_diffm} hold for the NPIV
model (\ref{npiv}).
\end{proposition}

Proposition \ref{pro:NPIV-suff} allows for irregular functionals of the NPIV
model with severely ill-posed case.

\begin{condition}
\label{eqn:NPQIV-holder} (i) $\{E[F_{Y_{1}|Y_{2}X}(h(Y_{2}),Y_{2},\cdot
)|\cdot ]:h\in \mathcal{H}\}\subseteq \Lambda _{c}^{\gamma }(\mathcal{X})$,
with $\gamma >0.5$; (ii) $\sup_{y_{1},y_{2},x}|\frac{%
df_{Y_{1}|Y_{2}X}(y_{1},y_{2},x)}{dy_{1}}|\leq C<\infty $.
\end{condition}

\begin{condition}
\label{con:NPQIV-rate1} $n(\log \log n)^{4}\delta _{s,n}^{4}=o(1)$
\end{condition}

\begin{proposition}
\label{pro:NPQIV-suff} Let all conditions for Remark \ref%
{remark-normality-gnpiv1} hold. Under conditions \ref{eqn:NPIV-holder}(ii)
and \ref{eqn:NPQIV-holder}-\ref{con:NPQIV-rate1}, Assumptions \ref%
{ass:anor-mtilde} and \ref{ass:cont_diffm} hold for the NPQIV model (\ref%
{npqiv}).
\end{proposition}

It is clear that Condition \ref{con:NPQIV-rate1} rules out severely
ill-posed case, and hence Proposition \ref{pro:NPQIV-suff} only allows for
irregular functionals of the NPQIV model with mildly ill-posed case.

\section{Simulation Studies and An Empirical Illustration}

\label{sec:sec_simulation}

\renewcommand{\thetable}{7.\arabic{table}}

\renewcommand{\thefigure}{7.\roman{figure}}

This section first presents simulation studies for SQLR and sieve t tests of
linear and nonlinear hypotheses for the NPQIV and NPIV models respectively.
It then provides an empirical illustration of the optimally weighted SQLR
inferences for a NPQIV Engel curve. In this section, we use the series LS
estimator (\ref{mhat}) of $m(x,h)$ with $p^{J_{n}}(x)$ as its basis, and $%
q^{k(n)}$ as the basis approximating the unknown structure function $h_{0}$.
We use $p^{J}=\mathrm{P-Spline}(r,k)$ to denote $r$th degree polynomial
spline with $k$ (quantile) equally spaced knots, hence $J=(r+1)+k$ is the
total number of sieve terms. We use $p^{J}=\mathrm{Pol}(J)$ to denote power
series up to $(J-1)$th degree. See, e.g., \cite{C_bookchp07} for definitions
of these and other sieve bases.

\subsection{Simulation Studies\label{sec:sec_simulation1}}

We run Monte Carlo (MC) studies to assess the finite sample performance of
SQLR and sieve t tests of linear and nonlinear hypotheses in two models: the
NPQIV (\ref{npqiv}) and the NPIV (\ref{npiv}).

For all cases, our design is based on the MC design of \cite{NP_ECMA03} and
\cite{Santos_ECMA} for a NPIV model, which we adapt to cover both NPIV and
NPQIV models. Specifically, we generate i.i.d. draws of $(Y_{2},X,U^{\ast })$
from
\begin{equation*}
\left[
\begin{array}{c}
Y_{2}^{\ast } \\
X^{\ast } \\
U^{\ast }%
\end{array}%
\right] \sim N\left( 0,%
\begin{bmatrix}
1 & 0.8 & 0.5 \\
0.8 & 1 & 0 \\
0.5 & 0 & 1%
\end{bmatrix}%
\right) ,
\end{equation*}%
and $Y_{2}=2(\Phi (Y_{2}^{\ast }/3)-0.5)$ and $X=2(\Phi (X^{\ast }/3)-0.5)$.
The true function $h_{0}$ is given by $h_{0}(\cdot )=2\sin (\pi \cdot )$. We
consider 5,000 MC repetitions and $n=750$ for each of the cases studied
below. We use $Pen(h)=||h||_{L^{2}}^{2}+||\nabla h||_{L^{2}}^{2}$ in all the
simulations, and have used a very small $\lambda _{n}=10^{-5}$ in most cases
(except for the cases where we study the sensitivity to the choice of $\lambda
_{n} $).

\textbf{Summary of sensitivity checks}: For NPQIV and NPIV models, for both
SQLR and sieve t tests of linear and nonlinear hypotheses, as long as $%
J_{n}>k(n)+1$ with not too large $k(n)$, the MC sizes of the tests are good
and insensitive to the choices of basis $q^{k(n)}$ and $p^{J_{n}}$ or the
very small penalty $\lambda _{n}$. This is consistent with previous MC
findings in \cite{BCK_Emetrica07} and \cite{CP_WP07} for PSMD estimation of
NPIV and NPQIV respectively.

\medskip

\noindent \textbf{NPQIV model: SQLR test for an irregular linear functional.}
We consider the NPQIV model $Y_{1}=h_{0}(Y_{2})+U=2\sin (\pi Y_{2})+U$ with $%
U=2(\Phi (U^{\ast })-\gamma )$. This last transformation is done to ensure
that $E[1\{U\leq 0\}|X]=\gamma $. To save space we only present the case
with $\gamma =0.5$. The parameter of interest is $\phi (h_{0})=h_{0}(0)$,
hence $\phi $ is an irregular linear functional. We study the finite sample
properties of the SQLR and bootstrap-SQLR tests. The SQLR-based confidence
intervals are specially well-suited for models like NPQIV where the
generalized residual function is non-smooth yet the optimal weighting matrix
is easy to compute.

\textbf{Size}. Table \ref{tab:NPQIV-SQLR} reports the simulated size of the
SQLR test of $H_{0}\colon \phi (h_{0})=0$ as a function of the nominal size
(NS), for different choices of $q^{k(n)}$ and $p^{J_{n}}$, and different
values of the tuning parameters $(\lambda _{n},k(n),J_{n})$.

\begin{table}[h]
\centering
\caption{Size of the SQLR test of $\protect\phi%
	(h_{0}) = 0$ for NPQIV model.}
\begin{tabular}{cccccc}
\hline\hline
$q^{k(n)}$ & $p^{J_{n}}$ & $\lambda_{n}$ & 10\% & 5\% & 1\% \\ \hline
\multirow{3}{*}{Pol(4)} & Pol(7) & {\small {$(1 \times 10^{-3})$ }} & 0.099 &
0.055 & 0.008 \\
& Pol(7) & {\small {$(2 \times 10^{-4})$ }} & 0.096 & 0.048 & 0.008 \\
& Pol(7) & {\small {$(4 \times 10^{-5})$ }} & 0.107 & 0.053 & 0.010 \\ \hline
\multirow{3}{*}{Pol(6)} & Pol(7) & {\small {$(1 \times 10^{-3})$ }} & 0.133 &
0.068 & 0.011 \\
& Pol(7) & {\small {$(2 \times 10^{-4})$ }} & 0.091 & 0.036 & 0.006 \\
& Pol(7) & {\small {$(4 \times 10^{-5})$ }} & 0.105 & 0.052 & 0.008 \\ \hline
\multirow{3}{*}{Pol(6)} & Pol(9) & {\small {$(1 \times 10^{-5})$ }} & 0.107 &
0.055 & 0.012 \\
& Pol(15) & {\small {$(1 \times 10^{-5})$ }} & 0.109 & 0.058 & 0.014 \\
& Pol(21) & {\small {$(1 \times 10^{-5})$ }} & 0.112 & 0.058 & 0.013 \\
\hline
\multirow{4}{*}{P-Spline(3,2)} & Pol(9) & {\small {$(1 \times 10^{-5})$ }} &
0.103 & 0.049 & 0.010 \\
& Pol(10) & {\small {$(1 \times 10^{-5})$ }} & 0.104 & 0.051 & 0.010 \\
& Pol(15) & {\small {$(1 \times 10^{-5})$ }} & 0.105 & 0.049 & 0.009 \\
& Pol(21) & {\small {$(1 \times 10^{-5})$ }} & 0.105 & 0.052 & 0.009 \\
\hline
\multirow{3}{*}{P-Spline(3,2)} & P-Spline(5,3) & {\small {$(1 \times 10^{-5})$
}} & 0.098 & 0.049 & 0.008 \\
& P-Spline(5,9) & {\small {$(1 \times 10^{-5})$ }} & 0.103 & 0.050 & 0.009
\\
& P-Spline(5,18) & {\small {$(1 \times 10^{-5})$ }} & 0.106 & 0.051 & 0.009
\\ \hline
\end{tabular}%
\label{tab:NPQIV-SQLR}
\end{table}

\noindent Table \ref{tab:NPQIV-SQLR} shows that for small value of $k(n)$,
say in $(k(n),J_{n})=(4,7)$ (i.e., rows 1-3), the SQLR test performs well
and is fairly insensitive to different choices of $\lambda _{n}$. For a
fixed relatively small $J_{n}=7$, rows 1-6 indicate that as $k(n)$
increases, the results become a bit more sensitive to the choice of $\lambda
_{n}$. For a fixed very small penalty $\lambda _{n}=10^{-5}$, rows 7-16 show
that the results are fairly insensitive to different choices of $J_{n}$ and
basis for $p^{J_{n}}$ and $q^{k(n)}$ as long as $J_{n}>k(n)+1$.

\begin{figure}[h]
\centering
\includegraphics[height=1.5in,width=4in]{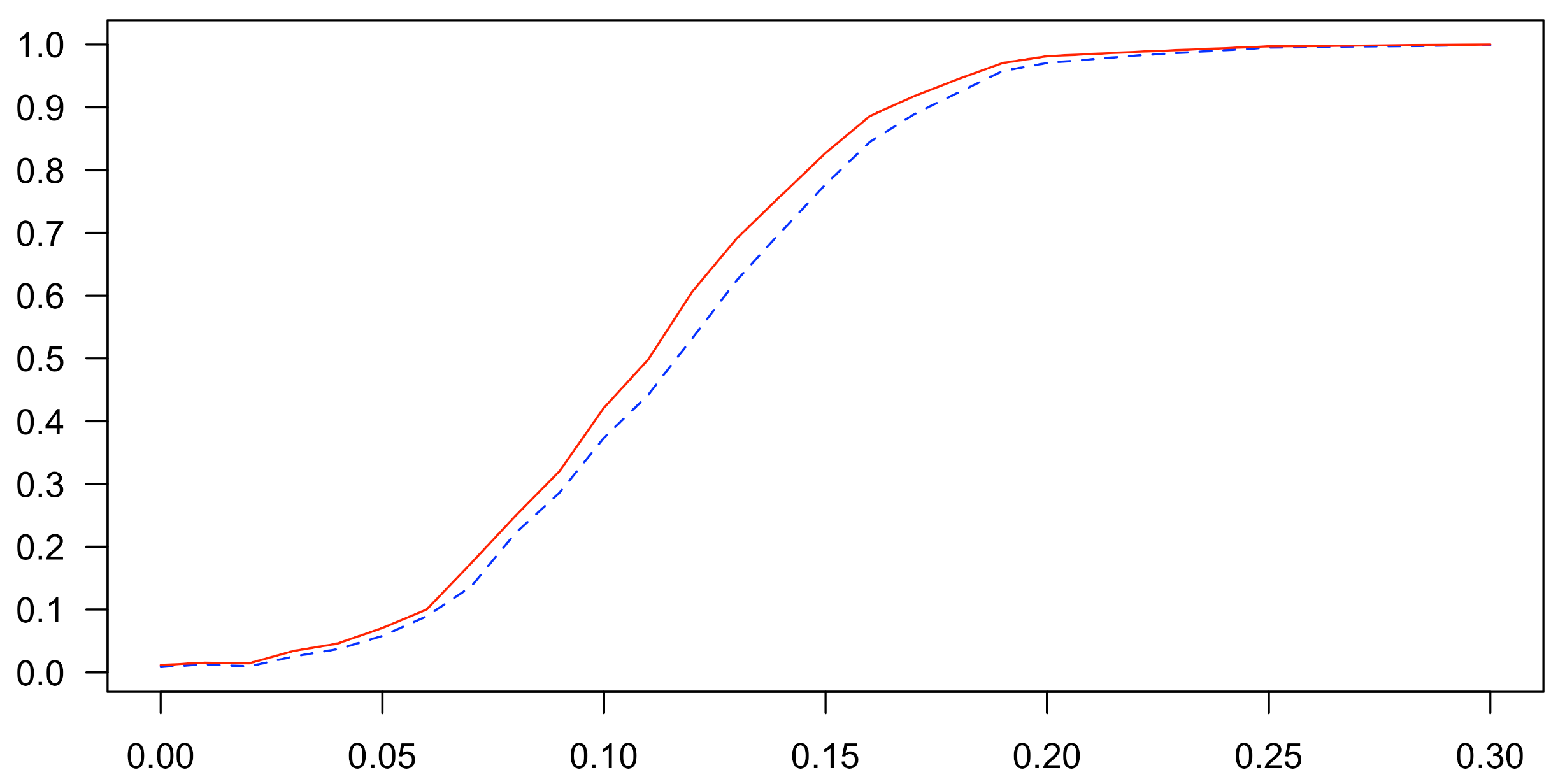}
\includegraphics[height=1.5in,width=4in]{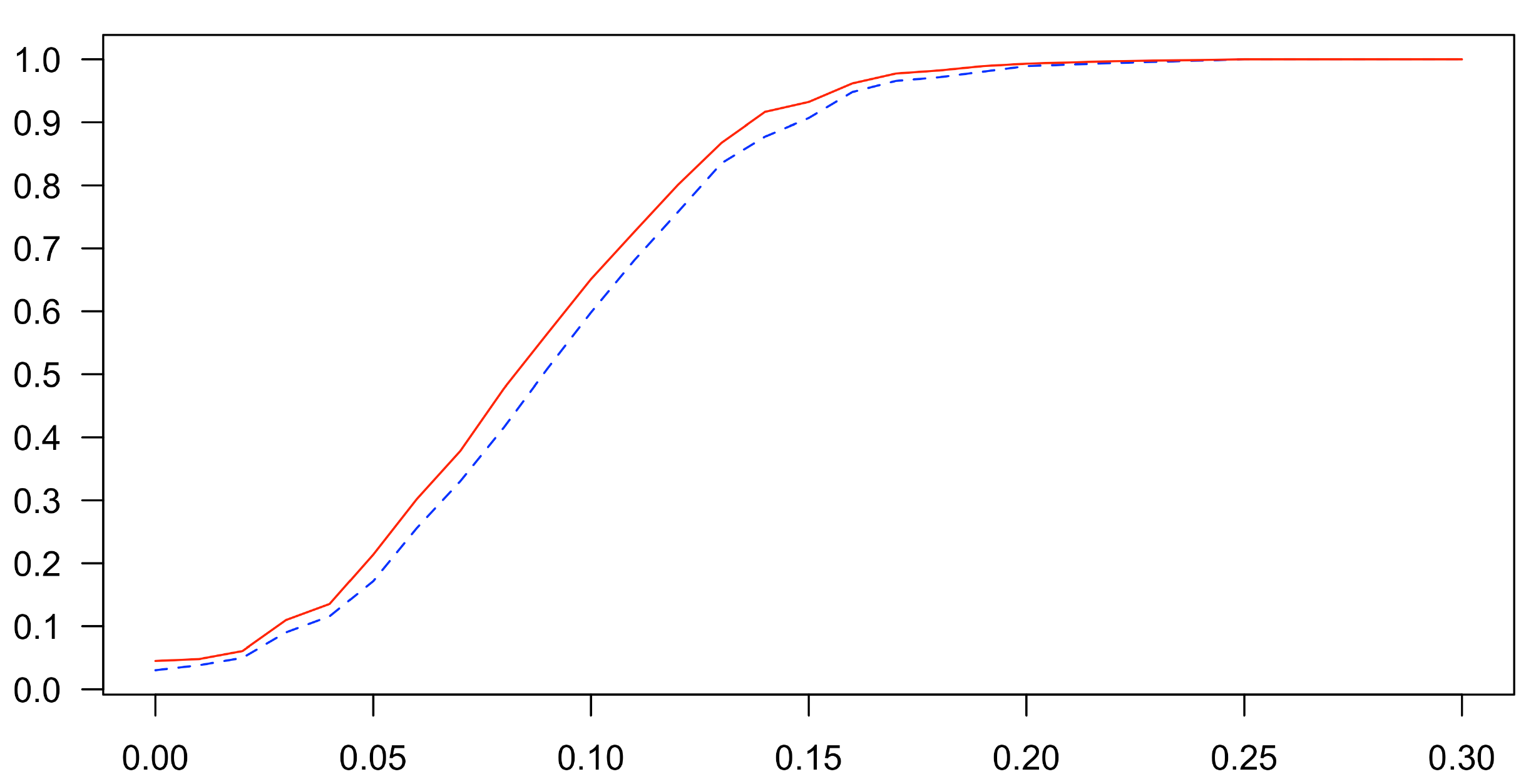}
\caption{Rejection probabilities at 1\% (upper panel) and 5\% level (lower panel) of the
null hypothesis as a function of $r = \protect\phi(h_{0})$ for the SQLR
(solid red line) and for the bootstrap SQLR (dashed blue line) for NPQIV.}
\label{fig:SQLR-NPQIV}
\end{figure}

\textbf{Local power}. Figure \ref{fig:SQLR-NPQIV}
shows the rejection probabilities at 5\% (lower panel) and 1\% (upper panel) level of the null hypothesis as a
function of $r$ where $r\colon \phi (h_{0})=r$ for the SQLR (solid red line) and the bootstrap SQLR (dashed blue line) with multinomial weights. To save space we only report the local power results corresponding to the case of P-Spline(3,2) for $q^{k(n)}$, Pol(10) for $p^{J_{n}} $ and $\lambda
_{n}= 10^{-5}$ in Table \ref{tab:NPQIV-SQLR}. We employ 500 bootstrap evaluations per MC replication,
and lower the number of MC repetitions to 1,000 to ease the computational
burden. We note that since our functional $\phi (h)=h(0)$ is
estimated at a slower than root-$n$ rate, the deviations considered for $r$
which are in the range of $[0,8/\sqrt{n}]$ are indeed
\textquotedblleft small\textquotedblright.
We can see from the figure that the
bootstrap SQLR performance is similar to its non-bootstrapped counterpart.
We expect that the performance will improve if we increase number of
bootstrap runs. (We also run simulation studies corresponding to the case of Pol(4) for $q^{k(n)}$, Pol(7) for $p^{J_{n}} $ and $\lambda _{n}=2\times 10^{-4}$ in Table \ref{tab:NPQIV-SQLR}, and the local power patterns are similar to the ones reported here.)

\medskip

\noindent \textbf{NPIV model: sieve variance estimators for an irregular
linear functional.} We now consider the NPIV model: $%
Y_{1}=h_{0}(Y_{2})+0.76U=2\sin (\pi Y_{2})+0.76U$, with $U=U^{\ast }$ so the
identifying condition of NPIV holds: $E[U|X]=0$. The parameter of interest
is $\phi (h_{0})=h_{0}(0)$, and the null hypothesis is $H_{0}\colon \phi
(h_{0})=0$. We focus on the finite sample performance of the sieve variance
estimators for irregular linear functionals. We compute two sieve variance
estimators:%
\begin{equation*}
\widehat{V}_{1}=q^{k(n)}(0)^{\prime }\widehat{D}_{n}^{-1}\widehat{\mho }_{n}%
\widehat{D}_{n}^{-1}q^{k(n)}(0)\text{\quad and\quad }\widehat{V}%
_{2}=q^{k(n)}(0)^{\prime }\widehat{D}_{n}^{-1}\widehat{\Omega }_{n}\widehat{D%
}_{n}^{-1}q^{k(n)}(0),
\end{equation*}%
where $\widehat{D}_{n}=n^{-1}\left( \widehat{C}_{n}(P^{\prime }P)^{-}%
\widehat{C}_{n}^{\prime }\right) $, $\widehat{C}_{n}\equiv
\sum_{i=1}^{n}q^{k(n)}(Y_{2i})p^{J_{n}}(X_{i})^{\prime }$, $\widehat{\mho }%
_{n}$ is given in equation (\ref{2sls-vhat}), and $\widehat{\Omega }_{n}=%
\frac{1}{n}\widehat{C}_{n}(P^{\prime }P)^{-}\left(
\sum_{i=1}^{n}p^{J_{n}}(X_{i})\widehat{\Sigma }_{0}(X_{i})p^{J_{n}}(X_{i})^{%
\prime }\right) (P^{\prime }P)^{-}\widehat{C}_{n}^{\prime }$ with $\widehat{U%
}_{j}=Y_{1j}-\widehat{h}(Y_{2j})$ and $\widehat{\Sigma }_{0}(x)=\left(
\sum_{j=1}^{n}\widehat{U}_{j}^{2}p^{J_{n}}(X_{j})^{\prime }\right)
(P^{\prime }P)^{-}p^{J_{n}}(x)$. (See Theorem \ref{thm:VE2} in Appendix \ref%
{app:appB} for the definition and consistency of $\widehat{V}_{2}$ as
another sieve variance estimator for any plug-in PSMD $\phi (\widehat{\alpha
})$.)

\begin{table}[h]
\centering
\caption{Relative performance of $\hat{V}_{1}$ and $%
	\hat{V}_{2}$: $Med_{MC} \left[\left|\frac{\hat{V}_{j}}{||v^{%
			\ast}_{n}||^{2}_{sd}}-1\right| \right]$, and Nominal size and MC rejection
	frequencies for t tests $\hat{t}_{j}$ for $j=1,2$ for a linear functional of
	NPIV.}
\begin{tabular}{cccccccc}
\hline\hline
&  & \multicolumn{2}{c}{$Med_{MC}$} & \multicolumn{2}{c}{5\%} &
\multicolumn{2}{c}{10\%} \\ \hline
$q^{k(n)}$ & $p^{J_{n}}$ & $\widehat{V}_{1}$ & $\widehat{V}_{2}$ & $\widehat{%
V}_{1}$ & $\widehat{V}_{2}$ & $\widehat{V}_{1}$ & $\widehat{V}_{2}$ \\ \hline
\multirow{4}{*}{Pol(4)} & Pol(6) & 0.0946 & 0.0937 & 0.0512 & 0.0514 & 0.0980
& 0.0974 \\
& Pol(10) & 0.0922 & 0.0920 & 0.0536 & 0.0532 & 0.0992 & 0.0990 \\
& Pol(12) & 0.0918 & 0.0917 & 0.0538 & 0.0532 & 0.1002 & 0.0998 \\
& Pol(16) & 0.0911 & 0.0912 & 0.0540 & 0.0538 & 0.1000 & 0.0998 \\ \hline
\multirow{4}{*}{Pol(4)} & P-Spline(3,2) & 0.0939 & 0.0942 & 0.051 & 0.0516 &
0.0984 & 0.0986 \\
& P-Spline(3,5) & 0.0939 & 0.0920 & 0.053 & 0.0532 & 0.0990 & 0.0984 \\
& P-Spline(3,11) & 0.0923 & 0.0925 & 0.055 & 0.0548 & 0.1014 & 0.1014 \\
& P-Spline(3,17) & 0.0922 & 0.0917 & 0.0542 & 0.0538 & 0.100 & 0.1008 \\
\hline
\multirow{4}{*}{P-Spline(3,2)} & Pol(12) & 0.0938 & 0.0930 & 0.0572 & 0.0564 &
0.1082 & 0.1074 \\
& Pol(16) & 0.0936 & 0.0936 & 0.0582 & 0.0578 & 0.1082 & 0.1082 \\
& Pol(18) & 0.0936 & 0.0935 & 0.0580 & 0.0578 & 0.1088 & 0.1086 \\
& Pol(20) & 0.0936 & 0.0937 & 0.0580 & 0.0574 & 0.1086 & 0.1092 \\ \hline
\multirow{4}{*}{P-Spline(3,2)} & P-Spline(3,2) & 0.1106 & 0.1116 & 0.0606 &
0.0598 & 0.1130 & 0.1120 \\
& P-Spline(3,5) & 0.1019 & 0.1023 & 0.0584 & 0.0574 & 0.1122 & 0.1116 \\
& P-Spline(3,11) & 0.0961 & 0.0960 & 0.0572 & 0.0566 & 0.1100 & 0.1094 \\
& P-Spline(3,17) & 0.0949 & 0.0944 & 0.0570 & 0.0566 & 0.1082 & 0.1080 \\
\hline
\multirow{4}{*}{P-Spline(3,2)} & P-Spline(5,3) & 0.1007 & 0.0998 & 0.0586 &
0.0576 & 0.1102 & 0.1088 \\
& P-Spline(5,6) & 0.1011 & 0.1009 & 0.0586 & 0.0578 & 0.1100 & 0.1092 \\
& P-Spline(5,12) & 0.1007 & 0.1009 & 0.0580 & 0.0572 & 0.1110 & 0.1096 \\
& P-Spline(5,18) & 0.1009 & 0.1010 & 0.0580 & 0.0570 & 0.1106 & 0.1092 \\
\hline
\end{tabular}%
\label{tab:ve}
\end{table}

Table \ref{tab:ve} reports the results for different choices of bases for $%
q^{k(n)}$ and $p^{J_{n}}$, and for different values of $k(n)$ and $J_{n}$;
in all cases we use a very small $\lambda _{n}=10^{-5}$. This table shows $%
Med_{MC}\left[ \left\vert \frac{\widehat{V}_{j}}{||v_{n}^{\ast }||_{sd}^{2}}%
-1\right\vert \right] $ for $j=1,2$, where $||v_{n}^{\ast }||_{sd}$ is
computed using the MC variance of $\sqrt{n}\widehat{h}_{n}(0)$ and $%
Med_{MC}[\cdot ]$ is the MC median. It also shows the nominal size and MC
rejection frequencies of the two sieve t tests $\widehat{t}_{j}=\sqrt{n}%
\frac{\widehat{h}_{n}(0)-0}{\sqrt{\widehat{V}_{j}}}$ for $j=1,2$.

\begin{figure}[tbp]
\centering
\includegraphics[height=2.5in,width=4in]{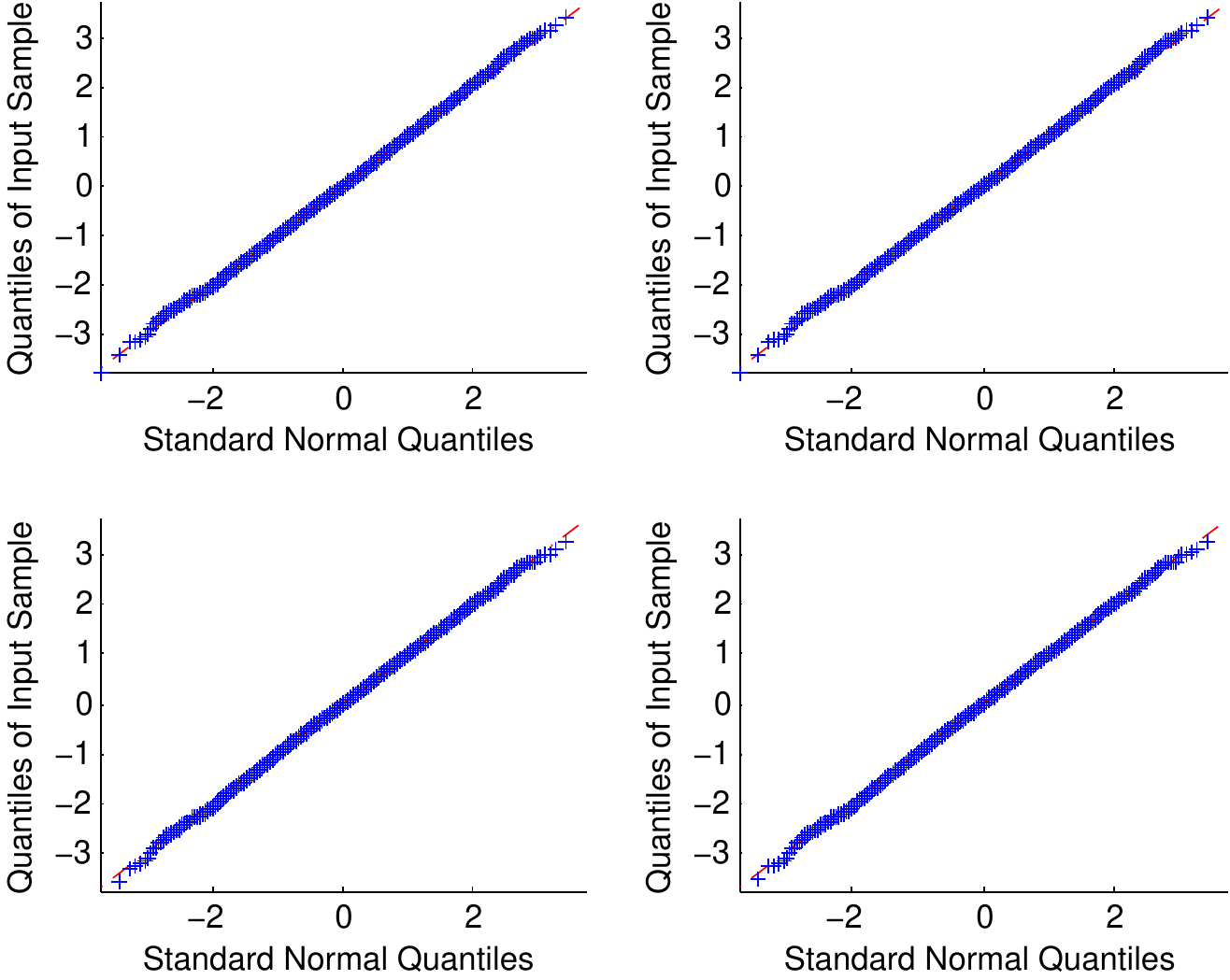}
\caption{QQ-Plot for t tests $\hat{t}_{j}$ for $j=1,2$ for a linear functional (first row) and a nonlinear functional (second row) of NPIV, with
$q^{k(n)}=Pol(4)$ and $p^{J_{n}}=Pol(16)$.}
\label{fig:QQplot-MC3-MC4}
\end{figure}

We note that the two sieve variance estimators have almost identical
performance and the associated sieve t tests have good rejection
probabilities. These results are fairly robust to different choices of basis
for $q^{k(n)}$ and $p^{J_{n}}$ and different values of $k(n)$ and $J_{n}$ as
long as $J_{n}>k(n)+1$. Figure \ref{fig:QQplot-MC3-MC4} (first row) shows
the QQ-Plot for the sieve t tests $\widehat{t}_{j}=\sqrt{n}\frac{\widehat{h}%
_{n}(0)-0}{\sqrt{\widehat{V}_{j}}}$ under the null for $j=1,2$ for the case
Pol(4)-Pol(16) in the table; the right panel in the first row corresponds to
$\hat{t}_{1}$ and the left panel in the first row to $\hat{t}_{2}$. Both
sieve t tests are almost identical to each other and to the standard normal.

\medskip

\noindent \textbf{NPIV model: sieve variance estimators for an irregular
\emph{nonlinear} functional.} This case is identical to the previous one for
the NPIV model, except that the functional of interest is $\phi (h_{0})=\exp
\{h_{0}(0)\}$, and the null hypothesis is $H_{0}\colon \phi (h_{0})=1$. This
choice of $\phi $ allows us to evaluate the finite sample performance of
sieve t statistics for a nonlinear functional.

\begin{table}[h]
\centering
\caption{Relative performance of $\hat{V}_{1}$ and $%
	\hat{V}_{2}$: $Med_{MC} \left[\left|\frac{\hat{V}_{j}}{||v^{%
			\ast}_{n}||^{2}_{sd}}-1\right| \right]$, and Nominal size and MC rejection
	frequencies for t tests $\hat{t}_{j}$ for $j=1,2$ for a nonlinear functional
	of NPIV.}
\begin{tabular}{cccccccc}
\hline\hline
&  & \multicolumn{2}{c}{$Med_{MC}$} & \multicolumn{2}{c}{5\%} &
\multicolumn{2}{c}{10\%} \\ \hline
$q^{k(n)}$ & $p^{J_{n}}$ & $\widehat{V}_{1}$ & $\widehat{V}_{2}$ & $\widehat{%
V}_{1}$ & $\widehat{V}_{2}$ & $\widehat{V}_{1}$ & $\widehat{V}_{2}$ \\ \hline
\multirow{4}{*}{Pol(4)} & Pol(6) & 0.0990 & 0.0985 & 0.0528 & 0.0530 & 0.0982
& 0.0988 \\
& Pol(10) & 0.0971 & 0.0958 & 0.0524 & 0.0522 & 0.1014 & 0.1012 \\
& Pol(12) & 0.0967 & 0.0959 & 0.0526 & 0.0526 & 0.1020 & 0.1018 \\
& Pol(16) & 0.0961 & 0.0958 & 0.0524 & 0.0528 & 0.1018 & 0.1014 \\ \hline
\multirow{4}{*}{Pol(4)} & P-Spline(3,2) & 0.0996 & 0.0983 & 0.0534 & 0.0530 &
0.0978 & 0.0976 \\
& P-Spline(3,5) & 0.0982 & 0.0969 & 0.0538 & 0.0542 & 0.0990 & 0.0992 \\
& P-Spline(3,11) & 0.0985 & 0.0984 & 0.0554 & 0.0552 & 0.1014 & 0.1010 \\
& P-Spline(3,17) & 0.0982 & 0.0978 & 0.0544 & 0.0546 & 0.1010 & 0.1008 \\
\hline
\multirow{4}{*}{P-Spline(3,2)} & Pol(12) & 0.1011 & 0.1009 & 0.0580 & 0.0568 &
0.1120 & 0.1122 \\
& Pol(16) & 0.1014 & 0.1005 & 0.0588 & 0.0574 & 0.1128 & 0.1126 \\
& Pol(18) & 0.1014 & 0.1007 & 0.0582 & 0.0568 & 0.1130 & 0.1122 \\
& Pol(20) & 0.1015 & 0.1006 & 0.0580 & 0.0568 & 0.1138 & 0.1128 \\ \hline
\multirow{4}{*}{P-Spline(3,2)} & P-Spline(3,2) & 0.1191 & 0.1192 & 0.0620 &
0.0612 & 0.1132 & 0.1120 \\
& P-Spline(3,5) & 0.1090 & 0.1103 & 0.0596 & 0.0594 & 0.1140 & 0.1134 \\
& P-Spline(3,11) & 0.1028 & 0.1032 & 0.0582 & 0.0572 & 0.1130 & 0.1126 \\
& P-Spline(3,17) & 0.1029 & 0.1029 & 0.0588 & 0.0580 & 0.1124 & 0.1112 \\
\hline
\multirow{4}{*}{P-Spline(3,2)} & P-Spline(5,3) & 0.1059 & 0.1064 & 0.0594 &
0.0592 & 0.1114 & 0.1104 \\
& P-Spline(5,6) & 0.1066 & 0.1076 & 0.0598 & 0.0586 & 0.1124 & 0.1118 \\
& P-Spline(5,12) & 0.1071 & 0.1079 & 0.0594 & 0.0586 & 0.1126 & 0.1120 \\
& P-Spline(5,18) & 0.1069 & 0.1079 & 0.0594 & 0.0586 & 0.1122 & 0.1120 \\
\hline
\end{tabular}%
\label{tab:ve1}
\end{table}

Table \ref{tab:ve1} shows $Med_{MC}$ and rejection probabilities for this
nonlinear case. By comparing the results with those in Table \ref{tab:ve} we
note that the results are very similar in both cases. Figure \ref%
{fig:QQplot-MC3-MC4} (second row) shows the QQ-Plot for the two sieve t
tests for the non-linear case; the right panel in the second row corresponds
to $\hat{t}_{1}$ whereas the left panel in the second row corresponds to $%
\hat{t}_{2}$. These results suggest that our sieve t tests perform equally
well for both functionals.

Finally we wish to point out that we have tried other bases such as Hermite
polynomials and cosine series and even larger $J_{n}$ in these two NPIV MC
studies, the results are all similar to the ones reported here and hence are
not presented due to the lack of space.

\subsection{An Empirical Application}

\label{sec:application}

We compute SQLR based confidence bands for nonparametric quantile IV Engel
curves using the British FES data set from \cite{BCK_Emetrica07}:
\begin{equation*}
E[1\{Y_{1,i}\leq h_{0}(Y_{2,i})\}\mid X_{i}]=0.5,
\end{equation*}%
where $Y_{1,i}$ is the budget share of the $i-$th household on a particular
non-durable goods, say food-in consumption; $Y_{2,i}$ is the log-total
expenditure of the household, which is endogenous, and hence we use $X_{i}$,
the gross earnings of the head of the household, to instrument it. We work
with the \textquotedblleft no kids\textquotedblright\ sub-sample of the data
set, which consists of $n=628$ observations. \cite{BCK_Emetrica07} estimated
NPIV Engel curves using this data set. But, as explained by \cite%
{Koenker_2005} and others, quantile Engel curves are more informative.

We estimate $h_{0}(\cdot )$ for food-in quantile Engel curve via the
optimally weighted PSMD procedure with $\widehat{\Sigma }=\Sigma _{0}=0.25$,
using a polynomial spline (P-spline) sieve $\mathcal{H}_{k(n)}$ with $k(n)=4$%
, $Pen(h)=||h||_{L^{2}}^{2}+||\nabla h||_{L^{2}}^{2}$ with $\lambda
_{n}=0.0005$, and a Hermite polynomial LS basis $p^{J_{n}}(X)$ with $J_{n}=6$%
. We also considered other bases such as P-splines as $p^{J_{n}}(X)$ and
results remained essentially the same. See \cite{CP_WP07a} for PSMD
estimates of NPQIV Engel curves for other non-durable goods.

\begin{figure}[tbp]
\centering
\includegraphics[height=3in,width=4.5in]{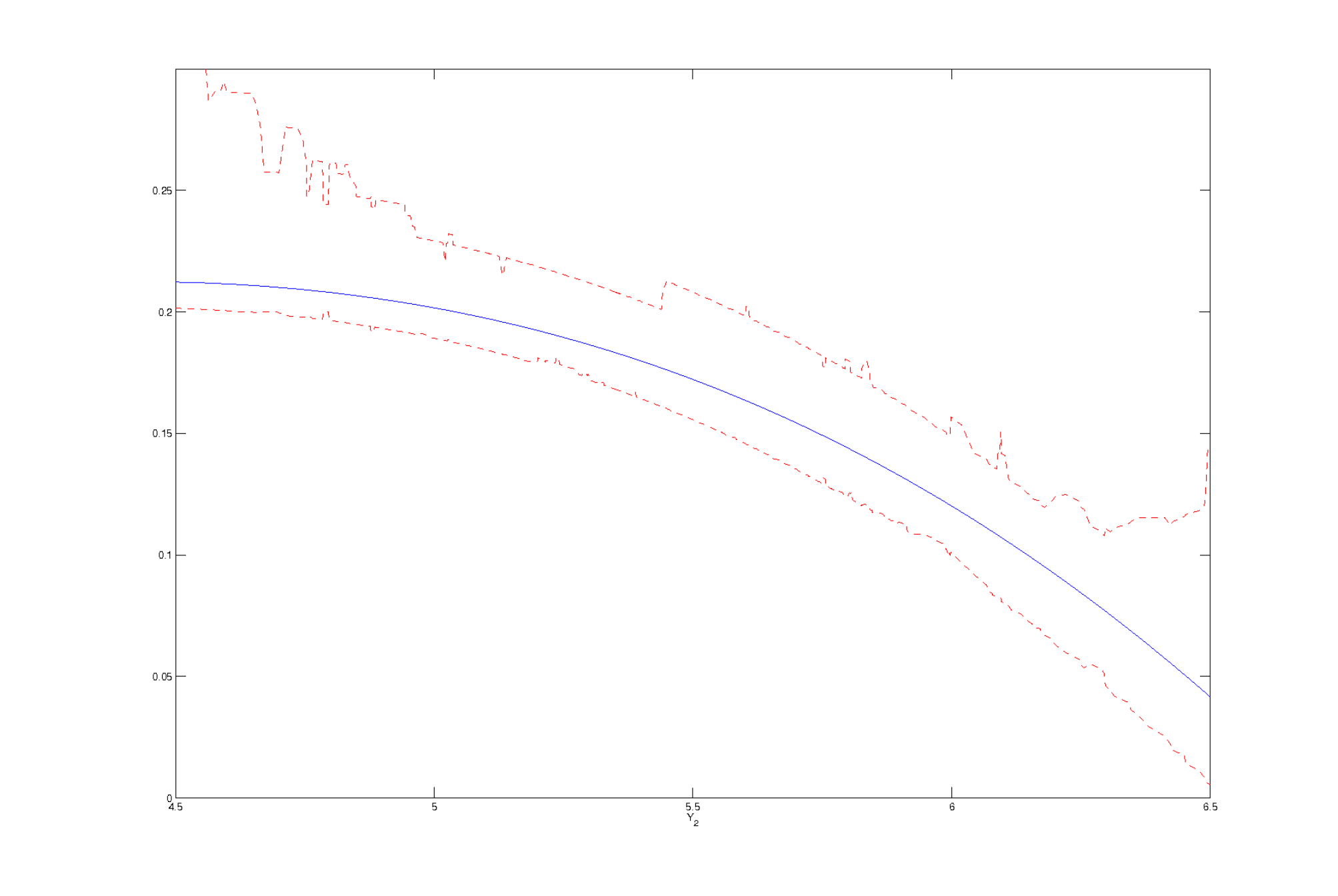}
\caption{{\protect\footnotesize {PSMD Estimate of the NPQIV food-in Engel
curve (blue solid line), with the 95\% pointwise confidence bands (red dash
lines).}}}
\label{fig:EC}
\end{figure}

We use the fact that the optimally weighted SQLR of testing $\phi
(h)=h(y_{2})$ (for any fixed $y_{2}$) is asymptotically $\chi _{1}^{2}$ to
construct pointwise confidence bands. That is, for each $y_{2}$ in the
sample we construct a grid, $(r_{i} )_{i=1}^{30}$. For each $i={1,...,30}$,
we compute the value of the SQLR test statistic under $h(y_{2})=r_{i}$ for $%
(r_{i})_{i=1}^{30}$. We then, take the smallest interval that included all
points $r_{i}$ that yield a corresponding value of the SQLR test below the
95\% percentile of $\chi _{1}^{2}$.\footnote{%
The grid $(r_{i})_{i=1}^{30}$ was constructed to have $r_{15}=\widehat{h}%
_{n}(y_{2})$, for all $i\leq 15$, $r_{i+1}\leq r_{i}\leq r_{15}$ decreasing
in steps of length $0.002$ (approx) and for all $i\geq 15$, $r_{i+1}\geq
r_{i}\geq r_{15}$ increasing in steps of length $0.008$ (approx); finally,
the extremes, $r_{1}$ and $r_{30}$, were chosen so the SQLR test at those
points was above the 95\% percentile of $\chi _{1}^{2}$. We tried different
lengths and step sizes and the results remain qualitatively unchanged. For
some observations, which only account for less than 4\% of the sample, the
confidence interval was degenerate at a point; this result was due to
numerical approximation issues, and these observations were excluded from the reported
results.} Figure \ref{fig:EC} presents the results, where the solid blue
line is the point estimate and the red dashed lines are the 95\% pointwise
confidence bands. We can see that the confidence bands get wider towards the
extremes of the sample, but are tighter in the middle.

To test whether the quantile IV Engel curve for food-in is linear or not,
one can test whether $\phi (h_{0})\equiv \int \left\vert \nabla
^{2}h(y_{2})\right\vert ^{2}w(y_{2})dy_{2}=0$ using our SQLR test. Let $%
w(\cdot )=(\sigma _{Y_{2}})^{-1}\exp \left( -\frac{1}{2}(\sigma
_{Y_{2}}^{-1}(\cdot -\mu _{Y_{2}}))^{2}\right) 1\{t_{0.01}\leq \cdot \leq
t_{0.99}\}$ where $\mu _{Y_{2}}$, $\sigma _{Y_{2}}$, $t_{0.01}$ and $%
t_{0.99} $ are the sample mean, standard deviation and the 1\% and 99\%
quantiles of $Y_{2}$. The value of the SQLR is (approx.) 38 and the p-value
is smaller than 0.0001, and hence we reject the null hypothesis of linearity.%
\footnote{%
We use the standard Riemann sum with 1000 terms to compute the integral. We
also considered other choices of $w$ such that $w(\cdot )=1\{t_{0.25}\leq
\cdot \leq t_{0.75}\}$ and $w(\cdot )=1\{t_{0.01}\leq \cdot \leq t_{0.99}\}$%
. Although the numerical value of the SQLR test changes, all produce
p-values below 0.0001.}

\section{Conclusion}

\label{sec:conclusion}

In this paper, we provide unified asymptotic theories for PSMD based
inferences on possibly irregular parameters $\phi (\alpha _{0})$ of the
general semi/nonparametric conditional moment restrictions $E[\rho
(Y,X;\alpha _{0})|X]=0$. Under regularity conditions that allow for any
consistent nonparametric estimator of the conditional mean function $%
m(X,\alpha )\equiv E[\rho (Y,X;\alpha )|X]$, we establish the asymptotic
normality of the plug-in PSMD estimator $\phi (\widehat{\alpha }_{n})$ of $%
\phi (\alpha _{0})$, as well as the asymptotically tight distribution of a
possibly non-optimally weighted SQLR statistic under the null hypothesis of $%
\phi (\alpha _{0})=\phi _{0}$. As a simple yet useful by-product, we
immediately obtain that an optimally weighted SQLR statistic is
asymptotically chi-square distributed under the null hypothesis. For
(pointwise) smooth residuals $\rho (Z;\alpha )$ (in $\alpha $), we propose
several simple consistent sieve variance estimators for $\phi (\widehat{%
\alpha }_{n})$ (in the text and in online Appendix \ref{app:appB}), and
establish the asymptotic chi-square distribution of sieve Wald statistics.
We also establish local power properties of SQLR and sieve Wald tests in
Appendix \ref{app:appA}. Under conditions that are virtually the same as
those for the limiting distributions of the original-sample sieve Wald and
SQLR statistics, we establish the consistency of the generalized residual
bootstrap sieve Wald and SQLR statistics. All these results are valid
regardless of whether $\phi (\alpha _{0})$ is regular or not. While SQLR and
bootstrap SQLR are useful for models with (pointwise) non-smooth $\rho
(Z;\alpha )$, sieve Wald statistic is computationally attractive for models
with smooth $\rho (Z;\alpha )$. Monte Carlo studies and an empirical
illustration of a nonparametric quantile IV regression demonstrate the good
finite sample performance of our inference procedures.

This paper assumes that the semi/nonparametric conditional moment
restrictions $E[\rho (Y,X;\alpha _{0})|X]=0$ uniquely identifies the unknown
true parameter value $\alpha _{0}\equiv (\theta _{0}^{\prime },h_{0})$, and
conduct inferences that are robust to whether a possibly nonlinear
functional of $\alpha _{0}$ is root-$n$ estimable or not. Recently, for the
NPIV model $E[Y_{1}-h_{0}(Y_{2})|X]=0$ without assuming point identification
of $h_{0}$, \cite{Santos_JOE} proposed a root-$n$ asymptotically normal
estimation of a regular linear functional of $h_{0}$ and \cite{Santos_ECMA}
considered Bierens' type test of the NPIV. \cite{CPT_WP12} is
extending the SQLR inference procedure to allow for partial identification
of the general model $E[\rho (Y,X;\alpha _{0})|X]=0$.

\baselineskip=15pt

\bigskip




\setstretch{1}

\bibliography{mybib_FUNC}

\appendix
\newpage

\section{Additional Results and Sufficient Conditions}

\label{app:appA}

\renewcommand{\theassumption}{\thesection.\arabic{assumption}}

Appendix \ref{app:appA} consists of several subsections. Subsection \ref%
{sec:phi-further} presents additional results on (sieve) Riesz
representation of the functional of interest. Subsection \ref%
{sec:bconsistency} derives the convergence rates of the bootstrap PSMD
estimator. Subsection \ref{subsec-A5} presents asymptotic properties under
local alternatives of the SQLR and the sieve Wald tests, and of their
bootstrap versions. Subsection \ref{subsec-A6} provides some inference
results for functionals of increasing dimension. Subsection \ref{subsec-A1}
provides some low level sufficient conditions for the high level LQA
Assumption \ref{ass:LAQ}(i) and the bootstrap LQA Assumption \ref{ass:LAQ_B}%
(i) with series LS estimated conditional mean functions $m(\cdot ,\alpha )$.
Subsection \ref{subsec-A2} states useful lemmas with series LS estimated
conditional mean functions $m(\cdot ,\alpha )$. See online supplemental
Appendix \ref{app:appC} for the proofs of all the results in this Appendix.

\subsection{Additional discussion on (sieve) Riesz representation}

\label{sec:phi-further}

The discussion in Subsection \ref{sec:phi} on Riesz representation seems to
depend on the weighting matrix $\Sigma $, but, under Assumption \ref%
{ass:sieve}(iv), we have $||\cdot ||\asymp ||\cdot ||_{0}$, (i.e., the norm $%
||\cdot ||$ (using $\Sigma $) is equivalent to the norm $||\cdot ||_{0}$
(using $\Sigma _{0}$) defined in (\ref{fmetric0})), and the space $\mathbf{V}
$ (or $\overline{\mathbf{V}}$) under $||\cdot ||$ is equivalent to that
under $||\cdot ||_{0}$. Therefore, under Assumption \ref{ass:sieve}(iv), $%
\frac{d\phi (\alpha _{0})}{d\alpha }[\cdot ]$ is bounded on $(\mathbf{V}%
,||\cdot ||)$ iff $\frac{d\phi (\alpha _{0})}{d\alpha }[\cdot ]$ is \textit{%
bounded} on $(\mathbf{V},||\cdot ||_{0})$, i.e.,
\begin{equation*}
\sup_{v\in \mathbf{V},v\neq 0}\frac{\left\vert \frac{d\phi (\alpha _{0})}{%
d\alpha }\left[ v\right] \right\vert }{\left\Vert v\right\Vert _{0}}<\infty ,
\end{equation*}%
and in this case we call $\phi \left( \cdot \right) $ \textit{regular} at $%
\alpha _{0}$. Likewise, $\frac{d\phi (\alpha _{0})}{d\alpha }[\cdot ]$ is
unbounded on $(\mathbf{V},||\cdot ||)$ iff $\frac{d\phi (\alpha _{0})}{%
d\alpha }[\cdot ]$ is \textit{unbounded} on $(\mathbf{V},||\cdot ||_{0})$
i.e., $\sup_{v\in \mathbf{V},v\neq 0}\left\{ \left\vert \frac{d\phi (\alpha
_{0})}{d\alpha }\left[ v\right] \right\vert /\left\Vert v\right\Vert
_{0}\right\} =\infty $, and in this case we call $\phi \left( \cdot \right) $
\textit{irregular} at $\alpha _{0}$.

For the specific Euclidean functional $\phi \left( \alpha \right) =\lambda
^{\prime }\theta $ of the semi/nonparametric model (\ref{semi00}), following
the proof in appendix E of \cite{AC_Emetrica03}, it is easy to see the
equivalence between $\sup_{v\in \overline{\mathbf{V}},v\neq 0}\left\{
\left\vert \frac{d\phi (\alpha _{0})}{d\alpha }\left[ v\right] \right\vert
/\left\Vert v\right\Vert _{0}\right\} <\infty $ (or $=\infty $) and the%
\textit{\ semiparametric efficiency bound} of $\theta _{0}$ being \textit{%
non-singular} (or \textit{singular}). Previously \cite{KhanTamer2010} and
\cite{GPowell2012} call $\theta _{0}$ irregular when the semiparametric
efficiency bound of $\theta _{0}$ is singular. This motivates us to call a
general functional $\phi \left( \cdot \right) $ \textit{irregular} at $%
\alpha _{0}$ whenever $\frac{d\phi (\alpha _{0})}{d\alpha }[\cdot ]$ is
unbounded on $(\mathbf{V},||\cdot ||_{0})$.

It is clear that $\phi \left( \cdot \right) $ being regular at $\alpha _{0}$
(i.e., $\frac{d\phi (\alpha _{0})}{d\alpha }[\cdot ]$ being bounded on $(%
\mathbf{V},||\cdot ||_{0})$) is a necessary condition for the root-$n$ rate
of convergence of $\phi (\widehat{\alpha }_{n})-\phi (\alpha _{0})$.
Unfortunately for a complicated semi/nonparametric model (\ref{semi00}), it
is difficult to compute $\sup_{v\in \mathbf{V},v\neq 0}\left\{ \left\vert
\frac{d\phi (\alpha _{0})}{d\alpha }\left[ v\right] \right\vert /\left\Vert
v\right\Vert _{0}\right\} $ explicitly; and hence difficult to verify its
root-$n$ estimableness. But, we could always compute $\left\Vert v_{n}^{\ast
}\right\Vert =\sup_{v\in \overline{\mathbf{V}}_{k(n)},v\neq 0}\left\{
\left\vert \frac{d\phi (\alpha _{0})}{d\alpha }\left[ v\right] \right\vert
/\left\Vert v\right\Vert \right\} $ in closed form (by Lemma \ref%
{lem:P-sieve}) and then check whether $\lim_{k(n)\rightarrow \infty
}\left\Vert v_{n}^{\ast }\right\Vert ^{2}<\infty $ or not (by Lemma \ref%
{lem:sieve-Riesz-property}).

\begin{remark}
\label{remark-norm} For a semi/nonparametric model (\ref{semi00}) with $%
\alpha _{0}=(\theta _{0}^{\prime },h_{0})$, it is convenient to rewrite $%
D_{n}$ and its inverse in Lemma \ref{lem:P-sieve} as
\begin{equation*}
D_{n}\equiv \left(
\begin{array}{cc}
I_{11} & I_{n,12} \\
I_{n,12}^{\prime } & I_{n,22}%
\end{array}%
\right) \text{\quad and\quad }D_{n}^{-1}=\left(
\begin{array}{cc}
I_{n}^{11} & -I_{11}^{-1}I_{n,12}I_{n}^{22} \\
-I_{n,22}^{-1}I_{n,12}^{\prime }I_{n}^{11} & I_{n}^{22}%
\end{array}%
\right) ,
\end{equation*}%
$I_{11}=E\left[ \left( \frac{dm(X,\alpha _{0})}{d\theta ^{\prime }}\right)
^{\prime }\Sigma (X)^{-1}\frac{dm(X,\alpha _{0})}{d\theta ^{\prime }}\right]
$, $I_{n,22}=E\left[ \left( \frac{dm(X,\alpha _{0})}{dh}[\psi ^{k(n)}(\cdot
)^{\prime }]\right) ^{\prime }\Sigma (X)^{-1}\left( \frac{dm(X,\alpha _{0})}{%
dh}[\psi ^{k(n)}(\cdot )^{\prime }]\right) \right] $, $I_{n,12}=E\left[
\left( \frac{dm(X,\alpha _{0})}{d\theta ^{\prime }}\right) ^{\prime }\Sigma
(X)^{-1}\left( \frac{dm(X,\alpha _{0})}{dh}[\psi ^{k(n)}(\cdot )^{\prime
}]\right) \right] $, $I_{n}^{11}=\left(
I_{11}-I_{n,12}I_{n,22}^{-1}I_{n,21}^{\prime }\right) ^{-1}$ and $%
I_{n}^{22}=\left( I_{n,22}-I_{n,21}^{\prime }I_{11}^{-1}I_{n,12}\right)
^{-1} $. For instance, for the Euclidean functional $\phi (\alpha )=\lambda
^{\prime }\theta $, we have $\digamma _{n}=(\lambda ^{\prime },\mathbf{0}%
_{k(n)}^{\prime })^{\prime }$ with $\mathbf{0}_{k(n)}^{\prime }=\left[
0,...,0\right] _{1\times k(n)}$, and hence $v_{n}^{\ast }=(v_{\theta
,n}^{\ast \prime },\psi ^{k(n)}(\cdot )^{\prime }\beta _{n}^{\ast })^{\prime
}\in \overline{\mathbf{V}}_{k(n)}$ with $v_{\theta ,n}^{\ast
}=I_{n}^{11}\lambda $, $\beta _{n}^{\ast }=-I_{n,22}^{-1}I_{n,21}^{\prime
}v_{\theta ,n}^{\ast }$, and $\left\Vert v_{n}^{\ast }\right\Vert
^{2}=\lambda ^{\prime }I_{n}^{11}\lambda $. Thus the functional $\phi
(\alpha )=\lambda ^{\prime }\theta $ is regular iff $\lim_{k(n)\rightarrow
\infty }\lambda ^{\prime }I_{n}^{11}\lambda <\infty $; in this case,%
\begin{equation*}
\lim_{k(n)\rightarrow \infty }\left\Vert v_{n}^{\ast }\right\Vert
^{2}=\lim_{k(n)\rightarrow \infty }\lambda ^{\prime }I_{n}^{11}\lambda
=\lambda ^{\prime }\mathcal{I}_{\ast }^{-1}\lambda =\left\Vert v^{\ast
}\right\Vert ^{2},
\end{equation*}%
where
\begin{equation}
\mathcal{I}_{\ast }=\inf_{\mathbf{w}}E\left[ \left\Vert \Sigma (X)^{-\frac{1%
}{2}}\left( \frac{dm(X,\alpha _{0})}{d\theta ^{\prime }}-\frac{dm(X,\alpha
_{0})}{dh}[\mathbf{w}]\right) \right\Vert _{e}^{2}\right] ,  \label{I*}
\end{equation}%
and $v^{\ast }=(v_{\theta }^{\ast \prime },v_{h}^{\ast }\left( \cdot \right)
)^{\prime }\in \overline{\mathbf{V}}$ where $v_{\theta }^{\ast }\equiv
\mathcal{I}_{\ast }^{-1}\lambda $, $v_{h}^{\ast }\equiv -\mathbf{w}^{\ast
}\times v_{\theta }^{\ast }$, and $\mathbf{w}^{\ast }$ solves (\ref{I*}).
That is, $v^{\ast }=(v_{\theta }^{\ast \prime },v_{h}^{\ast }\left( \cdot
\right) )^{\prime }$ becomes the Riesz representer for $\phi (\alpha
)=\lambda ^{\prime }\theta $ previously computed in \cite{AC_Emetrica03} and
\cite{CP_WP07a}. Moreover, if $\Sigma (X)=\Sigma _{0}(X)$, then $\mathcal{I}%
_{\ast }$ becomes the semiparametric efficiency bound for $\theta _{0}$ that
was derived in \cite{CHAMBERLAIN_ECMA92} and \cite{AC_Emetrica03} for the
model (\ref{semi00}). In this setup, Lemma \ref{lem:sieve-Riesz-property}
implies that one could check whether $\theta _{0}$ has non-singular
efficiency bound or not by checking if $\lim_{k(n)\rightarrow \infty
}\lambda ^{\prime }I_{n}^{11}\lambda <\infty $ or not.
\end{remark}

\subsection{Consistency and convergence rate of the bootstrap PSMD
estimators \label{sec:bconsistency}}

In this subsection we establish the consistency and the convergence rate of
the bootstrap PSMD estimator $\widehat{\alpha }_{n}^{B}$ (and the restricted
bootstrap PSMD estimator $\widehat{\alpha }_{n}^{R,B}$) under virtually the
same conditions as those imposed for the consistency and the convergence
rate of the original-sample PSMD estimator $\widehat{\alpha }_{n}$.

The next assumption is needed to control the difference of the bootstrap
criterion function $\widehat{Q}_{n}^{B}(\alpha )$ and the original-sample
criterion function $\widehat{Q}_{n}(\alpha )$. Let $\{\overline{\delta }%
_{m,n}^{\ast }\}_{n=1}^{\infty }$ be a sequence of real valued positive
numbers such that $\overline{\delta }_{m,n}^{\ast }=o(1)$ and $\overline{%
\delta }_{m,n}^{\ast }\geq \delta _{n}$. Let $c_{0}^{\ast }$ and $c^{\ast }$
be finite positive constants.

\begin{assumption}[Bootstrap sample criterion]
\label{ass:rates_B} (i) $\widehat{Q}_{n}^{B}(\widehat{\alpha }_{n})\leq
c_{0}^{\ast }\widehat{Q}_{n}(\widehat{\alpha }_{n})+o_{P_{V^{\infty
}|Z^{\infty }}}(\frac{1}{n})$ wpa1($P_{Z^{\infty }}$); (ii) $\widehat{Q}%
_{n}^{B}(\alpha )\geq c^{\ast }\widehat{Q}_{n}(\alpha )-O_{P_{V^{\infty
}|Z^{\infty }}}((\overline{\delta }_{m,n}^{\ast })^{2})$ uniformly over $%
\mathcal{A}_{k(n)}^{M_{0}}$ wpa1($P_{Z^{\infty }}$); (iii) $\widehat{Q}%
_{n}^{B}(\alpha )\geq c^{\ast }\widehat{Q}_{n}(\alpha )-O_{P_{V^{\infty
}|Z^{\infty }}}(\delta _{n}^{2})$ uniformly over $\mathcal{A}_{osn}$ wpa1($%
P_{Z^{\infty }}$).
\end{assumption}

Assumption \ref{ass:rates_B}(i)(ii) is analogous to Assumption \ref%
{ass:rates} for the original sample, while Assumption \ref{ass:rates_B}(iii)
is analogous to Assumption \ref{ass:weak_equiv}(iii) for the original
sample. Again, when $\widehat{m}^{B}(x,\alpha )$ is the bootstrap series LS
estimator (\ref{mhat_B}) of $m(x,\alpha )$, under virtually the same
sufficient conditions as those in \cite{CP_WP07} and \cite{CP_WP07a} for
their original-sample series LS estimator $\widehat{m}(x,\alpha )$,
Assumption \ref{ass:rates_B} can be verified.\footnote{%
The verification is amounts to follow the proof of Lemma C.2 of \cite%
{CP_WP07} except that the original-sample series LS estimator $\widehat{m}%
(x,\alpha )$ is replaced by its bootstrap version $\widehat{m}^{B}(x,\alpha
) $.}

\begin{lemma}
\label{lem:cons_boot} Let Assumption \ref{ass:rates_B}(i)(ii) and conditions
for Lemma \ref{thm:Thm-ill-suff_pencompact3} hold. Then:
\begin{equation*}
\text{(1)\quad }||\widehat{\alpha }_{n}^{B}-\alpha
_{0}||_{s}=o_{P_{V^{\infty }|Z^{\infty }}}(1)~wpa1(P_{Z^{\infty }})\quad
and\quad Pen\left( \widehat{h}_{n}^{B}\right) =O_{P_{V^{\infty }|Z^{\infty
}}}(1)~wpa1(P_{Z^{\infty }}).
\end{equation*}%
(2) In addition, let Assumption \ref{ass:weak_equiv}(i)(ii)(iii) and
Assumption \ref{ass:rates_B}(iii) hold, then:%
\begin{eqnarray*}
||\widehat{\alpha }_{n}^{B}-\alpha _{0}|| &=&O_{P_{V^{\infty }|Z^{\infty
}}}\left( \delta _{n}\right) ~wpa1(P_{Z^{\infty }}); \\
||\widehat{\alpha }_{n}^{B}-\alpha _{0}||_{s} &=&O_{P_{V^{\infty }|Z^{\infty
}}}(||\Pi _{n}\alpha _{0}-\alpha _{0}||_{s}+\tau _{n}\times \delta
_{n})~wpa1(P_{Z^{\infty }}).
\end{eqnarray*}%
(3) The above results remain true when $\widehat{\alpha }_{n}^{B}$ is
replaced by $\widehat{\alpha }_{n}^{R,B}$.
\end{lemma}

Lemma \ref{lem:cons_boot}(2) and (3) show that $\widehat{\alpha }_{n}^{B}\in
\mathcal{N}_{osn}$ wpa1 and $\widehat{\alpha }_{n}^{R,B}\in \mathcal{N}%
_{osn} $ wpa1 regardless of whether the null $H_{0}:$ $\phi (\alpha
_{0})=\phi _{0}$ is true or not.

\subsection{Asymptotic behaviors under local alternatives\label{subsec-A5}}

In this subsection we consider the behavior of SQLR, sieve Wald and their
bootstrap versions under local alternatives. That is, we consider local
alternatives along the curve $\{\boldsymbol{\alpha }_{n}\in \mathcal{N}%
_{osn}:n\in \{1,2,...\}\}$, where
\begin{equation}
\boldsymbol{\alpha }_{n}=\alpha _{0}+d_{n}\Delta _{n}\quad with\quad \frac{%
d\phi (\alpha _{0})}{d\alpha }[\Delta _{n}]=\kappa \times \left(
1+o(1)\right) \neq 0  \label{eq:locaalt-1}
\end{equation}%
for any $(d_{n},\Delta _{n})\in \mathbb{R}_{+}\times \overline{\mathbf{V}}%
_{k(n)}$ such that $d_{n}  ||\Delta _{n}|| \leq M_{n}\delta _{n}$, $%
d_{n}||\Delta _{n}||_{s}\leq M_{n}\delta _{s,n}$ for all $n$. The
restriction on the rates under both norms is to ensure that the required
assumptions for studying the asymptotic behavior under these alternatives
(Assumption \ref{ass:phi} in particular) hold. This choice of local
alternatives is to simplify the presentation and could be relaxed somewhat.

Since we are now interested in the behavior of the test statistics under
local alternatives, we need to be more explicit about the underlying
probability, in a.s. or in probability statements. Henceforth, we use $%
P_{n,Z^{\infty }}$ to denote the probability measure over sequences $%
Z^{\infty }$ induced by the model at $\boldsymbol{\alpha }_{n}$ (we leave $%
P_{Z^{\infty }}$ to denote the one associated to $\alpha _{0}$).

\subsubsection{SQLR and SQLR$^{B}$ under local alternatives}

In this subsection we consider the behavior of the SQLR and the bootstrap
SQLR, under local alternatives along the curve $\{\boldsymbol{\alpha }%
_{n}\in \mathcal{N}_{osn}:n\in \{1,2,...\}\}$ defined in (\ref{eq:locaalt-1}%
).

\begin{theorem}
\label{thm:chi2_localt} Let conditions for Lemma \ref{thm:ThmCONVRATEGRAL}
and Proposition \ref{pro:conv-rate-RPSMDE} and Assumption \ref{ass:LAQ}
(with $\left\vert B_{n}-||u_{n}^{\ast }||^{2}\right\vert =o_{P_{n,Z^{\infty
}}}(1)$) hold under the local alternatives $\boldsymbol{\alpha }_{n}$
defined in (\ref{eq:locaalt-1}). Let Assumption \ref{ass:phi} hold and $\hat{\alpha}^{R}_{n} \in \mathcal{N}_{osn}$ wpa1-$P_{n,Z^{\infty}}$. Then,
under the \emph{local alternatives} $\boldsymbol{\alpha }_{n}$,

(1) if $d_{n}=n^{-1/2}||v_{n}^{\ast }||_{sd}$, then $||u_{n}^{\ast
}||^{2}\times \widehat{QLR}_{n}(\phi _{0})\Rightarrow \chi _{1}^{2}(\kappa
^{2})$;

(2) if $n^{1/2}||v_{n}^{\ast }||_{sd}^{-1}d_{n}\rightarrow \infty $, then $%
\lim_{n\rightarrow \infty }\left( ||u_{n}^{\ast }||^{2}\times \widehat{QLR}%
_{n}(\phi _{0})\right) =\infty $ in probability.
\end{theorem}

The statement that assumptions hold under the local alternatives $%
\boldsymbol{\alpha }_{n}$ really means that the assumptions hold when the
true DGP model is indexed by $\boldsymbol{\alpha }_{n}$ (as opposed to $%
\alpha _{0}$). For instance, this change impacts on Assumption \ref{ass:LAQ}
by changing the \textquotedblleft centering\textquotedblright\ of the
expansion to $\boldsymbol{\alpha }_{n}$ and also changing \textquotedblleft
in probability\textquotedblright\ statements to hold under $P_{n,Z^{\infty
}} $ as opposed to $P_{Z^{\infty }}$.

If we had a likelihood function instead of our criterion function, we could
adapt Le Cam's 3rd Lemma to show that Assumption \ref{ass:LAQ} under local
alternatives holds directly. Since our criterion function is not a
likelihood we cannot proceed in this manner, and we directly assume it.
Also, if we only consider contiguous alternatives, i.e., curves $\{%
\boldsymbol{\alpha }_{n}\}_{n}$ that yield probability measures $%
P_{n,Z^{\infty }}$ that are contiguous to $P_{Z^{\infty }}$, then any
statement in a.s. or wpa1 under $P_{Z^{\infty }}$ holds automatically under $%
P_{n,Z^{\infty }}$.

The next proposition presents the relative efficiency under local
alternatives of tests based on the non- and optimally weighted SQLR
statistics. We show ---aligned with the literature for regular cases--- that
optimally weighted SQLR statistic is more efficient than the non-optimally
weighted one.

\begin{proposition}
\label{pro:ARE-localalt} Let all conditions for Theorem \ref{thm:chi2_localt}
hold. Then, under the local alternatives $\boldsymbol{\alpha }_{n}$ defined
in (\ref{eq:locaalt-1}) with $d_{n}=n^{-1/2}||v_{n}^{\ast }||_{sd}$, we
have: for any $t$,
\begin{equation*}
\lim_{n\rightarrow \infty }P_{n,Z^{\infty }}(||u_{n}^{\ast }||^{2}\times
\widehat{QLR}_{n}(\phi _{0})\geq t)\leq \liminf_{n\rightarrow \infty
}P_{n,Z^{\infty }}(\widehat{QLR}_{n}^{0}(\phi _{0})\geq t).
\end{equation*}
\end{proposition}

The next theorem shows the consistency of our bootstrap SQLR statistic under
the local alternatives $\boldsymbol{\alpha }_{n}$ in (\ref{eq:locaalt-1}).
This result completes that in Remark \ref{remark-QLR-B}.

\begin{theorem}
\label{thm:BSQLR-loc-alt} Let conditions for Theorem \ref{thm:bootstrap}
hold under local alternatives $\boldsymbol{\alpha }_{n}$ defined in (\ref%
{eq:locaalt-1}). Then: (1)%
\begin{equation*}
\frac{\widehat{QLR}_{n}^{B}(\widehat{\phi }_{n})}{\sigma _{\omega }^{2}}%
=\left( \sqrt{n}\frac{\mathbb{Z}_{n}^{\omega -1}(\boldsymbol{\alpha }_{n})}{%
\sigma _{\omega }||u_{n}^{\ast }||}\right) ^{2}+o_{P_{V^{\infty }|Z^{\infty
}}}(1)=O_{P_{V^{\infty }|Z^{\infty }}}(1)~wpa1(P_{n,Z^{\infty }});\quad
\text{and}
\end{equation*}%
\begin{equation*}
\sup_{t\in \mathbb{R}}\left\vert P_{V^{\infty }|Z^{\infty }}\left( \frac{%
\widehat{QLR}_{n}^{B}(\widehat{\phi }_{n})}{\sigma _{\omega }^{2}}\leq t\mid
Z^{n}\right) -P_{Z^{\infty }}\left( \widehat{QLR}_{n}(\phi _{0})\leq t\mid
H_{0}\right) \right\vert =o_{P_{V^{\infty }|Z^{\infty
}}}(1)~wpa1(P_{n,Z^{\infty }}).
\end{equation*}%
(2) In addition, let conditions for Theorem \ref{thm:chi2_localt} hold.
Then: for any $\tau \in (0,1)$,

$\tau <\lim_{n\rightarrow \infty }P_{n,Z^{\infty }}\left( \widehat{QLR}%
_{n}(\phi _{0})\geq \widehat{c}_{n}(1-\tau )\right) <1$ under $%
d_{n}=n^{-1/2}||v_{n}^{\ast }||_{sd}$;

$\lim_{n\rightarrow \infty }P_{n,Z^{\infty }}\left( \widehat{QLR}_{n}(\phi
_{0})\geq \widehat{c}_{n}(1-\tau )\right) =1$ under $n^{1/2}||v_{n}^{\ast
}||_{sd}^{-1}d_{n}\rightarrow \infty $,

\noindent where $\widehat{c}_{n}(a)$ is the $a-th$ quantile of the
distribution of $\frac{\widehat{QLR}_{n}^{B}(\widehat{\phi }_{n})}{\sigma
_{\omega }^{2}}$ (conditional on data $\{Z_{i}\}_{i=1}^{n}$).
\end{theorem}

\subsubsection{Sieve Wald and bootstrap sieve Wald tests under local
alternatives}

\label{sec:loc-alt}

The next result establishes the asymptotic behavior of the sieve Wald test
statistic $\mathcal{W}_{n}=\left( \sqrt{n}\frac{\phi (\widehat{\alpha }%
_{n})-\phi _{0}}{||\widehat{v}_{n}^{\ast }||_{n,sd}}\right) ^{2}$ under the
local alternative along the curve $\boldsymbol{\alpha }_{n}$ defined in (\ref%
{eq:locaalt-1}).

\begin{theorem}
\label{thm:t-localt} Let $\widehat{\alpha }_{n}$ be the PSMD estimator (\ref%
{psmd}), conditions for Lemma \ref{thm:ThmCONVRATEGRAL} and Theorem \ref%
{thm:VE} and Assumption \ref{ass:LAQ} hold under the local alternatives $%
\boldsymbol{\alpha }_{n}$ defined in (\ref{eq:locaalt-1}). Let Assumption %
\ref{ass:phi} hold. Then, under the \emph{local alternatives }$\boldsymbol{%
\alpha }_{n}$,

(1) if $d_{n}=n^{-1/2}||v_{n}^{\ast }||_{sd}$, then $\mathcal{W}%
_{n}\Rightarrow \chi _{1}^{2}(\kappa ^{2})$;

(2) if $n^{1/2}||v_{n}^{\ast }||_{sd}^{-1}d_{n}\rightarrow \infty $, then $%
\lim_{n\rightarrow \infty }\mathcal{W}_{n}=\infty $ in probability.
\end{theorem}

\begin{remark}
By the same proof as that of Proposition \ref{pro:ARE-localalt}, one can
establish the asymptotically relative efficiency results for the sieve Wald
test statistic.
\end{remark}

The next theorem shows the consistency of our bootstrap sieve Wald test
statistic under the local alternatives $\boldsymbol{\alpha }_{n}$ in (\ref%
{eq:locaalt-1}). This result completes that in Remark \ref{remark-Wald-B}.

\begin{theorem}
\label{thm:waldB_con} Let all conditions for Theorem \ref{thm:bootstrap_2}%
(1) hold under local alternatives $\boldsymbol{\alpha }_{n}$ defined in (\ref%
{eq:locaalt-1}). Then: (1) for $j=1,2,$
\begin{equation*}
\sup_{t\in \mathbb{R}}\left\vert P_{V^{\infty }|Z^{\infty }}\left( \widehat{W%
}_{j,n}^{B}\leq t\mid Z^{n}\right) -P_{Z^{\infty }}\left( \widehat{W}%
_{n}\leq t\right) \right\vert =o_{P_{V^{\infty }|Z^{\infty
}}}(1)~wpa1(P_{n,Z^{\infty }}).
\end{equation*}%
(2) In addition, let conditions for Theorem \ref{thm:t-localt} hold. Then:
for any $\tau \in (0,1)$,

(2a) If $d_{n}=n^{-1/2}||v_{n}^{\ast }||_{sd}$ then:%
\begin{eqnarray*}
P_{n,Z^{\infty }}\left( \mathcal{W}_{n}\geq \widehat{c}_{j,n}(1-\tau
)\right) &=&\tau +\Pr \left( \chi _{1}^{2}(\kappa ^{2})\geq \widehat{c}%
_{j,n}(1-\tau )\right) -\Pr \left( \chi _{1}^{2}\geq \widehat{c}%
_{j,n}(1-\tau )\right)\\
&& +o_{P_{V^{\infty }|Z^{\infty
}}}(1)~wpa1(P_{n,Z^{\infty }})
\end{eqnarray*}%
and $\tau <\lim_{n\rightarrow \infty }P_{n,Z^{\infty }}\left( \mathcal{W}%
_{n}\geq \widehat{c}_{j,n}(1-\tau )\right) <1$,

(2b) If $\sqrt{n}||v_{n}^{\ast }||_{sd}^{-1}d_{n}\rightarrow \infty $ then: $%
\lim_{n\rightarrow \infty }P_{n,Z^{\infty }}\left( \mathcal{W}_{n}\geq
\widehat{c}_{j,n}(1-\tau )\right) =1$.

\noindent where $\widehat{c}_{j,n}(a)$ be the $a-th$ quantile of the
distribution of $\mathcal{W}_{j,n}^{B}\equiv \left( \widehat{W}%
_{j,n}^{B}\right) ^{2}$ (conditional on the data $\{Z_{i}\}_{i=1}^{n}$).
\end{theorem}

\subsection{Local asymptotic theory under increasing dimension of $\protect%
\phi \label{subsec-A6}$}

In this section we extend some inference results to the case of
vector-valued functional $\phi $ (i.e., $d_{\phi }\equiv d(n)>1$), and in
fact $d(n)$ could grow with $n$.

We first introduce some notation. Let $v_{j,n}^{\ast }$ be the sieve Riesz
representer corresponding to $\phi _{j}$ for $j=1,...,d(n)$ and let $\mathbf{%
v}_{n}^{\ast }\equiv (v_{1,n}^{\ast },...,v_{d(n),n}^{\ast })$. For each $x$%
, we use $\frac{dm(x,\alpha _{0})}{d\alpha }[\mathbf{v}_{n}^{\ast }]$ to
denote a $d_{\rho }\times d(n)-$matrix with $\frac{dm(x,\alpha _{0})}{%
d\alpha }[v_{j,n}^{\ast }]$ as its $j-$th column for $j=1,...,d(n)$.
Finally, let
\begin{equation*}
\Omega _{sd,n}\equiv E\left[ \left( \frac{dm(X,\alpha _{0})}{d\alpha }[%
\mathbf{v}_{n}^{\ast }]\right) ^{\prime }\Sigma ^{-1}(X)\Sigma _{0}(X)\Sigma
^{-1}(X)\left( \frac{dm(X,\alpha _{0})}{d\alpha }[\mathbf{v}_{n}^{\ast
}]\right) \right] \in \mathbb{R}^{d(n)\times d(n)}
\end{equation*}%
and
\begin{equation*}
\Omega _{n}\equiv \langle \mathbf{v}_{n}^{\ast \prime },\mathbf{v}_{n}^{\ast
}\rangle \equiv E\left[ \left( \frac{dm(X,\alpha _{0})}{d\alpha }[\mathbf{v}%
_{n}^{\ast }]\right) ^{\prime }\Sigma ^{-1}(X)\left( \frac{dm(X,\alpha _{0})%
}{d\alpha }[\mathbf{v}_{n}^{\ast }]\right) \right] \in \mathbb{R}%
^{d(n)\times d(n)}.
\end{equation*}%
Observe that for $d(n)=1$, $\Omega _{sd,n}=||v_{n}^{\ast }||_{sd}^{2}$ and $%
\Omega _{n}=||v_{n}^{\ast }||^{2}$. Also, for the case $\Sigma =\Sigma _{0}$%
, we would have
\begin{equation*}
\Omega _{n}=\Omega _{sd,n}=\Omega _{0,n}\equiv E\left[ \left( \frac{%
dm(X,\alpha _{0})}{d\alpha }[\mathbf{v}_{n}^{\ast }]\right) ^{\prime }\Sigma
_{0}^{-1}(X)\left( \frac{dm(X,\alpha _{0})}{d\alpha }[\mathbf{v}_{n}^{\ast
}]\right) \right] .
\end{equation*}

Let
\begin{equation*}
\mathcal{T}_{n}^{M}\equiv \{t\in \mathbb{R}^{d(n)}:||t||_{e}\leq
M_{n}n^{-1/2}\sqrt{d(n)}\}\quad \text{and\quad }\alpha (t)\equiv \alpha +%
\mathbf{v}_{n}^{\ast }(\Omega _{sd,n})^{-1/2}t.
\end{equation*}%
Let $(c_{n})_{n}$ be a real-valued positive sequence that converges to zero
as $n\rightarrow \infty $. The following assumption is analogous to
Assumption \ref{ass:phi} but for vector-valued $\phi $. Under Assumption \ref%
{ass:sieve}(iv), we could use $\Omega _{n}$ instead of $\Omega _{sd,n}$ in
Assumption \ref{ass:phivv1}(ii)(iii) below.

\begin{assumption}
\label{ass:phivv1} (i) for each $j=1,...,d(n)$, $\frac{d\phi _{j}(\alpha
_{0})}{d\alpha }$ satisfies Assumption \ref{ass:phi}(i); and for each $v\neq
0$, $\frac{d\phi (\alpha _{0})}{d\alpha }[v]\equiv \left( \frac{d\phi
_{1}(\alpha _{0})}{d\alpha }[v],...,\frac{d\phi _{d(n)}(\alpha _{0})}{%
d\alpha }[v]\right) ^{\prime }$ is linearly independent;
\begin{equation*}
\text{(ii)\quad }\sup_{(\alpha ,t)\in \mathcal{N}_{osn}\times \mathcal{T}%
_{n}^{M}}\left\Vert (\Omega _{sd,n})^{-1/2}\left\{ \phi \left( \alpha
(t)\right) -\phi (\alpha _{0})-\frac{d\phi (\alpha _{0})}{d\alpha }[\alpha
(t)-\alpha _{0}]\right\} \right\Vert _{e}=O\left( c_{n}\right) ;
\end{equation*}%
(iii) $\left\Vert (\Omega _{sd,n})^{-1/2}\frac{d\phi (\alpha _{0})}{d\alpha }%
[\alpha _{0,n}-\alpha _{0}]\right\Vert _{e}=O\left( c_{n}\right) $; (iv) $%
c_{n}=o(n^{-1/2})$.
\end{assumption}

For any $v\in \overline{\mathbf{V}}_{k(n)}$, we use $\langle \mathbf{v}%
_{n}^{\ast \prime },v\rangle $ to denote a $d(n)\times 1$ vector with
components $\langle v_{j,n}^{\ast },v\rangle $ for $j=1,...,d(n)$. Then $%
\frac{d\phi (\alpha _{0})}{d\alpha }[v]=\langle \mathbf{v}_{n}^{\ast \prime
},v\rangle $ with $\frac{d\phi _{j}(\alpha _{0})}{d\alpha }[v]=\langle
v_{j,n}^{\ast },v\rangle $ for $j=1,...,d(n)$. Let $\mathbf{Z}_{n}\equiv (%
\mathbb{Z}_{1,n}||v_{1,n}^{\ast }||_{sd},...,\mathbb{Z}%
_{d(n),n}||v_{d(n),n}^{\ast }||_{sd})^{\prime }$, where $\mathbb{Z}_{j,n}$
is the notation for $\mathbb{Z}_{n}$ defined in (\ref{effscore})
corresponding to the $j-$th sieve Riesz representer.

The next assumption is analogous to Assumption \ref{ass:LAQ}(i) but for the
vector valued case. Let $(a_{n},b_{n},s_{n})_{n}$ be real-valued positive
sequences that converge to zero as $n\rightarrow \infty $.

\begin{assumption}
\label{ass:LAQ_vv} (i) For all $n$, for all $(\alpha ,t)\in \mathcal{N}%
_{osn}\times \mathcal{T}_{n}^{M}$ with $\alpha (t)\in \mathcal{A}_{k(n)}$,
\begin{equation*}
\sup_{(\alpha ,t_{n})\in \mathcal{N}_{osn}\times \mathcal{T}%
_{n}^{M}}r_{n}(t_{n})\left\vert \widehat{\Lambda }_{n}(\alpha (t_{n}),\alpha
)-t_{n}^{\prime }(\Omega _{sd,n})^{-1/2}\left\{ \mathbf{Z}_{n}+\langle
\mathbf{v}_{n}^{\ast \prime },\alpha -\alpha _{0}\rangle \right\}
-t_{n}^{\prime }\frac{\mathbb{B}_{n}}{2}t_{n}\right\vert =O_{P_{Z^{\infty
}}}(1)
\end{equation*}%
where $r_{n}(t_{n})=\left( \max
\{||t_{n}||_{e}^{2}b_{n},||t_{n}||_{e}a_{n},s_{n}\}\right) ^{-1}$ and $(%
\mathbb{B}_{n})_{n}$ is such that, for each $n$, $\mathbb{B}_{n}$ is a $%
Z^{n} $ measurable positive definite matrix in $\mathbb{R}^{d(n)\times d(n)}$
and $\mathbb{B}_{n}=O_{P_{Z^{\infty }}}(1)$; (ii) $s_{n}nd(n)=o(1)$, $b_{n}%
\sqrt{d(n)}=o(1)$, $\sqrt{nd(n)}\times a_{n}=o(1)$.
\end{assumption}

In the rest of this section as well as in its proofs, since there is no risk
of confusion, we use $o_{P}$ and $O_{P}$ to denote $o_{P_{Z^{\infty }}}$ and
$O_{P_{Z^{\infty }}}$ respectively.

The next theorem extends Theorem \ref{thm:theta_anorm} to the case of
vector-valued functionals $\phi $ (of increasing dimension). Let $\mu
_{3,n}\equiv E\left[ \left\Vert \Omega _{sd,n}^{-1/2}\left( \frac{%
dm(X,\alpha _{0})}{d\alpha }[\mathbf{v}_{n}^{\ast }]\right) ^{\prime }\rho
(Z,\alpha _{0})\right\Vert _{e}^{3}\right] $.

\begin{theorem}
\label{thm:wald_vv} Let Conditions for Lemma \ref{thm:ThmCONVRATEGRAL},
Assumptions \ref{ass:phivv1} and \ref{ass:LAQ_vv} hold. Then:

(1) $n(\phi (\widehat{\alpha }_{n})-\phi (\alpha _{0}))^{\prime }\Omega
_{sd,n}^{-1}(\phi (\widehat{\alpha }_{n})-\phi (\alpha _{0}))=n\mathbf{Z}%
_{n}^{\prime }\Omega _{sd,n}^{-1}\mathbf{Z}_{n}+o_{P}\left( \sqrt{d(n)}%
\right) ;$

(2) for a fixed $d(n)=d$, if $\sqrt{n}\Omega _{sd,n}^{-1/2}\mathbf{Z}%
_{n}\Rightarrow N(0,I_{d})$ then%
\begin{equation*}
n(\phi (\widehat{\alpha }_{n})-\phi (\alpha _{0}))^{\prime }\Omega
_{sd,n}^{-1}(\phi (\widehat{\alpha }_{n})-\phi (\alpha _{0}))\Rightarrow
\chi _{d}^{2};
\end{equation*}

(3) if $d(n)\rightarrow \infty $, $d(n)=o(\sqrt{n}\mu _{3,n}^{-1})$, then:
\begin{equation*}
\frac{n(\phi (\widehat{\alpha }_{n})-\phi (\alpha _{0}))^{\prime }\Omega
_{sd,n}^{-1}(\phi (\widehat{\alpha }_{n})-\phi (\alpha _{0}))-d(n)}{\sqrt{%
2d(n)}}\Rightarrow N(0,1).
\end{equation*}
\end{theorem}

Theorem \ref{thm:wald_vv}(3) essentially states that the asymptotic
distribution of $n(\phi (\widehat{\alpha }_{n})-\phi (\alpha _{0}))^{\prime
}\Omega _{sd,n}^{-1}(\phi (\widehat{\alpha }_{n})-\phi (\alpha _{0}))$ is
close to $\chi _{d(n)}^{2}$. Moreover, as $N(d(n),2d(n))$ is close to $\chi
_{d(n)}^{2}$ for large $d(n)$ one could simulate from either distribution.
However, since $d(n)$ grows slowly (depends on the rate of $\mu _{3,n}$),%
\footnote{%
The condition $d(n)=o(\sqrt{n}\mu _{3,n}^{-1})$ is used for a coupling
argument regarding $\Omega _{sd,n}^{-1/2}\sqrt{n}\mathbf{Z}_{n}$ and a
multivariate Gaussian $N(0,I_{d(n)})$. See, e.g., Section 10.4 of \cite%
{Pollard_90}.} it might be more convenient to use $\chi _{d(n)}^{2}$ in
finite samples.

Let%
\begin{equation*}
\mathbb{D}_{n}\equiv \Omega _{sd,n}^{1/2}\Omega _{n}^{-1}\Omega _{sd,n}^{1/2}
\end{equation*}%
which, under Assumption \ref{ass:sieve}(iv), is bounded in the sense that $%
\mathbb{D}_{n}\asymp I_{d(n)}$ (see Lemma \ref{lem:D_BDD} in Appendix \ref%
{app:appC}). It is obvious that if $\Sigma =\Sigma _{0}$ then $\mathbb{D}%
_{n}=I_{d(n)}$. Note that $\mathbb{D}_{n}$ becomes $||u_{n}^{\ast }||^{-2}$
for a scalar-valued functional $\phi $.

The next result extends Theorem \ref{thm:chi2} for the SQLR statistic to the
case of vector-valued functionals $\phi $ (of increasing dimension). Recall
that $\widehat{QLR}_{n}^{0}(\phi _{0})$ is the SQLR statistic $\widehat{QLR}%
_{n}(\phi _{0})$ when $\Sigma =\Sigma _{0}$.

\begin{theorem}
\label{thm:chi2_vv} Let Conditions for Lemma \ref{thm:ThmCONVRATEGRAL} and
Proposition \ref{pro:conv-rate-RPSMDE} (in Appendix \ref{app:appB}) hold.
Let Assumptions \ref{ass:phivv1} and \ref{ass:LAQ_vv} hold with $%
\max_{t:||t||_{e}=1}|t^{\prime }\{\mathbb{B}_{n}-\mathbb{D}%
_{n}^{-1}\}t|=O_{P}(b_{n})$. Then: under the null hypothesis of $\phi
(\alpha _{0})=\phi _{0}$,

(1) $\widehat{QLR}_{n}(\phi _{0})=(\sqrt{n}\Omega _{sd,n}^{-1/2}\mathbf{Z}%
_{n})^{\prime }\mathbb{D}_{n}(\sqrt{n}\Omega _{sd,n}^{-1/2}\mathbf{Z}%
_{n})+o_{P}(\sqrt{d(n)});$

(2) if $\Sigma =\Sigma _{0}$, then $\widehat{QLR}_{n}^{0}(\phi _{0})=n%
\mathbf{Z}_{n}^{\prime }\Omega _{0,n}^{-1}\mathbf{Z}_{n}+o_{P}\left( \sqrt{%
d(n)}\right) ;$ for a fixed $d(n)=d$ if $\sqrt{n}\Omega _{0,n}^{-1/2}\mathbf{%
Z}_{n}\Rightarrow N(0,I_{d})$ then $\widehat{QLR}_{n}^{0}(\phi
_{0})\Rightarrow \chi _{d}^{2}$;

(3) if $\Sigma =\Sigma _{0}$ and $d(n)\rightarrow \infty $, $d(n)=o(\sqrt{n}%
\mu _{3,n}^{-1})$, then: $\frac{\widehat{QLR}_{n}^{0}(\phi _{0})-d(n)}{\sqrt{%
2d(n)}}\Rightarrow N(0,1).$
\end{theorem}

Theorem \ref{thm:chi2_vv}(2) is a multivariate version of Theorem \ref%
{thm:chi2}(2). Theorem \ref{thm:chi2_vv}(3) shows that the optimally
weighted SQLR preserves the Wilks phenomenon that is previously shown for
the likelihood ratio statistic for semiparametric likelihood models. Again,
as $d(n)$ grows slowly with $n$, Theorem \ref{thm:chi2_vv}(3) essentially
states that the asymptotic null distribution of $\widehat{QLR}_{n}^{0}(\phi
_{0})$ is close to $\chi _{d(n)}^{2}$.

Given Theorems \ref{thm:wald_vv} and \ref{thm:chi2_vv} and their proofs, it
is obvious that we can repeat the results on the consistency of the
bootstrap SQLR and sieve Wald as well as the local power properties of SQLR
and sieve Wald tests to vector-valued $\phi $ (of increasing dimension). We
do not state these results here due to the length of the paper. We suspect
that one could slightly improve Assumptions \ref{ass:phivv1} and \ref%
{ass:LAQ_vv} and the coupling condition $d(n)=o(\sqrt{n}\mu _{3,n}^{-1})$ so
that the dimension $d(n)$ might grow faster with $n$, but this will be a
subject of future research.

\subsection{Sufficient conditions for LQA(i) and LQA$^{B}$(i) with series LS
estimator $\widehat{m}\label{subsec-A1}$}



\begin{assumption}
\label{ass:m_ls} (i) $\mathcal{X}$ is a compact connected subset of $\mathbb{%
R}^{d_{x}}$ with Lipschitz continuous boundary, and $f_{X}$ is bounded and
bounded away from zero over $\mathcal{X}$; (ii) The smallest and largest
eigenvalues of $E[p^{J_{n}}(X)p^{J_{n}}(X)^{\prime }]$ are bounded and
bounded away from zero for all $J_{n}$; (iii) $\sup_{x\in \mathcal{X}%
}|p_{j}(x)|\leq const.<\infty $ for all $j=1,...,J_{n}$ and $J_{n}\log
(J_{n})=o(n)$ for $p^{J_{n}}(X)$ a polynomial spline or
trigonometric polynomial sieve; (iv) There is $p^{J_{n}}(X)^{\prime }\pi $
such that $\sup_{x}|g(x)-p^{J_{n}}(x)^{\prime }\pi |=O(b_{m,J_{n}})=o(1)$
uniformly in $g\in \{m(\cdot ,\alpha ):\alpha \in \mathcal{A}%
_{k(n)}^{M_{0}}\}$.
\end{assumption}

Thanks to lemma 5.2 in \cite{CC_WP13}, Assumption \ref{ass:m_ls}(iii) now
allows $J_{n}\log (J_{n})=o(n)$ for $p^{J_{n}}(X)$ being a (tensor-product)
trigonometric polynomial in addition to a polynomial spline
sieve. Let $\mathcal{O}_{on}\equiv \{\rho (\cdot ,\alpha )-\rho (\cdot
,\alpha _{0}):\alpha \in \mathcal{N}_{osn}\}$. Denote
\begin{equation*}
1\leq \sqrt{C_{n}}\equiv \int_{0}^{1}\sqrt{1+\log (N_{[]}(w(M_{n}\delta
_{s,n})^{\kappa },\mathcal{O}_{on},||\cdot ||_{L^{2}(f_{Z})}))}dw<\infty .
\end{equation*}

\begin{assumption}
\label{ass:rho_Donsker} (i) There is a sequence $\left\{ \bar{\rho}%
_{n}(Z)\right\} _{n}$ of measurable functions such that $\sup_{\mathcal{A}%
_{k(n)}^{M_{0}}}|\rho (Z,\alpha )|\leq \bar{\rho}_{n}(Z)$ a.s.-$Z$ and $E[|%
\bar{\rho}_{n}(Z)|^{2}|X]\leq const.<\infty $; (ii) there exist some $\kappa
\in (0,1]$ and $K\colon \mathcal{X}\rightarrow \mathbb{R}$ measurable with $%
E[|K(X)|^{2}]\leq const.$ such that $\forall \delta >0$,
\begin{equation*}
E\left[ \sup_{\alpha \in \mathcal{N}_{0sn}\colon ||\alpha -\alpha ^{\prime
}||_{s}\leq \delta }\left\Vert \rho (Z,\alpha )-\rho (Z,\alpha ^{\prime
})\right\Vert _{e}^{2} | X=x\right] \leq K(x)^{2}\delta ^{2\kappa
},~\forall \alpha ^{\prime }\in \mathcal{N}_{osn}\cup \{\alpha
_{0}\}~and~all~n,
\end{equation*}%
and $\max \left\{ (M_{n}\delta _{n})^{2},(M_{n}\delta _{s,n})^{2\kappa
}\right\} =(M_{n}\delta _{s,n})^{2\kappa }$; (iii) $n\delta
_{n}^{2}(M_{n}\delta _{s,n})^{\kappa }\sqrt{C_{n}}\max \left\{ (M_{n}\delta
_{s,n})^{\kappa }\sqrt{C_{n}},M_{n}\right\} =o(1)$; (iv) $\sup_{\mathcal{X}%
}||\widehat{\Sigma }(x)-\Sigma (x)||\times (M_{n}\delta
_{n})=o_{P_{Z^{\infty }}}(n^{-1/2})$; $\delta _{n}\asymp \sqrt{\frac{J_{n}}{n%
}}=\max \{\sqrt{\frac{J_{n}}{n}},b_{m,J_{n}}\}=o(n^{-1/4})$.
\end{assumption}

Let $\widetilde{m}(X,\alpha )\equiv \left( \sum_{i=1}^{n}m(X_{i},\alpha
)p^{J_{n}}(X_{i})^{\prime }\right) (P^{\prime }P)^{-}p^{J_{n}}(X)$ be the LS
projection of $m(X,\alpha )$ onto $p^{J_{n}}(X)$, and let $g(X,u_{n}^{\ast
})\equiv \{\frac{dm(X,\alpha _{0})}{d\alpha }[u_{n}^{\ast }]\}^{\prime
}\Sigma (X)^{-1}$ and $\widetilde{g}(X,u_{n}^{\ast })$ be its LS projection
onto $p^{J_{n}}(X)$.

\begin{assumption}
\label{ass:anor-mtilde} (i) $E_{P_{Z^{\infty }}}\left[ \left\Vert \frac{d%
\widetilde{m}(X,\alpha _{0})}{d\alpha }[u_{n}^{\ast }]-\frac{dm(X,\alpha
_{0})}{d\alpha }[u_{n}^{\ast }]\right\Vert _{e}^{2}\right] (M_{n}\delta
_{n})^{2}=o(n^{-1})$;

\noindent(ii) $E_{P_{Z^{\infty }}}\left[ \left\Vert \widetilde{g}%
(X,u_{n}^{\ast })-g(X,u_{n}^{\ast })\right\Vert _{e}^{2}\right] (M_{n}\delta
_{n})^{2}=o(n^{-1})$;

\noindent(iii) $\sup_{\mathcal{N}_{osn}}n^{-1}\sum_{i=1}^{n}\{||m(X_{i},%
\alpha )||_{e}^{2}-E[||m(X_{1},\alpha )||_{e}^{2}]\}=o_{P}(n^{-1/2})$;

\noindent(iv) $\sup_{\mathcal{N}_{osn}}n^{-1}\sum_{i=1}^{n}\{g(X_{i},u_{n}^{%
\ast })m(X_{i},\alpha )-E[g(X_{1},u_{n}^{\ast })m(X_{1},\alpha
)]\}=o_{P}(n^{-1/2}).$
\end{assumption}

\begin{assumption}
\label{ass:cont_diffm} (i) $m(X,\alpha )$ is twice continuously pathwise
differentiable in $\alpha \in \mathcal{N}_{os}$, a.s.-X;
\begin{equation*}
\text{(ii)\quad }E\left[ \sup_{\alpha \in \mathcal{N}_{osn}}\left\Vert \frac{%
dm(X,\alpha )}{d\alpha }[u_{n}^{\ast }]-\frac{dm(X,\alpha _{0})}{d\alpha }%
[u_{n}^{\ast }]\right\Vert _{e}^{2}\right] \times (M_{n}\delta
_{n})^{2}=o(n^{-1});
\end{equation*}%
(iii) $E\left[ \sup_{\alpha \in \mathcal{N}_{osn}}\left\Vert \frac{%
d^{2}m(X,\alpha )}{d\alpha ^{2}}[u_{n}^{\ast },u_{n}^{\ast }]\right\Vert
_{e}^{2}\right] \times (M_{n}\delta _{n})^{2}=o(1)$; (iv) Uniformly over $%
\alpha _{1}\in \mathcal{N}_{os}$ and $\alpha _{2}\in \mathcal{N}_{osn}$,
\begin{equation*}
E\left[ g(X,u_{n}^{\ast })\left( \frac{dm(X,\alpha _{1})}{d\alpha }[\alpha
_{2}-\alpha _{0}]-\frac{dm(X,\alpha _{0})}{d\alpha }[\alpha _{2}-\alpha
_{0}]\right) \right] =o(n^{-1/2}).
\end{equation*}
\end{assumption}

Assumptions \ref{ass:m_ls} and \ref{ass:rho_Donsker} are comparable to those
imposed in \cite{CP_WP07a} for a non-smooth residual function $\rho
(Z,\alpha )$. These assumptions ensure that the sample criterion function $%
\widehat{Q}_{n}$ is well approximated by a \textquotedblleft
smooth\textquotedblright\ version of it. Assumptions \ref{ass:anor-mtilde}
and \ref{ass:cont_diffm} are similar to those imposed in \cite{AC_Emetrica03}%
, \cite{AC_JOE07} and \cite{CP_WP07a}, except that we use the scaled sieve
Riesz representer $u_{n}^{\ast }\equiv v_{n}^{\ast }/\left\Vert v_{n}^{\ast
}\right\Vert _{sd}$. This is because we allow for possibly irregular
functionals (i.e., possibly $\left\Vert v_{n}^{\ast }\right\Vert \rightarrow
\infty $), while the above mentioned papers only consider regular
functionals (i.e., $\left\Vert v_{n}^{\ast }\right\Vert \rightarrow
\left\Vert v^{\ast }\right\Vert <\infty $). We refer readers to these papers
for detailed discussions and verifications of these assumptions in examples
of the general model (\ref{semi00}).

\subsection{Lemmas for series LS estimator $\widehat{m}(x,\protect\alpha )$
and its bootstrap version$\label{subsec-A2}$}

The next lemma (Lemma \ref{lem:suff_mcond_boot}) extends Lemma C.3 of \cite%
{CP_WP07} and Lemma A.1 of \cite{CP_WP07a} to the bootstrap version. Denote%
\begin{equation*}
\ell _{n}(x,\alpha )\equiv \widetilde{m}(x,\alpha )+\widehat{m}(x,\alpha
_{0})\text{\quad and\quad }\ell _{n}^{B}(x,\alpha )\equiv \widetilde{m}%
(x,\alpha )+\widehat{m}^{B}(x,\alpha _{0}).
\end{equation*}

\begin{lemma}
\label{lem:suff_mcond_boot} Let $\widehat{m}^{B}(\cdot ,\alpha )$ be the
bootstrap series LS estimator (\ref{mhat_B}). Let Assumptions \ref{ass:sieve}%
(iv), \ref{ass:weak_equiv}(i)(ii), \ref{ass:VE}(iii), \ref{ass:m_ls}, \ref%
{ass:rho_Donsker}(i)(ii), and \ref{ass:Wboot} or \ref{ass:Wboot_e} hold.
Then:

(1) For all $\delta >0$, there is a $M(\delta )>0$ such that for all $M\geq
M(\delta )$,
\begin{equation*}
P_{Z^{\infty }}\left( P_{V^{\infty }|Z^{\infty }}\left( \sup_{\alpha \in
\mathcal{N}_{osn}}\frac{\overline{\tau }_{n}}{n}\sum_{i=1}^{n}\left\Vert
\widehat{m}^{B}(X_{i},\alpha )-\ell _{n}^{B}(X_{i},\alpha )\right\Vert
_{e}^{2}\geq M\mid Z^{n}\right) \geq \delta \right) <\delta
\end{equation*}%
eventually, with $\overline{\tau }_{n}^{-1}\equiv (\delta _{n})^{2}\left(
M_{n}\delta _{s,n}\right) ^{2\kappa }C_{n}$.

(2) For all $\delta >0$, there is a $M(\delta )>0$ such that for all $M\geq
M(\delta )$,
\begin{equation*}
P_{Z^{\infty }}\left( P_{V^{\infty }|Z^{\infty }}\left( \sup_{\alpha \in
\mathcal{N}_{osn}}\frac{\tau _{n}^{\prime }}{n}\sum_{i=1}^{n}\left\Vert \ell
_{n}^{B}(X_{i},\alpha )\right\Vert _{e}^{2}\geq M\mid Z^{n}\right) \geq
\delta \right) <\delta
\end{equation*}%
eventually, with
\begin{equation*}
(\tau _{n}^{\prime })^{-1}=\max \{\frac{J_{n}}{n},b_{m,J_{n}}^{2},(M_{n}%
\delta _{n})^{2}\}=const.\times (M_{n}\delta _{n})^{2}.
\end{equation*}

(3) Let Assumption \ref{ass:rho_Donsker}(iii) hold. For all $\delta >0$,
there is $N(\delta )$ such that, for all $n\geq N(\delta )$,
\begin{equation*}
P_{Z^{\infty }}\left( P_{V^{\infty }|Z^{\infty }}\left( \sup_{\mathcal{N}%
_{osn}}\frac{s_{n}}{n}\left\vert \sum_{i=1}^{n}\left\Vert \widehat{m}%
^{B}\left( X_{i},\alpha \right) \right\Vert _{\widehat{\Sigma }%
^{-1}}^{2}-\sum_{i=1}^{n}\left\Vert \ell _{n}^{B}\left( X_{i},\alpha \right)
\right\Vert _{\widehat{\Sigma }^{-1}}^{2}\right\vert \geq \delta \mid
Z^{n}\right) \geq \delta \right) <\delta
\end{equation*}%
with
\begin{equation*}
s_{n}^{-1}\leq (\delta _{n})^{2}(M_{n}\delta _{s,n})^{\kappa }\sqrt{C_{n}}%
\max \left\{ (M_{n}\delta _{s,n})^{\kappa }\sqrt{C_{n}},M_{n}\right\}
L_{n}=o(n^{-1}),
\end{equation*}%
where $\left\{ L_{n}\right\} _{n=1}^{\infty }$ is a slowly divergent
sequence of positive real numbers (such a choice of $L_{n}$ exists under
assumption \ref{ass:rho_Donsker}(iii)).
\end{lemma}

Recall that%
\begin{equation*}
\mathbb{Z}_{n}^{\omega }=\frac{1}{n}\sum_{i=1}^{n}\omega _{i,n}\left( \frac{%
dm(X_{i},\alpha _{0})}{d\alpha }[u_{n}^{\ast }]\right) ^{\prime }\Sigma
(X_{i})^{-1}\rho (Z_{i},\alpha _{0})=\frac{1}{n}%
\sum_{i=1}^{n}g(X_{i},u_{n}^{\ast })\omega _{i,n}\rho (Z_{i},\alpha _{0}).
\end{equation*}

\begin{lemma}
\label{lem:T1n} Let all of the conditions for Lemma \ref{lem:suff_mcond_boot}%
(2) hold. If Assumptions \ref{ass:rho_Donsker}(iv), \ref{ass:anor-mtilde}
and \ref{ass:cont_diffm}(i)(ii)(iv) hold, then: for all $\delta >0$, there
is a $N(\delta )$ such that for all $n\geq N(\delta )$, {\footnotesize {%
\begin{eqnarray*}
& P_{Z^{\infty }}\left( P_{V^{\infty }|Z^{\infty }}\left( \sup_{\mathcal{N}%
_{osn}}\sqrt{n}\left\vert \frac{1}{n}\sum_{i=1}^{n}\left( \frac{d\widetilde{m%
}(X_{i},\alpha )}{d\alpha }[u_{n}^{\ast }]\right) ^{\prime }(\widehat{\Sigma
}(X_{i}))^{-1}\ell _{n}^{B}(X_{i},\alpha )-\left\{ \mathbb{Z}_{n}^{\omega
}+\langle u_{n}^{\ast },\alpha -\alpha _{0}\rangle \right\} \right\vert \geq
\delta \mid Z^{n}\right) \geq \delta \right) \\
& <\delta .
\end{eqnarray*}%
}}
\end{lemma}

\begin{lemma}
\label{lem:T3n} Let all of the conditions for Lemma \ref{lem:suff_mcond_boot}%
(2) hold. If Assumption \ref{ass:cont_diffm}(i)(iii) holds, then: for all $%
\delta >0$, there is a $N(\delta )$ such that for all $n\geq N(\delta )$,
\begin{equation*}
P_{Z^{\infty }}\left( P_{V^{\infty }|Z^{\infty }}\left( \sup_{\mathcal{N}%
_{osn}}n^{-1}\sum_{i=1}^{n}\left( \frac{d^{2}\widetilde{m}(X_{i},\alpha )}{%
d\alpha ^{2}}[u_{n}^{\ast },u_{n}^{\ast }]\right) ^{\prime }(\widehat{\Sigma
}(X_{i}))^{-1}\ell _{n}^{B}(X_{i},\alpha )\geq \delta \mid Z^{n}\right) \geq
\delta \right) <\delta .
\end{equation*}
\end{lemma}

\begin{lemma}
\label{lem:T2n} Let Assumptions \ref{ass:sieve}(iv), \ref{ass:weak_equiv}%
(i), \ref{ass:VE}(iii), \ref{ass:m_ls}, \ref{ass:anor-mtilde}(i), \ref%
{ass:cont_diffm}(ii) hold. Then: (1) For all $\delta >0$ there is a $%
M(\delta )>0$, such that for all $M\geq M(\delta )$,
\begin{equation*}
P_{Z^{\infty }}\left( \sup_{\mathcal{N}_{osn}}\frac{1}{n}\sum_{i=1}^{n}%
\left( \frac{d\widetilde{m}(X_{i},\alpha )}{d\alpha }[u_{n}^{\ast }]\right)
^{\prime }\widehat{\Sigma }^{-1}(X_{i})\left( \frac{d\widetilde{m}%
(X_{i},\alpha )}{d\alpha }[u_{n}^{\ast }]\right) \geq M\right) <\delta
\end{equation*}%
eventually.

(2) If in addition, Assumption \ref{ass:LLN_triangular} holds, then: For all
$\delta >0$, there is a $N(\delta )$ such that for all $n\geq N(\delta )$,
{\small {%
\begin{equation*}
P_{Z^{\infty }}\left( \sup_{\mathcal{N}_{osn}}\left\vert \frac{1}{n}%
\sum_{i=1}^{n}\left( \frac{d\widetilde{m}(X_{i},\alpha )}{d\alpha }%
[u_{n}^{\ast }]\right) ^{\prime }\widehat{\Sigma }^{-1}(X_{i})\left( \frac{d%
\widetilde{m}(X_{i},\alpha )}{d\alpha }[u_{n}^{\ast }]\right) -||u_{n}^{\ast
}||^{2}\right\vert \geq \delta \right) <\delta .
\end{equation*}%
}}
\end{lemma}

\pagebreak

\setcounter{page}{1} \thispagestyle{empty}

\begin{center}
	{\large{Supplement to \textquotedblleft Sieve Wald and QLR Inference on Semi/Nonparametric Conditional Moment Models" by X. Chen and D. Pouzo}}
\end{center}

\medskip

This supplemental material consists of Appendices \ref{app:appB}, \ref{app:appC} and \ref{app:appD} to the paper {\it{\textquotedblleft Sieve Wald and QLR Inference on Semi/Nonparametric Conditional Moment Models"}} by X. Chen and D. Pouzo.

\section{Additional Results and Proofs of the Results in the Main Text}

\label{app:appB}

In Appendix \ref{app:appB}, we provide the proofs of all the lemmas,
theorems and propositions stated in the main text. Additional results on
consistent sieve variance estimators and bootstrap sieve t statistics are
also presented.

\subsection{Proofs for Section \protect\ref{sec:conditions} on basic
conditions}

\noindent \textbf{Proof of Lemma \ref{lem:sieve-Riesz-property}}: For
\textbf{Result (1).} Observe that $\frac{d\phi (\alpha _{0})}{d\alpha }%
[\cdot ]$ is bounded on $(\mathbf{V},||\cdot ||)$; and in this case equation
(\ref{RRT-0}) holds. By definitions of $v_{n}^{\ast }$ and $v^{\ast }$, we
have: $\frac{d\phi (\alpha _{0})}{d\alpha }[v]=\left\langle v_{n}^{\ast
},v\right\rangle $ and $\frac{d\phi (\alpha _{0})}{d\alpha }[v]=\left\langle
v^{\ast },v\right\rangle $ for all $v\in \overline{\mathbf{V}}_{k(n)}$. Thus%
\begin{equation*}
\left\langle v^{\ast }-v_{n}^{\ast },v\right\rangle =0\text{ for all }v\in
\overline{\mathbf{V}}_{k(n)}\text{\quad and\quad }\left\Vert v^{\ast
}\right\Vert ^{2}=\left\Vert v^{\ast }-v_{n}^{\ast }\right\Vert
^{2}+\left\Vert v_{n}^{\ast }\right\Vert ^{2}\text{.}
\end{equation*}%
Since $\overline{\mathbf{V}}_{k(n)}$ is a finite dimensional Hilbert space
we have $v_{n}^{\ast }=\arg \min_{v\in \overline{\mathbf{V}}%
_{k(n)}}\left\Vert v^{\ast }-v\right\Vert $. Since $\overline{\mathbf{V}}%
_{k(n)}$ is dense in $(\overline{\mathbf{V}},||\cdot ||)$ we have $%
\left\Vert v^{\ast }-v_{n}^{\ast }\right\Vert \rightarrow 0$ and $\left\Vert
v_{n}^{\ast }\right\Vert \rightarrow \left\Vert v^{\ast }\right\Vert <\infty
$ as $k(n)\rightarrow \infty $.

For \textbf{Result (2).} We show this part by contradiction. That is,
assume that $\lim_{k(n)\rightarrow \infty }\left\Vert v_{n}^{\ast
}\right\Vert =C^{\ast }<\infty $. Since $\frac{d\phi (\alpha _{0})}{d\alpha }
$ is unbounded under $||\cdot ||$ in $\mathbf{V}$, we have: for any $M>0$,
there exists a $v_{M}\in \mathbf{V}$ such that $\left\vert \frac{d\phi
(\alpha _{0})}{d\alpha }[v_{M}]\right\vert >M||v_{M}||$.

Since $v_{M}\in \mathbf{V}$, and $\{\mathbf{V}_{k}\}_{k}$ is dense (under $%
||\cdot ||_{s}$) in $\mathbf{V}$, there exists a sequence $(v_{n,M})_{n}$
such that $v_{n,M}\in \mathbf{V}_{k(n)}$ and $\lim_{n\rightarrow \infty
}||v_{n,M}-v_{M}||_{s}=0$. This result and the fact that $||\cdot ||\leq
C||\cdot ||_{s}$ for some finite $C>0$, imply that $\lim_{n\rightarrow
\infty }||v_{n,M}||=||v_{M}||$. Also, since $\frac{d\phi (\alpha _{0})}{%
d\alpha }[\cdot ]$ is continuous or bounded on $(\mathbf{V},||\cdot ||_{s})$%
, we have:
\begin{equation*}
\lim_{n\rightarrow \infty }\left\vert \frac{d\phi (\alpha _{0})}{d\alpha }%
[v_{n,M}-v_{M}]\right\vert =0.
\end{equation*}

Hence, there exists a $N(M)$ such that
\begin{equation*}
\left\vert \frac{d\phi (\alpha _{0})}{d\alpha }[v_{n,M}]\right\vert \geq
M||v_{n,M}||
\end{equation*}%
for all $n\geq N(M)$. Since $v_{n,M}\in \mathbf{V}_{k(n)}$, the previous
inequality implies that
\begin{equation*}
||v_{n}^{\ast }||=\sup_{v\in \overline{\mathbf{V}}_{k(n)}:\left\Vert
v\right\Vert \neq 0}\frac{\left\vert \frac{d\phi (\alpha _{0})}{d\alpha }%
[v]\right\vert }{\left\Vert v\right\Vert }\geq M
\end{equation*}%
for all $n\geq N(M)$. Since $M$ is arbitrary we have $\lim_{k(n)\rightarrow
\infty }\left\Vert v_{n}^{\ast }\right\Vert =\infty $. A contradiction.
\textit{Q.E.D.}

\subsection{Proofs for Section \protect\ref{sec:AsymDist} on sieve t (Wald)
and SQLR}

\begin{lemma}
\label{lem:theta_anorm(1)} Let $\widehat{\alpha }_{n}$ be the PSMD estimator
(\ref{psmd}) and conditions for Lemma \ref{thm:ThmCONVRATEGRAL} hold. Let
Assumptions \ref{ass:phi}(i) and \ref{ass:LAQ}(i) hold. Then:%
\begin{equation*}
\sqrt{n}\langle u_{n}^{\ast },\widehat{\alpha }_{n}-\alpha _{0}\rangle =-%
\sqrt{n}\mathbb{Z}_{n}+o_{P_{Z^{\infty }}}(1).
\end{equation*}
\end{lemma}

\noindent \textbf{Proof of Lemma \ref{lem:theta_anorm(1)}:} We note that $%
n^{-1}\sum_{i=1}^{n}\left\Vert \widehat{m}(X_{i},\alpha )\right\Vert _{%
\widehat{\Sigma }^{-1}}^{2}=\widehat{Q}_{n}(\alpha )$. By Assumption \ref%
{ass:LAQ}(i), we have: for any $\epsilon _{n}\in \mathcal{T}_{n}$,
\begin{align}
& n^{-1}\sum_{i=1}^{n}\left\Vert \widehat{m}(X_{i},\widehat{\alpha }%
_{n}+\epsilon _{n}u_{n}^{\ast })\right\Vert _{\widehat{\Sigma }%
^{-1}}^{2}-n^{-1}\sum_{i=1}^{n}\left\Vert \widehat{m}(X_{i},\widehat{\alpha }%
_{n})\right\Vert _{\widehat{\Sigma }^{-1}}^{2}  \notag \\
=& 2\epsilon _{n}\{\mathbb{Z}_{n}+\langle u_{n}^{\ast },\widehat{\alpha }%
_{n}-\alpha _{0}\rangle \}+\epsilon _{n}^{2}B_{n}+o_{P_{Z^{\infty
}}}(r_{n}^{-1}),  \label{eqn:proof_Thm4.1_1}
\end{align}%
where $r_{n}^{-1}=\max \{\epsilon _{n}^{2},\epsilon
_{n}n^{-1/2},s_{n}^{-1}\} $ with $s_{n}^{-1}=o(n^{-1})$, and
\begin{equation*}
\mathbb{Z}_{n}=n^{-1}\sum_{i=1}^{n}\left( \frac{dm(X_{i},\alpha _{0})}{%
d\alpha }[u_{n}^{\ast }]\right) ^{\prime }\Sigma (X_{i})^{-1}\rho
(Z_{i},\alpha _{0}).
\end{equation*}

By adding
\begin{equation*}
E_{n}(\widehat{\alpha }_{n},\epsilon _{n})\equiv o(n^{-1})+\lambda
_{n}\left( Pen\left( \widehat{h}_{n}+\epsilon _{n}\frac{v_{h,n}^{\ast }}{%
\left\Vert v_{n}^{\ast }\right\Vert _{sd}}\right) -Pen\left( \widehat{h}%
_{n}\right) \right)
\end{equation*}%
to both sides of equation (\ref{eqn:proof_Thm4.1_1}), we have, by the
definition of the approximate minimizer $\widehat{\alpha }_{n}$ and the fact
$\widehat{\alpha }_{n}+\epsilon _{n}u_{n}^{\ast }\in \mathcal{A}_{k(n)}$
that, for all $\epsilon _{n}\in \mathcal{T}_{n}$
\begin{equation*}
2\epsilon _{n}\{\mathbb{Z}_{n}+\langle u_{n}^{\ast },\widehat{\alpha }%
_{n}-\alpha _{0}\rangle \}+\epsilon _{n}^{2}B_{n}+E_{n}(\widehat{\alpha }%
_{n},\epsilon _{n})+o_{P_{Z^{\infty }}}(r_{n}^{-1})\geq 0.
\end{equation*}%
Or, equivalently, for any $\delta >0$ and some $N(\delta )$ {\small {%
\begin{equation}
P_{Z^{\infty }}\left( ~\forall \epsilon _{n}~:~\widehat{\alpha }%
_{n}+\epsilon _{n}u_{n}^{\ast }\in \mathcal{N}_{osn},~2\epsilon _{n}\{%
\mathbb{Z}_{n}+\langle u_{n}^{\ast },\widehat{\alpha }_{n}-\alpha
_{0}\rangle \}+\epsilon _{n}^{2}B_{n}+E_{n}(\widehat{\alpha }_{n},\epsilon
_{n})\geq -\delta r_{n}^{-1}\right) \geq 1-\delta  \label{eqn:ANORM_eqn1}
\end{equation}%
}} for all $n\geq N(\delta )$. In particular, this holds for $\epsilon
_{n}\equiv \pm \{s_{n}^{-1/2}+o(n^{-1/2})\}=\pm o(n^{-1/2})$ since $%
s_{n}^{-1/2}=o(n^{-1/2})$. Under this choice of $\epsilon _{n}$, $%
r_{n}^{-1}=\max \{s_{n}^{-1},s_{n}^{-1/2}n^{-1/2}\}$. Moreover Assumptions %
\ref{A_3.6}(i)(ii) and \ref{ass:weak_equiv}(iv) imply that $E(\widehat{%
\alpha }_{n},\epsilon _{n})=o_{P_{Z^{\infty }}}(n^{-1})$. Thus $\sqrt{n}%
\epsilon _{n}^{-1}E(\widehat{\alpha }_{n},\epsilon _{n})=o_{P_{Z^{\infty }}}(%
\sqrt{n}\epsilon _{n}^{-1}n^{-1})=o_{P_{Z^{\infty }}}(1)$. Thus, from
equation (\ref{eqn:ANORM_eqn1}), it follows, {\small {%
\begin{equation*}
P_{Z^{\infty }}\left( A_{n,\delta }\geq \sqrt{n}\{\mathbb{Z}_{n}+\langle
u_{n}^{\ast },\widehat{\alpha }_{n}-\alpha _{0}\rangle \}\geq B_{n,\delta
}\right) \geq 1-\delta
\end{equation*}%
}} eventually, where
\begin{equation*}
A_{n,\delta }\equiv -0.5\sqrt{n}\epsilon _{n}B_{n}-\delta \sqrt{n}\epsilon
_{n}^{-1}r_{n}^{-1}+0.5\delta
\end{equation*}%
and
\begin{equation*}
B_{n,\delta }\equiv -0.5\sqrt{n}\epsilon _{n}B_{n}-0.5\sqrt{n}\delta
\epsilon _{n}^{-1}r_{n}^{-1}-0.5\delta
\end{equation*}%
(here the $0.5\delta $ follows from the previous algebra regarding $\sqrt{n}%
\epsilon _{n}^{-1}E(\widehat{\alpha }_{n},\epsilon _{n})$). Note that $\sqrt{%
n}\epsilon _{n}=o(1)$, $B_{n}=O_{P_{Z^{\infty }}}(1)$, and $\sqrt{n}\epsilon
_{n}^{-1}r_{n}^{-1}=\pm \max \{s_{n}^{-1/2}\sqrt{n},1\}\asymp \pm 1$. Thus
\begin{equation*}
P_{Z^{\infty }}\left( 2\delta \geq \sqrt{n}\{\mathbb{Z}_{n}+\langle
u_{n}^{\ast },\widehat{\alpha }_{n}-\alpha _{0}\rangle \}\geq -2\delta
\right) \geq 1-\delta ,~eventually.
\end{equation*}%
Hence we have established $\sqrt{n}\langle u_{n}^{\ast },\widehat{\alpha }%
_{n}-\alpha _{0}\rangle =-\sqrt{n}\mathbb{Z}_{n}+o_{P_{Z^{\infty }}}(1)$.
\textit{Q.E.D.}

\medskip

\noindent \textbf{Proof of Theorem \ref{thm:theta_anorm}}\textsc{:} By Lemma %
\ref{lem:theta_anorm(1)} and Assumption \ref{ass:LAQ}(ii), we immediately
obtain: $\sqrt{n}\langle u_{n}^{\ast },\widehat{\alpha }_{n}-\alpha
_{0}\rangle \Rightarrow N(0,1)$. Hence, in order to show the result, it
suffices to prove that
\begin{equation*}
\sqrt{n}\frac{\phi (\widehat{\alpha }_{n})-\phi (\alpha _{0})}{||v_{n}^{\ast
}||_{sd}}=\sqrt{n}\langle u_{n}^{\ast },\widehat{\alpha }_{n}-\alpha
_{0}\rangle +o_{P_{Z^{\infty }}}(1).
\end{equation*}

By Riesz representation Theorem and the orthogonality property of $\alpha
_{0,n}$, it follows
\begin{equation*}
\frac{d\phi (\alpha _{0})}{d\alpha }[\widehat{\alpha }_{n}-\alpha
_{0,n}]=\left\langle v_{n}^{\ast },\widehat{\alpha }_{n}-\alpha
_{0,n}\right\rangle =\left\langle v_{n}^{\ast },\widehat{\alpha }_{n}-\alpha
_{0}\right\rangle .
\end{equation*}%
By Assumptions \ref{ass:sieve}(iv) and \ref{ass:phi}(i) we have $%
||v_{n}^{\ast }||_{sd}\asymp ||v_{n}^{\ast }||$. This and Assumption \ref%
{ass:phi} (ii)(iii) imply
\begin{align*}
\sqrt{n}\frac{\phi (\widehat{\alpha }_{n})-\phi (\alpha _{0})}{||v_{n}^{\ast
}||_{sd}}& =\sqrt{n}||v_{n}^{\ast }||_{sd}^{-1}\frac{d\phi (\alpha _{0})}{%
d\alpha }[\widehat{\alpha }_{n}-\alpha _{0}]+o_{P_{Z^{\infty }}}(1) \\
& =\sqrt{n}||v_{n}^{\ast }||_{sd}^{-1}\frac{d\phi (\alpha _{0})}{d\alpha }[%
\widehat{\alpha }_{n}-\alpha _{0,n}]+\sqrt{n}||v_{n}^{\ast }||_{sd}^{-1}%
\frac{d\phi (\alpha _{0})}{d\alpha }[\alpha _{0,n}-\alpha
_{0}]+o_{P_{Z^{\infty }}}(1) \\
& =\sqrt{n}||v_{n}^{\ast }||_{sd}^{-1}\frac{d\phi (\alpha _{0})}{d\alpha }[%
\widehat{\alpha }_{n}-\alpha _{0,n}]+o_{P_{Z^{\infty }}}(1) \\
& =\sqrt{n}||v_{n}^{\ast }||_{sd}^{-1}\left\langle v_{n}^{\ast },\widehat{%
\alpha }_{n}-\alpha _{0}\right\rangle +o_{P_{Z^{\infty }}}(1).
\end{align*}%
Thus
\begin{equation*}
\sqrt{n}\frac{\phi (\widehat{\alpha }_{n})-\phi (\alpha _{0})}{||v_{n}^{\ast
}||_{sd}}=\sqrt{n}\frac{\langle v_{n}^{\ast },\widehat{\alpha }_{n}-\alpha
_{0}\rangle }{||v_{n}^{\ast }||_{sd}}+o_{P_{Z^{\infty }}}(1),
\end{equation*}%
and the claimed result now follows from Lemma \ref{lem:theta_anorm(1)} and
Assumption \ref{ass:LAQ}(ii). \textit{Q.E.D.}

\medskip

\noindent \textbf{Proof of Lemma \ref{lem:P-sieve}:} By the definitions of $%
\overline{\mathbf{V}}_{k(n)}$ and the sieve Riesz representer $v_{n}^{\ast
}\in \overline{\mathbf{V}}_{k(n)}$ of $\frac{d\phi (\alpha _{0})}{d\alpha }%
[\cdot ]$ given in (\ref{RRT-1}), we know that $v_{n}^{\ast }=(v_{\theta
,n}^{\ast \prime },v_{h,n}^{\ast }\left( \cdot \right) )^{\prime
}=(v_{\theta ,n}^{\ast \prime },\psi ^{k(n)}(\cdot )^{\prime }\beta
_{n}^{\ast })^{\prime }\in \overline{\mathbf{V}}_{k(n)}$ solves the
following optimization problem:%
\begin{align}
\frac{d\phi (\alpha _{0})}{d\alpha }[v_{n}^{\ast }]& =\left\Vert v_{n}^{\ast
}\right\Vert ^{2}=\sup_{v=\left( v_{\theta }^{\prime },v_{h}\right) ^{\prime
}\in \overline{\mathbf{V}}_{k(n)},v\neq 0}\frac{\left\vert \frac{\partial
\phi (\alpha _{0})}{\partial \theta ^{\prime }}v_{\theta }+\frac{\partial
\phi (\alpha _{0})}{\partial h}[v_{h}(\cdot )]\right\vert ^{2}}{E\left[
\left( \frac{dm(X,\alpha _{0})}{d\alpha }[v]\right) ^{\prime }\Sigma
(X)^{-1}\left( \frac{dm(X,\alpha _{0})}{d\alpha }[v]\right) \right] }  \notag
\\
& =\sup_{\gamma =\left( v_{\theta }^{\prime },\beta ^{\prime }\right)
^{\prime }\in \mathbb{R}^{d_{\theta }+k(n)},\gamma \neq 0}\frac{\gamma
^{\prime }\digamma _{n}\digamma _{n}^{\prime }\gamma }{\gamma ^{\prime
}D_{n}\gamma },  \label{finite-riesz}
\end{align}%
where $D_{n}=E\left[ \left( \frac{dm(X,\alpha _{0})}{d\alpha }[\overline{%
\psi }^{k(n)}(\cdot )^{\prime }]\right) ^{\prime }\Sigma (X)^{-1}\left(
\frac{dm(X,\alpha _{0})}{d\alpha }[\overline{\psi }^{k(n)}(\cdot )^{\prime
}]\right) \right] $ is a $(d_{\theta }+k(n))\times (d_{\theta }+k(n))$
positive definite matrix such that%
\begin{equation*}
\gamma ^{\prime }D_{n}\gamma \equiv E\left[ \left( \frac{dm(X,\alpha _{0})}{%
d\alpha }[v]\right) ^{\prime }\Sigma (X)^{-1}\left( \frac{dm(X,\alpha _{0})}{%
d\alpha }[v]\right) \right] \text{\quad for all }v=\left( v_{\theta
}^{\prime },\psi ^{k(n)}(\cdot )^{\prime }\beta \right) ^{\prime }\in
\overline{\mathbf{V}}_{k(n)},
\end{equation*}%
and $\digamma _{n}\equiv \left( \frac{\partial \phi (\alpha _{0})}{\partial
\theta ^{\prime }},\frac{\partial \phi (\alpha _{0})}{\partial h}[\psi
^{k(n)}(\cdot )^{\prime }]\right) ^{\prime }=\frac{d\phi (\alpha _{0})}{%
d\alpha }[\overline{\psi }^{k(n)}(\cdot )]$ is a $(d_{\theta }+k(n))\times 1$
vector.

The sieve Riesz representation (\ref{RRT-1}) becomes: for all $v=\left(
v_{\theta }^{\prime },\psi ^{k(n)}(\cdot )^{\prime }\beta \right) ^{\prime
}\in \overline{\mathbf{V}}_{k(n)},$
\begin{equation}
\frac{d\phi (\alpha _{0})}{d\alpha }[v]=\digamma _{n}^{\prime }\gamma
=\left\langle v_{n}^{\ast },v\right\rangle =\gamma _{n}^{\ast \prime
}D_{n}\gamma \text{\quad for all }\gamma =(v_{\theta }^{\prime },\beta
^{\prime })^{\prime }\in \mathbb{R}^{d_{\theta }+k(n)}.  \label{finite-RRT}
\end{equation}

It is obvious that the optimal solution of $\gamma $ in (\ref{finite-riesz})
or in (\ref{finite-RRT}) has a closed-form expression:%
\begin{equation*}
\gamma _{n}^{\ast }=\left( v_{\theta ,n}^{\ast \prime },\beta _{n}^{\ast
\prime }\right) ^{\prime }=D_{n}^{-}\digamma _{n}.
\end{equation*}%
The sieve Riesz representer is then given by
\begin{equation*}
v_{n}^{\ast }=(v_{\theta ,n}^{\ast \prime },v_{h,n}^{\ast }\left( \cdot
\right) )^{\prime }=(v_{\theta ,n}^{\ast \prime },\psi ^{k(n)}(\cdot
)^{\prime }\beta _{n}^{\ast })^{\prime }\in \overline{\mathbf{V}}_{k(n)}.
\end{equation*}%
Consequently, $\left\Vert v_{n}^{\ast }\right\Vert ^{2}=\gamma _{n}^{\ast
\prime }D_{n}\gamma _{n}^{\ast }=\digamma _{n}^{\prime }D_{n}^{-}\digamma
_{n}$. \textit{Q.E.D.}

\medskip

\textbf{Another consistent variance estimator}. For $\left\Vert v_{n}^{\ast
}\right\Vert _{sd}^{2}=E\left( S_{n,i}^{\ast }S_{n,i}^{\ast \prime }\right) $
given in (\ref{svar}) and (\ref{P-svar}), by Lemma \ref{lem:P-sieve}, it has
an alternative closed form expression:
\begin{equation*}
||v_{n}^{\ast }||_{sd}^{2}=\digamma _{n}^{\prime }D_{n}^{-}\Omega
_{n}D_{n}^{-}\digamma _{n},
\end{equation*}%
\begin{equation*}
\Omega _{n}\equiv E\left[ \left( \frac{dm(X,\alpha _{0})}{d\alpha }[%
\overline{\psi }^{k(n)}(\cdot )^{\prime }]\right) ^{\prime }\Sigma
(X)^{-1}\Sigma _{0}(X)\Sigma (X)^{-1}\left( \frac{dm(X,\alpha _{0})}{d\alpha
}[\overline{\psi }^{k(n)}(\cdot )^{\prime }]\right) \right] =\mho _{n}.
\end{equation*}%
Therefore, in addition to the sieve variance estimator $||\widehat{v}%
_{n}^{\ast }||_{n,sd}$ given in (\ref{svar-hat1}), we can define another
simple plug-in sieve variance estimator:
\begin{equation}
||\widehat{v}_{n}^{\ast }||_{n,sd}^{2}=||\widehat{v}_{n}^{\ast }||_{n,%
\widehat{\Sigma }^{-1}\widehat{\Sigma }_{0}\widehat{\Sigma }^{-1}}^{2}=\frac{%
1}{n}\sum_{i=1}^{n}\left( \frac{d\widehat{m}(X_{i},\widehat{\alpha }_{n})}{%
d\alpha }[\widehat{v}_{n}^{\ast }]\right) ^{\prime }\widehat{\Sigma }%
_{i}^{-1}\widehat{\Sigma }_{0i}\widehat{\Sigma }_{i}^{-1}\left( \frac{d%
\widehat{m}(X_{i},\widehat{\alpha }_{n})}{d\alpha }[\widehat{v}_{n}^{\ast
}]\right)  \label{svar-hat}
\end{equation}%
with $\widehat{\Sigma }_{0i}=\widehat{\Sigma }_{0}(X_{i})$ where $\widehat{%
\Sigma }_{0}(x)$ is a consistent estimator of $\Sigma _{0}(x)$, e.g. $%
\widehat{E}_{n}[\rho (Z,\widehat{\alpha }_{n})\rho (Z,\widehat{\alpha }%
_{n})^{\prime } | X=x]$, where $\widehat{E}_{n}[\cdot | X=x]$ is some consistent estimator of a conditional mean
function of $X$, such as a series, kernel or local polynomial based
estimator.

The sieve variance estimator given in (\ref{svar-hat}) can also be expressed
as
\begin{equation}
||\widehat{v}_{n}^{\ast }||_{n,sd}^{2}=\widehat{V}_{2}\equiv \widehat{%
\digamma }_{n}^{\prime }\widehat{D}_{n}^{-}\widehat{\Omega }_{n}\widehat{D}%
_{n}^{-}\widehat{\digamma }_{n}\text{\quad with}  \label{P-svar-hat}
\end{equation}%
\begin{equation*}
\widehat{\Omega }_{n}=\frac{1}{n}\sum_{i=1}^{n}\left( \frac{d\widehat{m}%
(X_{i},\widehat{\alpha }_{n})}{d\alpha }[\overline{\psi }^{k(n)}(\cdot
)^{\prime }]\right) ^{\prime }\widehat{\Sigma }_{i}^{-1}\widehat{\Sigma }%
_{0i}\widehat{\Sigma }_{i}^{-1}\left( \frac{d\widehat{m}(X_{i},\widehat{%
\alpha }_{n})}{d\alpha }[\overline{\psi }^{k(n)}(\cdot )^{\prime }]\right) .
\end{equation*}

\begin{assumption}
\label{ass:sigma-smooth} (i) $\sup_{v\in \overline{\mathbf{V}}%
_{k(n)}^{1}}\left\vert \langle v,v\rangle _{n,\Sigma ^{-1}\Sigma _{0}\Sigma
^{-1}}-\langle v,v\rangle _{\Sigma ^{-1}\Sigma _{0}\Sigma ^{-1}}\right\vert
=o_{P_{Z^{\infty }}}(1)$; and

\noindent (ii) $\sup_{\alpha \in \mathcal{N}_{osn}}\sup_{x\in \mathcal{X}}||%
\widehat{E}_{n}[\rho (z,\alpha )\rho (z,\alpha )^{\prime }|X=x]-E[\rho
(z,\alpha )\rho (z,\alpha )^{\prime }|X=x]||_{e}=o_{P_{Z^{\infty }}}(1)$.
\end{assumption}

\begin{theorem}
\label{thm:VE2} Let Assumption \ref{ass:VE}(i)-(iv), Assumption \ref%
{ass:sigma-smooth} and assumptions for Lemma \ref{thm:ThmCONVRATEGRAL} hold.
Then: Results (1) and (2) of Theorem \ref{thm:VE} hold with $||\widehat{v}%
_{n}^{\ast }||_{n,sd}^{2}$ given in (\ref{svar-hat}).
\end{theorem}

Monte Carlo studies indicate that both sieve variance estimators perform
well and similarly in finite samples.

\medskip

\noindent \textbf{Proof of Theorems \ref{thm:VE} and \ref{thm:VE2}}: In the
proof we use simplified notation $o_{P_{Z^{\infty }}}(1)=o_{P}(1)$. Also,
Result (2) trivially follows from Result (1) and Theorem \ref%
{thm:theta_anorm}. So we only show Result (1). For \textbf{Result (1),} by
the triangle inequality, we have: that
\begin{eqnarray*}
\left\vert \frac{||\widehat{v}_{n}^{\ast }||_{n,sd}-||v_{n}^{\ast }||_{sd}}{%
||v_{n}^{\ast }||_{sd}}\right\vert &\leq &\left\vert \frac{||\widehat{v}%
_{n}^{\ast }||_{n,sd}-||\widehat{v}_{n}^{\ast }||_{sd}}{||v_{n}^{\ast
}||_{sd}}\right\vert +\left\vert \frac{||\widehat{v}_{n}^{\ast
}||_{sd}-||v_{n}^{\ast }||_{sd}}{||v_{n}^{\ast }||_{sd}}\right\vert \\
&\leq &\left\vert \frac{||\widehat{v}_{n}^{\ast }||_{n,sd}-||\widehat{v}%
_{n}^{\ast }||_{sd}}{||v_{n}^{\ast }||_{sd}}\right\vert +\frac{||\widehat{v}%
_{n}^{\ast }-v_{n}^{\ast }||_{sd}}{||v_{n}^{\ast }||_{sd}}.
\end{eqnarray*}%
This and the fact $\frac{||\widehat{v}_{n}^{\ast }-v_{n}^{\ast }||_{sd}}{%
||v_{n}^{\ast }||_{sd}}\asymp \frac{||\widehat{v}_{n}^{\ast }-v_{n}^{\ast }||%
}{||v_{n}^{\ast }||}$ (under Assumption \ref{ass:sieve}(iv)) imply that
Result (1) follows from:
\begin{equation}
\frac{||\widehat{v}_{n}^{\ast }-v_{n}^{\ast }||}{||v_{n}^{\ast }||}=o_{P}(1),
\label{step1}
\end{equation}%
and
\begin{equation}
\left\vert \frac{||\widehat{v}_{n}^{\ast }||_{n,sd}-||\widehat{v}_{n}^{\ast
}||_{sd}}{||v_{n}^{\ast }||_{sd}}\right\vert =o_{P}(1).  \label{step2}
\end{equation}

We will establish results (\ref{step1}) and (\ref{step2}) in Step 1 and Step
2 below.

\textsc{Step 1.} Observe that result (\ref{step1}) is about the consistency
of the empirical sieve Riesz representer $\widehat{v}_{n}^{\ast }$ in $%
||\cdot ||$ norm, which is the same whether we use $\hat{\rho}_{i}\hat{\rho}%
_{i}^{\prime }$ or $\hat{\Sigma}_{0i}$ to compute the sieve variance
estimators (\ref{svar-hat1}) or (\ref{svar-hat}). By the Riesz
representation theorem, we have for all $v\in \overline{\mathbf{V}}_{k(n)}$,
\begin{equation*}
\frac{d\phi (\widehat{\alpha }_{n})}{d\alpha }[v]=\langle \widehat{v}%
_{n}^{\ast },v\rangle _{n,\widehat{\Sigma }^{-1}}\quad \text{and\quad }\frac{%
d\phi (\alpha _{0})}{d\alpha }[v]=\langle v_{n}^{\ast },v\rangle =\langle
v_{n}^{\ast },v\rangle _{\Sigma ^{-1}}.
\end{equation*}%
Hence, by Assumption \ref{ass:VE}(i), we have:
\begin{eqnarray*}
o_{P}(1) &=&\sup_{v\in \overline{\mathbf{V}}_{k(n)}}\left\vert \frac{\langle
\widehat{v}_{n}^{\ast },v\rangle _{n,\widehat{\Sigma }^{-1}}-\langle
v_{n}^{\ast },v\rangle }{||v||}\right\vert \\
&=&\sup_{v\in \overline{\mathbf{V}}_{k(n)}}\left\vert \frac{\langle \widehat{%
v}_{n}^{\ast },v\rangle _{n,\widehat{\Sigma }^{-1}}-\langle \widehat{v}%
_{n}^{\ast },v\rangle }{||\widehat{v}_{n}^{\ast }||\times ||v||}||\widehat{v}%
_{n}^{\ast }||+\frac{\langle \widehat{v}_{n}^{\ast },v\rangle -\langle
v_{n}^{\ast },v\rangle }{||v||}\right\vert \\
&\geq &\sup_{v\in \overline{\mathbf{V}}_{k(n)}}\left\vert \frac{\langle
\widehat{v}_{n}^{\ast }-v_{n}^{\ast },v\rangle }{||v||}\right\vert
-\sup_{\varpi \in \overline{\mathbf{V}}_{k(n)}:||\varpi ||=1}\left\vert
\langle \widehat{\varpi }_{n}^{\ast },\varpi \rangle _{n,\widehat{\Sigma }%
^{-1}}-\langle \widehat{\varpi }_{n}^{\ast },\varpi \rangle \right\vert
\times ||\widehat{v}_{n}^{\ast }||,
\end{eqnarray*}%
where $\varpi \equiv v/||v||$ and $\widehat{\varpi }_{n}^{\ast }\equiv
\widehat{v}_{n}^{\ast }/||\widehat{v}_{n}^{\ast }||$. First note that
\begin{align*}
\left\vert \langle \widehat{\varpi }_{n}^{\ast },\varpi \rangle _{n,\widehat{%
\Sigma }^{-1}}-\langle \widehat{\varpi }_{n}^{\ast },\varpi \rangle
\right\vert & \leq \left\vert \langle \widehat{\varpi }_{n}^{\ast },\varpi
\rangle _{n,\widehat{\Sigma }^{-1}}-\langle \widehat{\varpi }_{n}^{\ast
},\varpi \rangle _{n,\Sigma ^{-1}}\right\vert +\left\vert \langle \widehat{%
\varpi }_{n}^{\ast },\varpi \rangle _{n,\Sigma ^{-1}}-\langle \widehat{%
\varpi }_{n}^{\ast },\varpi \rangle _{\Sigma ^{-1}}\right\vert \\
& \equiv |T_{1n}(\varpi )|+|T_{2n}(\varpi )|.
\end{align*}%
By Assumption \ref{ass:VE}(ii), we have: $\sup_{\varpi \in \overline{\mathbf{%
V}}_{k(n)}:||\varpi ||=1}|T_{2n}(\varpi )|=o_{P}(1)$. Note that
\begin{equation*}
T_{1n}(\varpi )=n^{-1}\sum_{i=1}^{n}\left( \frac{d\widehat{m}(X_{i},\widehat{%
\alpha }_{n})}{d\alpha }[\widehat{\varpi }_{n}^{\ast }]\right) ^{\prime }\{%
\widehat{\Sigma }^{-1}(X_{i})-\Sigma ^{-1}(X_{i})\}\left( \frac{d\widehat{m}%
(X_{i},\widehat{\alpha }_{n})}{d\alpha }[\varpi ]\right) .
\end{equation*}%
By the triangle inequality, Assumptions \ref{ass:sieve}(iv) and \ref{ass:VE}%
(ii)(iii), we obtain
\begin{eqnarray*}
|T_{1n}(\varpi )|& \leq & \sup_{x\in \mathcal{X}}||\widehat{\Sigma }%
^{-1}(x)-\Sigma ^{-1}(x)||_{e}\sqrt{n^{-1}\sum_{i=1}^{n}\left\Vert \frac{d%
\widehat{m}(X_{i},\widehat{\alpha }_{n})}{d\alpha }[\widehat{\varpi }%
_{n}^{\ast }]\right\Vert _{e}^{2}}\sqrt{n^{-1}\sum_{i=1}^{n}\left\Vert \frac{%
d\widehat{m}(X_{i},\widehat{\alpha }_{n})}{d\alpha }[\varpi ]\right\Vert
_{e}^{2}} \\
& \leq & o_{P}(1)\times O_{P}\left( \sqrt{\langle \widehat{\varpi }_{n}^{\ast
},\widehat{\varpi }_{n}^{\ast }\rangle _{n,\Sigma ^{-1}}}\times \sqrt{%
\langle \varpi ,\varpi \rangle _{n,\Sigma ^{-1}}}\right) =o_{P}(1)\times
O_{P}(1)=o_{P}(1).
\end{eqnarray*}

Hence%
\begin{equation*}
0<\sup_{v\in \overline{\mathbf{V}}_{k(n)},v\neq 0}\left\vert \frac{\langle
\widehat{v}_{n}^{\ast }-v_{n}^{\ast },v\rangle }{||v||}\right\vert
=o_{P}\left( 1+||\widehat{v}_{n}^{\ast }||\right) .
\end{equation*}%
In particular, for $v=\widehat{v}_{n}^{\ast }-v_{n}^{\ast }$, this implies
\begin{equation*}
\frac{||\widehat{v}_{n}^{\ast }-v_{n}^{\ast }||}{||v_{n}^{\ast }||}=\frac{%
o_{P}\left( 1+||\widehat{v}_{n}^{\ast }||\right) }{||v_{n}^{\ast }||}.
\end{equation*}%
Note that $||v_{n}^{\ast }||\geq const.>0$ and $\frac{||\widehat{v}%
_{n}^{\ast }||}{||v_{n}^{\ast }||}\leq \frac{||\widehat{v}_{n}^{\ast
}-v_{n}^{\ast }||}{||v_{n}^{\ast }||}+1$, and thus, the previous equation
implies
\begin{equation*}
\frac{||\widehat{v}_{n}^{\ast }-v_{n}^{\ast }||}{||v_{n}^{\ast }||}%
(1-o_{P}(1))=o_{P}(1)\text{\quad and\quad }\frac{||\widehat{v}_{n}^{\ast }||%
}{||v_{n}^{\ast }||}=O_{P}(1).
\end{equation*}

\textsc{Step 2.} We now show that result (\ref{step2}) holds for the sieve
variance estimators $||\hat{v}_{n}^{\ast }||_{n,sd}^{2}$ defined in (\ref%
{svar-hat1}) and (\ref{svar-hat}). By Assumption \ref{ass:sieve}(iv), we
have:
\begin{eqnarray*}
\left\vert \frac{||\widehat{v}_{n}^{\ast }||_{n,sd}-||\widehat{v}_{n}^{\ast
}||_{sd}}{||v_{n}^{\ast }||_{sd}}\right\vert &=&\left\vert \frac{||\widehat{v%
}_{n}^{\ast }||_{n,sd}-||\widehat{v}_{n}^{\ast }||_{sd}}{||\widehat{v}%
_{n}^{\ast }||_{sd}}\right\vert \times \frac{||\widehat{v}_{n}^{\ast }||_{sd}%
}{||v_{n}^{\ast }||_{sd}}\asymp \left\vert \frac{||\widehat{v}_{n}^{\ast
}||_{n,sd}}{||\widehat{v}_{n}^{\ast }||_{sd}}-1\right\vert \times \frac{||%
\widehat{v}_{n}^{\ast }||}{||v_{n}^{\ast }||} \\
&\leq &\left( \frac{||\widehat{v}_{n}^{\ast }||_{n,sd}}{||\widehat{v}%
_{n}^{\ast }||_{sd}}+1\right) \left\vert \frac{||\widehat{v}_{n}^{\ast
}||_{n,sd}}{||\widehat{v}_{n}^{\ast }||_{sd}}-1\right\vert \times \frac{||%
\widehat{v}_{n}^{\ast }||}{||v_{n}^{\ast }||}=\left\vert \frac{||\widehat{v}%
_{n}^{\ast }||_{n,sd}^{2}}{||\widehat{v}_{n}^{\ast }||_{sd}^{2}}%
-1\right\vert \times \frac{||\widehat{v}_{n}^{\ast }||}{||v_{n}^{\ast }||} \\
&=&\left\vert ||\widehat{\varpi }_{n}^{\ast }||_{n,sd}^{2}-||\widehat{\varpi
}_{n}^{\ast }||_{sd}^{2}\right\vert \times \frac{||\widehat{v}_{n}^{\ast
}||^{2}}{||\widehat{v}_{n}^{\ast }||_{sd}^{2}}\times \frac{||\widehat{v}%
_{n}^{\ast }||}{||v_{n}^{\ast }||} \\
&=&\left\vert ||\widehat{\varpi }_{n}^{\ast }||_{n,sd}^{2}-||\widehat{\varpi
}_{n}^{\ast }||_{sd}^{2}\right\vert \times O_{P}(1),
\end{eqnarray*}%
where $\widehat{\varpi }_{n}^{\ast }\equiv \widehat{v}_{n}^{\ast }/||%
\widehat{v}_{n}^{\ast }||$, $\frac{||\widehat{v}_{n}^{\ast }||}{%
||v_{n}^{\ast }||}=O_{P}(1)$ (by Step 1), and $\frac{||\widehat{v}_{n}^{\ast
}||^{2}}{||\widehat{v}_{n}^{\ast }||_{sd}^{2}}=O_{P}(1)$ (by Assumption \ref%
{ass:sieve}(iv) and i.i.d. data). Thus, it suffices to show that
\begin{equation}
\left\vert ||\widehat{\varpi }_{n}^{\ast }||_{n,sd}^{2}-||\widehat{\varpi }%
_{n}^{\ast }||_{sd}^{2}\right\vert =o_{P}(1).  \label{step2-w}
\end{equation}

\textsc{Step 2a for the Estimator $||\hat{v}_{n}^{\ast }||_{n,sd}^{2}$
defined in (\ref{svar-hat1}).} We now establish the result (\ref{step2-w})
when the sieve variance estimator is defined in (\ref{svar-hat1}).

Let $\widehat{M}(Z_{i},\alpha )=\hat{\Sigma}_{i}^{-1}\rho (Z_{i},\alpha
)\rho (Z_{i},\alpha )^{\prime }\hat{\Sigma}_{i}^{-1}$ and $M(z,\alpha
_{0})\equiv \Sigma ^{-1}(x)\rho (z,\alpha _{0})\rho (z,\alpha _{0})^{\prime
}\Sigma ^{-1}(x)$ and $M_{i}=M(Z_{i},\alpha _{0})$. Also let $\widehat{T}%
_{i}[v_{n}]\equiv \frac{d\hat{m}(X_{i},\hat{\alpha}_{n})}{d\alpha }[v_{n}]$,
$T_{i}[v_{n}]\equiv \frac{dm(X_{i},\alpha _{0})}{d\alpha }[v_{n}]$ and $%
\Sigma (x,\alpha )\equiv E[\rho (Z,\alpha )\rho (Z,\alpha )^{\prime }|x]$.

It turns out that $\left\vert ||\widehat{\varpi }_{n}^{\ast }||_{n,sd}^{2}-||%
\widehat{\varpi }_{n}^{\ast }||_{sd}^{2}\right\vert $ can be bounded above
by
\begin{align*}
& \sup_{v_{n}\in \overline{\mathbf{V}}_{k(n)}^{1}}\left\vert
n^{-1}\sum_{i=1}^{n}\widehat{T}_{i}[v_{n}]^{\prime }\widehat{M}(Z_{i},\hat{%
\alpha}_{n})\widehat{T}_{i}[v_{n}]-n^{-1}\sum_{i=1}^{n}\widehat{T}%
_{i}[v_{n}]^{\prime }M_{i}\widehat{T}_{i}[v_{n}]\right\vert \\
+& \sup_{v_{n}\in \overline{\mathbf{V}}_{k(n)}^{1}}\left\vert
n^{-1}\sum_{i=1}^{n}\widehat{T}_{i}[v_{n}]^{\prime }M_{i}\widehat{T}%
_{i}[v_{n}]-E[T_{i}[v_{n}]^{\prime }M_{i}T_{i}[v_{n}]]\right\vert \\
+& \sup_{v_{n}\in \overline{\mathbf{V}}_{k(n)}^{1}}\left\vert
E[T_{i}[v_{n}]^{\prime }M_{i}T_{i}[v_{n}]]-E[T_{i}[v_{n}]\prime \Sigma
^{-1}(X_{i})\Sigma (X_{i},\alpha _{0})\Sigma
^{-1}(X_{i})T_{i}[v_{n}]]\right\vert \\
\equiv & A_{1n}+A_{2n}+A_{3n}.
\end{align*}%
Note that $A_{3n}=0$ by the fact that $E[M_{i}|X_{i}]=\Sigma
^{-1}(X_{i})\Sigma (X_{i},\alpha _{0})\Sigma ^{-1}(X_{i})$, and that $%
A_{2n}=o_{P}(1)$ by Assumption \ref{ass:VE}(v). Thus it remains to show that
$A_{1n}=o_{P}(1)$. We note that
\begin{align*}
A_{1n}\leq & \sup_{z}\sup_{\alpha \in \mathcal{N}_{osn}}||\widehat{M}%
(z,\alpha )-M(z,\alpha _{0})||_{e}\sup_{v_{n}\in \overline{\mathbf{V}}%
_{n}^{1}}\left\vert n^{-1}\sum_{i=1}^{n}\widehat{T}_{i}[v_{n}]^{\prime }%
\widehat{T}_{i}[v_{n}]\right\vert \\
\leq & Const.\times \sup_{z}\sup_{\alpha \in \mathcal{N}_{osn}}||\widehat{M}%
(z,\alpha )-M(z,\alpha _{0})||_{e}\sup_{v_{n}\in \overline{\mathbf{V}}%
_{n}^{1}}\left\vert n^{-1}\sum_{i=1}^{n}\widehat{T}_{i}[v_{n}]^{\prime
}M(Z_{i},\alpha _{0})\widehat{T}_{i}[v_{n}]\right\vert
\end{align*}%
where the first inequality follows from the fact that for matrices $A$ and $%
B $, $|A^{\prime }BA|\leq ||A||^2_{e}||B||_{e}$ and Assumption \ref{ass:sieve}%
(iv). Observe that by Assumptions \ref{ass:VE}(iii)(iv) and \ref{ass:sieve}%
(iv),
\begin{align*}
& \sup_{z}\sup_{\alpha \in \mathcal{N}_{osn}}||\widehat{M}(z,\alpha
)-M(z,\alpha _{0})||_{e} \\
& \leq \sup_{z}\sup_{\alpha \in \mathcal{N}_{osn}}||\hat{\Sigma}%
^{-1}(x)\{\rho (z,\alpha )\rho (z,\alpha )^{\prime }-\rho (z,\alpha
_{0})\rho (z,\alpha _{0})^{\prime }\}\hat{\Sigma}^{-1}(x)||_{e} \\
& +\sup_{z}||\hat{\Sigma}^{-1}(x)\rho (z,\alpha _{0})\rho (z,\alpha
_{0})^{\prime }\hat{\Sigma}^{-1}(x)-\Sigma ^{-1}(x)\rho (z,\alpha _{0})\rho
(z,\alpha _{0})^{\prime }\Sigma ^{-1}(x)||_{e}.
\end{align*}%
The first term in the RHS is $o_{P}(1)$ by Assumptions \ref{ass:VE}(iii)(iv)
and \ref{ass:sieve}(iv); the second term in the RHS is also of order $%
o_{P}(1)$ by Assumptions \ref{ass:VE}(iii) and \ref{ass:sieve}(iv) and the
fact that $\rho (Z,\alpha _{0})\rho (Z,\alpha _{0})^{\prime }=O_{P}(1)$. By
Assumption \ref{ass:VE}(v), $\sup_{v_{n}\in \overline{\mathbf{V}}%
_{n}^{1}}\left\vert n^{-1}\sum_{i=1}^{n}\widehat{T}_{i}[v_{n}]^{\prime
}M(Z_{i},\alpha _{0})\widehat{T}_{i}[v_{n}]\right\vert =O_{P}(1)$. Hence $%
A_{1n}=o_{P}(1)$ and result (\ref{step2-w}) holds.

\textsc{Step 2b for the Estimator $||\hat{v}_{n}^{\ast }||_{n,sd}^{2}$
defined in (\ref{svar-hat}).} Since we already provide a detailed proof for
result (\ref{step2-w}) in Step 2a for the case of (\ref{svar-hat1}), here we
present a more succinct proof for the case of (\ref{svar-hat}).

By the triangle inequality,
\begin{equation*}
\left\vert ||\widehat{\varpi }_{n}^{\ast }||_{n,sd}^{2}-||\widehat{\varpi }%
_{n}^{\ast }||_{sd}^{2}\right\vert \leq \left\vert ||\widehat{\varpi }%
_{n}^{\ast }||_{n,sd}^{2}-||\widehat{\varpi }_{n}^{\ast }||_{n,\Sigma
^{-1}\Sigma _{0}\Sigma ^{-1}}^{2}\right\vert +\left\vert ||\widehat{\varpi }%
_{n}^{\ast }||_{n,\Sigma ^{-1}\Sigma _{0}\Sigma ^{-1}}^{2}-||\widehat{\varpi
}_{n}^{\ast }||_{sd}^{2}\right\vert \equiv R_{1n}+R_{2n}.
\end{equation*}%
By Assumptions \ref{ass:sieve}(iv), \ref{ass:VE}(iii)(iv) and \ref%
{ass:sigma-smooth}, we have:%
\begin{equation*}
\sup_{x\in \mathcal{X}}||\widehat{\Sigma }^{-1}(x)\widehat{\Sigma }_{0}(x)%
\widehat{\Sigma }^{-1}(x)-\Sigma ^{-1}(x)\Sigma _{0}(x)\Sigma
^{-1}(x)||_{e}=o_{P}(1),
\end{equation*}%
where $\widehat{\Sigma }_{0}(x)=\hat{E}_{n}[\rho (Z,\hat{\alpha}_{n})\rho (Z,%
\hat{\alpha}_{n})^{\prime }|x]$. Therefore, by Assumptions \ref{ass:sieve}%
(iv) and \ref{ass:VE}(ii) and similar algebra to the one used to bound $%
T_{1n}(\varpi )$, we have:
\begin{equation*}
R_{1n}\leq o_{P}(1)\times n^{-1}\sum_{i=1}^{n}\left\Vert \frac{d\widehat{m}%
(X_{i},\widehat{\alpha }_{n})}{d\alpha }[\widehat{\varpi }_{n}^{\ast
}]\right\Vert _{e}^{2}=o_{P}(1)\times O_{P}(1)=o_{P}(1).
\end{equation*}%
Also by Assumption \ref{ass:sigma-smooth}, $R_{2n}=o_{P}(1)$. Thus result (%
\ref{step2-w}) holds. \textit{Q.E.D.}

\medskip

Before we prove Theorem \ref{thm:chi2}, we introduce some notation that will
simplify the presentation of the proofs. For any $\bar{\phi}\in \mathbb{R}$
let $\mathcal{A}(\bar{\phi})\equiv \{\alpha \in \mathcal{A}:\phi (\alpha )=%
\bar{\phi}\}$, and $\mathcal{A}_{k(n)}(\bar{\phi})\equiv \mathcal{A}(\bar{%
\phi})\cap \mathcal{A}_{k(n)}$. In particular, let $\mathcal{A}^{0}\equiv
\mathcal{A}(\phi (\alpha _{0}))$ and $\mathcal{A}_{k(n)}^{0}\equiv \mathcal{A%
}_{k(n)}(\phi (\alpha _{0}))$.

Also, we need to show that for any deviation of $\alpha $ of the type $%
\alpha +tu_{n}^{\ast }$, there exists a $t$ such that $\phi (\alpha
+tu_{n}^{\ast })$ is \textquotedblleft close\textquotedblright\ to $\phi
(\alpha _{0})$. Formally,

\begin{lemma}
\label{pro:phi_solve} Let Assumption \ref{ass:phi} hold. (1) For any $n\in
\{1,2,....\}$, any $r\in \{r:|r|\leq 2M_{n}||v_{n}^{\ast }||\delta _{n}\}$,
and any $\alpha \in \mathcal{N}_{osn}$, there exists a $t\in \mathcal{T}_{n}$
such that $\phi (\alpha +tu_{n}^{\ast })-\phi (\alpha _{0})=r$ and $\alpha
+tu_{n}^{\ast }\in \mathcal{A}_{k(n)}$. (2) For any $r\in \{r:|r|\leq
||v_{n}^{\ast }||\tau _{n}\}$ and any $\alpha \in \{\alpha \in \mathcal{A}%
_{k(n)}:||\alpha -\alpha _{0}||\leq \tau _{n}\}$ with some positive sequence
$(\tau _{n})_{n}$ such that $\tau _{n}=O(\delta _{n})$, the $t$ in part
(1) also satisfies $|t|\leq \max \{C\tau _{n},o(n^{-1/2})\}$ for some
constant $C>0$.
\end{lemma}

\smallskip

\noindent \textbf{Proof of Lemma \ref{pro:phi_solve}}: For Part (1), we
first show that there exists a $t\in \mathcal{T}_{n}$ such that $\phi
(\alpha +tu_{n}^{\ast })-\phi (\alpha _{0})=r$. By Assumption \ref{ass:phi},
there exists a $(F_{n})_{n}$ such that $F_{n}>0$ and $%
F_{n}=o(n^{-1/2}||v_{n}^{\ast }||)$ and, for any $\alpha \in \mathcal{N}%
_{osn}$ and $t\in \mathcal{T}_{n}$,
\begin{equation}
\left\vert \phi (\alpha +tu_{n}^{\ast })-\phi (\alpha _{0})-\langle
v_{n}^{\ast },\alpha -\alpha _{0}\rangle -t\frac{||v_{n}^{\ast }||^{2}}{%
||v_{n}^{\ast }||_{sd}}\right\vert \leq F_{n}.  \label{eqn:phi_solve1}
\end{equation}%
(note that by assumption \ref{ass:phi}, $F_{n}$ does not depend on $\alpha $
nor $t$).

For any $r\in \{|r|\leq 2M_{n}||v_{n}^{\ast }||\delta _{n}\}$, we define $%
(t_{l})_{l=1,2}$ as
\begin{equation*}
t_{l}||u_{n}^{\ast }||^{2}=-\langle u_{n}^{\ast },\alpha -\alpha _{0}\rangle
+a_{l,n}F_{n}||v_{n}^{\ast }||_{sd}^{-1}+r||v_{n}^{\ast }||_{sd}^{-1}.
\end{equation*}%
where $a_{l}=(-1)^{l}2$. Note that, by assumption \ref{ass:phi}(i) (the
second part), $||u_{n}^{\ast }||^{-2}\leq c^{-2}$, and thus
\begin{equation*}
|t_{l}|\leq c^{-2}\left( ||u_{n}^{\ast }||\times ||\alpha -\alpha
_{0}||+2|F_{n}|\times ||v_{n}^{\ast }||_{sd}^{-1}+|r|\times ||v_{n}^{\ast
}||_{sd}^{-1}\right) .
\end{equation*}%
Without loss of generality, we can re-normalize $M_{n}$ so that $%
c^{-2}C<M_{n}$ and $C\geq 1$. Hence,
\begin{align*}
|t_{l}|& \leq c^{-2}\left( ||u_{n}^{\ast }||\times ||\alpha -\alpha
_{0}||+2|F_{n}|\times ||v_{n}^{\ast }||_{sd}^{-1}+|r|\times ||v_{n}^{\ast
}||_{sd}^{-1}\right) \\
& =c^{-2}\left( ||u_{n}^{\ast }||\times ||\alpha -\alpha
_{0}||+2|F_{n}|\times ||v_{n}^{\ast }||_{sd}^{-1}+|r|\times ||v_{n}^{\ast
}||^{-1}||u_{n}^{\ast }||\right) \\
& \leq c^{-2}C\left( ||u_{n}^{\ast }||\times ||\alpha -\alpha
_{0}||+2|F_{n}|\times ||v_{n}^{\ast }||_{sd}^{-1}+|r|\times ||v_{n}^{\ast
}||^{-1}\right) \leq 4M_{n}^{2}\delta _{n},
\end{align*}%
where the third inequality follows from Assumption \ref{ass:phi}(i) (the
second part), and the last inequality follows from the facts that $\alpha
\in \mathcal{N}_{osn}$, $c^{-2}C2|F_{n}|\times ||v_{n}^{\ast
}||_{sd}^{-1}=o(n^{-1/2})\leq M_{n}^{2}\delta _{n}$, $r\in \{|r|\leq
2M_{n}||v_{n}^{\ast }||\delta _{n}\}$. Thus, $t_{l}$ is a valid choice in
the sense that $t_{l}\in \mathcal{T}_{n}$ for $l=1,2$.

Thus, this result and equation (\ref{eqn:phi_solve1}) imply
\begin{align*}
\phi (\alpha +t_{1}u_{n}^{\ast })-\phi (\alpha _{0})\leq & \langle
v_{n}^{\ast },\alpha -\alpha _{0}\rangle +t_{1}\frac{||v_{n}^{\ast }||^{2}}{%
||v_{n}^{\ast }||_{sd}}+F_{n} \\
=& ||v_{n}^{\ast }||_{sd}\left( \langle u_{n}^{\ast },\alpha -\alpha
_{0}\rangle +t_{1}||u_{n}^{\ast }||^{2}+F_{n}||v_{n}^{\ast
}||_{sd}^{-1}\right) \\
=& r-F_{n}<r.
\end{align*}%
Hence, $\phi (\alpha +t_{1}u_{n}^{\ast })-\phi (\alpha _{0})<r$. Similarly,
\begin{align*}
\phi (\alpha +t_{2}u_{n}^{\ast })-\phi (\alpha _{0})\geq & \langle
v_{n}^{\ast },\alpha -\alpha _{0}\rangle +t_{2}\frac{||v_{n}^{\ast }||^{2}}{%
||v_{n}^{\ast }||_{sd}}-F_{n} \\
=& ||v_{n}^{\ast }||_{sd}\left( \langle u_{n}^{\ast },\alpha -\alpha
_{0}\rangle +t_{2}||u_{n}^{\ast }||^{2}-F_{n}||v_{n}^{\ast
}||_{sd}^{-1}\right) \\
=& r+F_{n}>r
\end{align*}%
and thus $\phi (\alpha +t_{2}u_{n}^{\ast })-\phi (\alpha _{0})>r$. Since $%
t\mapsto \phi (\alpha +tu_{n}^{\ast })$ is continuous, there exists a $t\in
\lbrack t_{1},t_{2}]$ such that $\phi (\alpha +tu_{n}^{\ast })-\phi (\alpha
_{0})=r$. Clearly, $t\in \mathcal{T}_{n}$.

The fact that $\alpha (t)\equiv \alpha +tu_{n}^{\ast }\in \mathcal{A}_{k(n)}$
for $\alpha \in \mathcal{N}_{osn}$ and $t\in \mathcal{T}_{n}$ follows from
the fact that the sieve space $\mathcal{A}_{k(n)}$ is assumed to be convex with non-empty interior.
Part (2) can be proved in the same way as that for Part (1). \textit{Q.E.D.}

\medskip

\noindent \textbf{Proof of Theorem \ref{thm:chi2}}: Result (2) directly
follows from Result (1) with $\Sigma =\Sigma _{0}$ and $||u_{n}^{\ast }||=1$%
. The proof of Result (1) consists of several steps.

\textsc{Step 1.} For any $t_{n}\in \mathcal{T}_{n}$ wpa1., by Assumption \ref%
{ass:LAQ} and Lemma \ref{lem:theta_anorm(1)}, we have:
\begin{align}
0.5\left( \widehat{Q}_{n}(\widehat{\alpha }_{n}(-t_{n}))-\widehat{Q}_{n}(%
\widehat{\alpha }_{n})\right) =& -t_{n}\{\mathbb{Z}_{n}+\langle u_{n}^{\ast
},\widehat{\alpha }_{n}-\alpha _{0}\rangle \}+\frac{B_{n}}{2}%
t_{n}^{2}+o_{P_{Z^{\infty }}}(r_{n}^{-1})  \notag \\
=& \frac{B_{n}}{2}t_{n}^{2}+o_{P_{Z^{\infty }}}(r_{n}^{-1}),
\label{eqn:proof_chi2_2}
\end{align}%
where $r_{n}^{-1}=\max \{t_{n}^{2},t_{n}n^{-1/2},s_{n}^{-1}\}$ and $%
s_{n}^{-1}=o(n^{-1})$.

And under the null hypothesis, $\widehat{\alpha }_{n}^{R}\in \mathcal{N}%
_{osn}\cap \mathcal{A}_{k(n)}^{0}$ wpa1,
\begin{align}
0.5\left( \widehat{Q}_{n}(\widehat{\alpha }_{n}^{R}(t_{n}))-\widehat{Q}_{n}(%
\widehat{\alpha }_{n}^{R})\right) =& t_{n}\{\mathbb{Z}_{n}+\langle
u_{n}^{\ast },\widehat{\alpha }_{n}^{R}-\alpha _{0}\rangle \}+\frac{B_{n}}{2}%
t_{n}^{2}+o_{P_{Z^{\infty }}}(r_{n}^{-1})  \notag \\
=& t_{n}\mathbb{Z}_{n}+\frac{B_{n}}{2}t_{n}^{2}+o_{P_{Z^{\infty
}}}(r_{n}^{-1}),  \label{eqn:proof_chi2_1}
\end{align}%
where the last line follows from the fact that $t_{n}\langle u_{n}^{\ast },%
\widehat{\alpha }_{n}^{R}-\alpha _{0}\rangle =o_{P_{Z^{\infty
}}}(r_{n}^{-1}) $. To show this, note that under the null hypothesis, $%
\widehat{\alpha }_{n}^{R}\in \mathcal{N}_{osn}\cap \mathcal{A}_{k(n)}^{0}$
wpa1. This and Assumption \ref{ass:phi}(ii) imply that
\begin{equation*}
\left\vert \underbrace{\phi (\widehat{\alpha }_{n}^{R})-\phi (\alpha _{0})}%
_{=0}-\frac{d\phi (\alpha _{0})}{d\alpha }[\widehat{\alpha }_{n}^{R}-\alpha
_{0}]\right\vert =o_{P_{Z^{\infty }}}(n^{-1/2}||v_{n}^{\ast }||).
\end{equation*}%
Thus
\begin{equation*}
P_{Z^{\infty }}\left( \frac{\sqrt{n}}{||v_{n}^{\ast }||}\left\vert \frac{%
d\phi (\alpha _{0})}{d\alpha }[\widehat{\alpha }_{n}^{R}-\alpha
_{0}]\right\vert <\delta \right) \geq 1-\delta
\end{equation*}%
eventually. By similar calculations to those in the proof of Theorem \ref%
{thm:theta_anorm}, we have
\begin{equation*}
P_{Z^{\infty }}\left( \sqrt{n}\left\vert \langle u_{n}^{\ast },\widehat{%
\alpha }_{n}^{R}-\alpha _{0}\rangle \right\vert <\delta \right) \geq
1-\delta ,~eventually.
\end{equation*}%
Hence, $\langle u_{n}^{\ast },\widehat{\alpha }_{n}^{R}-\alpha _{0}\rangle
=o_{P_{Z^{\infty }}}(n^{-1/2})$, and thus $t_{n}\langle u_{n}^{\ast },%
\widehat{\alpha }_{n}^{R}-\alpha _{0}\rangle =o_{P_{Z^{\infty
}}}(n^{-1/2}t_{n})=o_{P_{Z^{\infty }}}(r_{n}^{-1})$.\newline

\textsc{Step 2.} We choose $t_{n}=-\mathbb{Z}_{n}B_{n}^{-1}$. Note that
under assumption \ref{ass:LAQ}, $t_{n}\in \mathcal{T}_{n}$ wpa1. By the
definition of $\widehat{\alpha }_{n}$, we have, under the null hypothesis,
\begin{align*}
0.5\left( \widehat{Q}_{n}(\widehat{\alpha }_{n}^{R})-\widehat{Q}_{n}(%
\widehat{\alpha }_{n})\right) \geq & 0.5\left( \widehat{Q}_{n}(\widehat{%
\alpha }_{n}^{R})-\widehat{Q}_{n}(\widehat{\alpha }_{n}^{R}(t_{n}))\right)
-o_{P_{Z^{\infty }}}(n^{-1}) \\
=& \frac{1}{2}\mathbb{Z}_{n}^{2}B_{n}^{-1}-o_{P_{Z^{\infty }}}(\max
\{B_{n}^{-2}\mathbb{Z}_{n}^{2},-B_{n}^{-1}\mathbb{Z}_{n}n^{-1/2},s_{n}^{-1}%
\})-o_{P_{Z^{\infty }}}(n^{-1}) \\
=& \frac{1}{2}\mathbb{Z}_{n}^{2}B_{n}^{-1}+o_{P_{Z^{\infty }}}(n^{-1}),
\end{align*}%
where the first inequality follows from the fact that, since $t_{n}\in
\mathcal{T}_{n}$ and $\widehat{\alpha }_{n}^{R}\in \mathcal{N}_{osn}$ wpa1,
then $\widehat{\alpha }_{n}^{R}(t_{n})\in \mathcal{A}_{k(n)}$ wpa1; and the
second line follows from equation (\ref{eqn:proof_chi2_1}) with $t_{n}=-%
\mathbb{Z}_{n}B_{n}^{-1}$.\newline

\textsc{Step 3.} We choose $t_{n}^{\ast }\in \mathcal{T}_{n}$ wpa1 such that
(a) $\phi (\widehat{\alpha }_{n}(t_{n}^{\ast }))=\phi (\alpha _{0})$, $%
\widehat{\alpha }_{n}(t_{n}^{\ast })\in \mathcal{A}_{k(n)}$, and (b) $%
t_{n}^{\ast }=\mathbb{Z}_{n}\frac{||v_{n}^{\ast }||_{sd}^{2}}{||v_{n}^{\ast
}||^{2}}+o_{P_{Z^{\infty }}}(n^{-1/2})=O_{P_{Z^{\infty }}}(n^{-1/2})$.

Suppose such a $t_{n}^{\ast }$ exists, then $\left[ r_{n}(t_{n}^{\ast })%
\right] ^{-1}=\max \{(t_{n}^{\ast })^{2},t_{n}^{\ast
}n^{-1/2},o(n^{-1})\}=O_{P_{Z^{\infty }}}(n^{-1})$. By the definition of $%
\widehat{\alpha }_{n}^{R}$, we have, under the null hypothesis,%
\begin{align*}
0.5\left( \widehat{Q}_{n}(\widehat{\alpha }_{n}^{R})-\widehat{Q}_{n}(%
\widehat{\alpha }_{n})\right) \leq & 0.5\left( \widehat{Q}_{n}(\widehat{%
\alpha }_{n}(t_{n}^{\ast }))-\widehat{Q}_{n}(\widehat{\alpha }_{n})\right)
+o_{P_{Z^{\infty }}}(n^{-1}) \\
=& t_{n}^{\ast }\{\mathbb{Z}_{n}+\langle u_{n}^{\ast },\widehat{\alpha }%
_{n}-\alpha _{0}\rangle \}+\frac{B_{n}}{2}\left( t_{n}^{\ast }\right)
^{2}+o_{P_{Z^{\infty }}}(n^{-1}) \\
=& \frac{B_{n}}{2}\left( \mathbb{Z}_{n}\frac{||v_{n}^{\ast }||_{sd}^{2}}{%
||v_{n}^{\ast }||^{2}}+o_{P_{Z^{\infty }}}(n^{-1/2})\right)
^{2}+o_{P_{Z^{\infty }}}(n^{-1}) \\
=& \frac{1}{2}\mathbb{Z}_{n}^{2}B_{n}^{-1}+o_{P_{Z^{\infty }}}(n^{-1})=\frac{%
1}{2}\mathbb{Z}_{n}^{2}\frac{||v_{n}^{\ast }||_{sd}^{2}}{||v_{n}^{\ast
}||^{2}}+o_{P_{Z^{\infty }}}(n^{-1}),
\end{align*}%
where the second line follows from Assumption \ref{ass:LAQ}(i) and the fact
that $t_{n}^{\ast }$ satisfying (b), $\left[ r_{n}(t_{n}^{\ast })\right]
^{-1}=O_{P_{Z^{\infty }}}(n^{-1})$; the third line follows from equation (%
\ref{eqn:proof_chi2_2}) and the fact that $t_{n}^{\ast }$ satisfying (b);
and the last line follows from Assumptions \ref{ass:phi}(i) and \ref{ass:LAQ}%
(ii), $\left\vert B_{n}-||u_{n}^{\ast }||^{2}\right\vert =o_{P_{Z^{\infty
}}}(1)$ and $u_{n}^{\ast }=v_{n}^{\ast }/\left\Vert v_{n}^{\ast }\right\Vert
_{sd}$.

We now show that there is a $t_{n}^{\ast }\in \mathcal{T}_{n}$ wpa1 such
that (a) and (b) hold. Denote $r\equiv \phi (\widehat{\alpha }_{n})-\phi
(\alpha _{0})$. Since $\widehat{\alpha }_{n}\in \mathcal{N}_{osn}$ wpa1 and $%
\phi (\widehat{\alpha }_{n})-\phi (\alpha _{0})=O_{P_{Z^{\infty
}}}(||v_{n}^{\ast }||/\sqrt{n})$ (see the proof of Theorem \ref%
{thm:theta_anorm}), we have $|r|\leq 2M_{n}||v_{n}^{\ast }||\delta _{n}$.
Thus, by Lemma \ref{pro:phi_solve}, there is a $t_{n}^{\ast }\in \mathcal{T}%
_{n}$ wpa1 such that $\widehat{\alpha }_{n}(t_{n}^{\ast })=\widehat{\alpha }%
_{n}+t_{n}^{\ast }u_{n}^{\ast }\in \mathcal{A}_{k(n)}$ and $\phi (\widehat{%
\alpha }_{n}(t_{n}^{\ast }))=\phi (\alpha _{0})$, so (a) holds. Moreover, by
Assumption \ref{ass:phi}(ii), such a choice of $t_{n}^{\ast }$ also
satisfies
\begin{equation*}
\left\vert \underbrace{\phi (\widehat{\alpha }_{n}(t_{n}^{\ast }))-\phi
(\alpha _{0})}_{=0}-\frac{d\phi (\alpha _{0})}{d\alpha }[\widehat{\alpha }%
_{n}-\alpha _{0}+t_{n}^{\ast }u_{n}^{\ast }]\right\vert =o_{P_{Z^{\infty
}}}(||v_{n}^{\ast }||/\sqrt{n}).
\end{equation*}%
By Assumption \ref{ass:phi}(i) and the definition of $u_{n}^{\ast
}=v_{n}^{\ast }/\left\Vert v_{n}^{\ast }\right\Vert _{sd}$ we have: $\frac{%
d\phi (\alpha _{0})}{d\alpha }[t_{n}^{\ast }u_{n}^{\ast }]=t_{n}^{\ast }%
\frac{||v_{n}^{\ast }||^{2}}{||v_{n}^{\ast }||_{sd}}$. Thus
\begin{equation*}
P_{Z^{\infty }}\left( \frac{\sqrt{n}}{||v_{n}^{\ast }||}\left\vert \frac{%
d\phi (\alpha _{0})}{d\alpha }[\widehat{\alpha }_{n}-\alpha
_{0}]+t_{n}^{\ast }\frac{||v_{n}^{\ast }||^{2}}{||v_{n}^{\ast }||_{sd}}%
\right\vert <\delta \right) \geq 1-\delta
\end{equation*}%
eventually. By similar algebra to that in the proof of Theorem \ref%
{thm:theta_anorm} it follows that the LHS of the equation above is majorized
by 
\begin{align*}
& P_{Z^{\infty }}\left( \frac{\sqrt{n}}{||v_{n}^{\ast }||}\left\vert \langle
v_{n}^{\ast },\widehat{\alpha }_{n}-\alpha _{0}\rangle +t_{n}^{\ast }\frac{%
||v_{n}^{\ast }||^{2}}{||v_{n}^{\ast }||_{sd}}\right\vert <\delta \right)
+\delta \\
& =P_{Z^{\infty }}\left( \frac{\sqrt{n}}{||v_{n}^{\ast }||}\left\vert -%
\mathbb{Z}_{n}||v_{n}^{\ast }||_{sd}+t_{n}^{\ast }\frac{||v_{n}^{\ast }||^{2}%
}{||v_{n}^{\ast }||_{sd}}\right\vert <\delta \right) +\delta \\
& =P_{Z^{\infty }}\left( \sqrt{n}\frac{||v_{n}^{\ast }||_{sd}}{||v_{n}^{\ast
}||}\left\vert -\mathbb{Z}_{n}+t_{n}^{\ast }\frac{||v_{n}^{\ast }||^{2}}{%
||v_{n}^{\ast }||_{sd}^{2}}\right\vert <\delta \right) +\delta ,
\end{align*}%
where the second line follows from the proof of Lemma \ref%
{lem:theta_anorm(1)}. Since $\frac{||v_{n}^{\ast }||_{sd}}{||v_{n}^{\ast }||}%
\asymp const.$ (by Assumption \ref{ass:phi}(i)), we obtain:
\begin{equation*}
P_{Z^{\infty }}\left( \sqrt{n}\left\vert t_{n}^{\ast }-\mathbb{Z}_{n}\frac{%
||v_{n}^{\ast }||_{sd}^{2}}{||v_{n}^{\ast }||^{2}}\right\vert <\delta
\right) \geq 1-\delta ,~eventually.
\end{equation*}%
Since $\sqrt{n}\mathbb{Z}_{n}=O_{P_{Z^{\infty }}}(1)$ (Assumption \ref%
{ass:LAQ}(ii)), we have: $t_{n}^{\ast }=O_{P_{Z^{\infty }}}(n^{-1/2})$, and
in fact, $\sqrt{n}t_{n}^{\ast }=\sqrt{n}\mathbb{Z}_{n}\frac{||v_{n}^{\ast
}||_{sd}^{2}}{||v_{n}^{\ast }||^{2}}+o_{P_{Z^{\infty }}}(1)$ and hence (b)
holds. \textit{Q.E.D.}

\medskip

Let $\mathcal{A}^{R}\equiv \{\alpha \in \mathcal{A}:\phi (\alpha )=\phi
_{0}\}$ be the restricted parameter space. Then $\alpha _{0}\in \mathcal{A}%
^{R}$ iff the null hypothesis $H_{0}:\phi (\alpha _{0})=\phi _{0}$ holds.
Also, $\mathcal{A}_{k(n)}^{R}\equiv \{\alpha \in \mathcal{A}_{k(n)}\colon
\phi (\alpha )=\phi _{0}\}$ is a sieve space for $\mathcal{A}^{R}$. Let $\{%
\bar{\alpha}_{0,n}\in \mathcal{A}_{k(n)}^{R}\}$ be a sequence such that $||%
\bar{\alpha}_{0,n}-\alpha _{0}||_{s}\leq \inf_{\alpha \in \mathcal{A}%
_{k(n)}^{R}}||\alpha -\alpha _{0}||_{s}+o(n^{-1})$.\footnote{%
Sufficient conditions for $\overline{\alpha }_{0,n}\in \mathcal{A}%
_{k(n)}^{R} $ to solve $\inf_{\alpha \in \mathcal{A}_{k(n)}^{R}}\left\Vert
\alpha -\alpha _{0}\right\Vert _{s}$ under the null include either (a) $%
\mathcal{A}_{k(n)}$ is compact (in $||\cdot ||_{s}$) and $\phi $ is
continuous (in $||\cdot ||_{s}$), or (b) $\mathcal{A}_{k(n)}$ is convex and $%
\phi $ is linear.}

\begin{assumption}
\label{ass:conv-rate-RPSMDE} (i) $|Pen(\bar{h}_{0,n})-Pen(h_{0})|=O(1)$ and $%
Pen(h_{0})<\infty $; (ii) $\widehat{Q}_{n}(\bar{\alpha}_{0,n})\leq c_{0}Q(%
\bar{\alpha}_{0,n})+o_{P_{Z^{\infty }}}(n^{-1})$.
\end{assumption}

This assumption on $\bar{\alpha}_{0,n}\in \mathcal{A}_{k(n)}^{R}$ is the
same as Assumptions \ref{A_3.6}(ii) and \ref{ass:rates}(i) imposed on $\Pi
_{n}\alpha _{0}\in \mathcal{A}_{k(n)}$, and can be verified in the same way
provided that $\alpha _{0}\in \mathcal{A}^{R}$.

\begin{proposition}
\label{pro:conv-rate-RPSMDE} Let $\widehat{\alpha }_{n}^{R}\in \mathcal{A}%
_{k(n)}^{R}$ be the restricted PSMD estimator (\ref{Rpsmd}) and $\alpha
_{0}\in \mathcal{A}^{R}$. Let Assumptions \ref{ass:sieve}, \ref{A_3.6}(iii), %
\ref{ass:rates}(ii), \ref{ass:conv-rate-RPSMDE} and $Q(\bar{\alpha}%
_{0,n})+o(n^{-1})=O(\lambda _{n})=o(1)$ hold. Then:

(1) $Pen(\widehat{h}_{n}^{R})=O_{P_{Z^{\infty }}}(1)$ and $||\widehat{\alpha
}_{n}^{R}-\alpha _{0}||_{s}=o_{P_{Z^{\infty }}}(1)$;

(2) Further, let $Q(\bar{\alpha}_{0,n})\asymp Q(\Pi _{n}\alpha _{0})$ and
Assumptions \ref{A_3.6}(ii), \ref{ass:rates}(i) and \ref{ass:weak_equiv}%
(i)(ii)(iii) hold. Then: $||\widehat{\alpha }_{n}^{R}-\alpha
_{0}||=O_{P_{Z^{\infty }}}(\delta _{n})$ and $||\widehat{\alpha }%
_{n}^{R}-\alpha _{0}||_{s}=O_{P_{Z^{\infty }}}(||\alpha _{0}-\Pi _{n}\alpha
_{0}||_{s}+\tau _{n}\delta _{n})$.
\end{proposition}

\medskip

\noindent \textbf{Proof of Proposition \ref{pro:conv-rate-RPSMDE}.} The
proof is very similar to those for theorem 3.2 and remark 4.1 in \cite%
{CP_WP07} by recognizing that $\mathcal{A}_{k(n)}^{R}$ is a sieve for $%
\alpha _{0}\in \mathcal{A}^{R}$.

For \textbf{Result (1)}, we first want to show that $\widehat{\alpha }%
_{n}^{R}\in \mathcal{A}_{k(n)}^{R}\cap \{Pen(h)\leq M\}$ for some $M>0$ wpa1-%
$P_{Z^{\infty }}$. By definitions of $\widehat{\alpha }_{n}^{R}$ and $\bar{%
\alpha}_{0,n}$, Assumption \ref{ass:conv-rate-RPSMDE}(i)(ii) and the
condition that $Q(\bar{\alpha}_{0,n})+o(n^{-1})=O(\lambda _{n})$, we have:
\begin{equation*}
Pen(\widehat{h}_{n}^{R})\leq \frac{\widehat{Q}_{n}(\bar{\alpha}_{0,n})}{%
\lambda _{n}}+Pen(\bar{h}_{0,n})+\frac{o(n^{-1})}{\lambda _{n}}\leq \frac{Q(%
\bar{\alpha}_{0,n})+o(n^{-1})}{\lambda _{n}}+O_{P_{Z^{\infty
}}}(1)=O_{P_{Z^{\infty }}}(1).
\end{equation*}%
Therefore, for any $\epsilon >0$, $\Pr (Pen(\widehat{h}_{n}^{R})\geq
M)<\epsilon $ for some $M$, eventually.

We now show that $\Pr (||\widehat{\alpha }_{n}^{R}-\alpha _{0}||_{s}\geq
\epsilon )=o(1)$ for any $\epsilon >0$. Let $\mathcal{A}_{k(n)}^{R,M}\equiv
\mathcal{A}_{k(n)}^{R}\cap \{Pen(h)\leq M\}$ and $\mathcal{A}^{R,M}\equiv
\mathcal{A}^{R}\cap \{Pen(h)\leq M\}$. These sets are compact under $||\cdot
||_{s}$ (by Assumption \ref{A_3.6}(iii) and the $||\cdot ||_{s}-$ continuity
of $\phi $). Assumptions \ref{ass:sieve}(i)(iv) and \ref%
{ass:conv-rate-RPSMDE}(i) imply that $\alpha _{0}\in \mathcal{A}^{R,M}$ and $%
\bar{\alpha}_{0,n}\in \mathcal{A}_{k(n)}^{R,M}$. Under assumption \ref%
{ass:sieve}(ii), $cl\left( \cup _{k}\mathcal{A}_{k}\right) \supseteq
\mathcal{A}$ and thus $cl\left( \cup _{k}\mathcal{A}_{k}^{R,M}\right)
\supseteq \mathcal{A}^{R,M}$. Therefore $||\bar{\alpha}_{0,n}-\alpha
_{0}||_{s}=o(1)$ by the definition of $\bar{\alpha}_{0,n}$ and the fact that
$\mathcal{A}_{k(n)}^{R,M}$ being dense in $\mathcal{A}^{R,M}$.

By standard calculations, it follows that, for any $\epsilon >0$,
\begin{eqnarray*}
&&\Pr (||\widehat{\alpha }_{n}^{R}-\alpha _{0}||_{s}\geq \epsilon )\\
&\leq& \Pr
\left( \inf_{\mathcal{A}_{k(n)}^{R,M}:||\alpha -\alpha _{0}||_{s}\geq
\epsilon }\{\widehat{Q}_{n}(\alpha )+\lambda _{n}Pen(h)\}\leq \widehat{Q}%
_{n}(\bar{\alpha}_{0,n})+\lambda _{n}Pen(\bar{h}_{0,n})+o_{P}(n^{-1})\right)\\
&& +0.5\epsilon.
\end{eqnarray*}

Moreover (up to omitted constants)
\begin{eqnarray*}
& & \Pr \left( ||\widehat{\alpha }_{n}^{R}-\alpha _{0}||_{s}\geq \epsilon
\right) \\
& \leq & \Pr \left( \inf_{\mathcal{A}_{k(n)}^{R,M}:||\alpha -\alpha
_{0}||_{s}\geq \epsilon }\{Q(\alpha )+\lambda _{n}Pen(h)\}\leq Q(\bar{\alpha}%
_{0,n})+\lambda _{n}Pen(\bar{h}_{0,n})+O_{P}(\bar{\delta}%
_{m,n}^{2})+o_{P}(n^{-1})\right) +\epsilon \\
& \leq  & \Pr \left( \inf_{\mathcal{A}^{R,M}:||\alpha -\alpha _{0}||_{s}\geq
\epsilon }\{Q(\alpha )+\lambda _{n}Pen(h)\}\leq Q(\bar{\alpha}%
_{0,n})+\lambda _{n}Pen(\bar{h}_{0,n})+O_{P}(\bar{\delta}%
_{m,n}^{2})+o_{P}(n^{-1})\right) +\epsilon ,
\end{eqnarray*}%
where the first line follows by Assumptions \ref{ass:rates}(ii) and \ref%
{ass:conv-rate-RPSMDE} and the second by $\mathcal{A}_{k(n)}^{R,M}\subseteq
\mathcal{A}^{R,M}$. Since $\mathcal{A}^{R,M}$ is compact under $||\cdot
||_{s}$, $\alpha _{0}\in \mathcal{A}^{R,M}$ is unique and $Q$ is continuous
(Assumption \ref{ass:sieve}), then $\inf_{\mathcal{A}^{R,M}:||\alpha -\alpha
_{0}||_{s}\geq \epsilon }\{Q(\alpha )+\lambda _{n}Pen(h)\}\geq c(\epsilon
)>0 $; however, the term $Q(\bar{\alpha}_{0,n})+\lambda _{n}Pen(\bar{h}%
_{0,n})+O_{P}(\bar{\delta}_{m,n}^{2})+o_{P}(n^{-1})=o_{P}(1)$ and thus the
desired result follows.

For \textbf{Result (2)}, we now show that $||\widehat{\alpha }%
_{n}^{R}-\alpha _{0}||=O_{P_{Z^{\infty }}}(\kappa _{n})$ where
$\kappa _{n}^{2}\equiv \max \{\delta _{n}^{2},||\bar{\alpha}_{0,n}-\alpha
_{0}||^{2},\lambda _{n},o(n^{-1})\}$. Let $\mathcal{A}_{osn}^{R}=\{\alpha
\in \mathcal{A}_{osn}:\phi (\alpha )=\phi (\alpha _{0})\}$ and $\mathcal{A}%
_{os}^{R}=\{\alpha \in \mathcal{A}_{os}:\phi (\alpha )=\phi (\alpha _{0})\}$%
. Result (1) implies that $\widehat{\alpha }_{n}^{R}\in \mathcal{A}%
_{osn}^{R} $ wpa1. To show Result (2), we employ analogous arguments to
those for Result (1) and obtain that for all large $K>0$,%
\begin{eqnarray*}
& &\Pr \left( ||\widehat{\alpha }_{n}^{R}-\alpha _{0}||\geq K\kappa
_{n}\right) \\
& \leq &\Pr \left( \inf_{\mathcal{A}_{osn}^{R}:||\alpha -\alpha _{0}||\geq
K\kappa _{n}}Q(\alpha )+\lambda _{n}Pen(h)\leq Q(\bar{\alpha}_{0,n})+\lambda
_{n}Pen(\bar{h}_{0,n})+O_{P}(\delta _{n}^{2})+o_{P}(n^{-1})\right) +\epsilon
\\
& \leq & \Pr \left( \inf_{\mathcal{A}_{os}^{R}:||\alpha -\alpha _{0}||\geq
K\kappa _{n}}||\alpha -\alpha _{0}||^{2}\leq Const.\{||\bar{\alpha}%
_{0,n}-\alpha _{0}||^{2}+\lambda _{n}Pen(\bar{h}_{0,n})+O_{P}(\delta
_{n}^{2})+o_{P}(n^{-1})\right) +\epsilon \\
& \leq  & \Pr \left( K^{2}\kappa _{n}^{2}\leq Const.||\bar{\alpha}_{0,n}-\alpha
_{0}||^{2}+O(\lambda _{n})+O_{P}(\delta _{n}^{2})+o_{P}(n^{-1})\right)
+\epsilon ,
\end{eqnarray*}%
where the first inequality is due to Assumption \ref{ass:conv-rate-RPSMDE}%
(ii) and the assumption that $\widehat{Q}_{n}(\alpha )\geq cQ(\alpha
)-O_{P_{Z^{\infty }}}(\delta _{n}^{2})$ uniformly over $\mathcal{A}_{osn}$;
the second inequality is due to Assumption \ref{ass:weak_equiv}. By our
choice of $\kappa _{n}$ the first term in the RHS is zero for large $K$. So
the desired result follows. The fact that $\kappa _{n}$ coincides with $%
\delta _{n}$ follows from the fact that $||\bar{\alpha}_{0,n}-\alpha
_{0}||^{2}\asymp Q(\bar{\alpha}_{0,n})\asymp Q(\Pi _{n}\alpha _{0})$ by
assumption in the Proposition.

Finally, the convergence rate under $||\cdot ||_{s}$ is obtain by applying
the previous result and the definition of $\tau _{n}$. \textit{Q.E.D.}

\medskip

\noindent \textbf{Proof of Theorem \ref{thm:QLR-H1}}: Since $\sup_{h\in
\mathcal{H}}Pen(h)<\infty $, the relevant parameter set is $\mathcal{A}%
^{M}\equiv \{\alpha \in \mathcal{A}:Pen(h)\leq M\}$ with $M=\sup_{h\in
\mathcal{H}}Pen(h)$, which is non-empty and compact (in $||\cdot ||_{s}$)
under Assumptions \ref{ass:sieve}(i)(ii) and \ref{A_3.6}(iii). Let $\mathcal{%
A}^{R,M}=\mathcal{A}^{M}\cap \{\alpha \in \mathcal{A}:\phi (\alpha )=\phi
_{0}\}$. Since $\phi $ is continuous in $||\cdot ||_{s}$, $\mathcal{A}^{R,M}$
is also compact (in $||\cdot ||_{s}$). Note that $\alpha _{0}\in \mathcal{A}%
^{R,M}$ iff the null $H_{0}:\phi (\alpha _{0})=\phi _{0}$ holds.

If $\mathcal{A}^{R,M}$ is empty, then there does not exist any $\alpha \in
\mathcal{A}^{M}$ such that $\phi (\alpha )=\phi _{0}$, and hence it holds
trivially that $\widehat{QLR}_{n}(\phi _{0})\geq nC$ for some $C>0$ wpa1.

If $\mathcal{A}^{R,M}$ is non-empty, under Assumption \ref{ass:sieve}(iii)
we have: $\min_{\alpha \in \mathcal{A}^{R,M}}Q(\alpha )$ is achieved at some
point within $\mathcal{A}^{R,M}$, say, $\overline{\alpha }\in \mathcal{A}%
^{R,M}$. This and Assumption \ref{ass:sieve}(i)(iv) imply that $Q(\overline{%
\alpha })=\min_{\alpha \in \mathcal{A}^{R,M}}Q(\alpha )>0=Q(\alpha _{0})$
under the fixed alternatives $H_{1}:$ $\phi (\alpha _{0})\neq \phi _{0}$.

By definitions of $\widehat{\alpha }_{n}$ and $\Pi _{n}\alpha _{0}$ and
Assumption \ref{ass:rates}(i), we have:%
\begin{equation*}
\widehat{Q}_{n}(\widehat{\alpha }_{n})\leq \widehat{Q}_{n}(\Pi _{n}\alpha
_{0})\leq c_{0}Q(\Pi _{n}\alpha _{0})+o_{P_{Z^{\infty }}}(n^{-1}).
\end{equation*}

Since $M=\sup_{h\in \mathcal{H}}Pen(h)<\infty $, we also have that $\widehat{%
\alpha }_{n}^{R}\in \mathcal{A}_{k(n)}^{R,M}\subseteq \mathcal{A}_{k(n)}^{M}$
wpa1, so by Assumption \ref{ass:rates}(ii), we have:%
\begin{equation*}
\widehat{Q}_{n}(\widehat{\alpha }_{n}^{R})\geq cQ(\widehat{\alpha }%
_{n}^{R})-O_{P_{Z^{\infty }}}(\bar{\delta}_{m,n}^{2})\geq c\times
\min_{\alpha \in \mathcal{A}^{R,M}}Q(\alpha )-O_{P_{Z^{\infty }}}(\bar{\delta%
}_{m,n}^{2}).
\end{equation*}%
Thus%
\begin{equation*}
\widehat{Q}_{n}(\widehat{\alpha }_{n}^{R})-\widehat{Q}_{n}(\widehat{\alpha }%
_{n})\geq c\times \min_{\alpha \in \mathcal{A}^{R,M}}Q(\alpha )-c_{0}Q(\Pi
_{n}\alpha _{0})-o_{P_{Z^{\infty }}}(n^{-1})-O_{P_{Z^{\infty }}}(\bar{\delta}%
_{m,n}^{2})=cQ(\overline{\alpha })+o_{P_{Z^{\infty }}}(1).
\end{equation*}%
Thus under the fixed alternatives $H_{1}:$ $\phi (\alpha _{0})\neq \phi _{0}$%
,%
\begin{equation*}
\frac{\widehat{QLR}_{n}(\phi _{0})}{n}\geq cQ(\overline{\alpha })>0\quad
\text{wpa1.}
\end{equation*}%
\textit{Q.E.D.}

\medskip

\textbf{A consistent variance estimator for optimally weighted PSMD
estimator. }To stress the fact that we consider the optimally weighted PSMD
procedure, we use $v_{n}^{0}$ and $||v_{n}^{0}||_{0}$ to denote the
corresponding $v_{n}^{\ast }$ and $||v_{n}^{\ast }||$ computed using the
optimal weighting matrix $\Sigma =\Sigma _{0}$. That is,
\begin{equation*}
||v_{n}^{0}||_{0}^{2}=E\left[ \left( \frac{dm(X,\alpha _{0})}{d\alpha }%
[v_{n}^{0}]\right) ^{\prime }\Sigma _{0}(X)^{-1}\left( \frac{dm(X,\alpha
_{0})}{d\alpha }[v_{n}^{0}]\right) \right] .
\end{equation*}%
We call the corresponding sieve score, $S_{n,i}^{0}\equiv \left( \frac{%
dm(X_{i},\alpha _{0})}{d\alpha }[v_{n}^{0}]\right) ^{\prime }\Sigma
_{0}(X_{i})^{-1}\rho (Z_{i},\alpha _{0})$, the optimal sieve score. Note
that $||v_{n}^{0}||_{sd}^{2}=Var(S_{n,i}^{0})=||v_{n}^{0}||_{0}^{2}$. By
Theorem \ref{thm:theta_anorm}, $||v_{n}^{0}||_{sd}^{2}=||v_{n}^{0}||_{0}^{2}$
is the variance of the optimally weighted PSMD estimator $\phi (\widehat{%
\alpha }_{n})$. We could compute a consistent estimator $\widehat{%
||v_{n}^{0}||_{0}^{2}}$ of the variance $||v_{n}^{0}||_{0}^{2}$ by looking
at the \textquotedblleft slope\textquotedblright\ of the optimally weighted
criterion $\widehat{Q}_{n}^{0}$:
\begin{equation}
\widehat{||v_{n}^{0}||_{0}^{2}}\equiv \left( \frac{\widehat{Q}_{n}^{0}(%
\widetilde{\alpha }_{n})-\widehat{Q}_{n}^{0}(\widehat{\alpha }_{n})}{%
\varepsilon _{n}^{2}}\right) ^{-1},  \label{ovar-est}
\end{equation}%
where $\widetilde{\alpha }_{n}$ is an approximate minimizer of $\widehat{Q}%
_{n}^{0}(\alpha )$ over $\{\alpha \in \mathcal{A}_{k(n)}\colon \phi (\alpha
)=\phi (\widehat{\alpha }_{n})-\varepsilon _{n}\}$.

\begin{theorem}
\label{thm:est_avar} Let $\widehat{\alpha }_{n}$ be the optimally weighted
PSMD estimator (\ref{psmd}) with $\Sigma =\Sigma _{0}$, and conditions for
Lemma \ref{thm:ThmCONVRATEGRAL}, Assumptions \ref{ass:phi} and \ref{ass:LAQ}
hold with $||v_{n}^{0}||_{sd}=||v_{n}^{0}||_{0}$ and $\left\vert
B_{n}-1\right\vert =o_{P_{Z^{\infty }}}(1)$. Let $cn^{-1/2}\leq \frac{%
\varepsilon _{n}}{||v_{n}^{0}||_{0}}\leq C\delta _{n}$ for finite constants $%
c,C>0$. Then: $\widetilde{\alpha }_{n}\in \mathcal{N}_{osn}$ wpa1-$%
P_{Z^{\infty }}$, and
\begin{equation*}
\frac{\widehat{||v_{n}^{0}||_{0}^{2}}}{||v_{n}^{0}||_{0}^{2}}%
=1+o_{P_{Z^{\infty }}}(1).
\end{equation*}
\end{theorem}

When $\widehat{\alpha }_{n}$ is the optimally weighted PSMD estimator of $%
\alpha _{0}$, Theorem \ref{thm:est_avar} suggests $\widehat{%
||v_{n}^{0}||_{0}^{2}}$ defined in (\ref{ovar-est}) as an alternative
consistent variance estimator for $\phi (\widehat{\alpha }_{n})$. Compared
to Theorems \ref{thm:VE} and \ref{thm:VE2}, this alternative variance
estimator $\widehat{||v_{n}^{0}||_{0}^{2}}$ allows for a non-smooth residual
function $\rho (Z,\alpha )$ (such as the one in NPQIV), but is only valid
for an optimally weighted PSMD estimator.

\medskip

\noindent \textbf{Proof of Theorem \ref{thm:est_avar}} Recall that for the
optimally weighted criterion case $u_{n}^{\ast }=v_{n}^{0}/||v_{n}^{0}||_{0}$%
, and hence $||u_{n}^{\ast }||=1$, $B_{n}=1+o_{P_{Z^{\infty }}}(1)$. To
simplify notation, in this proof we use $\langle \cdot ,\cdot \rangle $, $%
||\cdot ||$ and $\widehat{Q}_{n}(\cdot )$ for the ones corresponding to the
optimal weighting matrix $\Sigma =\Sigma _{0}$.

\textbf{We first show that} $\widetilde{\alpha }_{n}\in \mathcal{N}_{osn}$
\textbf{wpa1}. Recall that $\widetilde{\alpha }_{n}$ is defined as an
approximate optimally weighted PSMD estimator constrained to $\{\alpha \in
\mathcal{A}_{k(n)}\colon \phi (\alpha )=\phi (\widehat{\alpha }%
_{n})-\varepsilon _{n}\}$. In the following since there is no risk of
confusion, we use $_{P}$ instead of $_{P_{Z^{\infty }}}$.

Let $r=\phi (\widehat{\alpha }_{n})-\phi (\alpha _{0})-\varepsilon _{n}$.
Since $\varepsilon _{n}\leq C||v_{n}^{0}||_{0}\delta _{n}$ (by assumption),
and $\widehat{\alpha }_{n}\in \mathcal{N}_{osn}$ wpa1, $\phi (\widehat{%
\alpha }_{n})-\phi (\alpha _{0})=O_{P}(||v_{n}^{0}||_{0}/\sqrt{n})$ (by
Theorem \ref{thm:theta_anorm}), we have $|r|\leq C||v_{n}^{0}||_{0}O(\delta
_{n}+n^{-1/2})\leq C||v_{n}^{0}||_{0}\delta _{n}$ for some $C>0$. Also note
that $||\widehat{\alpha }_{n}-\alpha _{0}||\leq C\delta _{n}$ wpa1. Thus, by
Lemma \ref{pro:phi_solve}(2), there exists a $t_{n}^{\ast }\in \mathcal{T}%
_{n}$ such that $\phi (\widehat{\alpha }_{n}(t_{n}^{\ast }))=\phi (\widehat{%
\alpha }_{n})-\varepsilon _{n}$ and $\widehat{\alpha }_{n}(t_{n}^{\ast })=%
\widehat{\alpha }_{n}+t_{n}^{\ast }u_{n}^{\ast }\in \mathcal{A}_{k(n)}$ and $%
t_{n}^{\ast }=O(\delta _{n})$. Henceforth, let $\bar{\alpha}_{n}\equiv
\widehat{\alpha }_{n}(t_{n}^{\ast })$. Observe that
\begin{equation*}
||\bar{\alpha}_{n}-\alpha _{0}||\leq \delta _{n}+t_{n}^{\ast }=O(\delta _{n})%
\text{,}
\end{equation*}%
and
\begin{equation*}
||\bar{\alpha}_{n}-\alpha _{0}||_{s}\leq ||\widehat{\alpha }_{n}-\alpha
_{0}||_{s}+t_{n}^{\ast }||u_{n}^{\ast }||_{s}\leq \delta _{s,n}+t_{n}^{\ast
}\tau _{n}
\end{equation*}%
which is of order $\delta _{s,n}$. Therefore, $\bar{\alpha}_{n}$ satifies:
(a) $\bar{\alpha}_{n}\in \mathcal{N}_{osn}$ wpa1., and $t_{n}^{\ast }\in
\mathcal{T}_{n}$ with $t_{n}^{\ast }=O(\delta _{n})$; and (b) $\bar{\alpha}%
_{n}\in \{\alpha \in \mathcal{A}_{k(n)}\colon \phi (\alpha )=\phi (\widehat{%
\alpha }_{n})-\varepsilon _{n}\}$.

We now establish the consistency of $\widetilde{\alpha }_{n}$ using the
properties of $\bar{\alpha}_{n}$. We observe that, for any $\epsilon >0$,
\begin{equation*}
\Pr (||\widetilde{\alpha }_{n}-\alpha _{0}||_{s}\geq \epsilon )\leq \Pr
\left( \inf_{\mathcal{B}_{n}:||\alpha -\alpha _{0}||_{s}\geq \epsilon }%
\widehat{Q}_{n}(\alpha )\leq \widehat{Q}_{n}(\bar{\alpha}_{n})+o(n^{-1})+%
\lambda _{n}Pen(\bar{h}_{n})\right)
\end{equation*}%
where $\mathcal{B}_{n}\equiv \{\alpha \in \mathcal{A}_{k(n)}^{M_{0}}:\phi
(\alpha )=\phi (\widehat{\alpha }_{n})-\varepsilon _{n}\}$ and the
inequality is valid because $\bar{\alpha}_{n}\in \mathcal{B}_{n}$ by (a) and
(b). Under (a) and Lemma \ref{thm:ThmCONVRATEGRAL}, $\lambda _{n}Pen(\bar{h}%
_{n})=O_{P}(\lambda _{n})=o(n^{-1})$.

By (a), under assumption \ref{ass:LAQ}(i)
\begin{equation*}
\widehat{Q}_{n}(\bar{\alpha}_{n})=\widehat{Q}_{n}(\widehat{\alpha }%
_{n})+t_{n}^{\ast }\{\mathbb{Z}_{n}+\langle u_{n}^{\ast },\widehat{\alpha }%
_{n}-\alpha _{0}\rangle \}+0.5(t_{n}^{\ast })^{2}+o_{P}\left( t_{n}^{\ast
}n^{-1/2}+(t_{n}^{\ast })^{2}+o(n^{-1})\right) .
\end{equation*}

By Lemma \ref{lem:theta_anorm(1)}, $\mathbb{Z}_{n}+\langle u_{n}^{\ast },%
\widehat{\alpha }_{n}-\alpha _{0}\rangle =o_{P}(n^{-1/2})$ and thus, given
that $t_{n}^{\ast }=O(\delta _{n})$, the previous display implies that
\begin{equation*}
\widehat{Q}_{n}(\bar{\alpha}_{n})\leq \widehat{Q}_{n}(\widehat{\alpha }%
_{n})+o_{P}(n^{-1/2}\delta _{n}+\delta _{n}^{2}+o(n^{-1}))\leq O_{P}(\delta
_{n}^{2})
\end{equation*}%
Therefore,
\begin{equation*}
\Pr (||\widetilde{\alpha }_{n}-\alpha _{0}||_{s}\geq \epsilon )\leq \Pr
\left( \inf_{\mathcal{B}_{n}:||\alpha -\alpha _{0}||_{s}\geq \epsilon }%
\widehat{Q}_{n}(\alpha )\leq \widehat{Q}_{n}(\widehat{\alpha }%
_{n})+O(\lambda _{n}+\delta _{n}^{2})\right) .
\end{equation*}%
Since $\widehat{Q}_{n}(\widehat{\alpha }_{n})\leq \widehat{Q}_{n}(\Pi
_{n}\alpha _{0})+O(\lambda _{n})$ by definition of $\widehat{\alpha }_{n}$
and from the fact that $\mathcal{B}_{n}\subseteq \mathcal{A}_{k(n)}^{M_{0}}$%
, it follows that
\begin{equation*}
\Pr (||\widetilde{\alpha }_{n}-\alpha _{0}||_{s}\geq \epsilon )\leq \Pr
\left( \inf_{\mathcal{A}_{n}^{M_{0}}:||\alpha -\alpha _{0}||_{s}\geq
\epsilon }\widehat{Q}_{n}(\alpha )\leq \widehat{Q}_{n}(\Pi _{n}\alpha
_{0})+O(\lambda _{n}+\delta _{n}^{2})\right) .
\end{equation*}%
The rest of the consistency proof follows from identical steps to the
standard one; see \cite{CP_WP07a}.

In order to show the rate, by similar arguments to the previous ones
\begin{equation*}
\widehat{Q}_{n}(\widetilde{\alpha }_{n})\leq \widehat{Q}_{n}(\Pi _{n}\alpha
_{0})+O(\lambda _{n}+\delta _{n}^{2}),
\end{equation*}%
under our assumptions $\widehat{Q}_{n}(\widetilde{\alpha }_{n})\geq c||%
\widetilde{\alpha }_{n}-\alpha _{0}||^{2}-O_{P}(\delta _{n}^{2})$ and $%
\widehat{Q}_{n}(\Pi _{n}\alpha _{0})\leq c_{0}Q(\Pi _{n}\alpha
_{0})+o_{P}(n^{-1})$, so the desired rate under $||\cdot ||$ follows. The
rate under $||\cdot ||_{s}$ immediately follows using the definition of
sieve measure of local ill-posedness $\tau _{n}$. Thus $\widetilde{\alpha }%
_{n}\in \mathcal{N}_{osn}$ wpa1.

\textbf{We now show that} $\frac{\widehat{||v_{n}^{0}||_{0}^{2}}}{%
||v_{n}^{0}||_{0}^{2}}=1+o_{P_{Z^{\infty }}}(1)$. This part of proof
consists of several steps that are similar to those in the proof of Theorem %
\ref{thm:chi2}, and hence we omit some details. We first provide an
asymptotic expansion for $n(\widehat{Q}_{n}(\widetilde{\alpha }_{n})-%
\widehat{Q}_{n}(\widehat{\alpha }_{n}))$ using Assumption \ref{ass:LAQ}(i)
(with $B_{n}=1+o_{P_{Z^{\infty }}}(1)$), and then show that this is enough
to establish the desired result.

In the following we let $t_{n}\equiv \varepsilon _{n}/||v_{n}^{0}||_{0}$. By
the assumption on $\varepsilon _{n}$ we have: $cn^{-1/2}\leq t_{n}\leq
C\delta _{n}$. Therefore, $t_{n}\in \mathcal{T}_{n}$, $t_{n}=o_{P_{Z^{\infty
}}}(1)$ and $o_{P_{Z^{\infty }}}\left( \frac{1}{t_{n}}n^{-1/2}\right)
=o_{P_{Z^{\infty }}}(1)$.

\textsc{Step 1:} First, we note that $\widehat{\alpha }_{n}\in \mathcal{N}%
_{osn}$ wpa1, that $-t_{n}\in \mathcal{T}_{n}$ and $\widehat{\alpha }%
_{n}\left( -t_{n}\right) \in \mathcal{A}_{k(n)}$. So we can apply Assumption %
\ref{ass:LAQ}(i) with $\alpha =\widehat{\alpha }_{n}$ and $-t_{n}$ as the
direction, and obtain:
\begin{eqnarray}\notag
\frac{(\widehat{Q}_{n}(\widehat{\alpha }_{n}(-t_{n}))-\widehat{Q}_{n}(%
\widehat{\alpha }_{n}))}{t_{n}^{2}}& = &\frac{-2}{t_{n}}\{\mathbb{Z}%
_{n}+\langle u_{n}^{\ast },\widehat{\alpha }_{n}-\alpha _{0}\rangle
\}+1\\
&& +o_{P}\left( \max \left\{ 1,\frac{n^{-1/2}}{t_{n}},\frac{o(n^{-1})}{%
t_{n}^{2}}\right\} \right)  \notag \\
& = & 1+o_{P_{Z^{\infty }}}\left( 1\right) ,  \label{eqn:estvar_eq1}
\end{eqnarray}%
where the last equality follows from the fact that $\langle u_{n}^{\ast },%
\widehat{\alpha }_{n}-\alpha _{0}\rangle +\mathbb{Z}_{n}=o_{P_{Z^{\infty
}}}(n^{-1/2})$ (by Lemma \ref{lem:theta_anorm(1)}), and that $%
o_{P_{Z^{\infty }}}\left( \frac{1}{t_{n}}n^{-1/2}\right) =o_{P_{Z^{\infty
}}}(1)$ (by our choice of $t_{n}$).

\textsc{Step 2:} Since $\widetilde{\alpha }_{n}\in \mathcal{N}_{osn}$ wpa1, $%
t_{n}\in \mathcal{T}_{n}$ and $\widetilde{\alpha }_{n}(t_{n})\in \mathcal{A}%
_{k(n)}$, we can apply Assumption \ref{ass:LAQ}(i) with $\alpha =\widetilde{%
\alpha }_{n}$ and $t_{n}$ as the direction, and obtain:
\begin{eqnarray}\notag
\frac{(\widehat{Q}_{n}(\widetilde{\alpha }_{n}(t_{n}))-\widehat{Q}_{n}(%
\widetilde{\alpha }_{n}))}{t_{n}^{2}}& = &\frac{2}{t_{n}}\{\mathbb{Z}%
_{n}+\langle u_{n}^{\ast },\widetilde{\alpha }_{n}-\alpha _{0}\rangle
\}+1\\
&& +o_{P}\left( \max \left\{ 1,\frac{n^{-1/2}}{t_{n}},\frac{o(n^{-1})}{%
t_{n}^{2}}\right\} \right)  \notag \\
& = &-1+o_{P_{Z^{\infty }}}\left( 1\right) ,  \label{eqn:estvar_eq2}
\end{eqnarray}%
where the last line follows from the definition of the restricted estimator $%
\widetilde{\alpha }_{n}$. This is because $\phi (\widetilde{\alpha }%
_{n})=\phi (\widehat{\alpha }_{n})-\varepsilon _{n}$, by Assumptions \ref%
{ass:phi}(i)(ii),
\begin{equation*}
\left\vert -\varepsilon _{n}-\frac{d\phi (\alpha _{0})}{d\alpha }[\widehat{%
\alpha }_{n}-\widetilde{\alpha }_{n}]\right\vert =o_{P_{Z^{\infty
}}}(||v_{n}^{0}||_{0}/\sqrt{n}).
\end{equation*}%
Hence $\langle v_{n}^{0},\widetilde{\alpha }_{n}-\alpha _{0}\rangle =\langle
v_{n}^{0},\widehat{\alpha }_{n}-\alpha _{0}\rangle -\varepsilon
_{n}+o_{P_{Z^{\infty }}}(||v_{n}^{0}||_{0}/\sqrt{n})$. This implies that $%
\mathbb{Z}_{n}+\langle u_{n}^{\ast },\widetilde{\alpha }_{n}-\alpha
_{0}\rangle =-\frac{\varepsilon _{n}}{||v_{n}^{0}||_{0}}+o_{P_{Z^{\infty
}}}(n^{-1/2})=-t_{n}+o_{P_{Z^{\infty }}}(n^{-1/2})$.

\textsc{Step 3:} It is easy to see that, from equation (\ref{eqn:estvar_eq2}%
) and by the definition of $\hat{\alpha}_{n}$,
\begin{equation*}
\frac{(\widehat{Q}_{n}(\widetilde{\alpha }_{n})-\widehat{Q}_{n}(\widehat{%
\alpha }_{n}))}{t_{n}^{2}}\geq \frac{(\widehat{Q}_{n}(\widetilde{\alpha }%
_{n}))-\widehat{Q}_{n}(\widetilde{\alpha }_{n}(t_{n}))}{t_{n}^{2}}%
-o_{P_{Z^{\infty }}}(1)=1+o_{P_{Z^{\infty }}}(1).
\end{equation*}

Also, from equation (\ref{eqn:estvar_eq1}), Assumption \ref{ass:LAQ}(i) and
by the definition of $\tilde{\alpha}_{n}$,
\begin{eqnarray*}
\frac{(\widehat{Q}_{n}(\widetilde{\alpha }_{n})-\widehat{Q}_{n}(\widehat{%
\alpha }_{n}))}{t_{n}^{2}} &\leq &\frac{(\widehat{Q}_{n}(\widehat{\alpha }%
_{n}(t_{n}^{\ast }))-\widehat{Q}_{n}(\widehat{\alpha }_{n}))}{t_{n}^{2}}%
+o_{P_{Z^{\infty }}}(1) \\
&=&\frac{2t_{n}^{\ast }\{\mathbb{Z}_{n}+\langle u_{n}^{\ast },\widehat{%
\alpha }_{n}-\alpha _{0}\rangle \} + \left( t_{n}^{\ast }\right)
^{2}}{t_{n}^{2}}\\
&& + t^{-2}_{n} o_{P}\left( \max \left\{ \left( t_{n}^{\ast }\right) ^{2},t_{n}^{\ast
}n^{-\frac{1}{2}},o(n^{-1})\right\} \right) +o_{P}(1) \\
&=&\frac{-2}{t_{n}}\{\mathbb{Z}_{n}+\langle u_{n}^{\ast },\widehat{\alpha }%
_{n}-\alpha _{0}\rangle \}+1+o_{P_{Z^{\infty }}}(1) \\
&=&1+o_{P_{Z^{\infty }}}(1),
\end{eqnarray*}%
provided that there is a $t_{n}^{\ast }\in \mathcal{T}_{n}$ such that (3a) $%
\phi (\widehat{\alpha }_{n}(t_{n}^{\ast }))=\phi (\widehat{\alpha }%
_{n})-\varepsilon _{n}$ and (3b) $t_{n}^{\ast }/t_{n}=-1+o_{P_{Z^{\infty
}}}(1)$. In Step 5 we verify that such a $t_{n}^{\ast }$ exists.

By putting these inequalities together, it follows
\begin{equation}
||v_{n}^{0}||_{0}^{2}\frac{\widehat{Q}_{n}(\widetilde{\alpha }_{n})-\widehat{%
Q}_{n}(\widehat{\alpha }_{n})}{\varepsilon _{n}^{2}}=\frac{(\widehat{Q}_{n}(%
\widetilde{\alpha }_{n})-\widehat{Q}_{n}(\widehat{\alpha }_{n}))}{t_{n}^{2}}%
=1+o_{P_{Z^{\infty }}}(1).  \label{eqn:estvar_1}
\end{equation}

\textsc{Step 4:} By equation (\ref{eqn:estvar_1}) we have:%
\begin{equation*}
\frac{||v_{n}^{0}||_{0}^{2}}{\widehat{||v_{n}^{0}||_{0}^{2}}}%
=1+o_{P_{Z^{\infty }}}(1)\text{,\quad with\quad }\widehat{%
||v_{n}^{0}||_{0}^{2}}\equiv \left( \frac{\widehat{Q}_{n}(\widetilde{\alpha }%
_{n})-\widehat{Q}_{n}(\widehat{\alpha }_{n})}{\varepsilon _{n}^{2}}\right)
^{-1},
\end{equation*}%
which implies that $0.5\leq \frac{||v_{n}^{0}||_{0}^{2}}{\widehat{%
||v_{n}^{0}||_{0}^{2}}}\leq 1.5$ with probability $P_{Z^{\infty }}$
approaching one. By continuous mapping theorem, we obtain:%
\begin{equation*}
\frac{\widehat{||v_{n}^{0}||_{0}^{2}}}{||v_{n}^{0}||_{0}^{2}}%
=1+o_{P_{Z^{\infty }}}(1).
\end{equation*}

\textsc{Step 5:} We now show that there is a $t_{n}^{\ast }\in \mathcal{T}%
_{n}$ such that (3a) and (3b) in Step 3 hold. Denote $r\equiv \phi (\widehat{%
\alpha }_{n})-\phi (\alpha _{0})-\varepsilon _{n}$. Since $\varepsilon
_{n}\leq C||v_{n}^{0}||_{0}\delta _{n}$, and $\widehat{\alpha }_{n}\in
\mathcal{N}_{osn}$ wpa1, $\phi (\widehat{\alpha }_{n})-\phi (\alpha
_{0})=O_{P}(||v_{n}^{0}||_{0}/\sqrt{n})$ (by Theorem \ref{thm:theta_anorm}),
we have $|r|\leq ||v_{n}^{0}||_{0}\delta _{n}(M_{n}+C)\leq
2M_{n}||v_{n}^{0}||_{0}\delta _{n}$ (since $C<M_{n}$ eventually). Thus, by
Lemma \ref{pro:phi_solve}, there exists a $t_{n}^{\ast }\in \mathcal{T}_{n}$
such that $\phi (\widehat{\alpha }_{n}(t_{n}^{\ast }))=\phi (\widehat{\alpha
}_{n})-\varepsilon _{n}$ and $\widehat{\alpha }_{n}(t_{n}^{\ast })=\widehat{%
\alpha }_{n}+t_{n}^{\ast }u_{n}^{\ast }\in \mathcal{A}_{k(n)}$, and hence
(3a) holds. Moreover, by Assumption \ref{ass:phi}(i)(ii), such a choice of $%
t_{n}^{\ast }$ also satisfies

\begin{equation*}
\left\vert \underbrace{\phi (\widehat{\alpha }_{n}(t_{n}^{\ast }))-\phi (%
\widehat{\alpha }_{n})}_{=-\varepsilon _{n}}-\frac{d\phi (\alpha _{0})}{%
d\alpha }[t_{n}^{\ast }u_{n}^{\ast }]\right\vert =o_{P_{Z^{\infty
}}}(||v_{n}^{0}||_{0}n^{-1/2}).
\end{equation*}%
Since $u_{n}^{\ast }=v_{n}^{0}/||v_{n}^{0}||_{0}$ for optimally weighted
criterion case, we have: $\frac{d\phi (\alpha _{0})}{d\alpha }[u_{n}^{\ast
}]=||v_{n}^{0}||_{0}$. Thus%
\begin{equation*}
\left\vert -\varepsilon _{n}-t_{n}^{\ast }||v_{n}^{0}||_{0}\right\vert
=o_{P_{Z^{\infty }}}(||v_{n}^{0}||_{0}n^{-1/2}).
\end{equation*}%
Since $t_{n}\equiv \varepsilon _{n}/||v_{n}^{0}||_{0}$, we obtain: $%
\left\vert -t_{n}-t_{n}^{\ast }\right\vert =o_{P_{Z^{\infty }}}(n^{-1/2})$,
and hence%
\begin{equation*}
\left\vert (t_{n}^{\ast }/t_{n})+1\right\vert =o_{P_{Z^{\infty
}}}(n^{-1/2}/t_{n})=o_{P_{Z^{\infty }}}(1)
\end{equation*}%
due to the fact that $cn^{-1/2}\leq t_{n}\leq C\delta _{n}$. Thus (3b)
holds. \textit{Q.E.D.}

\subsection{Proofs for Section \protect\ref{sec:bootstrap} on bootstrap
inference}

Throughout the Appendices, we sometimes use the simplified term
\textquotedblleft wpa1\textquotedblright\ in the bootstrap world while its
precise meaning is given in Section \ref{sec:bootstrap}.

Recall that $\mathbb{Z}_{n}^{\omega }\equiv \frac{1}{n}\sum_{i=1}^{n}\omega
_{i}g(X_{i},u_{n}^{\ast })\rho (Z_{i},\alpha _{0})$ with $%
g(X_{i},u_{n}^{\ast })\equiv \left( \frac{dm(X_{i},\alpha _{0})}{d\alpha }%
[u_{n}^{\ast }]\right) ^{\prime }\Sigma (X_{i})^{-1}$.

\begin{lemma}
\label{lem:LAR_boot} Let $\widehat{\alpha }_{n}^{B}$ be the bootstrap PSMD
estimator and conditions for Lemma \ref{thm:ThmCONVRATEGRAL} and Lemma \ref%
{lem:cons_boot} hold. Let Assumption \ref{ass:LAQ_B}(i) hold. Then: (1) for
all $\delta >0$, there exists a $N(\delta )$ such that for all $n\geq
N(\delta )$,
\begin{equation*}
P_{Z^{\infty }}\left( P_{V^{\infty }|Z^{\infty }}\left( \sqrt{n}\left\vert
\langle u_{n}^{\ast },\widehat{\alpha }_{n}^{B}-\alpha _{0}\rangle +\mathbb{Z%
}_{n}^{\omega }\right\vert \geq \delta |Z^{n}\right) <\delta \right) \geq
1-\delta \text{.}
\end{equation*}%
(2) If, in addition, assumptions of Lemma \ref{lem:theta_anorm(1)} hold,
then
\begin{equation*}
\sqrt{n}\langle u_{n}^{\ast },\widehat{\alpha }_{n}^{B}-\widehat{\alpha }%
_{n}\rangle =-\sqrt{n}\mathbb{Z}_{n}^{\omega -1}+o_{P_{V^{\infty }|Z^{\infty
}}}(1)~wpa1(P_{Z^{\infty }}).
\end{equation*}
\end{lemma}

\smallskip

\noindent \textbf{Proof of Lemma \ref{lem:LAR_boot}:} The proof is very
similar to that of Lemma \ref{lem:theta_anorm(1)}, so we only present the
main steps.

For \textbf{Result (1)}. Under Assumption \ref{ass:LAQ_B}(i) and using the
fact that $\widehat{\alpha }_{n}^{B}$ is an approximate minimizer of $%
\widehat{Q}_{n}^{B}(\alpha )+\lambda _{n}Pen(h)$ on $\mathcal{A}_{k(n)}$ and
$\widehat{\alpha }_{n}^{B}\in \mathcal{N}_{osn}$ wpa1, it follows (see the
proof of Lemma \ref{lem:theta_anorm(1)} for details), for sufficiently large
$n$,
\begin{equation*}
P_{Z^{\infty }}\left( P_{V^{\infty }|Z^{\infty }}\left( 2\epsilon _{n}\{%
\mathbb{Z}_{n}^{\omega }+\langle u_{n}^{\ast },\widehat{\alpha }%
_{n}^{B}-\alpha _{0}\rangle \}+\epsilon _{n}^{2}B_{n}^{\omega }+E_{n}(%
\widehat{\alpha }_{n}^{B},\epsilon _{n})\geq -\delta r_{n}^{-1}|Z^{n}\right)
\geq 1-\delta \right) >1-\delta ,
\end{equation*}%
where $r_{n}$ and $E_{n}$ are defined as in the proof of Lemma \ref%
{lem:theta_anorm(1)}, and $\epsilon _{n}=\pm \{s_{n}^{-1/2}+o(n^{-1/2})\}$.
Dividing by $2\epsilon _{n}$ and multiplying by $\sqrt{n}$, it follows that
\begin{equation*}
P_{Z^{\infty }}\left( P_{V^{\infty }|Z^{\infty }}\left( A_{n,\delta
}^{\omega }\geq \sqrt{n}\{\mathbb{Z}_{n}^{\omega }+\langle u_{n}^{\ast },%
\widehat{\alpha }_{n}^{B}-\alpha _{0}\rangle \}\geq B_{n,\delta }^{\omega
}|Z^{n}\right) \geq 1-\delta \right) >1-\delta
\end{equation*}%
eventually, where
\begin{align*}
A_{n,\delta }^{\omega }& \equiv -0.5\sqrt{n}\epsilon _{n}B_{n}^{\omega
}-\delta \sqrt{n}\epsilon _{n}^{-1}r_{n}^{-1}+0.5\delta \\
B_{n,\delta }^{\omega }& \equiv -0.5\sqrt{n}\epsilon _{n}B_{n}^{\omega
}-\delta \sqrt{n}\epsilon _{n}^{-1}r_{n}^{-1}-0.5\delta .
\end{align*}

Since $\sqrt{n}\epsilon _{n}=o(1)$ and $B_{n}^{\omega }=O_{P_{V^{\infty
}|Z^{\infty }}}(1)$ wpa1($P_{Z^{\infty }}$) and $|\sqrt{n}\epsilon
^{-1}_{n}r_{n}^{-1}|\asymp 1$, it follows, for sufficiently large $n$,
\begin{equation*}
P_{Z^{\infty }}\left( P_{V^{\infty }|Z^{\infty }}\left( 2\delta \geq \sqrt{n}%
\{\mathbb{Z}_{n}^{\omega }+\langle u_{n}^{\ast },\widehat{\alpha }%
_{n}^{B}-\alpha _{0}\rangle \}\geq -2\delta |Z^{n}\right) \geq 1-\delta
\right) >1-\delta .
\end{equation*}%
Or equivalently, for sufficiently large $n$,
\begin{equation*}
P_{Z^{\infty }}\left( P_{V^{\infty }|Z^{\infty }}\left( \left\vert \sqrt{n}\{%
\mathbb{Z}_{n}^{\omega }+\langle u_{n}^{\ast },\widehat{\alpha }%
_{n}^{B}-\alpha _{0}\rangle \}\right\vert \geq 2\delta |Z^{n}\right) <\delta
\right) \geq 1-\delta .
\end{equation*}

\textbf{Result (2)} directly follows from Result (1) and Lemma \ref%
{lem:theta_anorm(1)}. \textit{Q.E.D.}

\medskip

\noindent \textbf{Proof of Theorem \ref{thm:VE-boot2}} We note that
Assumption \ref{ass:VE_boot2} implies that $|n^{-1}\sum_{i=1}^{n}\widehat{T}%
_{i}[v_{n}]^{\prime }\hat{M}_{i}^{B}\widehat{T}_{i}[v_{n}]-\sigma _{\omega
}^{2}n^{-1}\sum_{i=1}^{n}\hat{T}_{i}[v_{n}]^{\prime }\hat{M}_{i}\hat{T}%
_{i}[v_{n}]|=o_{P_{V^{\infty }|Z^{\infty }}}(1)$, uniformly over $v_{n}\in
\overline{\mathbf{V}}_{k(n)}^{1}$ with $\hat{M}_{i}=\widehat{M}(Z_{i},%
\widehat{\alpha }_{n})$ and $\widehat{T}_{i}[v_{n}]\equiv \frac{d\hat{m}%
(X_{i},\hat{\alpha}_{n})}{d\alpha }[v_{n}]$. The rest of the proof follows
directly from that of Theorem \ref{thm:VE}(1) for the sieve variance defined
in (\ref{svar-hat1}) case. \textit{Q.E.D.}

\medskip

\noindent \textbf{Proof of Theorem \ref{thm:bootstrap_2}} By Lemma \ref%
{lem:LAR_boot} and steps analogous to those used to show Theorem \ref%
{thm:theta_anorm}, it follows%
\begin{equation}
\sqrt{n}\frac{\phi (\widehat{\alpha }_{n}^{B})-\phi (\widehat{\alpha }_{n})}{%
\sigma _{\omega }||v_{n}^{\ast }||_{sd}}=-\sqrt{n}\frac{\mathbb{Z}%
_{n}^{\omega -1}}{\sigma _{\omega }}+o_{P_{V^{\infty }|Z^{\infty
}}}(1)~wpa1(P_{Z^{\infty }}).  \label{b-regular}
\end{equation}

For \textbf{Result (1)}, we note that the result for $\widehat{W}_{2,n}^{B}$
follows directly from Theorem \ref{thm:VE-boot2} and the proof of the Result
(1) for $\widehat{W}_{1,n}^{B}\equiv \sqrt{n}\frac{\phi (\widehat{\alpha }%
_{n}^{B})-\phi (\widehat{\alpha }_{n})}{\sigma _{\omega }||\widehat{v}%
_{n}^{\ast }||_{n,sd}}$.

In fact, for both $j=1,2$, Theorem \ref{thm:VE}(1), equation (\ref{b-regular}%
) and Theorem \ref{thm:VE-boot2} imply that%
\begin{equation}
\widehat{W}_{j,n}^{B}=-\sqrt{n}\frac{\mathbb{Z}_{n}^{\omega -1}}{\sigma
_{\omega }}+o_{P_{V^{\infty }|Z^{\infty }}}(1)~wpa1(P_{Z^{\infty }});
\label{old-result1}
\end{equation}%
Equation (\ref{old-result1}) and Assumptions \ref{ass:LAQ}(ii) and \ref%
{ass:LAQ_B}(ii) imply that:
\begin{equation*}
\left\vert \mathcal{L}_{V^{\infty }|Z^{\infty }}\left( \widehat{W}%
_{j,n}^{B}\mid Z^{n}\right) -\mathcal{L}\left( \widehat{W}_{n}\right)
\right\vert =o_{P_{Z^{\infty }}}(1).
\end{equation*}%
Result (1) now follows from the following two equations:
\begin{equation}
\sup_{t\in \mathbb{R}}|P_{V^{\infty }|Z^{\infty }}(\widehat{W}_{j,n}^{B}\leq
t|Z^{n})-\Phi (t)|=o_{P_{V^{\infty }|Z^{\infty }}}(1)~wpa1(P_{Z^{\infty }}),
\label{Polya1}
\end{equation}%
and
\begin{equation}
\sup_{t\in \mathbb{R}}|P_{Z^{\infty }}(\widehat{W}_{n}\leq t)-\Phi
(t)|=o_{P_{Z^{\infty }}}(1),  \label{Polya2}
\end{equation}%
where $\Phi ()$ is the cdf of a standard normal. Equation (\ref{Polya2})
follows directly from Theorem \ref{thm:VE}(2) and Polya's theorem (see e.g.,
\cite{Bickel_SS1992}). Equation (\ref{Polya1}) follows by the same arguments
in Lemma 10.11 in \cite{Kosorok} (which are in turn analogous to those used
in the proof of Polya's theorem).

\textbf{Result (2)} follows from equation (\ref{b-regular}) and the fact
that $||v_{n}^{\ast }||_{sd}\rightarrow ||v^{\ast }||_{sd}\in (0,\infty )$
for regular functionals. \textit{Q.E.D.}

\noindent \textbf{Proof of Theorem \ref{thm:bootstrap}}: For \textbf{Result
(1),} denote
\begin{eqnarray*}
\mathcal{F}_{n}& \equiv & n\frac{\inf_{\mathcal{A}_{k(n)}(\widehat{\phi }_{n})}%
\widehat{Q}_{n}^{B}(\alpha )-\widehat{Q}_{n}^{B}(\widehat{\alpha }_{n}^{B})}{%
\sigma _{\omega }^{2}}=\frac{\widehat{QLR}_{n}^{B}(\widehat{\phi }_{n})}{%
\sigma _{\omega }^{2}}\\
&=& n\frac{\widehat{Q}_{n}^{B}(\widehat{\alpha }%
_{n}^{R,B})-\widehat{Q}_{n}^{B}(\widehat{\alpha }_{n}^{B})}{\sigma _{\omega
}^{2}}+o_{P_{V^{\infty }|Z^{\infty }}}(1)~wpa1(P_{Z^{\infty }})
\end{eqnarray*}%
where $\mathcal{A}_{k(n)}(\widehat{\phi }_{n})\equiv \{\alpha \in \mathcal{A}%
_{k(n)}\colon \phi (\alpha )=\phi (\widehat{\alpha }_{n})\}$. Since $%
o_{P_{V^{\infty }|Z^{\infty }}}(1)~wpa1(P_{Z^{\infty }})$ will not affect
the asymptotic results we omit it from the rest of the proof to ease the notational
burden. We want to show that for all $\delta >0$, there exists a $N(\delta )$
such that
\begin{equation*}
P_{Z^{\infty }}\left( P_{V^{\infty }|Z^{\infty }}\left( \left\vert \mathcal{F%
}_{n}-\left( \sqrt{n}\frac{\mathbb{Z}_{n}^{\omega -1}}{\sigma _{\omega
}||u_{n}^{\ast }||}\right) ^{2}\right\vert \geq \delta \mid Z^{n}\right)
<\delta \right) \geq 1-\delta
\end{equation*}%
for all $n\geq N(\delta )$. We divide the proof in several steps.\newline

\textsc{Step 1.} By assumption $|B_{n}^{\omega }-||u_{n}^{\ast
}||^{2}|=o_{P_{V^{\infty }|Z^{\infty }}}(1)~wpa1(P_{Z^{\infty }})$ and $%
||u_{n}^{\ast }||\in (c,C)$, we have: $\left\vert \frac{||u_{n}^{\ast }||^{2}%
}{B_{n}^{\omega }}-1\right\vert =o_{P_{V^{\infty }|Z^{\infty
}}}(1)~wpa1(P_{Z^{\infty }})$. Therefore, it suffices to show that
\begin{equation}
P_{Z^{\infty }}\left( P_{V^{\infty }|Z^{\infty }}\left( \left\vert \mathcal{F%
}_{n}-\left( \sqrt{n}\frac{\mathbb{Z}_{n}^{\omega -1}}{\sigma _{\omega }%
\sqrt{B_{n}^{\omega }}}\right) ^{2}\right\vert \geq \delta \mid Z^{n}\right)
<\delta \right) \geq 1-\delta  \label{eqn:boot_00}
\end{equation}%
eventually.

\textsc{Step 2.} By Assumption \ref{ass:LAQ_B}(i), for all $\delta >0$,
there is a $M>0$ such that
\begin{equation*}
P_{Z^{\infty }}\left( P_{V^{\infty }|Z^{\infty }}\left( \sqrt{n}|\mathbb{Z}%
_{n}^{\omega -1}/B_{n}^{\omega }|\geq M\mid Z^{n}\right) <\delta \right)
\geq 1-\delta
\end{equation*}%
eventually. Thus $t_{n}=-\mathbb{Z}_{n}^{\omega -1}/B_{n}^{\omega }\in
\mathcal{T}_{n}$ wpa1. By the definition of $\widehat{\alpha }_{n}^{B}$, and
the fact that $\widehat{\alpha }_{n}^{R,B}\in \mathcal{N}_{osn}$ wpa1 (by
Lemma \ref{lem:cons_boot}(3)),
\begin{equation*}
\mathcal{F}_{n}\geq n\frac{\widehat{Q}_{n}^{B}(\widehat{\alpha }_{n}^{R,B})-%
\widehat{Q}_{n}^{B}(\widehat{\alpha }_{n}^{R,B}(t_{n}))}{\sigma _{\omega
}^{2}}-o_{P_{V^{\infty }|Z^{\infty }}}(1)~wpa1(P_{Z^{\infty }}).
\end{equation*}

By specializing Assumption \ref{ass:LAQ_B}(i) to $\alpha =\widehat{\alpha }%
_{n}^{R,B}$ and $t_{n}=-\mathbb{Z}_{n}^{\omega -1}/B_{n}^{\omega }$, it
follows
\begin{eqnarray}
&&0.5(\widehat{Q}_{n}^{B}(\widehat{\alpha }_{n}^{R,B}(-\frac{\mathbb{Z}%
_{n}^{\omega -1}}{B_{n}^{\omega }}))-\widehat{Q}_{n}^{B}(\widehat{\alpha }%
_{n}^{R,B}))  \label{eqn:boot_1} \\
&=&-\frac{\mathbb{Z}_{n}^{\omega -1}}{B_{n}^{\omega }}\{\mathbb{Z}%
_{n}^{\omega }+\langle u_{n}^{\ast },\widehat{\alpha }_{n}^{R,B}-\alpha
_{0}\rangle \}+\frac{(\mathbb{Z}_{n}^{\omega -1})^{2}}{2B_{n}^{\omega }}%
+o_{P_{V^{\infty }|Z^{\infty }}}(r_{n}^{-1})~wpa1(P_{Z^{\infty }}).  \notag
\end{eqnarray}

By Assumption \ref{ass:phi}(i)(ii), and the fact that $\widehat{\alpha }%
_{n}^{R,B}\in \mathcal{N}_{osn}$ wpa1,
\begin{equation*}
P_{Z^{\infty }}\left( P_{V^{\infty }|Z^{\infty }}\left( \frac{\sqrt{n}}{%
||v_{n}^{\ast }||}\left\vert \underbrace{\phi (\widehat{\alpha }%
_{n}^{R,B})-\phi (\widehat{\alpha }_{n})}_{=0}-\frac{d\phi (\alpha _{0})}{%
d\alpha }[\widehat{\alpha }_{n}^{R,B}-\widehat{\alpha }_{n}]\right\vert \geq
\delta \mid Z^{n}\right) <\delta \right) \geq 1-\delta
\end{equation*}%
eventually. Also by definition $\frac{d\phi (\alpha _{0})}{d\alpha }[%
\widehat{\alpha }_{n}^{R,B}-\widehat{\alpha }_{n}]=\langle v_{n}^{\ast },%
\widehat{\alpha }_{n}^{R,B}-\widehat{\alpha }_{n}\rangle $. This and
Assumption \ref{ass:phi}(i) imply that%
\begin{equation}
\sqrt{n}\langle u_{n}^{\ast },\widehat{\alpha }_{n}^{R,B}-\widehat{\alpha }%
_{n}\rangle =o_{P_{V^{\infty }|Z^{\infty }}}(1)~wpa1(P_{Z^{\infty }}).
\label{boot-correct-center}
\end{equation}%
Equation (\ref{boot-correct-center}) and $\sqrt{n}\langle u_{n}^{\ast },%
\widehat{\alpha }_{n}-\alpha _{0}\rangle =-\sqrt{n}\mathbb{Z}%
_{n}+o_{P_{Z^{\infty }}}(1)$ (Lemma \ref{lem:theta_anorm(1)}) imply that%
\begin{equation*}
\sqrt{n}\langle u_{n}^{\ast },\widehat{\alpha }_{n}^{R,B}-\alpha _{0}\rangle
=-\sqrt{n}\mathbb{Z}_{n}+o_{P_{V^{\infty }|Z^{\infty
}}}(1)~wpa1(P_{Z^{\infty }}).
\end{equation*}%
Thus we can infer from equation (\ref{eqn:boot_1}) that
\begin{equation}
0.5(\widehat{Q}_{n}^{B}(\widehat{\alpha }_{n}^{R,B}(-\frac{\mathbb{Z}%
_{n}^{\omega -1}}{B_{n}^{\omega }}))-\widehat{Q}_{n}^{B}(\widehat{\alpha }%
_{n}^{R,B}))=-\frac{(\mathbb{Z}_{n}^{\omega -1})^{2}}{2B_{n}^{\omega }}%
+o_{P_{V^{\infty }|Z^{\infty }}}(r_{n}^{-1})~wpa1(P_{Z^{\infty }}).
\label{eqn:boot_03}
\end{equation}

Since $nr_{n}^{-1}=O(1)$, multiplying both sides by $-2n\sigma _{\omega
}^{-2}$, we obtain:%
\begin{equation*}
\mathcal{F}_{n}\geq \left( \sqrt{n}\frac{\mathbb{Z}_{n}^{\omega -1}}{\sigma
_{\omega }\sqrt{B_{n}^{\omega }}}\right) ^{2}-o_{P_{V^{\infty }|Z^{\infty
}}}(1)~wpa1(P_{Z^{\infty }}).
\end{equation*}

\textsc{Step 3.} In order to show
\begin{equation}
\mathcal{F}_{n}\leq \left( \sqrt{n}\frac{\mathbb{Z}_{n}^{\omega -1}}{\sigma
_{\omega }\sqrt{B_{n}^{\omega }}}\right) ^{2}+o_{P_{V^{\infty }|Z^{\infty
}}}(1)~wpa1(P_{Z^{\infty }}),  \label{eqn:boot_02}
\end{equation}%
we can repeat the same calculations as in Step 2, provided there exists a $%
t_{n}^{\ast }\in \mathcal{T}_{n}$ wpa1 such that \textbf{(a)} $\phi (%
\widehat{\alpha }_{n}^{B}(t_{n}^{\ast }))=\phi (\widehat{\alpha }_{n})$ with
$\widehat{\alpha }_{n}^{B}(t_{n}^{\ast })\in \mathcal{A}_{k(n)}$, and
\textbf{(b)} $t_{n}^{\ast }=Z_{n}^{\omega -1}/||u_{n}^{\ast
}||^{2}+o_{P_{V^{\infty }|Z^{\infty }}}(n^{-1/2})=O_{P_{V^{\infty
}|Z^{\infty }}}(n^{-1/2})~wpa1(P_{Z^{\infty }}).$

Because, by (a) and the definition of $\widehat{\alpha }_{n}^{R,B}$,
\begin{equation*}
n\frac{\widehat{Q}_{n}^{B}(\widehat{\alpha }_{n}^{R,B})-\widehat{Q}_{n}^{B}(%
\widehat{\alpha }_{n}^{B})}{\sigma _{\omega }^{2}}\leq n\frac{\widehat{Q}%
_{n}^{B}(\widehat{\alpha }_{n}^{B}(t_{n}^{\ast }))-\widehat{Q}_{n}^{B}(%
\widehat{\alpha }_{n}^{B})}{\sigma _{\omega }^{2}}+o_{P_{V^{\infty
}|Z^{\infty }}}(1)~wpa1(P_{Z^{\infty }}).
\end{equation*}%
By specializing Assumption \ref{ass:LAQ_B}(i) to $\alpha =\widehat{\alpha }%
_{n}^{B}\in \mathcal{N}_{osn}$ wpa1 (by Lemma \ref{lem:cons_boot}(2)), and $%
t_{n}^{\ast }$ as the direction, it follows
\begin{eqnarray*}
&&0.5(\widehat{Q}_{n}^{B}(\widehat{\alpha }_{n}^{B}(t_{n}^{\ast }))-\widehat{%
Q}_{n}^{B}(\widehat{\alpha }_{n}^{B})) \\
&=&t_{n}^{\ast }\{\mathbb{Z}_{n}^{\omega }+\langle u_{n}^{\ast },\widehat{%
\alpha }_{n}^{B}-\alpha _{0}\rangle \}+\frac{B_{n}^{\omega }}{2}(t_{n}^{\ast
})^{2}+o_{P_{V^{\infty }|Z^{\infty }}}(r_{n}^{-1})~wpa1(P_{Z^{\infty }}) \\
&=&\frac{B_{n}^{\omega }}{2}\left( \frac{Z_{n}^{\omega -1}}{||u_{n}^{\ast
}||^{2}}+o_{P_{V^{\infty }|Z^{\infty }}}(n^{-1/2})\right)
^{2}+o_{P_{V^{\infty }|Z^{\infty }}}(r_{n}^{-1})~wpa1(P_{Z^{\infty }}) \\
&=&\frac{1}{2}\left( \frac{\mathbb{Z}_{n}^{\omega -1}}{\sqrt{B_{n}^{\omega }}%
}\right) ^{2}+o_{P_{V^{\infty }|Z^{\infty }}}(r_{n}^{-1})~wpa1(P_{Z^{\infty
}}),
\end{eqnarray*}%
where the second equality is due to Lemma \ref{lem:LAR_boot}(2) and (b), the
third equality is due to the assumption $|B_{n}^{\omega }-||u_{n}^{\ast
}||^{2}|=o_{P_{V^{\infty }|Z^{\infty }}}(1)~wpa1(P_{Z^{\infty }})$ and $%
||u_{n}^{\ast }||\in (c,C)$. Thus equation (\ref{eqn:boot_02}) holds.

\textsc{Step 4.} We now show that there exists a $t_{n}^{\ast }$ such that
\textbf{(a)} and \textbf{(b)} hold in Step 3.

Let $r\equiv \phi (\hat{\alpha}_{n})-\phi (\alpha _{0})$. Since $\widehat{%
\alpha }_{n}^{B}\in \mathcal{N}_{osn}$ wpa1, and $\phi (\hat{\alpha}%
_{n})-\phi (\alpha _{0})=O_{P_{Z^{\infty }}}(||v_{n}^{\ast }||/\sqrt{n})$,
by Lemma \ref{pro:phi_solve}, there is a $t_{n}^{\ast }\in \mathcal{T}_{n}$
wpa1 satisfying (a) with $\widehat{\alpha }_{n}^{B}(t_{n}^{\ast })=\widehat{%
\alpha }_{n}^{B}+t_{n}^{\ast }u_{n}^{\ast }\in \mathcal{A}_{k(n)}$ and $\phi
(\widehat{\alpha }_{n}^{B}(t_{n}^{\ast }))-\phi (\alpha _{0})=r$. Moreover,
by assumption \ref{ass:phi}(i)(ii), such a choice of $t_{n}^{\ast }$ also
satisfies
\begin{equation*}
\left\vert \underbrace{\phi (\widehat{\alpha }_{n}^{B}(t_{n}^{\ast }))-\phi (%
\widehat{\alpha }_{n})}_{=0}-\frac{d\phi (\alpha _{0})}{d\alpha }[\widehat{%
\alpha }_{n}^{B}-\widehat{\alpha }_{n}+t_{n}^{\ast }u_{n}^{\ast
}]\right\vert =o_{P_{V^{\infty }|Z^{\infty }}}(||v_{n}^{\ast }||/\sqrt{n}%
)~wpa1(P_{Z^{\infty }}).
\end{equation*}%
Thus, for sufficiently large $n$,
\begin{equation*}
P_{Z^{\infty }}\left( P_{V^{\infty }|Z^{\infty }}\left( \frac{\sqrt{n}}{%
||v_{n}^{\ast }||}\left\vert \frac{d\phi (\alpha _{0})}{d\alpha }[\widehat{%
\alpha }_{n}^{B}-\widehat{\alpha }_{n}]+t_{n}^{\ast }\frac{||v_{n}^{\ast
}||^{2}}{||v_{n}^{\ast }||_{sd}}\right\vert \geq \delta \mid Z^{n}\right)
<\delta \right) \geq 1-\delta .
\end{equation*}%
By Assumption \ref{ass:phi}(i) and Lemma \ref{lem:LAR_boot}(2), it follows
that the LHS of the above equation is majorized by
\begin{align*}
& P_{Z^{\infty }}\left( P_{V^{\infty }|Z^{\infty }}\left( \frac{\sqrt{n}}{%
||v_{n}^{\ast }||}\left\vert \langle v_{n}^{\ast },\widehat{\alpha }_{n}^{B}-%
\widehat{\alpha }_{n}\rangle +t_{n}^{\ast }\frac{||v_{n}^{\ast }||^{2}}{%
||v_{n}^{\ast }||_{sd}}\right\vert \geq 2\delta \mid Z^{n}\right) <\delta
\right) +\delta \\
& =P_{Z^{\infty }}\left( P_{V^{\infty }|Z^{\infty }}\left( \frac{\sqrt{n}}{%
||v_{n}^{\ast }||}\left\vert -\mathbb{Z}_{n}^{\omega -1}||v_{n}^{\ast
}||_{sd}+t_{n}^{\ast }\frac{||v_{n}^{\ast }||^{2}}{||v_{n}^{\ast }||_{sd}}%
\right\vert \geq 2\delta \mid Z^{n}\right) <\delta \right) +\delta ,
\end{align*}%
Therefore,%
\begin{equation*}
\sqrt{n}t_{n}^{\ast }=\sqrt{n}Z_{n}^{\omega -1}/||u_{n}^{\ast
}||^{2}+o_{P_{V^{\infty }|Z^{\infty }}}(1)~wpa1(P_{Z^{\infty }}).
\end{equation*}%
Since $\sqrt{n}Z_{n}^{\omega -1}=O_{P_{V^{\infty }|Z^{\infty }}}(1)$ with
probability $P_{Z^{\infty }}$ approaching one (assumption \ref{ass:LAQ_B}%
(ii)) and $||u_{n}^{\ast }||^{2}=O(1)$, we have $t_{n}^{\ast
}=O_{P_{V^{\infty }|Z^{\infty }}}(n^{-1/2})$ with probability $P_{Z^{\infty
}}$ approaching one. Thus (b) holds.

\medskip

Before we prove \textbf{Result (2),} we wish to establish the following
\textbf{equation $(\ref{old-result2})$}:%
\begin{equation}
\left\vert \mathcal{L}_{V^{\infty }|Z^{\infty }}\left( \frac{\widehat{QLR}%
_{n}^{B}(\widehat{\phi }_{n})}{\sigma _{\omega }^{2}}\mid Z^{n}\right) -%
\mathcal{L}\left( \widehat{QLR}_{n}(\phi _{0})\mid H_{0}\right) \right\vert
=o_{P_{Z^{\infty }}}(1),  \label{old-result2}
\end{equation}%
where $\mathcal{L}\left( \widehat{QLR}_{n}(\phi _{0})\mid H_{0}\right) $
denotes the law of $\widehat{QLR}_{n}(\phi _{0})$ under the null $H_{0}:\phi
(\alpha )=\phi _{0}$, which will be simply denoted as $\mathcal{L}\left(
\widehat{QLR}_{n}(\phi _{0})\right) $ in the rest of the proof. By Result
(1), it suffices to show that for any $\delta >0$, there exists a $N(\delta
) $ such that
\begin{equation*}
P_{Z^{\infty }}\left( \sup_{f\in BL_{1}}\left\vert E\left[ f\left( \left[
\frac{\sqrt{n}\mathbb{Z}_{n}^{\omega -1}}{\sigma _{\omega }||u_{n}^{\ast }||}%
\right] ^{2}\right) \mid Z^{n}\right] -E[f(\widehat{QLR}_{n}(\phi
_{0}))]\right\vert \leq \delta \right) \geq 1-\delta
\end{equation*}%
for all $n\geq N(\delta )$. Let $\mathbb{Z}$ denote a standard normal random
variable (i.e., $\mathbb{Z}\sim N(0,1)$). If the following equation (\ref%
{eqn:proof-bootstrap_z}) holds, which will be shown at the end of the proof
of equation (\ref{old-result2}),
\begin{equation}
T_{n}\equiv \sup_{f\in BL_{1}}\left\vert E\left[ f\left( \left[ \frac{%
\mathbb{Z}}{||u_{n}^{\ast }||}\right] ^{2}\right) \right] -E[f(\widehat{QLR}%
_{n}(\phi _{0}))]\right\vert =o(1),  \label{eqn:proof-bootstrap_z}
\end{equation}

then, it suffices to show that
\begin{equation}
P_{Z^{\infty }}\left( \sup_{f\in BL_{1}}\left\vert E\left[ f\left( \left[
\frac{\sqrt{n}\mathbb{Z}_{n}^{\omega -1}}{\sigma _{\omega }||u_{n}^{\ast }||}%
\right] ^{2}\right) \mid Z^{n}\right] -E\left[ f\left( \left[ \frac{\mathbb{Z%
}}{||u_{n}^{\ast }||}\right] ^{2}\right) \right] \right\vert \leq \delta
\right) \geq 1-\delta  \label{eqn:proof_bootstrap_0}
\end{equation}%
for all $n\geq N(\delta )$.

Suppose we could show that
\begin{equation}
\sup_{f\in BL_{1}}\left\vert E\left[ f\left( \sqrt{n}\frac{\mathbb{Z}%
_{n}^{\omega -1}}{\sigma _{\omega }||u_{n}^{\ast }||}\right) \mid Z^{n}%
\right] -E\left[ f\left( \mathbb{Z}||u_{n}^{\ast }||^{-1}\right) \right]
\right\vert \rightarrow 0,~wpa1(P_{Z^{\infty }}),
\label{eqn:proof_bootstrap_1}
\end{equation}%
or equivalently,
\begin{equation*}
P_{Z^{\infty }}\left( \left\vert \mathcal{L}_{V^{\infty }|Z^{\infty }}\left(
\sqrt{n}\frac{\mathbb{Z}_{n}^{\omega -1}}{\sigma _{\omega }||u_{n}^{\ast }||}%
|Z^{n}\right) -\mathcal{L}(\mathbb{Z}||u_{n}^{\ast }||^{-1})\right\vert \leq
\delta \right) \geq 1-\delta ,~eventually.
\end{equation*}%
Then, by the continuous mapping theorem (see \cite{Kosorok} Theorem 10.8 and
the discussion in section 10.1.4), we have:%
\begin{equation*}
P_{Z^{\infty }}\left( \left\vert \mathcal{L}_{V^{\infty }|Z^{\infty }}\left(
\left( \sqrt{n}\frac{\mathbb{Z}_{n}^{\omega -1}}{\sigma _{\omega
}||u_{n}^{\ast }||}\right) ^{2}\mid Z^{n}\right) -\mathcal{L}\left( \left(
\mathbb{Z}||u_{n}^{\ast }||^{-1}\right) ^{2}\right) \right\vert \leq \delta
\right) \geq 1-\delta ,~eventually,
\end{equation*}%
and hence equation (\ref{eqn:proof_bootstrap_0}) follows.

It remains to show equation (\ref{eqn:proof_bootstrap_1}). By Assumption \ref%
{ass:LAQ_B}(ii), and the fact that if a sequence converges in probability,
for all subsequence, there exists a subsubsequence that converges almost
surely, it follows for all subsequence $(n_{k})_{k}$, there exists a
subsubsequence $(n_{k(j)})_{j}$ such that
\begin{equation*}
\left\vert \mathcal{L}_{V^{\infty }|Z^{\infty }}\left( \sqrt{n_{k(j)}}\frac{%
\mathbb{Z}_{n_{k(j)}}^{\omega -1}}{\sigma _{\omega }}\mid
Z^{n_{k(j)}}\right) -\mathcal{L}(\mathbb{Z})\right\vert \rightarrow
0,~a.s.-P_{Z^{\infty }}.
\end{equation*}

Since $||u_{n_{k(j)}}^{\ast }||\in (c,C)$, then there exists a further
subsequence (which we still denote as $n_{k(j)}$), such that $%
\lim_{j\rightarrow \infty }||u_{n_{k(j)}}^{\ast }||=d_{\infty }\in \lbrack
c,C]$. Also, since $\sqrt{n}\frac{\mathbb{Z}_{n}^{\omega -1}}{\sigma
_{\omega }}$ is a real valued sequence, by Helly's theorem, convergence in
distribution also holds for $(n_{k(j)})_{j}$. Therefore, by Slutsky theorem,
\begin{equation*}
\mathcal{L}_{V^{\infty }|Z^{\infty }}\left( \sqrt{n_{k(j)}}\frac{\mathbb{Z}%
_{n_{k(j)}}^{\omega -1}}{\sigma _{\omega }||u_{n_{k(j)}}^{\ast }||}\mid
Z^{n_{k(j)}}\right) -\mathcal{L}\left( \mathbb{Z}d_{\infty }^{-1}\right)
\rightarrow 0,~a.s.-P_{Z^{\infty }}.
\end{equation*}

Since $\lim_{j\rightarrow \infty }||u_{n_{k(j)}}^{\ast }||=d_{\infty }\in
\lbrack c,C]$ and $\mathbb{Z}$ is bounded in probability, this readily
implies
\begin{equation*}
\mathcal{L}_{V^{\infty }|Z^{\infty }}\left( \sqrt{n_{k(j)}}\frac{\mathbb{Z}%
_{n_{k(j)}}^{\omega -1}}{\sigma _{\omega }||u_{n_{k(j)}}^{\ast }||}\mid
Z^{n_{k(j)}}\right) -\mathcal{L}\left( \mathbb{Z}||u_{n_{k(j)}}^{\ast
}||^{-1}\right) \rightarrow 0,~a.s.-P_{Z^{\infty }}.
\end{equation*}

Therefore, it follows that
\begin{equation*}
\sup_{f\in BL_{1}}\left\vert E\left[ f\left( \sqrt{n_{k(j)}}\frac{\mathbb{Z}%
_{n_{k(j)}}^{\omega -1}}{\sigma _{\omega }||u_{n_{k(j)}}^{\ast }||}\right)
\mid Z^{n_{k(j)}}\right] -E\left[ f\left( \mathbb{Z}||u_{n_{k(j)}}^{\ast
}||^{-1}\right) \right] \right\vert \rightarrow 0,~a.s.-P_{Z^{\infty }}.
\end{equation*}%
Since the argument started with an arbitrary subsequence $n_{k}$, equation (%
\ref{eqn:proof_bootstrap_1}) holds.

To conclude the proof of equation (\ref{old-result2}), we now show that
equation (\ref{eqn:proof-bootstrap_z}) in fact holds (i.e., $T_{n}=o(1)$).
Again, it suffices to show that for any sub-sequence, there exists a
sub-sub-sequence such that $T_{n(j)}=o(1)$. For any sub-sequence, since $%
(||u_{n}^{\ast }||)_{n}$ is a bounded sequence (under Assumption \ref%
{ass:sieve}(iv)), there exists a further sub-sub-sequence (which we denote
as $(n(j))_{j}$) such that $\lim_{j\rightarrow \infty }||u_{n(j)}^{\ast
}||=d_{\infty }\in \lbrack c,C]$ for finite $c,C>0$. Observe that
\begin{align*}
T_{n(j)}\leq & \sup_{f\in BL_{1}}\left\vert E\left[ f\left( \left[ \frac{%
\mathbb{Z}}{||u_{n(j)}^{\ast }||}\right] ^{2}\right) \right] -E\left[
f\left( \left[ \frac{\mathbb{Z}}{d_{\infty }}\right] ^{2}\right) \right]
\right\vert \\
& +\sup_{f\in BL_{1}}\left\vert E\left[ f\left( \left[ \frac{\mathbb{Z}}{%
d_{\infty }}\right] ^{2}\right) \right] -E\left[ f\left( \left( \frac{%
||u_{n(j)}^{\ast }||}{d_{\infty }}\right) ^{2}\widehat{QLR}_{n(j)}(\phi
_{0})\right) \right] \right\vert \\
& +\sup_{f\in BL_{1}}\left\vert E\left[ f\left( \widehat{QLR}_{n(j)}(\phi
_{0})\right) \right] -E\left[ f\left( \left( \frac{||u_{n(j)}^{\ast }||}{%
d_{\infty }}\right) ^{2}\widehat{QLR}_{n(j)}(\phi _{0})\right) \right]
\right\vert .
\end{align*}%
The first term vanishes because $\mathbb{Z}$ is bounded in probability and $%
\lim_{j\rightarrow \infty }||u_{n(j)}^{\ast }||=d_{\infty }>0$; the third
term follows by the same reason (by Theorem \ref{thm:chi2} and Assumption %
\ref{ass:LAQ}(ii), $\widehat{QLR}_{n}(\phi _{0})$ is bounded in probability).

Finally, for any $f\in BL_{1}$, let $f(d_{\infty }^{-1}\cdot )\equiv f\circ
d_{\infty }^{-2}(\cdot )$. Since $f\circ d_{\infty }^{-2}$ is bounded and $%
|f\circ d_{\infty }^{-2}(t)-f\circ d_{\infty }^{-2}(s)|\leq d_{\infty
}^{-2}|t-s|\leq c^{-2}|t-s|$, we have $\{f\circ d_{\infty }^{-2}:f\in
BL_{1}\}\subseteq BL_{c^{-2}}$. Therefore, the second term in the previous
display is majorized by $\sup_{f\in BL_{c^{-2}}}\left\vert E\left[ f\left( %
\left[ \mathbb{Z}\right] ^{2}\right) \right] -E\left[ f\left(
||u_{n(j)}^{\ast }||^{2}\times \widehat{QLR}_{n(j)}(\phi _{0})\right) \right]
\right\vert $. Hence, to conclude the proof we need to show that
\begin{equation}
\lim_{j\rightarrow \infty }\sup_{f\in BL_{c^{-2}}}\left\vert E\left[ f\left(
\mathbb{Z}^{2}\right) \right] -E\left[ f\left( ||u_{n(j)}^{\ast
}||^{2}\times \widehat{QLR}_{n(j)}(\phi _{0})\right) \right] \right\vert =0.
\label{sub-thm3.2}
\end{equation}%
Theorem \ref{thm:chi2} (i.e., $||u_{n}^{\ast }||^{2}\times \widehat{QLR}%
_{n}(\phi _{0})=[\sqrt{n}\mathbb{Z}_{n}]^{2}+o_{P}(1)$) and Assumption \ref%
{ass:LAQ}(ii) directly imply that the above equation (\ref{sub-thm3.2})
actually holds for the whole sequence, which readily implies that for any
sub-sequence $(n(j))_{j}$ there is a sub-sub-sequence (which we still denote
as $(n(j))_{j}$) for which the previous display holds. \newline

Finally for \textbf{Result (2)}, we want to show that
\begin{equation*}
\sup_{t\in \mathbb{R}}\left\vert P_{V^{\infty }|Z^{\infty }}\left( \frac{%
\widehat{QLR}_{n}^{B}(\widehat{\phi }_{n})}{\sigma _{\omega }^{2}}\leq t\mid
Z^{n}\right) -P_{Z^{\infty }}\left( \widehat{QLR}_{n}(\phi _{0})\leq t\mid
H_{0}\right) \right\vert =o_{P_{Z^{\infty }}}(1).
\end{equation*}

Let $f_{t}(\cdot )\equiv 1\{\cdot \leq t\}$ for $t\in \mathbb{R}$. Under
this notation, the previous display can be cast as
\begin{equation*}
A_{n}\equiv \sup_{t\in \mathbb{R}}\left\vert E_{P_{V^{\infty }|Z^{\infty }}}%
\left[ f_{t}\left( \frac{\widehat{QLR}_{n}^{B}(\widehat{\phi }_{n})}{\sigma
_{\omega }^{2}}\right) \mid Z^{n}\right] -E_{P_{Z^{\infty }}}\left[
f_{t}\left( \widehat{QLR}_{n}(\phi _{0})\right) \right] \right\vert
=o_{P_{Z^{\infty }}}(1).
\end{equation*}%
Denote $\mathbb{Z}^{2}\sim \chi _{1}^{2}$ and
\begin{eqnarray*}
A_{1,n} &\equiv &\sup_{t^{\prime }\in \mathbb{R}}\left\vert E_{P_{V^{\infty
}|Z^{\infty }}}\left[ f_{t^{\prime }}\left( ||u_{n}^{\ast }||^{2}\times
\frac{\widehat{QLR}_{n}^{B}(\widehat{\phi }_{n})}{\sigma _{\omega }^{2}}%
\right) \mid Z^{n}\right] -E\left[ f_{t^{\prime }}\left( \mathbb{Z}%
^{2}\right) \right] \right\vert , \\
A_{2,n} &\equiv &\sup_{t^{\prime }\in \mathbb{R}}\left\vert E_{P_{Z^{\infty
}}}\left[ f_{t^{\prime }}\left( ||u_{n}^{\ast }||^{2}\times \widehat{QLR}%
_{n}(\phi _{0})\right) \right] -E\left[ f_{t^{\prime }}\left( \mathbb{Z}%
^{2}\right) \right] \right\vert .
\end{eqnarray*}%
Notice that%
\begin{eqnarray*}
A_{n} &=&\sup_{t\in \mathbb{R}}\left\vert E_{P_{V^{\infty }|Z^{\infty }}}%
\left[ f_{t||u_{n}^{\ast }||^{2}}\left( ||u_{n}^{\ast }||^{2}\frac{\widehat{%
QLR}_{n}^{B}(\widehat{\phi }_{n})}{\sigma _{\omega }^{2}}\right) \mid Z^{n}%
\right] -E_{P_{Z^{\infty }}}\left[ f_{t||u_{n}^{\ast }||^{2}}\left(
||u_{n}^{\ast }||^{2}\widehat{QLR}_{n}(\phi _{0})\right) \right] \right\vert
\\
&\leq &\sup_{t\in \mathbb{R}}\sup_{d\in \lbrack c,C]}\left\vert
E_{P_{V^{\infty }|Z^{\infty }}}\left[ f_{td^{2}}\left( ||u_{n}^{\ast }||^{2}%
\frac{\widehat{QLR}_{n}^{B}(\widehat{\phi }_{n})}{\sigma _{\omega }^{2}}%
\right) \mid Z^{n}\right] -E_{P_{Z^{\infty }}}\left[ f_{td^{2}}\left(
||u_{n}^{\ast }||^{2}\widehat{QLR}_{n}(\phi _{0})\right) \right] \right\vert
\\
&\leq &\sup_{t^{\prime }\in \mathbb{R}}\left\vert E_{P_{V^{\infty
}|Z^{\infty }}}\left[ f_{t^{\prime }}\left( ||u_{n}^{\ast }||^{2}\times
\frac{\widehat{QLR}_{n}^{B}(\widehat{\phi }_{n})}{\sigma _{\omega }^{2}}%
\right) \mid Z^{n}\right] -E_{P_{Z^{\infty }}}\left[ f_{t^{\prime }}\left(
||u_{n}^{\ast }||^{2}\times \widehat{QLR}_{n}(\phi _{0})\right) \right]
\right\vert \\
&\leq &A_{1,n}+A_{2,n}
\end{eqnarray*}%
where the first line follows from the property that $f_{t}(\cdot
)=f_{t\lambda }(\lambda \times \cdot )$ for any $\lambda \in \mathbb{R}_{+}$
; the second line follows because by assumption, $||u_{n}^{\ast }||^{2}\in
\lbrack c^{2},C^{2}]$; the third line follows simply because $\{1\{\cdot
\leq t\lambda \}:t\in \mathbb{R}~and~\lambda \in \mathbb{R}_{+}\}\subseteq
\{1\{\cdot \leq t\}:t\in \mathbb{R}\}$. Finally, the last line is due to the
triangle inequality and the definitions of $A_{1,n}$ and $A_{2,n}$.

By Theorem \ref{thm:chi2}, under the null, $||u_{n}^{\ast }||^{2}\times
\widehat{QLR}_{n}(\phi _{0})$ converges weakly to $\mathbb{Z}^{2}\sim \chi
_{1}^{2}$, whose distribution is continuous. Therefore, by Polya's theorem, $%
A_{2,n}=o(1)$. Similarly,
\begin{equation*}
A_{1,n}=\sup_{t^{\prime }\in \mathbb{R}}\left\vert P_{V^{\infty }|Z^{\infty
}}\left( ||u_{n}^{\ast }||^{2}\times \frac{\widehat{QLR}_{n}^{B}(\widehat{%
\phi }_{n})}{\sigma _{\omega }^{2}}\leq t^{\prime }\mid Z^{n}\right)
-P\left( \mathbb{Z}^{2}\leq t^{\prime }\right) \right\vert =o_{P_{Z^{\infty
}}}(1)
\end{equation*}%
by equation (\ref{old-result2}) and by the same arguments in Lemma 10.11 in
\cite{Kosorok}. \textit{Q.E.D.}

\medskip

We first recall some notation introduced in the main text. Let $\mathcal{T}%
_{n}\equiv \{t\in \mathbb{R}\colon |t|\leq 4M_{n}^{2}\delta _{n}\}$. For $%
t_{n}\in \mathcal{T}_{n}$, $\alpha (t_{n})\equiv \alpha +t_{n}u_{n}^{\ast }$
where $u_{n}^{\ast }=v_{n}^{\ast }/\left\Vert v_{n}^{\ast }\right\Vert _{sd}$
and $v_{n}^{\ast }=(v_{\theta ,n}^{\ast \prime },v_{h,n}^{\ast }\left( \cdot
\right) )^{\prime }$. To simplify presentation we use $r_{n}=r_{n}(t_{n})%
\equiv \left( \max \{t_{n}^{2},t_{n}n^{-1/2},o(n^{-1})\}\right) ^{-1}$.

\noindent \textbf{Proof of Lemma \ref{lem:Qdiff_B}}: For \textbf{Result (1)}%
, if $\omega \equiv 1$, then Assumption \ref{ass:LAQ_B}(i) simplifies to
\begin{equation*}
P_{Z^{\infty }}\left( P_{V^{\infty }|Z^{\infty }}\left( \sup_{(\alpha
,t_{n})\in \mathcal{N}_{osn}\times \mathcal{T}_{n}}r_{n}\left\vert \widehat{%
\Lambda }_{n}(\alpha (t_{n}),\alpha )-t_{n}\left\{ \mathbb{Z}_{n}+\langle
u_{n}^{\ast },\alpha -\alpha _{0}\rangle \right\} -\frac{B_{n}}{2}%
t_{n}^{2}\right\vert \geq \delta \mid Z^{n}\right) \leq \delta \right) \geq
1-\delta ;
\end{equation*}%
iff
\begin{equation*}
P_{Z^{\infty }}\left( \sup_{(\alpha ,t_{n})\in \mathcal{N}_{osn}\times
\mathcal{T}_{n}}r_{n}\left\vert \widehat{\Lambda }_{n}(\alpha (t_{n}),\alpha
)-t_{n}\left\{ \mathbb{Z}_{n}+\langle u_{n}^{\ast },\alpha -\alpha
_{0}\rangle \right\} -\frac{B_{n}}{2}t_{n}^{2}\right\vert \leq \delta
\right) \geq 1-\delta ,
\end{equation*}%
where $\widehat{\Lambda }_{n}(\alpha (t_{n}),\alpha )\equiv 0.5(\widehat{Q}%
_{n}(\alpha (t_{n}))-\widehat{Q}_{n}(\alpha ))$ and $B_{n}$ is a $Z^{n}$
measurable random variable with $B_{n}=O_{P_{Z^{\infty }}}(1)$. Therefore,
if we could verify Assumption \ref{ass:LAQ_B}(i) in Result (2), we also
verify Assumption \ref{ass:LAQ}(i).

For \textbf{Result (2)}, we divide its proof in several steps.

\textsc{Step 1:} We first introduce some notation. Let
\begin{equation*}
P_{n}(Z^{n})\equiv P_{V^{\infty }|Z^{\infty }}\left( \sup_{(\alpha
,t_{n})\in \mathcal{N}_{osn}\times \mathcal{T}_{n}}r_{n}\left\vert \widehat{%
\Lambda }_{n}^{B}(\alpha (t_{n}),\alpha )-t_{n}\left\{ \mathbb{Z}%
_{n}^{\omega }+\langle u_{n}^{\ast },\alpha -\alpha _{0}\rangle \right\} -%
\frac{B_{n}^{\omega }}{2}t_{n}^{2}\right\vert \geq \delta \mid Z^{n}\right) .
\end{equation*}%
Recall that $\ell _{n}^{B}(x,\alpha )\equiv \widetilde{m}(x,\alpha )+%
\widehat{m}^{B}(x,\alpha _{0})$. Let
\begin{equation*}
\widehat{L}_{n}^{B}(\alpha (t_{n}),\alpha )\equiv \frac{1}{2n}%
\sum_{i=1}^{n}\left\{ \ell _{n}^{B}(X_{i},\alpha (t_{n}))^{\prime }\widehat{%
\Sigma }(X_{i})^{-1}\ell _{n}^{B}(X_{i},\alpha (t_{n}))-\ell
_{n}^{B}(X_{i},\alpha )^{\prime }\widehat{\Sigma }(X_{i})^{-1}\ell
_{n}^{B}(X_{i},\alpha )\right\} .
\end{equation*}

We need to show that $P_{Z^{\infty }}(P_{n}(Z^{n})<\delta )\geq 1-\delta $
eventually which is equivalent to show that $P_{Z^{\infty
}}(P_{n}(Z^{n})>\delta )\leq \delta $ eventually. Hence, it suffices to show
that
\begin{equation*}
P_{Z^{\infty }}(\{P_{n}^{\prime }(Z^{n})>\delta \}\cap S_{n})+P_{Z^{\infty
}}(S_{n}^{C})\leq \delta ,~eventually,
\end{equation*}%
for some event $S_{n}$ that is measurable with respect to $Z^{n}$, and some $%
P_{n}^{\prime }(Z^{n})\geq P_{n}(Z^{n})$ \textit{a.s.}, here $S_{n}^{C}$
denotes the complement of $S_{n}$. In the following we take
\begin{equation*}
S_{n}\equiv \left\{ Z^{n}:P_{V^{\infty }|Z^{\infty }}\left( \sup_{(\alpha
,t_{n})\in \mathcal{N}_{osn}\times \mathcal{T}_{n}}r_{n}\left\vert \widehat{%
\Lambda }_{n}^{B}(\alpha (t_{n}),\alpha )-\widehat{L}_{n}^{B}(\alpha
(t_{n}),\alpha )\right\vert \geq 0.5\delta \mid Z^{n}\right) <0.5\delta
\right\} ,
\end{equation*}%
and
\begin{eqnarray*}
P_{n}^{\prime }(Z^{n}) &\equiv & P_{V^{\infty }|Z^{\infty }}\left(
\sup_{(\alpha ,t_{n})\in \mathcal{N}_{osn}\times \mathcal{T}%
_{n}}r_{n}\left\vert \widehat{L}_{n}^{B}(\alpha (t_{n}),\alpha
)-t_{n}\left\{ \mathbb{Z}_{n}^{\omega }+\langle u_{n}^{\ast },\alpha -\alpha
_{0}\rangle \right\} -\frac{B_{n}^{\omega }}{2}t_{n}^{2}\right\vert \geq
0.5\delta \mid Z^{n}\right) \\
& & +P_{V^{\infty }|Z^{\infty }}\left( \sup_{(\alpha ,t_{n})\in \mathcal{N}%
_{osn}\times \mathcal{T}_{n}}r_{n}\left\vert \widehat{\Lambda }%
_{n}^{B}(\alpha (t_{n}),\alpha )-\widehat{L}_{n}^{B}(\alpha (t_{n}),\alpha
)\right\vert \geq 0.5\delta \mid Z^{n}\right) .
\end{eqnarray*}

It follows that we \textquotedblleft only\textquotedblright\ need to show
that
\begin{equation*}
P_{Z^{\infty }}(S_{n}^{C})\leq 0.5\delta \text{\quad and\quad }P_{Z^{\infty
}}(\{P_{n}^{\prime }(Z^{n})>\delta \}\cap S_{n})\leq 0.5\delta ,~eventually.
\end{equation*}%
Since $P_{Z^{\infty }}(S_{n}^{C})$ can be expressed as%
\begin{equation*}
P_{Z^{\infty }}\left( P_{V^{\infty }|Z^{\infty }}\left( \sup_{(\alpha
,t_{n})\in \mathcal{N}_{osn}\times \mathcal{T}_{n}}r_{n}\left\vert \widehat{%
\Lambda }_{n}^{B}(\alpha (t_{n}),\alpha )-\widehat{L}_{n}^{B}(\alpha
(t_{n}),\alpha )\right\vert \geq 0.5\delta \mid Z^{n}\right) \geq 0.5\delta
\right) ,
\end{equation*}%
which, by Lemma \ref{lem:suff_mcond_boot}(3), is in fact less than $%
0.5\delta $. We only need to verify
\begin{equation*}
P_{Z^{\infty }}(\{P_{n}^{\prime }(Z^{n})>\delta \}\cap S_{n})\leq 0.5\delta
,~eventually.
\end{equation*}

It is easy to see that
\begin{eqnarray*}
& & P_{Z^{\infty }}(\{P_{n}^{\prime }(Z^{n})>\delta \}\cap S_{n}) \\
& \leq & P_{Z^{\infty }}\left( P_{V^{\infty }|Z^{\infty }}\left( \sup_{(\alpha
,t_{n})\in \mathcal{N}_{osn}\times \mathcal{T}_{n}}r_{n}\left\vert \widehat{L%
}_{n}^{B}(\alpha (t_{n}),\alpha )-t_{n}\left\{ \mathbb{Z}_{n}^{\omega
}+\langle u_{n}^{\ast },\alpha -\alpha _{0}\rangle \right\} -\frac{%
B_{n}^{\omega }}{2}t_{n}^{2}\right\vert \geq 0.5\delta \mid Z^{n}\right)
>0.5\delta \right) .
\end{eqnarray*}

Hence, in order to prove the desired result, it suffices to show that
\begin{equation}
P_{Z^{\infty }}\left( P_{V^{\infty }|Z^{\infty }}\left( \sup_{(\alpha
,t_{n})\in \mathcal{N}_{osn}\times \mathcal{T}_{n}}r_{n}\left\vert \widehat{L%
}_{n}^{B}(\alpha (t_{n}),\alpha )-t_{n}\left\{ \mathbb{Z}_{n}^{\omega
}+\langle u_{n}^{\ast },\alpha -\alpha _{0}\rangle \right\} -\frac{%
B_{n}^{\omega }}{2}t_{n}^{2}\right\vert \geq \delta \mid Z^{n}\right)
>\delta \right) <\delta  \label{eqn:main_lem_proof1}
\end{equation}%
eventually.\newline

\textsc{Step 2:} For any $\alpha \in \mathcal{N}_{osn}$ and $t_{n}\in
\mathcal{T}_{n}$, $\alpha (t_{n})=\alpha +t_{n}u_{n}^{\ast }$, under
Assumption \ref{ass:cont_diffm}(i), we can apply the mean value theorem (wrt
$t_{n}$) and obtain
\begin{align*}
\widehat{L}_{n}^{B}(\alpha (t_{n}),\alpha )=& \frac{t_{n}}{n}%
\sum_{i=1}^{n}\left( \frac{d\widetilde{m}(X_{i},\alpha )}{d\alpha }%
[u_{n}^{\ast }]\right) ^{\prime }\widehat{\Sigma }(X_{i})^{-1}\ell
_{n}^{B}(X_{i},\alpha ) \\
& +\frac{t_{n}^{2}}{2n}\int_{0}^{1}\sum_{i=1}^{n}\left( \frac{d\widetilde{m}%
(X_{i},\alpha (s))}{d\alpha }[u_{n}^{\ast }]\right) ^{\prime }\widehat{%
\Sigma }(X_{i})^{-1}\left( \frac{d\widetilde{m}(x,\alpha (s))}{d\alpha }%
[u_{n}^{\ast }]\right) ds \\
& +\frac{t_{n}^{2}}{2n}\int_{0}^{1}\sum_{i=1}^{n}\left( \frac{d^{2}%
\widetilde{m}(X_{i},\alpha (s))}{d\alpha ^{2}}[u_{n}^{\ast },u_{n}^{\ast
}]\right) ^{\prime }\widehat{\Sigma }(X_{i})^{-1}\ell _{n}^{B}(Z_{i},\alpha
(s))ds \\
\equiv & t_{n}T_{1n}^{B}(\alpha )+\frac{t_{n}^{2}}{2}\{T_{2n}(\alpha
)+T_{3n}^{B}(\alpha )\},
\end{align*}%
where $\alpha (s)\equiv \alpha +st_{n}u_{n}^{\ast }\in \mathcal{N}_{osn}$.

From these calculations and the fact that $P_{V^{\infty }|Z^{\infty
}}(a_{n}+b_{n}\geq d|Z^{n})\leq P_{V^{\infty }|Z^{\infty }}(a_{n}\geq
0.5d|Z^{n})+P_{V^{\infty }|Z^{\infty }}(b_{n}\geq 0.5d|Z^{n})$ a.s. for any
two measurable random variables $a_{n}~and~b_{n}$, it follows that
\begin{align*}
& P_{V^{\infty }|Z^{\infty }}\left( \sup_{(\alpha ,t_{n})\in \mathcal{N}%
_{osn}\times \mathcal{T}_{n}}r_{n}\left\vert \widehat{L}_{n}^{B}(\alpha
(t_{n}),\alpha )-t_{n}\left\{ \mathbb{Z}_{n}^{\omega }+\langle u_{n}^{\ast
},\alpha -\alpha _{0}\rangle \right\} -\frac{B_{n}^{\omega }}{2}%
t_{n}^{2}\right\vert \geq 0.5\delta \mid Z^{n}\right) \\
\leq & P_{V^{\infty }|Z^{\infty }}\left( \sup_{(\alpha ,t_{n})\in \mathcal{N}%
_{osn}\times \mathcal{T}_{n}}r_{n}t_{n}\left\vert T_{1n}^{B}(\alpha
)-\left\{ \mathbb{Z}_{n}^{\omega }+\langle u_{n}^{\ast },\alpha -\alpha
_{0}\rangle \right\} \right\vert \geq 0.25\delta \mid Z^{n}\right) \\
& +P_{V^{\infty }|Z^{\infty }}\left( \sup_{(\alpha ,t_{n})\in \mathcal{N}%
_{osn}\times \mathcal{T}_{n}}r_{n}\frac{t_{n}^{2}}{2}\left\vert
\{T_{2n}(\alpha )+T_{3n}^{B}(\alpha )\}-B_{n}^{\omega }\right\vert \geq
0.25\delta \mid Z^{n}\right) .
\end{align*}

Hence, in order to show equation (\ref{eqn:main_lem_proof1}), it suffices to
show that
\begin{equation*}
P_{Z^{\infty }}\left( P_{V^{\infty }|Z^{\infty }}\left( \sup_{(\alpha
,t_{n})\in \mathcal{N}_{osn}\times \mathcal{T}_{n}}r_{n}t_{n}\left\vert
T_{1n}^{B}(\alpha )-\left\{ \mathbb{Z}_{n}^{\omega }+\langle u_{n}^{\ast
},\alpha -\alpha _{0}\rangle \right\} \right\vert \geq \delta \mid
Z^{n}\right) \geq \delta \right) <\delta
\end{equation*}%
and
\begin{equation*}
P_{Z^{\infty }}\left( P_{V^{\infty }|Z^{\infty }}\left( \sup_{(\alpha
,t_{n})\in \mathcal{N}_{osn}\times \mathcal{T}_{n}}\frac{r_{n}t_{n}^{2}}{2}%
\left\vert \{T_{2n}(\alpha )+T_{3n}^{B}(\alpha )\}-B_{n}^{\omega
}\right\vert \geq \delta \mid Z^{n}\right) \geq \delta \right) <\delta
\end{equation*}%
eventually.

Since $r_{n}t_{n}\leq n^{1/2}$, by Lemma \ref{lem:T1n}, the first equation
holds. Since $r_{n}t_{n}^{2}\leq 1$, then in order to verify the second
equation it suffices to verify that, for any $\delta >0$,
\begin{equation*}
P_{Z^{\infty }}\left( \sup_{\alpha \in \mathcal{N}_{osn}}\left\vert
T_{2n}(\alpha )-B_{n}^{\omega }\right\vert \geq \delta \right) <\delta
,~\forall n\geq N(\delta ),
\end{equation*}%
and
\begin{equation*}
P_{Z^{\infty }}\left( P_{V^{\infty }|Z^{\infty }}\left( \sup_{\alpha \in
\mathcal{N}_{osn}}\left\vert T_{3n}^{B}(\alpha )\right\vert \geq \delta \mid
Z^{n}\right) \geq \delta \right) <\delta ,~\forall n\geq N(\delta ).
\end{equation*}%
By Lemmas \ref{lem:T2n}(1) and \ref{lem:T3n}, these two equations hold.

By our choice of $\ell _{n}^{B}()$ (in particular the fact that $\widetilde{m%
}$ is measurable with respect to $Z^{n}$), it follows that $B_{n}^{\omega
}=B_{n}=O_{P_{V^{\infty }|Z^{\infty }}}(1)~wpa1(P_{Z^{\infty }})$. Thus we
verified Assumption \ref{ass:LAQ_B}(i).

Finally, Lemma \ref{lem:T2n}(2) implies $\left\vert B_{n}^{\omega
}-||u_{n}^{\ast }||^{2}\right\vert =o_{P_{V^{\infty }|Z^{\infty
}}}(1)~wpa1(P_{Z^{\infty }})$ and $\left\vert B_{n}-||u_{n}^{\ast
}||^{2}\right\vert =o_{P_{Z^{\infty }}}(1)$. \textit{Q.E.D.}

\medskip

The following lemma is a LLN for triangular arrays.

\begin{lemma}
\label{lem:WLLN_TA} Let $((X_{i,n})_{i=1}^{n})_{n=1}^{\infty }$ be a
triangular array of real valued random variables such that (a) $%
X_{1,n},...,X_{n,n}$ are independent and $X_{i,n}\sim P_{i,n}$, for all $n$,
(b) $E[X_{i,n}]=0$ for all $i$ and $n$, and (c) there is a sequence of
non-negative real numbers $(b_{n})_{n}$ such that $b_{n}=o(\sqrt{n})$ and
\begin{equation*}
\limsup_{n\rightarrow \infty
}n^{-1}\sum_{i=1}^{n}E[|X_{i,n}|1\{|X_{i,n}|\geq b_{n}\}]=0.
\end{equation*}%
Then: for all $\epsilon >0$, there is a $N(\epsilon )$ such that
\begin{equation*}
\Pr \left( \left\vert n^{-1}\sum_{i=1}^{n}X_{i,n}\right\vert \geq \epsilon
\right) <\epsilon \quad \text{for all }n\geq N(\epsilon ).
\end{equation*}
\end{lemma}

\smallskip

\noindent \textbf{Proof of Lemma \ref{lem:WLLN_TA}}: We obtain the result by
modifying the proofs of \cite{Billingsley_book95} theorem 22.1 and of \cite%
{Feller_book} (p. 248). For any $\epsilon >0$, let
\begin{equation*}
X_{i,n}=X_{i,n}1\{|X_{i,n}|\leq b_{n}\}+X_{i,n}1\{|X_{i,n}|>b_{n}\}\equiv
X_{i,n}^{B}+X_{i,n}^{U}.
\end{equation*}%
Thus,
\begin{align*}
\Pr \left( \left\vert n^{-1}\sum_{i=1}^{n}X_{i,n}\right\vert \geq \epsilon
\right) & \leq \Pr \left( \left\vert
n^{-1}\sum_{i=1}^{n}X_{i,n}^{B}\right\vert \geq 0.5\epsilon \right) +\Pr
\left( \left\vert n^{-1}\sum_{i=1}^{n}X_{i,n}^{U}\right\vert \geq
0.5\epsilon \right) \\
& \equiv T_{1,\epsilon }+T_{2,\epsilon }.
\end{align*}

By conditions (b) and (c), it is easy to see that, for large enough $n$,
\begin{align*}
T_{1,\epsilon }& \leq \Pr \left( \left\vert
n^{-1}\sum_{i=1}^{n}\{X_{i,n}^{B}-E[X_{i,n}^{B}]\}\right\vert \geq
0.25\epsilon \right) +1\{n^{-1}\sum_{i=1}^{n}E[X_{i,n}^{B}]\geq 0.25\epsilon \} \\
& =\Pr \left( \left\vert
n^{-1}\sum_{i=1}^{n}\{X_{i,n}^{B}-E[X_{i,n}^{B}]\}\right\vert \geq
0.25\epsilon \right) \leq 2\exp \left( -const.\frac{\epsilon ^{2}n}{b_{n}^{2}%
}\right) ,
\end{align*}%
for some finite constant $const>0$, where the last inequality is due to
Hoeffding inequality (cf. \cite{VdV-W_book96} Appendix A.6). Thus, there is
a $N(\epsilon )$ such that for all $n\geq N(\epsilon )$, $T_{1,\epsilon
}<0.5\epsilon $.

For $T_{2,\epsilon }$, by Markov inequality and then by condition (c), we
have:
\begin{align*}
T_{2,\epsilon }& \leq (\epsilon /2)^{-1}n^{-1}\sum_{i=1}^{n}\int_{\{|x|\geq
b_{n}\}}|x|P_{i,n}(dx) \\
& =(\epsilon /2)^{-1}n^{-1}\sum_{i=1}^{n}\int |x|1\{|x|\geq
b_{n}\}P_{i,n}(dx)<0.5\epsilon
\end{align*}%
eventually. \textit{Q.E.D.}

\medskip

\noindent \textbf{Proof of Lemma \ref{lem:clt_B}}: We divide the proof into
several steps. \newline

\textsc{Step 1.} We first show that the event
\begin{equation*}
S_{n}\equiv \left\{ Z^{n}:\left\vert
n^{-1}\sum_{i=1}^{n}(g(X_{i},u_{n}^{\ast })\rho (Z_{i},\alpha
_{0}))^{2}-E[g(X,u_{n}^{\ast })\Sigma _{0}(X)g(X,u_{n}^{\ast })^{\prime
}]\right\vert \leq \delta \right\}
\end{equation*}%
\textit{occurs wpa1}($P_{Z^{\infty }}$). For this we apply Lemma \ref%
{lem:WLLN_TA}. Using the notation in the lemma, we let $X_{i,n}\equiv
(g(X_{i},u_{n}^{\ast })\rho (Z_{i},\alpha _{0}))^{2}-E[g(X,u_{n}^{\ast
})\Sigma _{0}(X)g(X,u_{n}^{\ast })^{\prime }]$, and thus conditions (a) and
(b) of Lemma \ref{lem:WLLN_TA} immediately follow (note that $%
E[g(X,u_{n}^{\ast })\Sigma _{0}(X)g(X,u_{n}^{\ast })^{\prime }]=1$). In
order to check condition (c), note first that for any generic random
variable $X$ with mean $\mu <\infty $, it follows
\begin{equation*}
E[|X-\mu |1\{|X-\mu |\geq b_{n}\}]\leq E[|X|1\{|X|\geq b_{n}-|\mu |\}]+|\mu
|\Pr \{|X|\geq b_{n}-|\mu |\}.
\end{equation*}%
Since $b_{n}$ is taken to diverge, we can \textquotedblleft
redefine\textquotedblright\ $b_{n}$ as $b_{n}-|\mu |$. Moreover,
\begin{equation*}
\Pr \{|X|\geq b_{n}-|\mu |\}\leq E[\max \{|X|,1\}1\{|X|\geq b_{n}-|\mu |\}].
\end{equation*}%
Again, since $b_{n}$ is taken to diverge the only relevant case is $|X|\geq
1 $. Therefore, it suffices to study $E[|X|1\{|X|\geq b_{n}\}]$ in order to
bound $E[|X-\mu |1\{|X-\mu |\geq b_{n}\}]$. Thus, applied to our case, it is
sufficient to verify that
\begin{equation*}
\limsup_{n\rightarrow \infty }n^{-1}\sum_{i=1}^{n}E\left[ (g(X_{i},u_{n}^{%
\ast })\rho (Z_{i},\alpha _{0}))^{2}1\left\{ (g(X_{i},u_{n}^{\ast })\rho
(Z_{i},\alpha _{0}))^{2}\geq b_{n}\right\} \right] =0,
\end{equation*}%
which holds under our assumption equation (\ref{CLT_triangular}). \smallskip

\textsc{Step 2.} Let $\sqrt{n}\frac{\mathbb{Z}_{n}^{\omega -1}}{\sigma
_{\omega }}=\frac{1}{\sqrt{n}}\sum_{i=1}^{n}\zeta _{i}\mathbf{s}_{i,n}$
where $\mathbf{s}_{i,n}\equiv g(X_{i},u_{n}^{\ast })\rho (Z_{i},\alpha _{0})$%
, and either$\ \left\{ \zeta _{i}\right\} _{i=1}^{n}$ is IID with $\zeta
_{i}=(\omega _{i}-1)\sigma _{\omega }^{-1}$ (under assumption \ref{ass:Wboot}%
) or$\ \left\{ \zeta _{i}\right\} _{i=1}^{n}$ is multinomial with $\zeta
_{i}=(\omega _{i,n}-1)$ (under assumption \ref{ass:Wboot_e}). In the
following we let $P_{\Omega }$ denote the conditional distribution of $%
\left\{ \zeta _{i}\right\} _{i=1}^{n}$ given the data $Z^{n}$, which is also
the unconditional distribution of $\left\{ \zeta _{i}\right\} _{i=1}^{n}$
since $\left\{ \zeta _{i}\right\} _{i=1}^{n}$ is independent of $Z^{n}$. We
want to establish that
\begin{equation*}
\sup_{f\in BL_{1}}\left\vert E\left[ f\left( \sqrt{n}\frac{\mathbb{Z}%
_{n}^{\omega -1}}{\sigma _{\omega }}\right) \mid Z^{n}\right] -E\left[
f\left( \mathbb{Z}\right) \right] \right\vert =o_{P_{Z^{\infty }}}(1),
\end{equation*}%
where $\mathbb{Z}\sim N(0,1)$. Which is equivalent to show that
\begin{equation*}
\frac{1}{\sqrt{n}}\sum_{i=1}^{n}\zeta _{i}\mathbf{s}_{i,n}\Rightarrow
\mathbb{Z},~wpa1(P_{Z^{\infty }}).
\end{equation*}%
Which, by \cite{Billingsley_book95} (Theorem 20.5, p. 268), in turn suffices
to show that any sub-sequence, contains a further sub-sequence, $(n_{k})_{k}$%
, such that
\begin{equation}
\frac{1}{\sqrt{n_{k}}}\sum_{i=1}^{n_{k}}\zeta _{i}\mathbf{s}%
_{i,n_{k}}\Rightarrow \mathbb{Z},~a.s.-(P_{Z^{\infty }}).
\label{subsq-clt-B}
\end{equation}%
Step 3 below establishes (\ref{subsq-clt-B}) under assumption \ref{ass:Wboot}%
, while Step 4 below establishes (\ref{subsq-clt-B}) under assumption \ref%
{ass:Wboot}.\smallskip

\textsc{Step 3.} (under assumption \ref{ass:Wboot}) Since the event $S_{n}$
occurs $wpa1(P_{Z^{\infty }})$ (Step 1), it follows that any sub-sequence,
contains a further sub-sequence such that $n_{k}^{-1}\sum_{i=1}^{n_{k}}(%
\mathbf{s}_{i,n_{k}})^{2}\rightarrow 1$ , $a.s.-(P_{Z^{\infty }})$.
Moreover, $\max_{i\leq n_{k}}|\mathbf{s}_{i,n_{k}}|/\sqrt{n_{k}}=o(1)$, $%
a.s.-(P_{Z^{\infty }})$. This follows since, for any $\epsilon >0$,
\begin{align*}
P_{Z^{\infty }}\left( \max_{i\leq n}|\mathbf{s}_{i,n}|\geq \epsilon \sqrt{n}%
\right) & \leq \sum_{i=1}^{n}\int_{|s|\geq \epsilon \sqrt{n}}P_{i,n}(ds)\leq
\epsilon ^{-2}n^{-1}\sum_{i=1}^{n}\int_{|s|\geq \epsilon \sqrt{n}%
}s^{2}P_{i,n}(ds) \\
& =\epsilon ^{-2}n^{-1}\sum_{i=1}^{n}E[\mathbf{s}_{i,n}^{2}1\{|\mathbf{s}%
_{i,n}|\geq \epsilon \sqrt{n}\}].
\end{align*}%
We note that $1\{|\mathbf{s}_{i,n}|\geq \epsilon \sqrt{n}\}\leq 1\{|\mathbf{s%
}_{i,n}|^{2}\geq b_{n}\}$ (provided that $|\mathbf{s}_{i,n}|\geq 1$, but if
it is not, then the proof is trivial). Hence by equation (\ref%
{CLT_triangular}) and the fact that $\mathbf{s}_{i,n}$ are row-wise iid, the
RHS is of order $o(1)$. Going to a sub-sequence establishes the result.

Under assumption \ref{ass:Wboot}, $\zeta _{i}=(\omega _{i}-1)\sigma _{\omega
}^{-1}$ is IID with mean zero, variance 1, 
 hence conditional on the event $S_{n}$,
for any $\epsilon >0,$
\begin{eqnarray*}
&&n_{k}^{-1}\sum_{i=1}^{n_{k}}E_{P_{\Omega }}\left[ (\zeta _{i}\mathbf{s}%
_{i,n_{k}})^{2}1\left\{ |\zeta _{i}\mathbf{s}_{i,n_{k}}|>\epsilon \sqrt{n_{k}%
}\right\} \right] \\
&\leq &\left( n_{k}^{-1}\sum_{i=1}^{n_{k}}|\mathbf{s}_{i,n_{k}}|^{2}\right)
\times E_{P_{\Omega }}\left[ \zeta _{1}^{2}\times 1\left\{ |\zeta
_{1}|\times \max_{1\leq i\leq n}|\mathbf{s}_{i,n_{k}}|>\epsilon \sqrt{n_{k}}%
\right\} \right] \\
&\leq &\left( n_{k}^{-1}\sum_{i=1}^{n_{k}}|\mathbf{s}_{i,n_{k}}|^{2}\right)
\times E_{P_{\Omega }}\left[ \zeta _{1}^{2}\times 1\left\{ |\zeta
_{1}|>\epsilon /\epsilon ^{\prime }\right\} \right] \rightarrow
0,~a.s.-(P_{Z^{\infty }}).
\end{eqnarray*}%
where the second inequality follows from the fact that $\max_{i\leq n_{k}}|%
\mathbf{s}_{i,n_{k}}|/\sqrt{n_{k}}<\epsilon ^{\prime }$, $a.s.-(P_{Z^{\infty
}})$ eventually. Since 
$\zeta_{1}$ are IID, by choosing the $\epsilon ^{\prime }$ (small relative to $%
\epsilon $), one can make the term $E_{P_{\Omega }}[\zeta _{1}^{2}1\left\{
|\zeta _{1}|>\epsilon /\epsilon ^{\prime }\right\} ]$ arbitrarily small. The
Lindeberg-Feller CLT then implies that $\frac{1}{\sqrt{n_{k}}}%
\sum_{i=1}^{n_{k}}\zeta _{i}\mathbf{s}_{i,n_{k}}\Rightarrow \mathbb{Z}$,~$%
a.s.-(P_{Z^{\infty }})$ where $\mathbb{Z}\sim N(0,1)$.

We have thus showed that any sub-sequence, contains a further sub-sequence
such that the above equation holds; therefore
\begin{align*}
\sup_{f\in BL_{1}}\left\vert E\left[ f\left( \frac{1}{\sqrt{n}}%
\sum_{i=1}^{n}\zeta _{i}\mathbf{s}_{i,n}\right) \mid Z^{n}\right] -E\left[
f\left( \mathbb{Z}\right) \right] \right\vert & =o_{P_{Z^{\infty }}}(1).
\end{align*}

\smallskip

\textsc{Step 4.} (under assumption \ref{ass:Wboot_e}) We proceed as in Step
3 to establish equation (\ref{subsq-clt-B}). The difference is that now $%
\left\{ \zeta _{i}\right\} _{i=1}^{n}$ is not iid, but exchangeable with $%
\zeta _{i}=\zeta _{i,n}\equiv (\omega _{i,n}-1)$. To overcome this, we
follow lemma 3.6.15 (or really proposition A.5.3) in VdV-W for a given
sub-sequence $(n_{k})_{k}$. To simplify notation we let $n=n_{k}$ and $%
\mathbf{s}_{i,n}=\mathbf{s}_{i,n_{k}}$.

Under assumption \ref{ass:Wboot_e} we have: $n^{-1}\sum_{i=1}^{n}\zeta
_{i,n}=0$, $n^{-1}\sum_{i=1}^{n}\zeta _{i,n}^{2}\rightarrow 1$, $%
n^{-1}\max_{1\leq i\leq n}\zeta _{i,n}^{2}=o_{P_{\Omega }}(1)$ and $%
\max_{1\leq i\leq n}E[\zeta _{i,n}^{4}]\leq c<\infty $. Conditional on the
event $S_{n}$, we also have $n^{-1}\sum_{i=1}^{n}\mathbf{s}_{i,n}\rightarrow
0$, $n^{-1}\sum_{i=1}^{n}\mathbf{s}_{i,n}^{2}\rightarrow 1$ and $%
n^{-1}\max_{1\leq i\leq n}\mathbf{s}_{i,n}^{2}=o(1)$ (this has already been
established in Step 3), and finally we need:

\begin{equation}
\limsup_{n\rightarrow \infty }n^{-2}\sum_{i=1}^{n}\sum_{j=1}^{n}(\mathbf{s}%
_{i,n}\zeta _{j,n})^{2}1\left\{ |\mathbf{s}_{i,n}\zeta _{j,n}|>\epsilon
\sqrt{n}\right\} =0,~a.s.-P_{Z^{\infty }}.  \label{eqn:clt_B-1}
\end{equation}%
To show equation (\ref{eqn:clt_B-1}), we note that
\begin{eqnarray*}
&&\limsup_{n\rightarrow \infty
}n^{-2}\sum_{i=1}^{n}\sum_{j=1}^{n}(\textbf{s}_{i,n}\zeta_{j,n})^{2}1\left\{
|\textbf{s}_{i,n}\zeta_{j,n}|>\epsilon \sqrt{n}\right\} \\
&\leq &\limsup_{n\rightarrow \infty }\left( \frac{\sum_{j=1}^{n}\zeta_{j,n}^{2}}{n%
}\times \frac{\sum_{i=1}^{n}(\textbf{s}_{i,n})^{2}1\left\{ |\textbf{s}_{i,n}|\times \max_{1\leq
j\leq n}|\zeta_{j,n}|>\epsilon \sqrt{n}\right\} }{n}\right) .
\end{eqnarray*}%
Since $n^{-1}\sum_{i=1}^{n}\zeta_{i,n}^{2}-1=o_{P_{\Omega }}(1)$, we have (with
possibly going to a subsequence) $\limsup_{n\rightarrow \infty
}n^{-1}\sum_{i=1}^{n}\zeta_{i,n}^{2}=1$ a.s-$P_{\Omega }~a.s.-P_{Z^{\infty }}$.
Hence, in order to show equation (\ref{eqn:clt_B-1}), it suffices to show
that
\begin{equation*}
\limsup_{n\rightarrow \infty }\frac{1}{n}\sum_{i=1}^{n}(\textbf{s}_{i,n})^{2}1\left\{
|\textbf{s}_{i,n}|\times \max_{1\leq j\leq n}|\zeta_{j,n}|>\epsilon \sqrt{n}\right\}
=0~a.s.-P_{Z^{\infty }}.
\end{equation*}
Let $L_{n}=\sqrt{\frac{n}{b_{n}}}\rightarrow \infty $ (such a choice exists
since $b_{n}=o(\sqrt{n})$ in equation (\ref{CLT_triangular})). It then follows from the properties of $\zeta_{j,n}$ that
\begin{equation*}
P_{\Omega }\left( \max_{1\leq j\leq n}|\zeta_{j,n}|>L_{n}\right) =o(1)
\end{equation*}%
Then, by possibly going to a subsequence, it
implies that $\max_{1\leq j\leq n}|\zeta_{j,n}|\leq L_{n}$ a.s.-$P_{Z^{\infty }}$
and $1\left\{ |\textbf{s}_{i,n}|\times \max_{1\leq j\leq n}|\zeta_{j,n}|>\epsilon \sqrt{n}%
\right\} \leq 1\{|\textbf{s}_{i,n}|L_{n}>\epsilon \sqrt{n}\}=1\{|\textbf{s}_{i,n}|^{2}>\epsilon ^{2}b_{n}\}$ a.s.-$P_{Z^{\infty
}}$. Thus
\begin{equation*}
\frac{1}{n}\sum_{i=1}^{n}(\textbf{s}_{i,n})^{2}1\left\{ |\textbf{s}_{i,n}|\times \max_{1\leq
j\leq n}|\zeta_{j,n}|>\epsilon \sqrt{n}\right\} \leq \frac{1}{n}%
\sum_{i=1}^{n}(\textbf{s}_{i,n})^{2}1\{|\textbf{s}_{i,n}|^{2}>\epsilon ^{2}b_{n}\}\rightarrow
0,~a.s.-(P_{Z^{\infty }}).
\end{equation*}%
Hence (by possibly going to subsequences) equation (\ref{eqn:clt_B-1}) follows from equation (\ref{CLT_triangular}%
). So, by lemma 3.6.15 (or proposition A.5.3) in \cite{VdV-W_book96},
\begin{equation*}
\frac{1}{\sqrt{n_{k}}}\sum_{i=1}^{n_{k}}\zeta _{i}\mathbf{s}%
_{i,n_{k}}\Rightarrow \mathbb{Z},~a.s.-(P_{Z^{\infty }}).
\end{equation*}%
The rest of the steps are analogous to those in Step 3 and will not be
repeated here. \textit{Q.E.D.}

\subsubsection{Alternative bootstrap sieve t statistics}

In this subsection we present additional bootstrap sieve t statistics.
Recall that $\widehat{W}_{n}\equiv \sqrt{n}\frac{\phi (\widehat{\alpha }%
_{n})-\phi (\alpha _{0})}{||\widehat{v}_{n}^{\ast }||_{n,sd}}$ is the
original sample sieve t statistic. The first one is $\widehat{W}%
_{1,n}^{B}\equiv \sqrt{n}\frac{\phi (\widehat{\alpha }_{n}^{B})-\phi (%
\widehat{\alpha }_{n})}{\sigma _{\omega }||\widehat{v}_{n}^{\ast }||_{n,sd}}$%
. In the definition of $\widehat{W}_{2,n}^{B}$ one could also define $||%
\widehat{v}_{n}^{\ast }||_{B,sd}^{2}$ using $\hat{\Sigma}_{0i}^{B}=\widehat{E%
}_{n}[\varrho (V,\widehat{\alpha }_{n})\varrho (V,\widehat{\alpha }%
_{n})^{\prime }|X=X_{i}]$ instead of $\varrho (V_{i},\widehat{\alpha }%
_{n})\varrho (V_{i},\widehat{\alpha }_{n})^{\prime }$, which will be a
bootstrap analog to $||\widehat{v}_{n}^{\ast }||_{n,sd}^{2}$ defined in
equation (\ref{svar-hat}).

Let $\widehat{W}_{3,n}^{B}\equiv \sqrt{n}\frac{\phi (\widehat{\alpha }%
_{n}^{B})-\phi (\widehat{\alpha }_{n})}{||\widehat{v}_{n}^{B}||_{B,sd}}$
where $||\widehat{v}_{n}^{B}||_{B,sd}^{2}$ is a bootstrap sieve variance
estimator that is constructed as follows. First, we define
\begin{equation*}
||\cdot ||_{B,M}^{2}\equiv n^{-1}\sum_{i=1}^{n}\left( \frac{d\widehat{m}%
^{B}(X_{i},\widehat{\alpha }_{n}^{B})}{d\alpha }[\cdot ]\right) ^{\prime
}M_{n,i}\left( \frac{d\widehat{m}^{B}(X_{i},\widehat{\alpha }_{n}^{B})}{%
d\alpha }[\cdot ]\right) ,
\end{equation*}%
where $M_{n,i}$ is some (almost surely) positive definite weighting matrix.
Let $\widehat{v}_{n}^{B}$ be a \emph{bootstrapped empirical Riesz representer%
} of the linear functional $\frac{d\phi (\widehat{\alpha }_{n}^{B})}{d\alpha
}[\cdot ]$ under $||\cdot ||_{B,\hat{\Sigma}^{-1}}$. We compute a bootstrap
sieve variance estimator as:
\begin{equation}
||\widehat{v}_{n}^{B}||_{B,sd}^{2}\equiv \frac{1}{n}\sum_{i=1}^{n}\left(
\frac{d\widehat{m}^{B}(X_{i},\widehat{\alpha }_{n}^{B})}{d\alpha }[\widehat{v%
}_{n}^{B}]\right) ^{\prime }\hat{\Sigma}_{i}^{-1}\varrho (V_{i},\widehat{%
\alpha }_{n}^{B})\varrho (V_{i},\widehat{\alpha }_{n}^{B})^{\prime }\hat{%
\Sigma}_{i}^{-1}\left( \frac{d\widehat{m}^{B}(X_{i},\widehat{\alpha }%
_{n}^{B})}{d\alpha }[\widehat{v}_{n}^{B}]\right)  \label{svar-hat-boot}
\end{equation}%
with $\varrho (V_{i},\alpha )\equiv (\omega _{i,n}-1)\rho (Z_{i},\alpha
)\equiv \rho ^{B}(V_{i},\alpha )-\rho (Z_{i},\alpha )$ for any $\alpha $.
That is, $||\widehat{v}_{n}^{B}||_{B,sd}^{2}$ is a bootstrap analog to $||%
\widehat{v}_{n}^{\ast }||_{n,sd}^{2}$ defined in equation (\ref{svar-hat1}).
One could also define $||\widehat{v}_{n}^{B}||_{B,sd}^{2}$ using $\widehat{E}%
_{n}[\varrho (V,\widehat{\alpha }_{n}^{B})\varrho (V,\widehat{\alpha }%
_{n}^{B})^{\prime }|X=X_{i}]$ instead of $\varrho (V_{i},\widehat{\alpha }%
_{n}^{B})\varrho (V_{i},\widehat{\alpha }_{n}^{B})^{\prime }$, which will be
a bootstrap analog to $||\widehat{v}_{n}^{\ast }||_{n,sd}^{2}$ defined in
equation (\ref{svar-hat}). In addition, one could also define $||\widehat{v}%
_{n}^{B}||_{B,sd}^{2}$ using $\widehat{\alpha }_{n}$ instead of $\widehat{%
\alpha }_{n}^{B}$. In terms of the first order asymptotic approximation,
this alternative definition yields the same asymptotic results. Due to space
considerations, we omit these alternative bootstrap sieve variance
estimators.

The bootstrap sieve variance estimator $||\widehat{v}_{n}^{B}||_{B,sd}^{2}$
also has a closed form expression: $||\widehat{v}_{n}^{B}||_{B,sd}^{2}=(%
\widehat{\digamma }_{n}^{B})^{\prime }(\widehat{D}_{n}^{B})^{-1}\widehat{%
\mho }_{3,n}^{B}(\widehat{D}_{n}^{B})^{-1}\widehat{\digamma }_{n}^{B}$ with%
\begin{eqnarray*}
\widehat{\digamma }_{n}^{B} &=&\frac{d\phi (\widehat{\alpha }_{n}^{B})}{%
d\alpha }[\overline{\psi }^{k(n)}(\cdot )^{\prime }],\text{ }\widehat{D}%
_{n}^{B}=\frac{1}{n}\sum_{i=1}^{n}\left( \frac{d\widehat{m}^{B}(X_{i},%
\widehat{\alpha }_{n}^{B})}{d\alpha }[\overline{\psi }^{k(n)}(\cdot
)^{\prime }]\right) ^{\prime }\widehat{\Sigma }_{i}^{-1}\left( \frac{d%
\widehat{m}^{B}(X_{i},\widehat{\alpha }_{n}^{B})}{d\alpha }[\overline{\psi }%
^{k(n)}(\cdot )^{\prime }]\right) , \\
\widehat{\mho }_{3,n}^{B} &=&\frac{1}{n}\sum_{i=1}^{n}\left( \frac{d\widehat{%
m}^{B}(X_{i},\widehat{\alpha }_{n}^{B})}{d\alpha }[\overline{\psi }%
^{k(n)}(\cdot )^{\prime }]\right) ^{\prime }\widehat{\Sigma }%
_{i}^{-1}(\omega _{i,n}-1)^{2}\rho (Z_{i},\widehat{\alpha }_{n}^{B})\rho
(Z_{i},\widehat{\alpha }_{n}^{B})^{\prime }\widehat{\Sigma }_{i}^{-1}\left(
\frac{d\widehat{m}^{B}(X_{i},\widehat{\alpha }_{n}^{B})}{d\alpha }[\overline{%
\psi }^{k(n)}(\cdot )^{\prime }]\right) .
\end{eqnarray*}%
This expression is computed in the same way as $||\widehat{v}_{n}^{\ast
}||_{n,sd}^{2}=\widehat{\digamma }_{n}^{\prime }\widehat{D}_{n}^{-1}\widehat{%
\mho }_{n}\widehat{D}_{n}^{-1}\widehat{\digamma }_{n}$ given in (\ref%
{P-svar-hat1}) but using bootstrap analogs. Note that this bootstrap sieve
variance only uses $\widehat{\alpha }_{n}^{B}$, and is easy to compute.

When specialized to the NPIV model (\ref{npiv}) in subsection \ref%
{sec:NPIVex1}, the expression $||\widehat{v}_{n}^{B}||_{B,sd}^{2}$
simplifies further, with $\widehat{\digamma }_{n}^{B}=\frac{d\phi (\widehat{h%
}_{n}^{B})}{d\alpha }[q^{k(n)}(\cdot )^{\prime }]$, $\widehat{D}_{n}^{B}=%
\frac{1}{n}\widehat{C}_{n}^{B}(P^{\prime }P)^{-}(\widehat{C}%
_{n}^{B})^{\prime }$, $\widehat{C}_{n}^{B}=\sum_{j=1}^{n}\omega
_{j,n}q^{k(n)}(Y_{2j})p^{J_{n}}(X_{j})^{\prime }$,%
\begin{equation*}
\widehat{\mho }_{3,n}^{B}=\frac{1}{n}\widehat{C}_{n}^{B}(P^{\prime
}P)^{-}\left( \sum_{i=1}^{n}p^{J_{n}}(X_{i})[(\omega _{i,n}-1)\widehat{U}%
_{i}^{B}]^{2}p^{J_{n}}(X_{i})^{\prime }\right) (P^{\prime }P)^{-}(\widehat{C}%
_{n}^{B})^{\prime },\quad \text{with }\widehat{U}_{i}^{B}=Y_{1i}-\widehat{h}%
_{n}^{B}(Y_{2i}).
\end{equation*}%
This expression is analogous to that for a 2SLS t-bootstrap test; see \cite%
{DM}. We leave it to further work to study whether this bootstrap sieve t
statistic might have second order refinement by choice of some IID bootstrap
weights.

Recall that $\hat{M}_{i}^{B}=(\omega _{i,n}-1)^{2}\hat{M}_{i}$ and $\hat{M}%
_{i}=\hat{\Sigma}_{i}^{-1}\rho (Z_{i},\widehat{\alpha }_{n})\rho (Z_{i},%
\widehat{\alpha }_{n})^{\prime }\hat{\Sigma}_{i}^{-1}$.

\begin{assumption}
\label{ass:VE_boot} (i) $\sup_{v_{1},v_{2}\in \overline{\mathbf{V}}%
_{k(n)}^{1}}|\langle v_{1},v_{2}\rangle _{B,\Sigma ^{-1}}-\langle
v_{1},v_{2}\rangle _{n,\Sigma ^{-1}}|=o_{P_{V^{\infty }|Z^{\infty
}}}(1)~wpa1(P_{Z^{\infty }})$;

\noindent(ii) $\sup_{v\in \overline{\mathbf{V}}_{k(n)}^{1}}|\langle
v,v\rangle _{B,\hat{M}^{B}}-\sigma _{\omega }^{2}\langle v,v\rangle _{n,\hat{%
M}}|=o_{P_{V^{\infty }|Z^{\infty }}}(1)~wpa1(P_{Z^{\infty }})$;

\noindent(iii) $\sup_{v\in \overline{\mathbf{V}}_{k(n)}^{1}}n^{-1}%
\sum_{i=1}^{n}(\omega _{i,n}-1)^{2}\left\Vert \frac{d\widehat{m}^{B}(X_{i},%
\hat{\alpha}_{n}^{B})}{d\alpha }[v]\right\Vert _{e}^{2}=O_{P_{V^{\infty
}|Z^{\infty }}}(1)~wpa1(P_{Z^{\infty }})$.
\end{assumption}

Assumption \ref{ass:VE_boot}(i)(ii) is analogous to Assumption \ref{ass:VE}%
(ii)(v). Assumption \ref{ass:VE_boot}(iii) is a mild one, for example, it is
implied by Assumptions for Lemma \ref{lem:cons_boot} and uniformly bounded
bootstrap weights (i.e., $|\omega _{i,n}|\leq C<\infty $ for all $i$).

The following result is a bootstrap version of Theorem \ref{thm:VE}.

\begin{theorem}
\label{thm:VE-boot} Let Conditions for Theorem \ref{thm:VE}(1) and Lemma \ref%
{lem:cons_boot}, Assumption \ref{ass:VE_boot} hold. Then:%
\begin{equation*}
\text{(1)}\quad \left\vert \frac{||\widehat{v}_{n}^{B}||_{B,sd}}{\sigma
_{\omega }||v_{n}^{\ast }||_{sd}}-1\right\vert =o_{P_{V^{\infty }|Z^{\infty
}}}(1)~wpa1(P_{Z^{\infty }}).
\end{equation*}%
(2) If further, conditions for Theorem \ref{thm:bootstrap_2}(1) hold, then:
\begin{equation*}
\widehat{W}_{3,n}^{B}=-\sqrt{n}\frac{\mathbb{Z}_{n}^{\omega -1}}{\sigma
_{\omega }}+o_{P_{V^{\infty }|Z^{\infty }}}(1)~wpa1(P_{Z^{\infty }}),
\end{equation*}%
\begin{equation*}
\left\vert \mathcal{L}_{V^{\infty }|Z^{\infty }}\left( \widehat{W}%
_{3,n}^{B}\mid Z^{n}\right) -\mathcal{L}\left( \widehat{W}_{n}\right)
\right\vert =o_{P_{Z^{\infty }}}(1),\quad \text{and}
\end{equation*}%
\begin{equation*}
\sup_{t\in \mathbb{R}}\left\vert P_{V^{\infty }|Z^{\infty }}(\widehat{W}%
_{3,n}^{B}\leq t|Z^{n})-P_{Z^{\infty }}(\widehat{W}_{n}\leq t)\right\vert
=o_{P_{V^{\infty }|Z^{\infty }}}(1)~wpa1(P_{Z^{\infty }}).
\end{equation*}
\end{theorem}

\smallskip

\noindent \textbf{Proof of Theorem \ref{thm:VE-boot}.} For \textbf{Result (1)%
}, the proof is analogous to the one for Theorem \ref{thm:VE}(1). As in the
proof of Theorem \ref{thm:VE}(1), it suffices to show that
\begin{equation}
\frac{||\hat{v}_{n}^{B}-v_{n}^{\ast }||}{||v_{n}^{\ast }||}=o_{P_{V^{\infty
}|Z^{\infty }}}(1)~wpa1(P_{Z^{\infty }}),  \label{eqn:VE_B_1}
\end{equation}%
and
\begin{equation}
\left\vert \frac{||\hat{v}_{n}^{B}||_{B,sd}-||\hat{v}_{n}^{B}||_{sd}}{%
||v_{n}^{\ast }||}\right\vert =o_{P_{V^{\infty }|Z^{\infty
}}}(1)~wpa1(P_{Z^{\infty }}).  \label{eqn:VE_B_2}
\end{equation}

Following the same derivations as in the proof of theorem \ref{thm:VE}(1)
step 1, for equation (\ref{eqn:VE_B_1}), it suffices to show
\begin{equation*}
|\langle \hat{\varpi}_{n}^{B},\varpi \rangle _{B,\hat{\Sigma}^{-1}}-\langle
\hat{\varpi}_{n}^{B},\varpi \rangle _{B,\Sigma ^{-1}}|=o_{P_{V^{\infty
}|Z^{\infty }}}(1)~and~|\langle \hat{\varpi}_{n}^{B},\varpi \rangle
_{B,\Sigma ^{-1}}-\langle \hat{\varpi}_{n}^{B},\varpi \rangle _{\Sigma
^{-1}}|=o_{P_{V^{\infty }|Z^{\infty }}}(1)
\end{equation*}%
$wpa1(P_{Z^{\infty }})$, uniformly over $\varpi \in \overline{\mathbf{V}}%
_{k(n)}^{1}$; where $\hat{\varpi}_{n}^{B}=\frac{\hat{v}_{n}^{B}}{||\hat{v}%
_{n}^{B}||}$. The first term follows by Assumptions \ref{ass:VE}(iii) and %
\ref{ass:sieve}(iv) and the fact that $\langle \varpi ,\varpi \rangle
_{B,\Sigma ^{-1}}=O_{P_{V^{\infty }|Z^{\infty }}}(1)$ $wpa1(P_{Z^{\infty }})$
(by Assumptions \ref{ass:VE_boot}(i) and \ref{ass:VE}(ii)). The second term
follows directly from these two assumptions.

Regarding equation (\ref{eqn:VE_B_2}), following the same derivations as in
the proof of Theorem \ref{thm:VE} step 2, it suffices to show that $%
\left\vert ||\hat{\varpi}_{n}^{B}||_{B,sd}^{2}-||\hat{\varpi}%
_{n}^{B}||_{sd}^{2}\right\vert =o_{P_{V^{\infty }|Z^{\infty }}}(1)$ $%
wpa1(P_{Z^{\infty }}).$ By the triangle inequality,
\begin{eqnarray*}
\sup_{v\in \overline{\mathbf{V}}_{k(n)}^{1}}|\langle v,v\rangle _{B,\hat{W}%
^{B}}-\sigma _{\omega }^{2}\langle v,v\rangle _{n,\hat{M}}|& \leq & \sup_{v\in
\overline{\mathbf{V}}_{k(n)}^{1}}\left\vert \langle v,v\rangle _{B,\hat{W}%
^{B}}-\langle v,v\rangle _{B,\hat{M}^{B}}\right\vert \\
&& +\sup_{v\in \overline{%
\mathbf{V}}_{k(n)}^{1}}\left\vert \langle v,v\rangle _{B,\hat{M}^{B}}-\sigma
_{\omega }^{2}\langle v,v\rangle _{n,\hat{M}}\right\vert \\
& \equiv & A_{1n}^{B}+A_{2n}^{B}
\end{eqnarray*}%
with $\hat{W}_{i}^{B}\equiv \hat{\Sigma}_{i}^{-1}\varrho (V_{i},\hat{\alpha}%
_{n}^{B})\varrho (V_{i},\hat{\alpha}_{n}^{B})^{\prime }\hat{\Sigma}%
_{i}^{-1}=(\omega _{i,n}-1)^{2}\hat{\Sigma}_{i}^{-1}\rho (Z_{i},\hat{\alpha}%
_{n}^{B})\rho (Z_{i},\hat{\alpha}_{n}^{B})^{\prime }\hat{\Sigma}_{i}^{-1}$
and $\hat{M}_{i}^{B}=(\omega _{i,n}-1)^{2}\hat{M}_{i}$ and $\hat{M}_{i}=\hat{%
\Sigma}_{i}^{-1}\rho (Z_{i},\widehat{\alpha }_{n})\rho (Z_{i},\widehat{%
\alpha }_{n})^{\prime }\hat{\Sigma}_{i}^{-1}$.

It is easy to see that $A_{1n}^{B}$ is bounded above by
\begin{eqnarray*}
& & \sup_{x}||\hat{\Sigma}^{-1}(x)\{\rho (z,\hat{\alpha}_{n}^{B})\rho (z,\hat{%
\alpha}_{n}^{B})^{\prime }-\rho (z,\hat{\alpha}_{n})\rho (z,\hat{\alpha}%
_{n})^{\prime }\}\hat{\Sigma}^{-1}(x)||_{e}n^{-1}\sum_{i=1}^{n}(\omega
_{i,n}-1)^{2}\left\Vert \hat{T}_{i}^{B}[v]\right\Vert _{e}^{2} \\
& \leq & 2\sup_{x}\sup_{\alpha \in \mathcal{N}_{osn}}||\hat{\Sigma}%
^{-1}(x)\{\rho (z,\alpha )\rho (z,\alpha )^{\prime }-\rho (z,\alpha
_{0})\rho (z,\alpha _{0})^{\prime }\}\hat{\Sigma}^{-1}(x)||_{e}n^{-1}%
\sum_{i=1}^{n}(\omega _{i,n}-1)^{2}\left\Vert \hat{T}_{i}^{B}[v]\right\Vert
_{e}^{2}
\end{eqnarray*}%
where $\hat{T}_{i}^{B}[v]\equiv \frac{d\widehat{m}^{B}(X_{i},\hat{\alpha}%
_{n}^{B})}{d\alpha }[v]$. The second line follows because $\hat{\alpha}%
^{B}\in \mathcal{N}_{osn}$ wpa1. The first term in the RHS is of order $%
o_{P_{Z^{\infty }}}(1)$ by Assumption \ref{ass:VE}(iv). The second term is $%
O_{P_{V^{\infty }|Z^{\infty }}}(1)$ by Assumption \ref{ass:VE_boot}(iii).

$A_{2n}^{B}$ is of order $o_{P_{V^{\infty }|Z^{\infty }}}(1)$ $%
wpa1(P_{Z^{\infty }})$ by Assumption \ref{ass:VE_boot}(ii).

Result (1) now follows from the same derivations as in the proof of Theorem %
\ref{thm:VE}(1) step 2a.

Given Result (1), \textbf{Result (2)} follows from exactly the same proof as
that of Theorem \ref{thm:bootstrap_2}(1), and is omitted. \textit{Q.E.D.}

\subsection{Proofs for Section \protect\ref{sec-ex} on examples}

\noindent \textbf{Proof of Proposition \ref{pro:NPIV-suff}.} By our
assumption over $clsp\{p_{j}:j=1,...,J\}$, $\frac{dm(x,\alpha _{0})}{d\alpha
}[u_{n}^{\ast }]\in clsp\{p_{j}:j=1,...,J_{n}\}$ provided $k(n)\leq J_{n}$,
and thus Assumption \ref{ass:anor-mtilde} (i) trivially holds. Since $\Sigma
=1$, Assumption \ref{ass:anor-mtilde} (ii) is the same as Assumption \ref%
{ass:anor-mtilde} (i).

We now show that Assumption \ref{ass:anor-mtilde}(iii)(iv) holds under
condition \ref{eqn:NPIV-holder}. First, condition \ref{eqn:NPIV-holder}(i)
implies that $\{(E[h(Y_{2})-h_{0}(Y_{2})|\cdot ])^{2}:h\in \mathcal{H}\}$ is
a P-Donsker class and, moreover,
\begin{equation*}
E[(E[h(Y_{2})-h_{0}(Y_{2})|X])^{4}]\leq 2c\times ||h-h_{0}||^{2}\rightarrow 0
\end{equation*}%
as $||h-h_{0}||_{L^{2}(f_{Y_{2}})}\rightarrow 0$. So by Lemma 1 in \cite%
{CLvK_Emetrica03}, Assumption \ref{ass:anor-mtilde}(iii) holds. Regarding
Assumption \ref{ass:anor-mtilde}(iv). By Theorem 2.14.2 in VdV-W, (up to
omitted constants)
\begin{equation*}
E\left[ \left\vert \sup_{f\in \mathcal{F}_{n}}n^{-1/2}\sum_{i=1}^{n}%
\{f(X_{i})-E[f(X_{i})]\}\right\vert \right] \leq
\int_{0}^{||F_{n}||_{L^{2}(f_{X})}}\sqrt{1+\log N_{[]}(u,\mathcal{F}%
_{n},||\cdot ||_{L^{2}(f_{X})})}du
\end{equation*}%
where $\mathcal{F}_{n}\equiv \{f:f=g(\cdot ,u_{n}^{\ast })(m(\cdot ,\alpha
)-m(\cdot ,\alpha _{0})),~some~\alpha \in \mathcal{N}_{osn}\}$ and
\begin{equation*}
F_{n}(x)\equiv \sup_{\mathcal{F}_{n}}|f(x)|=\sup_{\alpha \in \mathcal{N}%
_{osn}}|g(x,u_{n}^{\ast })\{m(x,\alpha )-m(x,\alpha _{0})\}|.
\end{equation*}

We claim that, under our assumptions,
\begin{align*}
N_{[]}(u,\mathcal{F}_{n}, ||\cdot||_{L^{2}(f_{X})}) \leq
N_{[]}(u,\Lambda_{c}^{\gamma}(\mathcal{X}), ||\cdot||_{L^{\infty}}).
\end{align*}
To show this claim, it suffices to show that given a radius $\delta>0$, if
we take $\{ [l_{j},u_{j}] \}_{j=1}^{N(\delta)}$ to be brackets of $%
\Lambda_{c}^{\gamma}(\mathcal{X})$ under $||\cdot||_{L^{\infty}}$, then we
can construct $\{ [l_{n,j},u_{n,j}] \}_{j=1}^{N(\delta)}$ such that: they
are valid brackets of $\mathcal{F}_{n}$, under $||\cdot||_{L^{2}(f_{X})}$.
To show this, observe that, for any $f_{n} \in \mathcal{F}_{n}$, there
exists a $\alpha \in \mathcal{N}_{osn}$, such that $f_{n} =
g(\cdot,u_{n}^{\ast })\{m(\cdot,\alpha ) - m(\cdot,\alpha_{0})\}$, and under
condition \ref{eqn:NPIV-holder}, it follows that there exists a $j \in
\{1,...,N(\delta)\}$ such that
\begin{align}
l_{j} \leq m(\cdot,\alpha ) - m(\cdot,\alpha_{0}) \leq u_{j},
\end{align}
hence, there exists a $[l_{n,j},u_{n,j}]$ such that, for all $x$,
\begin{align*}
l_{n,j}(x) = ( 1\{ g(x,u^{\ast}_{n}) > 0\} l_{j}(x) + 1\{ g(x,u^{\ast}_{n})
< 0\} u_{j}(x) ) g(x,u^{\ast}_{n}),
\end{align*}
and
\begin{align*}
u_{n,j}(x) = ( 1\{ g(x,u^{\ast}_{n}) > 0\} u_{j}(x) + 1\{ g(x,u^{\ast}_{n})
< 0\} l_{j}(x) ) g(x,u^{\ast}_{n}).
\end{align*}
such that $l_{n,j} \leq f_{n} \leq u_{n,j}$. Also, observe that
\begin{align*}
||l_{n,j} - u_{n,j}||_{L^{2}(f_{X})} = \sqrt{E[(g(X,u_{n}^{\ast }))^{2} (
u_{j}(X) - l_{j}(X) )^{2} ]} \leq ||u_{j} - l_{j}||_{L^{\infty}} \leq \delta
\end{align*}
because $E[(g(X,u_{n}^{\ast }))^{2}] = ||u^{\ast}_{n}||^{2} = 1$ and $%
||u_{j} - l_{j}||_{L^{\infty}} \leq \delta$ by construction.

Therefore,
\begin{align*}
E \left[ \left| \sup_{f \in \mathcal{F}_{n}} n^{-1/2} \sum_{i=1}^{n}\{
f(X_{i}) -E[f(X_{i})] \} \right| \right] \leq
\int_{0}^{||F_{n}||_{L^{2}(f_{X})}} \sqrt{1 + \log
N_{[]}(u,\Lambda_{c}^{\gamma}(\mathcal{X}), ||\cdot||_{L^{\infty}}) } du.
\end{align*}

Since by assumption $\gamma >0.5$, it is well-known that $\sqrt{1+\log
N_{[]}(u,\Lambda _{c}^{\gamma }(\mathcal{X}),||\cdot ||_{L^{\infty }})}$ is
integrable, so in order to show that $E\left[ \left\vert \sup_{f\in \mathcal{%
F}_{n}}n^{-1/2}\sum_{i=1}^{n}\{f(X_{i})-E[f(X_{i})]\}\right\vert \right]
=o(1)$, it suffices to show that $||F_{n}||_{L^{2}(f_{X})}=o(1)$. In order
to show this,
\begin{align*}
||F_{n}||_{L^{2}(f_{X})}\leq & \sqrt{E[(g(X,u_{n}^{\ast }))^{2}(\sup_{%
\mathcal{N}_{osn}}|m(X,\alpha )-m(X,\alpha _{0})|)^{2}]} \\
=& \sqrt{E[(g(X,u_{n}^{\ast }))^{2}(\sup_{\mathcal{N}%
_{osn}}|E[h(Y_{2})-h_{0}(Y_{2})|X]|)^{2}]} \\
=& \sqrt{E[(g(X,u_{n}^{\ast }))^{2}\sup_{\mathcal{N}_{osn}}\int
(h(y_{2})-h_{0}(y_{2}))^{2}f_{Y_{2}|X}(y_{2},X)dy_{2}]} \\
=& \sqrt{E[(g(X,u_{n}^{\ast }))^{2}\sup_{\mathcal{N}_{osn}}\int
(h(y_{2})-h_{0}(y_{2}))^{2}\frac{f_{Y_{2}X}(y_{2},X)}{%
f_{Y_{2}}(y_{2})f_{X}(X)}f_{Y_{2}}(y_{2})dy_{2}]} \\
\leq & \sup_{x,y_{2}}\frac{f_{Y_{2}X}(y_{2},x)}{f_{Y_{2}}(y_{2})f_{X}(x)}%
\sup_{\mathcal{N}_{osn}}||h-h_{0}||_{L^{2}(f_{Y_{2}})}\sqrt{%
E[(g(X,u_{n}^{\ast }))^{2}]} \\
\leq & Const.\times M_{n}\delta _{s,n}\rightarrow 0
\end{align*}%
where the last expression follows from the fact that $E[(g(X,u_{n}^{\ast
}))^{2}]=||u_{n}^{\ast }||^{2}=1$ and condition \ref{eqn:NPIV-holder}(ii),
that states that
\begin{equation*}
\sup_{x,y_{2}}\frac{f_{Y_{2}X}(y_{2},x)}{f_{Y_{2}}(y_{2})f_{X}(x)}\leq
Const.<\infty .
\end{equation*}

Hence, $E\left[ \left\vert \sup_{f\in \mathcal{F}_{n}}n^{-1/2}\sum_{i=1}^{n}%
\{f(X_{i})-E[f(X_{i})]\}\right\vert \right] =o(1)$ which implies assumption %
\ref{ass:anor-mtilde}(iv). Finally, Assumption \ref{ass:cont_diffm} is
automatically satisfied with the NPIV model. \textit{Q.E.D.}

\medskip

\noindent \textbf{Proof of Proposition \ref{pro:NPQIV-suff}.} Assumptions %
\ref{ass:anor-mtilde}(i) and (ii) hold by the same calculations as those in
the proof of Proposition \ref{pro:NPIV-suff} (for the NPIV model). Also,
under Condition \ref{eqn:NPQIV-holder}(i), $%
\{E[F_{Y_{1}|Y_{2}X}(h(Y_{2}),Y_{2},\cdot )|\cdot ]:h\in \mathcal{H}%
\}\subseteq \Lambda _{c}^{\gamma }(\mathcal{X})$ with $\gamma >0.5$,
Assumptions \ref{ass:anor-mtilde}(iii) and (iv) hold by similar calculations
to those in the proof of Proposition \ref{pro:NPIV-suff}.

Assumption \ref{ass:cont_diffm}(i) is standard in the literature. Regarding
Assumption \ref{ass:cont_diffm}(ii), observe that for any $h\in \mathcal{N}%
_{osn}$,%
\begin{eqnarray*}
& &\left\vert \frac{dm(x,h)}{dh}[u_{n}^{\ast }]-\frac{dm(x,h_{0})}{dh}%
[u_{n}^{\ast }]\right\vert \\
&=& \left\vert E\left[
\{f_{Y_{1}|Y_{2}X}(h(Y_{2}),Y_{2},x)-f_{Y_{1}|Y_{2}X}(h_{0}(Y_{2}),Y_{2},x)%
\}u_{n}^{\ast }(Y_{2})\mid X=x\right] \right\vert \\
&=& \left\vert \int \left\{ \int_{0}^{1}\frac{%
df_{Y_{1}|Y_{2}X}(h_{0}(t)(y_{2}),y_{2},x)}{dy_{1}}%
(h(y_{2})-h_{0}(y_{2}))u_{n}^{\ast }(y_{2})dt\right\}
f_{Y_{2}|X}(y_{2},x)dy_{2}\right\vert \\
&=& \left\vert \int \left( \int_{0}^{1}\frac{%
df_{Y_{1}|Y_{2}X}(h_{0}(t)(y_{2}),y_{2},x)}{dy_{1}}dt\right)
(h(y_{2})-h_{0}(y_{2}))u_{n}^{\ast }(y_{2})f_{Y_{2}}(y_{2})\left( \frac{%
f_{Y_{2}X}(y_{2},x)}{f_{Y_{2}}(y_{2})f_{X}(x)}\right) dy_{2}\right\vert \\
&=& \left\vert \int \Gamma _{1}(y_{2},x)\Gamma
_{2}(y_{2},x)(h(y_{2})-h_{0}(y_{2}))u_{n}^{\ast
}(y_{2})f_{Y_{2}}(y_{2})dy_{2}\right\vert \\
&\leq & ||\Gamma _{1}(\cdot ,x)\Gamma _{2}(\cdot ,x)||_{L^{\infty }}\times
||h-h_{0}||_{L^{2}(f_{Y_{2}})}||u_{n}^{\ast }||_{L^{2}(f_{Y_{2}})}
\end{eqnarray*}%
where $h_{0}(t)\equiv h_{0}+t\{h-h_{0}\}$ and $\Gamma _{1}(y_{2},x)\equiv
\left( \int_{0}^{1}\frac{df_{Y_{1}|Y_{2}X}(h_{0}(t)(y_{2}),y_{2},x)}{dy_{1}}%
dt\right) $ and $\Gamma _{2}(y_{2},x)\equiv \frac{f_{Y_{2}X}(y_{2},x)}{%
f_{Y_{2}}(y_{2})f_{X}(x)}$; the last line follows from Cauchy-Swarchz
inequality.

Under Condition \ref{eqn:NPQIV-holder}(ii), it follows that
\begin{equation*}
\sup_{y_{1},y_{2},x}|\frac{df_{Y_{1}|Y_{2}X}(y_{1},y_{2},x)}{dy_{1}}|\leq
C<\infty
\end{equation*}%
and, under Condition \ref{eqn:NPIV-holder}(ii), it follows that
\begin{equation*}
\sup_{x,y_{2}}\left\vert \frac{f_{Y_{2}X}(y_{2},x)}{f_{Y_{2}}(y_{2})f_{X}(x)}%
\right\vert \leq C<\infty .
\end{equation*}%
Then it is easy to see that $||\Gamma _{j}(\cdot ,x)||_{L^{\infty
}(f_{Y_{2}})}\leq C<\infty $ for both $j=1,2$. Thus%
\begin{equation*}
\left\vert \frac{dm(x,h)}{dh}[u_{n}^{\ast }]-\frac{dm(x,h_{0})}{dh}%
[u_{n}^{\ast }]\right\vert \leq C^{2}\times
||h-h_{0}||_{L^{2}(f_{Y_{2}})}||u_{n}^{\ast }||_{L^{2}(f_{Y_{2}})}
\end{equation*}

and thus, Assumption \ref{ass:cont_diffm}(ii) is satisfied provided that $%
n\times M_{n}^{2}\delta _{n}^{2}\sup_{h\in \mathcal{N}%
_{osn}}||h-h_{0}||_{L^{2}(f_{Y_{2}})}^{2}||u_{n}^{\ast
}||_{L^{2}(f_{Y_{2}})}^{2}=o(1)$. Since $||u_{n}^{\ast
}||_{L^{2}(f_{Y_{2}})}\leq c\mu _{k(n)}^{-1}$ it suffices to show that
\begin{equation*}
nM_{n}^{4}\delta _{n}^{2}(||\Pi _{n}h_{0}-h_{0}||_{L^{2}(f_{Y_{2}})}+\mu
_{k(n)}^{-1}\delta _{n})^{2}\mu _{k(n)}^{-2}=o(1).
\end{equation*}%
By assumption, $||\Pi _{n}h_{0}-h_{0}||_{L^{2}(f_{Y_{2}})}\leq Const.\times
\mu _{k(n)}^{-1}\delta _{n}=O(\delta _{s,n})$ and $\delta _{n}^{2}\asymp
Const.k(n)/n$ , then it suffices to show that
\begin{equation*}
nM_{n}^{4}\delta _{s,n}^{4}=o(1),
\end{equation*}%
which holds by Condition \ref{con:NPQIV-rate1}.

Regarding Assumption \ref{ass:cont_diffm}(iii), observe that for any $h\in
\mathcal{N}_{osn}$,
\begin{equation*}
\frac{d^{2}m(x,h)}{dh^{2}}[u_{n}^{\ast },u_{n}^{\ast }]=\int \frac{%
df_{Y_{1}|Y_{2}X}(h(y_{2}),y_{2},x)}{dy_{1}}(u_{n}^{\ast
}(y_{2}))^{2}f_{Y_{2}|X}(y_{2},x)dy_{2}.
\end{equation*}%
Again by Conditions \ref{eqn:NPQIV-holder}(ii) and \ref{eqn:NPIV-holder}%
(ii), it follows that $\left\vert \frac{d^{2}m(x,h)}{dh^{2}}[u_{n}^{\ast
},u_{n}^{\ast }]\right\vert \leq C^{2}\times ||u_{n}^{\ast
}||_{L^{2}(f_{Y_{2}})}^{2}$. Since $||u_{n}^{\ast }||_{L^{2}(f_{Y_{2}})}\leq
const\times \mu _{k(n)}^{-1}$, Assumption \ref{ass:cont_diffm}(iii) holds
because
\begin{equation*}
\mu _{k(n)}^{-2}\times (M_{n}\delta _{n})^{2}=o(1),~or~M_{n}^{2}\delta
_{s,n}^{2}=o(1).
\end{equation*}

Finally, we verify Assumption \ref{ass:cont_diffm}(iv). By our previous
calculations
\begin{eqnarray*}
& &\left\vert \frac{dm(x,h_{1})}{dh}[h_{2}-h_{0}]-\frac{dm(x,h_{0})}{dh}%
[h_{2}-h_{0}]\right\vert \\
&=& \left\vert \int \left( \int \frac{%
df_{Y_{1}|Y_{2}X}(h_{0}(y_{2})+t[h_{1}(y_{2})-h_{0}(y_{2})],y_{2},x)}{dy_{1}}%
dt\right)
(h_{1}(y_{2})-h_{0}(y_{2}))(h_{2}(y_{2})-h_{0}(y_{2}))f_{Y_{2}|X}(y_{2},x)dy_{2}\right\vert
\\
&\leq & C^{2}\times \int
|(h_{1}(y_{2})-h_{0}(y_{2}))(h_{2}(y_{2})-h_{0}(y_{2}))|f_{Y_{2}}(y_{2})dy_{2}
\\
&\leq & C^{2}\times
||h_{1}-h_{0}||_{L^{2}(f_{Y_{2}})}||h_{2}-h_{0}||_{L^{2}(f_{Y_{2}})},
\end{eqnarray*}%
where the first inequality follows from Conditions \ref{eqn:NPQIV-holder}%
(ii) and \ref{eqn:NPIV-holder}(ii), and the last one from Cauchy-Swarchz
inequality. This result and Cauchy-Swarchz inequality together imply that
\begin{align*}
& \left\vert E\left[ g(X,u_{n}^{\ast })\left( \frac{dm(X,h_{1})}{dh}%
[h_{2}-h_{0}]-\frac{dm(X,h_{0})}{dh}[h_{2}-h_{0}]\right) \right] \right\vert
\\
\leq & C^{2}\sqrt{E[(g(X,u_{n}^{\ast }))^{2}]}%
||h_{1}-h_{0}||_{L^{2}(f_{Y_{2}})}||h_{2}-h_{0}||_{L^{2}(f_{Y_{2}})} \\
\leq & const\times
||h_{1}-h_{0}||_{L^{2}(f_{Y_{2}})}||h_{2}-h_{0}||_{L^{2}(f_{Y_{2}})},
\end{align*}%
where the last line follows from $E[(g(X,u_{n}^{\ast }))^{2}]=||u^{\ast
}||^{2}\asymp 1$. Thus, Assumption \ref{ass:cont_diffm}(iv) follows if
\begin{equation*}
\delta _{s,n}^{2}=(||\Pi _{n}h_{0}-h_{0}||_{L^{2}(f_{Y_{2}})}+\mu
_{k(n)}^{-1}\delta _{n})^{2}=o(n^{-1/2})
\end{equation*}%
which holds by Condition \ref{con:NPQIV-rate1}. \textit{Q.E.D.}

\pagebreak

\section{Proofs of the Results in Appendix A}

\label{app:appC}

In Appendix \ref{app:appC}, we provide the proofs of all the lemmas,
theorems and propositions stated in Appendix \ref{app:appA}.

\subsection{Proofs for Section \protect\ref{sec:bconsistency} on convergence
rates of bootstrap PSMD estimators}

\noindent \textbf{Proof of Lemma \ref{lem:cons_boot}}: For \textbf{Result (1)%
}, we prove this result in two steps. First, we show that $\widehat{\alpha }%
_{n}^{B}\in \mathcal{A}_{k(n)}^{M_{0}}$ wpa1-$P_{V^{\infty }|Z^{\infty }}$
for any $Z^{\infty }$ in a set that occurs with $P_{Z^{\infty }}$ probability approaching one, where $\mathcal{A}_{k(n)}^{M_{0}}$ is defined in the
text. Second, we establish consistency, using the fact that we are in the $%
\mathcal{A}_{k(n)}^{M_{0}}$ set.\newline

\textsc{Step 1.} We show that for any $\delta >0$, there exists a $N(\delta
) $ such that
\begin{equation*}
P_{Z^{\infty }}\left( P_{V^{\infty }|Z^{\infty }}\left( \hat{\alpha}%
_{n}^{B}\notin \mathcal{A}_{k(n)}^{M_{0}} | Z^{n}\right) <\delta
\right) \geq 1-\delta ,~\forall n\geq N(\delta ).
\end{equation*}

To show this, note that, by definition of $\widehat{\alpha }_{n}^{B}$,
\begin{equation*}
\lambda _{n}Pen(\widehat{h}_{n}^{B})\leq \widehat{Q}_{n}^{B}(\widehat{\alpha
}_{n})+\lambda _{n}Pen(\widehat{h}_{n})+o_{P_{V^{\infty }|Z^{\infty }}}(%
\frac{1}{n}),~wpa1(P_{Z^{\infty }}).
\end{equation*}

By Assumption \ref{ass:rates_B}(i) and the definition of $\widehat{\alpha }%
_{n}\in \mathcal{A}_{k(n)}$,
\begin{eqnarray*}
\lambda _{n}Pen(\widehat{h}_{n}^{B}) &\leq &c_{0}^{\ast }\left( \widehat{Q}%
_{n}(\widehat{\alpha }_{n})+\lambda _{n}Pen(\widehat{h}_{n})\right)
+o_{P_{V^{\infty }|Z^{\infty }}}(\frac{1}{n}),~wpa1(P_{Z^{\infty }}) \\
&\leq &c_{0}^{\ast }\left( \widehat{Q}_{n}(\Pi _{n}\alpha _{0})+\lambda
_{n}Pen(\Pi _{n}h_{0})\right) +o_{P_{V^{\infty }|Z^{\infty }}}(\frac{1}{n}%
),~wpa1(P_{Z^{\infty }}).
\end{eqnarray*}%
By Assumptions \ref{A_3.6}(i)(ii) and \ref{ass:rates}(i),
\begin{equation*}
\lambda _{n}Pen(\widehat{h}_{n}^{B})\leq c_{0}^{\ast }c_{0}Q(\Pi _{n}\alpha
_{0})+\lambda _{n}Pen(h_{0})+O_{P_{V^{\infty }|Z^{\infty }}}(\lambda _{n}+o(%
\frac{1}{n})),~wpa1(P_{Z^{\infty }}).
\end{equation*}%
By the fact that $Q(\Pi _{n}\alpha _{0})+o(\frac{1}{n})=O(\lambda _{n})$,
the desired result follows.\newline

\textsc{Step 2.} We want to show that for any $\delta >0$, there exists a $%
N(\delta )$ such that
\begin{equation*}
P_{Z^{\infty }}\left( P_{V^{\infty }|Z^{\infty }}\left( ||\hat{\alpha}%
_{n}^{B}-\alpha _{0}||_{s}\geq \delta | Z^{n}\right) <\delta \right)
\geq 1-\delta ,~\forall n\geq N(\delta ),
\end{equation*}%
which is equivalent to show that $P_{Z^{\infty }}(P_{V^{\infty }|Z^{\infty
}}\left( ||\hat{\alpha}_{n}^{B}-\alpha _{0}||_{s}\geq \delta |
Z^{n}\right) >\delta )\leq \delta $ eventually. Note that
\begin{align*}
& P_{Z^{\infty }}\left( P_{V^{\infty }|Z^{\infty }}\left( ||\hat{\alpha}%
_{n}^{B}-\alpha _{0}||_{s}\geq \delta | Z^{n}\right) >\delta \right)
\\
\leq & P_{Z^{\infty }}\left( P_{V^{\infty }|Z^{\infty }}\left( \left\{ ||%
\hat{\alpha}_{n}^{B}-\alpha _{0}||_{s}\geq \delta \right\} \cap \{\hat{\alpha%
}_{n}^{B}\in \mathcal{A}_{k(n)}^{M_{0}}\}|Z^{n}\right) >0.5\delta
\right) \\
& +P_{Z^{\infty }}\left( P_{V^{\infty }|Z^{\infty }}\left( \hat{\alpha}%
_{n}^{B}\notin \mathcal{A}_{k(n)}^{M_{0}}|Z^{n}\right) >0.5\delta
\right) .
\end{align*}%
By step 1, the second summand in the RHS is negligible. Thus, it suffices to
show that
\begin{equation*}
P_{Z^{\infty }}\left( P_{V^{\infty }|Z^{\infty }}\left( \hat{\alpha}%
_{n}^{B}\in \mathcal{A}_{k(n)}^{M_{0}}\colon ||\hat{\alpha}_{n}^{B}-\alpha
_{0}||_{s}\geq \delta | Z^{n}\right) <\delta \right) \geq 1-\delta
,~\forall n\geq N(\delta ).
\end{equation*}%
(henceforth, we omit $\hat{\alpha}_{n}^{B}\in \mathcal{A}_{k(n)}^{M_{0}}$).
Note that, conditioning on $Z^{n}$, by Assumption \ref{ass:rates_B}(i)(ii),
the definition of $\widehat{\alpha }_{n}\in \mathcal{A}_{k(n)}^{M_{0}}$,
Assumption \ref{A_3.6}(i)(ii) and $\max \{\lambda _{n},o(\frac{1}{n}%
)\}=O(\lambda _{n})$, we have:
\begin{eqnarray*}
& & P_{V^{\infty }|Z^{\infty }}\left( ||\hat{\alpha}_{n}^{B}-\alpha
_{0}||_{s}\geq \delta | Z^{n}\right) \\
& \leq & P_{V^{\infty }|Z^{\infty }}\left( \inf_{\{\mathcal{A}%
_{k(n)}^{M_{0}}\colon ||\alpha -\alpha _{0}||_{s}\geq \delta \}}\left\{
\widehat{Q}_{n}^{B}(\alpha )+\lambda _{n}Pen(h)\right\} \leq \widehat{Q}%
_{n}^{B}(\widehat{\alpha }_{n})+\lambda _{n}Pen(\widehat{h}_{n})+o(\frac{1}{n%
})| Z^{n}\right) \\
& \leq & P_{V^{\infty }|Z^{\infty }}\left( \inf_{\{\mathcal{A}%
_{k(n)}^{M_{0}}\colon ||\alpha -\alpha _{0}||_{s}\geq \delta \}}\left\{
c^{\ast }\widehat{Q}_{n}(\alpha )+\lambda _{n}Pen(h)\right\} \leq
c_{0}^{\ast }\left[ \widehat{Q}_{n}(\widehat{\alpha }_{n})+\lambda _{n}Pen(%
\widehat{h}_{n})\right] +O(\lambda _{n})+\overline{\delta }_{m,n}^{\ast 2}%
| Z^{n}\right) \\
& \leq & P_{V^{\infty }|Z^{\infty }}\left( \inf_{\{\mathcal{A}%
_{k(n)}^{M_{0}}\colon ||\alpha -\alpha _{0}||_{s}\geq \delta \}}\left\{
c^{\ast }\widehat{Q}_{n}(\alpha )\right\} \leq c_{0}^{\ast }\left[ \widehat{Q%
}_{n}(\Pi _{n}\alpha _{0})+\lambda _{n}Pen(\Pi _{n}h_{0})\right] +O(\lambda
_{n})+\overline{\delta }_{m,n}^{\ast 2} | Z^{n}\right) .
\end{eqnarray*}

Thus, wpa1($P_{Z^{\infty }}$),%
\begin{equation*}
P_{V^{\infty }|Z^{\infty }}\left( ||\hat{\alpha}_{n}^{B}-\alpha
_{0}||_{s}\geq \delta | Z^{n}\right) \leq P_{V^{\infty }|Z^{\infty
}}\left( \inf_{\{\mathcal{A}_{k(n)}^{M_{0}}\colon ||\alpha -\alpha
_{0}||_{s}\geq \delta \}}c^{\ast }\widehat{Q}_{n}(\alpha )\leq c_{0}^{\ast }%
\widehat{Q}_{n}(\Pi _{n}\alpha _{0})+M(\lambda _{n}+\overline{\delta }%
_{m,n}^{\ast 2}) | Z^{n}\right) ,
\end{equation*}%
which can be bounded above by
\begin{align*}
& P_{V^{\infty }|Z^{\infty }}\left( \inf_{\{\mathcal{A}_{k(n)}^{M_{0}}\colon
||\alpha -\alpha _{0}||_{s}\geq \delta \}}c^{\ast }cQ(\alpha )\leq
c_{0}^{\ast }c_{0}Q(\Pi _{n}\alpha _{0})+M(\lambda _{n}+(\overline{\delta }%
_{m,n}+\overline{\delta }_{m,n}^{\ast })^{2})|Z^{n}\right) \\
& +P_{V^{\infty }|Z^{\infty }}\left( \sup_{\{\mathcal{A}_{k(n)}^{M_{0}}%
\colon ||\alpha -\alpha _{0}||_{s}\geq \delta \}}\widehat{Q}_{n}(\alpha
)-cQ(\alpha )<-M\overline{\delta }_{m,n}^{2}|Z^{n}\right) \\
& +P_{V^{\infty }|Z^{\infty }}\left( \widehat{Q}_{n}(\Pi _{n}\alpha
_{0})-c_{0}Q(\Pi _{n}\alpha _{0})>-o(\frac{1}{n})|Z^{n}\right) .
\end{align*}%
Therefore, for sufficiently large $n$,
\begin{eqnarray*}
& & P_{Z^{\infty }}\left( P_{V^{\infty }|Z^{\infty }}\left( ||\hat{\alpha}%
_{n}^{B}-\alpha _{0}||_{s}\geq \delta | Z^{n}\right) <\delta \right)
\leq 0.25\delta \\
& & +P_{Z^{\infty }}\left( \inf_{\{\mathcal{A}_{k(n)}^{M_{0}}\colon ||\alpha
-\alpha _{0}||_{s}\geq \delta \}}c^{\ast }cQ(\alpha )\leq c_{0}^{\ast
}c_{0}Q(\Pi _{n}\alpha _{0})+M(\lambda _{n}+(\overline{\delta }_{m,n}+%
\overline{\delta }_{m,n}^{\ast })^{2})\right) \\
& & +P_{Z^{\infty }}\left( \sup_{\{\mathcal{A}_{k(n)}^{M_{0}}\colon ||\alpha
-\alpha _{0}||_{s}\geq \delta \}}\widehat{Q}_{n}(\alpha )-cQ(\alpha )<-M%
\overline{\delta }_{m,n}^{2}\right) \\
&& +P_{Z^{\infty }}\left( \widehat{Q}%
_{n}(\Pi _{n}\alpha _{0})-c_{0}Q(\Pi _{n}\alpha _{0})>-o(\frac{1}{n})\right)
.
\end{eqnarray*}%
By Assumption \ref{ass:rates}, the third and fourth terms in the RHS are
less than $0.5\delta $. The second term in the RHS is not random. By
Assumptions \ref{ass:sieve}(ii) and \ref{A_3.6}(iii), $\mathcal{A}%
_{k(n)}^{M_{0}}$ is compact, and so is $\mathcal{A}^{M_{0}}\equiv \{\alpha
=(\theta ^{\prime },h)\in \mathcal{A}:\lambda _{n}Pen(h)\leq \lambda
_{n}M_{0}\}$. This fact, and Assumption \ref{ass:sieve}(iii) imply that $%
\inf_{\{\mathcal{A}_{k(n)}^{M_{0}}\colon ||\alpha -\alpha _{0}||_{s}\geq
\delta \}}c^{\ast }cQ(\alpha )\geq Q(\alpha (\delta ))$ some $\alpha (\delta
)\in \mathcal{A}^{M_{0}}\cap \{||\alpha -\alpha _{0}||_{s}\geq \delta \}$.
By Assumption \ref{ass:sieve}(i), $Q(\alpha (\delta ))>0$, so eventually,
since $c_{0}^{\ast }c_{0}Q(\Pi _{n}\alpha _{0})+M(\lambda _{n}+(\overline{%
\delta }_{m,n}+\overline{\delta }_{m,n}^{\ast })^{2})=o(1)$,
\begin{align*}
	P_{Z^{\infty
		}}\left( \inf_{\{\mathcal{A}_{k(n)}^{M_{0}}\colon ||\alpha -\alpha
		_{0}||_{s}\geq \delta \}}c^{\ast }cQ(\alpha )\leq c_{0}^{\ast }c_{0}Q(\Pi
	_{n}\alpha _{0})+M(\lambda _{n}+(\overline{\delta }_{m,n}+\overline{\delta }%
	_{m,n}^{\ast })^{2})\right) =0.
\end{align*}
\newline

For \textbf{Result (2)}, we want to show that for any $\delta >0$, there
exists a $M(\delta )$ such that
\begin{equation*}
P_{Z^{\infty }}\left( P_{V^{\infty }|Z^{\infty }}\left( \delta _{n}^{-1}||%
\hat{\alpha}_{n}^{B}-\alpha _{0}||\geq M^{\prime }\text{ }| Z^{n}\right) <\delta \right) \geq 1-\delta ,~\forall M^{\prime }\geq
M(\delta )
\end{equation*}%
eventually. By Assumptions \ref{ass:weak_equiv}(iii) and \ref{ass:rates_B}%
(iii), following the similar algebra as before, we have: for $M^{\prime }$
large enough,
\begin{align*}
& P_{V^{\infty }|Z^{\infty }}\left( \delta _{n}^{-1}||\hat{\alpha}%
_{n}^{B}-\alpha _{0}||\geq M^{\prime }\text{ }|Z^{n}\right) \\
\leq & P_{V^{\infty }|Z^{\infty }}\left( \inf_{\{\mathcal{A}_{osn}:\text{ }%
\delta _{n}^{-1}||\alpha -\alpha _{0}||\geq M^{\prime }\}}c^{\ast }cQ(\alpha
)\leq M(\lambda _{n}+\delta _{n}^{2})| Z^{n}\right) +\delta .
\end{align*}%
By Assumption \ref{ass:weak_equiv}(i)(ii) and $\delta _{n}=\sqrt{\max
\{\lambda _{n},\delta _{n}^{2}\}}$, we have:
\begin{align*}
& P_{V^{\infty }|Z^{\infty }}\left( \inf_{\{\mathcal{A}_{osn}:\text{ }\delta
_{n}^{-1}||\alpha -\alpha _{0}||\geq M^{\prime }\}}c^{\ast }cQ(\alpha )\leq
M(\lambda _{n}+\delta _{n}^{2})|Z^{n}\right) \\
\leq & 1\{c^{\ast }cc_{1}\left( M^{\prime }\delta _{n}\right) ^{2}\leq
M(\lambda _{n}+\delta _{n}^{2})\},
\end{align*}%
which is eventually naught, because $M^{\prime }$ can be chosen to be large.
The rate under $||\cdot ||_{s}$ immediately follows from this result and the
definition of the sieve measure of local ill-posedness $\tau _{n}$.

For \textbf{Result (3)}, we note that both $\hat{\alpha}_{n}^{R,B},\widehat{%
\alpha }_{n}\in \{\alpha \in \mathcal{A}_{k(n)}\colon \phi (\alpha )=\phi (%
\widehat{\alpha }_{n})\}$, and hence all the above proofs go through with $%
\hat{\alpha}_{n}^{R,B}$ replacing $\hat{\alpha}_{n}^{B}$. In particular, let
$\mathcal{A}_{k(n)}^{M_{0}}(\widehat{\phi })\equiv \{\alpha \in \mathcal{A}%
_{k(n)}^{M_{0}}\colon \phi (\alpha )=\phi (\widehat{\alpha }_{n})\}\subseteq
\mathcal{A}_{k(n)}^{M_{0}}$. Then: for any $\delta >0$,%
\begin{eqnarray*}
& & P_{V^{\infty }|Z^{\infty }}\left( \hat{\alpha}_{n}^{R,B}\in \mathcal{A}%
_{k(n)}^{M_{0}}(\widehat{\phi }):||\hat{\alpha}_{n}^{R,B}-\alpha
_{0}||_{s}\geq \delta | Z^{n}\right) \\
& \leq & P_{V^{\infty }|Z^{\infty }}\left( \inf_{\{\mathcal{A}_{k(n)}^{M_{0}}(%
\widehat{\phi })\colon ||\alpha -\alpha _{0}||_{s}\geq \delta \}}\left\{
\widehat{Q}_{n}^{B}(\alpha )+\lambda _{n}Pen(h)\right\} \leq \widehat{Q}%
_{n}^{B}(\widehat{\alpha }_{n})+\lambda _{n}Pen(\widehat{h}_{n})+o(\frac{1}{n%
})|Z^{n}\right) \\
& \leq & P_{V^{\infty }|Z^{\infty }}\left( \inf_{\{\mathcal{A}_{k(n)}^{M_{0}}(%
\widehat{\phi })\colon ||\alpha -\alpha _{0}||_{s}\geq \delta \}}\left\{
c^{\ast }\widehat{Q}_{n}(\alpha )+\lambda _{n}Pen(h)\right\} \leq A_{n}
| Z^{n}\right) \\
& \leq & P_{V^{\infty }|Z^{\infty }}\left( \inf_{\{\mathcal{A}%
_{k(n)}^{M_{0}}\colon ||\alpha -\alpha _{0}||_{s}\geq \delta \}}\left\{
c^{\ast }\widehat{Q}_{n}(\alpha )\right\} \leq c_{0}^{\ast }\left[ \widehat{Q%
}_{n}(\Pi _{n}\alpha _{0})+\lambda _{n}Pen(\Pi _{n}h_{0})\right] +O(\lambda
_{n})+\overline{\delta }_{m,n}^{\ast 2}|Z^{n}\right) .
\end{eqnarray*}
where $A_{n} \equiv c_{0}^{\ast }\left[ \widehat{Q}_{n}(\widehat{\alpha }_{n})+\lambda _{n}Pen(%
\widehat{h}_{n})\right] +O(\lambda _{n})+\overline{\delta }_{m,n}^{\ast 2}$.
The rest follows from the proof of Results (1) and (2). \textit{Q.E.D.}

\subsection{Proofs for Section \protect\ref{subsec-A5} on behaviors under
local alternatives}

\noindent \textbf{Proof of Theorem \ref{thm:chi2_localt}}: The proof is
analogous to that of Theorem \ref{thm:chi2}, hence we only present the main
steps. Let $\boldsymbol{\alpha }_{n}=\alpha _{0}+d_{n}\Delta _{n}$ with $%
\frac{d\phi (\alpha _{0})}{d\alpha }[\Delta _{n}]=\langle v_{n}^{\ast
},\Delta _{n}\rangle =\kappa _{n}=\kappa \times \left( 1+o(1)\right) \neq 0$.\\

\textsc{Step 1.} By assumption \ref{ass:LAQ}(i) under the local
alternatives, for any $t_{n}\in \mathcal{T}_{n}$,
\begin{equation}
0\leq 0.5\left( \widehat{Q}_{n}(\widehat{\alpha }_{n}(t_{n}))-\widehat{Q}%
_{n}(\widehat{\alpha }_{n})\right) =t_{n}\left\{ \mathbb{Z}_{n}(\boldsymbol{%
\alpha }_{n})+\langle u_{n}^{\ast },\widehat{\alpha }_{n}-\boldsymbol{\alpha
}_{n}\rangle \right\} +\frac{B_{n}}{2}t_{n}^{2}+o_{P_{n,Z^{\infty
}}}([r_{n}(t_{n})]^{-1})  \label{alt-u}
\end{equation}%
where $[r_{n}(t_{n})]^{-1}=\max \{t_{n}^{2},t_{n}n^{-1/2},s_{n}^{-1}\}$ and $%
s_{n}^{-1}=o(n^{-1})$. The LHS is always positive (up to possibly a
negligible terms given by the penalty function, see the proof of Theorem \ref%
{thm:theta_anorm}(1) for details) by definition of $\widehat{\alpha }_{n}$.
Hence, by choosing $t_{n}=\pm \{s_{n}^{-1/2}+o(n^{-1/2})\}$, it follows that
$\{\mathbb{Z}_{n}(\boldsymbol{\alpha }_{n})+\langle u_{n}^{\ast },\widehat{%
\alpha }_{n}-\boldsymbol{\alpha }_{n}\rangle \}=o_{P_{n,Z^{\infty
}}}(n^{-1/2})$. Since $\langle u_{n}^{\ast },\boldsymbol{\alpha }_{n}-\alpha
_{0}\rangle =\frac{d_{n}\kappa _{n}}{||v_{n}^{\ast }||_{sd}}$ by the
definition of local alternatives $\boldsymbol{\alpha }_{n}$, we obtain
equation (\ref{eqn:LA-LAR}):
\begin{equation}
\left\{ \mathbb{Z}_{n}(\boldsymbol{\alpha }_{n})+\langle u_{n}^{\ast },%
\widehat{\alpha }_{n}-\alpha _{0}\rangle -\frac{d_{n}\kappa _{n}}{||v_{n}^{\ast }||_{sd}}\right\} =\mathbb{Z}_{n}(\boldsymbol{\alpha }%
_{n})+\langle u_{n}^{\ast },\widehat{\alpha }_{n}-\boldsymbol{\alpha }%
_{n}\rangle =o_{P_{n,Z^{\infty }}}(n^{-1/2}),  \label{eqn:LA-LAR}
\end{equation}%
where $\mathbb{Z}_{n}(\boldsymbol{\alpha }_{n})$ is defined as that of $%
\mathbb{Z}_{n}$ but using $\rho (z,\boldsymbol{\alpha }_{n})$ instead of $%
\rho (z,\alpha _{0})$ (since $m(X,\boldsymbol{\alpha }_{n})=0$ a.s.-$X$
under the local alternative).

Next, by Assumption \ref{ass:LAQ}(i) under the local alternative, we have:
for any $t_{n}\in \mathcal{T}_{n}$,
\begin{equation}
0.5\left( \widehat{Q}_{n}(\widehat{\alpha }_{n}^{R}(t_{n}))-\widehat{Q}_{n}(%
\widehat{\alpha }_{n}^{R})\right) =t_{n}\left\{ \mathbb{Z}_{n}(\boldsymbol{%
\alpha }_{n})+\langle u_{n}^{\ast },\widehat{\alpha }_{n}^{R}-\boldsymbol{%
\alpha }_{n}\rangle \right\} +\frac{B_{n}}{2}t_{n}^{2}+o_{P_{n,Z^{\infty
}}}([r_{n}(t_{n})]^{-1}).  \label{alt-r}
\end{equation}%
By Assumption \ref{ass:phi}(ii)
\begin{equation*}
\sup_{\alpha \in \mathcal{N}_{0n}}\left\vert \phi (\alpha )-\phi (\alpha
_{0})-\frac{d\phi (\alpha _{0})}{d\alpha }[\alpha -\alpha _{0}]\right\vert
=o(n^{-1/2}||v_{n}^{\ast }||),
\end{equation*}%
and assumption $\widehat{\alpha }_{n}^{R}\in \mathcal{N}_{osn}$ wpa1-$%
P_{n,Z^{\infty }}$, and the fact that $\phi (\widehat{\alpha }_{n}^{R})-\phi
(\alpha _{0})=0$, following the same calculations as those in Step 1 of the
proof of Theorem \ref{thm:chi2}, we have:
\begin{equation*}
\langle u_{n}^{\ast },\widehat{\alpha }_{n}^{R}-\alpha _{0}\rangle
=o_{P_{n,Z^{\infty }}}(n^{-1/2}).
\end{equation*}%

Since $\boldsymbol{\alpha }_{n}=\alpha _{0}+d_{n}\Delta _{n}\in \mathcal{N}%
_{osn}$ with $\frac{d\phi (\alpha _{0})}{d\alpha }[\Delta _{n}]=\langle
v_{n}^{\ast },\Delta _{n}\rangle =\kappa _{n}$, we have:
\begin{equation*}
\langle u_{n}^{\ast },\widehat{\alpha }_{n}^{R}-\boldsymbol{\alpha }%
_{n}\rangle =\langle u_{n}^{\ast },\widehat{\alpha }_{n}^{R}-\alpha
_{0}\rangle -\frac{d_{n}\kappa _{n}}{||v_{n}^{\ast }||_{sd}}%
+o_{P_{n,Z^{\infty }}}(n^{-1/2})=-\frac{d_{n}\kappa _{n}}{||v_{n}^{\ast
}||_{sd}}+o_{P_{n,Z^{\infty }}}(n^{-1/2}).
\end{equation*}%
Therefore, by choosing $t_{n}\equiv -(\mathbb{Z}_{n}(\boldsymbol{\alpha }%
_{n})-\frac{d_{n}\kappa _{n}}{||v_{n}^{\ast }||_{sd}})B_{n}^{-1}$ in (\ref%
{alt-r}) with $[r_{n}(t_{n})]^{-1}=\max
\{t_{n}^{2},t_{n}n^{-1/2},o(n^{-1})\} $ (which is a valid choice), we
obtain:
\begin{eqnarray*}
0.5\left( \widehat{Q}_{n}(\widehat{\alpha }_{n})-\widehat{Q}_{n}(\widehat{%
\alpha }_{n}^{R})\right) &\leq &0.5\left( \widehat{Q}_{n}(\widehat{\alpha }%
_{n}^{R}(t_{n}))-\widehat{Q}_{n}(\widehat{\alpha }_{n}^{R})\right)
+o_{P_{n,Z^{\infty }}}(n^{-1}) \\
&=&-\frac{1}{2}\left( \frac{(\mathbb{Z}_{n}(\boldsymbol{\alpha }_{n})-\frac{%
d_{n}\kappa _{n}}{||v_{n}^{\ast }||_{sd}})}{\sqrt{B_{n}}}\right)
^{2}+o_{P_{n,Z^{\infty }}}([r_{n}(t_{n})]^{-1}).
\end{eqnarray*}%
By our assumption and the fact that $||u_{n}^{\ast }||\geq c>0$ for all $n$,
it follows that $B_{n}\geq c>0$ eventually, so
\begin{equation*}
0.5\left( \widehat{Q}_{n}(\widehat{\alpha }_{n})-\widehat{Q}_{n}(\widehat{%
\alpha }_{n}^{R})\right) \leq -\frac{1}{2}\left( \frac{(\mathbb{Z}_{n}(%
\boldsymbol{\alpha }_{n})-\frac{d_{n}\kappa _{n}}{||v_{n}^{\ast }||_{sd}})}{%
||u_{n}^{\ast }||}\right) ^{2}\times \left( 1+o_{P_{n,Z^{\infty
}}}(1)\right) .
\end{equation*}

\textsc{Step 2.} On the other hand, suppose there exists a $t_{n}^{\ast }$,
such that (a) $\phi (\widehat{\alpha }_{n}(t_{n}^{\ast }))=\phi (\alpha
_{0}) $, $\widehat{\alpha }_{n}(t_{n}^{\ast })\in \mathcal{A}_{k(n)}$, and
(b) $t_{n}^{\ast }=(\mathbb{Z}_{n}(\boldsymbol{\alpha }_{n})-\frac{%
d_{n}\kappa _{n}}{||v_{n}^{\ast }||_{sd}})\left( ||u_{n}^{\ast }||\right)
^{-2}+o_{P_{n,Z^{\infty }}}(n^{-1/2})$. Substituting this into (\ref{alt-u})
with $[r_{n}(t_{n}^{\ast })]^{-1}=\max \{(t_{n}^{\ast })^{2},t_{n}^{\ast
}n^{-1/2},o(n^{-1})\}$, we obtain:
\begin{align*}
0.5\left( \widehat{Q}_{n}(\widehat{\alpha }_{n})-\widehat{Q}_{n}(\widehat{%
\alpha }_{n}^{R})\right) & \geq 0.5\left( \widehat{Q}_{n}(\widehat{\alpha }%
_{n})-\widehat{Q}_{n}(\widehat{\alpha }_{n}(t_{n}^{\ast }))\right)
-o_{P_{n,Z^{\infty }}}(n^{-1}) \\
& =-\frac{B_{n}}{2}(t_{n}^{\ast })^{2}+o_{P_{n,Z^{\infty
}}}([r_{n}(t_{n}^{\ast })]^{-1}) \\
& =-\frac{B_{n}}{2}(\mathbb{Z}_{n}(\boldsymbol{\alpha }_{n})-\frac{%
d_{n}\kappa _{n}}{||v_{n}^{\ast }||_{sd}})^{2}\left( ||u_{n}^{\ast
}||\right) ^{-4}+o_{P_{n,Z^{\infty }}}([r_{n}(t_{n}^{\ast })]^{-1}) \\
& =-\frac{1}{2}\left( \frac{\mathbb{Z}_{n}(\boldsymbol{\alpha }_{n})-\frac{%
d_{n}\kappa _{n}}{||v_{n}^{\ast }||_{sd}}}{||u_{n}^{\ast }||}\right)
^{2}\times \left( 1+o_{P_{n,Z^{\infty }}}(1)\right)
\end{align*}%
where the second line follows from equation (\ref{eqn:LA-LAR}). Finally, we
observe that point (a) follows from Lemma \ref{pro:phi_solve}, with $r=0$.
Point (b) follows by analogous calculations to those in Step 3 of the proof
of Theorem \ref{thm:chi2}, except that now with $\widehat{\alpha }%
(t_{n}^{\ast })=\widehat{\alpha }_{n}+t_{n}^{\ast }u_{n}^{\ast }$,
\begin{eqnarray*}
\phi (\widehat{\alpha }(t_{n}^{\ast }))-\phi (\alpha _{0})& = &\frac{d\phi
(\alpha _{0})}{d\alpha }[\widehat{\alpha }_{n}-\alpha _{0}]+t_{n}^{\ast }%
\frac{||v_{n}^{\ast }||^{2}}{||v_{n}^{\ast }||_{sd}}+o_{P_{n,Z^{\infty
}}}(n^{-1/2}||v_{n}^{\ast }||) \\
& = &-\mathbb{Z}_{n}(\boldsymbol{\alpha }_{n})||v_{n}^{\ast }||_{sd}+\frac{%
d_{n}\kappa _{n}}{||v_{n}^{\ast }||_{sd}}||v_{n}^{\ast }||_{sd}\\
&& +\left( (%
\mathbb{Z}_{n}(\boldsymbol{\alpha }_{n})-\frac{d_{n}\kappa _{n}}{%
||v_{n}^{\ast }||_{sd}})\frac{||v_{n}^{\ast }||_{sd}^{2}}{||v_{n}^{\ast
}||^{2}}\right) \frac{||v_{n}^{\ast }||^{2}}{||v_{n}^{\ast }||_{sd}} \\
& + &o_{P_{n,Z^{\infty }}}(n^{-1/2}||v_{n}^{\ast }||) \\
& = & o_{P_{n,Z^{\infty }}}(n^{-1/2}||v_{n}^{\ast }||)
\end{eqnarray*}%
where the second line follows from equation (\ref{eqn:LA-LAR}) and some
straightforward algebra.

\textsc{Step 3.} Finally, the above calculations and $\kappa _{n}=\kappa
\left( 1+o(1)\right) $ imply that%
\begin{equation}
||u_{n}^{\ast }||^{2}\times \left( \widehat{Q}_{n}(\widehat{\alpha }%
_{n}^{R})-\widehat{Q}_{n}(\widehat{\alpha }_{n})\right) =\left( \mathbb{Z}%
_{n}(\boldsymbol{\alpha }_{n})-\frac{d_{n}\kappa \left( 1+o(1)\right) }{%
||v_{n}^{\ast }||_{sd}}\right) ^{2}\times \left( 1+o_{P_{n,Z^{\infty
}}}(1)\right) .  \label{alt-SQLR}
\end{equation}

For \textbf{Result (1), }equation (\ref{alt-SQLR}) with $%
d_{n}=n^{-1/2}||v_{n}^{\ast }||_{sd}$ implies that
\begin{equation*}
||u_{n}^{\ast }||^{2}\times \widehat{QLR}_{n}(\phi _{0})=\left( \sqrt{n}%
\mathbb{Z}_{n}(\boldsymbol{\alpha }_{n})-\kappa \left( 1+o(1)\right) \right)
^{2}\times \left( 1+o_{P_{n,Z^{\infty }}}(1)\right) \Rightarrow \chi
_{1}^{2}(\kappa ^{2}),
\end{equation*}%
which is due to $\sqrt{n}\mathbb{Z}_{n}(\boldsymbol{\alpha }_{n})\Rightarrow
N(0,1)$ under the local alternatives.

For \textbf{Result (2), }equation (\ref{alt-SQLR}) with $\sqrt{n}\frac{d_{n}%
}{||v_{n}^{\ast }||_{sd}}\rightarrow \infty $ implies that
\begin{eqnarray*}
||u_{n}^{\ast }||^{2}\times \widehat{QLR}_{n}(\phi _{0}) &=&\left( \sqrt{n}%
\mathbb{Z}_{n}(\boldsymbol{\alpha }_{n})-\sqrt{n}\frac{d_{n}\kappa \left(
1+o(1)\right) }{||v_{n}^{\ast }||_{sd}}\right) ^{2}\times \left(
1+o_{P_{n,Z^{\infty }}}(1)\right) \\
&=&\left( O_{P_{n,Z^{\infty }}}(1)-\sqrt{n}\frac{d_{n}\kappa \left(
1+o(1)\right) }{||v_{n}^{\ast }||_{sd}}\right) ^{2}\times \left(
1+o_{P_{n,Z^{\infty }}}(1)\right) ,
\end{eqnarray*}%
where the second line is due to $\sqrt{n}\mathbb{Z}_{n}(\boldsymbol{\alpha }%
_{n})\Rightarrow N(0,1)$ under the local alternatives. Since $\sqrt{n}\frac{%
d_{n}\kappa \left( 1+o(1)\right) }{||v_{n}^{\ast }||_{sd}}\rightarrow \infty
$ $($or $-\infty )$ if $\kappa >0$ (or $\kappa <0$), we have that $%
\lim_{n\rightarrow \infty }\left( ||u_{n}^{\ast }||^{2}\times \widehat{QLR}%
_{n}(\phi _{0})\right) =\infty $ in probability (under the alternative).
\textit{Q.E.D.}

\medskip

\noindent \textbf{Proof of Proposition \ref{pro:ARE-localalt}.} Recall that $%
\widehat{QLR}_{n}^{0}(\phi _{0})$ denotes the optimally-weighted SQLR
statistic. By inspection of the proof of Theorem \ref{thm:chi2_localt}, it
is easy to see that
\begin{equation*}
||u_{n}^{\ast }||^{2}\times \widehat{QLR}_{n}(\phi _{0})=\left( \sqrt{n}%
\mathbb{Z}_{n}(\boldsymbol{\alpha }_{n})-\kappa \right)
^{2}+o_{P_{n,Z^{\infty }}}(1)
\end{equation*}%
and
\begin{equation*}
\widehat{QLR}_{n}^{0}(\phi _{0})=\left( \sqrt{n}\mathbb{Z}_{n}(\boldsymbol{%
\alpha }_{n})-\kappa \frac{||v_{n}^{\ast }||_{sd}}{||v_{n}^{0}||_{0}}\right)
^{2}+o_{P_{n,Z^{\infty }}}(1)
\end{equation*}%
for local alternatives of the form described in equation (\ref{eq:locaalt-1}%
) with $d_{n}=n^{-1/2}||v_{n}^{\ast }||_{sd}$. Hence, the distribution of $%
||u_{n}^{\ast }||^{2}\times \widehat{QLR}_{n}(\phi _{0})$ is, asymptotically
close to $\chi _{1}^{2}(\kappa ^{2})$ and the distribution of $\widehat{QLR}%
_{n}^{0}(\phi _{0})$ is, asymptotically close to $\chi _{1}^{2}\left( \frac{%
||v_{n}^{\ast }||_{sd}^{2}}{||v_{n}^{0}||_{0}^{2}}\kappa ^{2}\right) $.

Let $A_{n}^{0}(z)\equiv \left( \frac{dm(x,\alpha _{0})}{d\alpha }%
[v_{n}^{0}]\right) ^{\prime }\left( \Sigma _{0}(x)\right) ^{-1}\rho
(z,\alpha _{0})$ and $A_{n}(z)\equiv \left( \frac{dm(x,\alpha _{0})}{d\alpha
}[v_{n}^{\ast }]\right) ^{\prime }\left( \Sigma (x)\right) ^{-1}\rho
(z,\alpha _{0})$ where $v_{n}^{0}$ is the Riesz representer under $||\cdot
||_{0}$. Since $\frac{\langle v_{n}^{\ast },v_{n}^{0}\rangle }{\langle
v_{n}^{0},v_{n}^{0}\rangle _{0}}E\left[ \left( A_{n}^{0}(Z)\right) \left(
A_{n}(Z)\right) ^{\prime }\right] =\frac{\langle v_{n}^{\ast
},v_{n}^{0}\rangle }{\langle v_{n}^{0},v_{n}^{0}\rangle _{0}}E\left[ \left(
\frac{dm(X,\alpha _{0})}{d\alpha }[v_{n}^{0}]\right) ^{\prime }\left( \Sigma
(X)\right) ^{-1}\left( \frac{dm(X,\alpha _{0})}{d\alpha }[v_{n}^{\ast
}]\right) \right] =\frac{(\langle v_{n}^{\ast },v_{n}^{0}\rangle )^{2}}{%
\langle v_{n}^{0},v_{n}^{0}\rangle _{0}}$ and $E\left[ \left(
A_{n}^{0}(Z)\right) \left( A_{n}^{0}(Z)\right) ^{\prime }\right] =\langle
v_{n}^{0},v_{n}^{0}\rangle _{0}$, we have:
\begin{align*}
& E\left[ \left( A_{n}(Z)-\frac{\langle v_{n}^{\ast },v_{n}^{0}\rangle }{%
\langle v_{n}^{0},v_{n}^{0}\rangle _{0}}A_{n}^{0}(Z)\right) \left( A_{n}(Z)-%
\frac{\langle v_{n}^{\ast },v_{n}^{0}\rangle }{\langle
v_{n}^{0},v_{n}^{0}\rangle _{0}}A_{n}^{0}(Z)\right) ^{\prime }\right] \\
& =E\left[ \left( A_{n}(Z)\right) \left( A_{n}(Z)\right) ^{\prime }\right] -%
\frac{(\langle v_{n}^{\ast },v_{n}^{0}\rangle )^{2}}{\langle
v_{n}^{0},v_{n}^{0}\rangle _{0}}=\langle v_{n}^{\ast },v_{n}^{\ast }\rangle
_{sd}-\frac{(\langle v_{n}^{\ast },v_{n}^{0}\rangle )^{2}}{\langle
v_{n}^{0},v_{n}^{0}\rangle _{0}}
\end{align*}%
Since the LHS is non-negative, the previous equation implies that $%
||v_{n}^{\ast }||_{sd}^{2}-\frac{(\langle v_{n}^{\ast },v_{n}^{0}\rangle
)^{2}}{\langle v_{n}^{0},v_{n}^{0}\rangle _{0}}\geq 0$. By definition of $%
v_{n}^{\ast }$ and $v_{n}^{0}$, it follows that
\begin{equation*}
\langle v_{n}^{\ast },v_{n}^{0}\rangle =\frac{d\phi (\alpha _{0})}{d\alpha }%
[v_{n}^{0}]=||v_{n}^{0}||_{0}^{2},
\end{equation*}%
and thus $||v_{n}^{\ast }||_{sd}^{2}\geq ||v_{n}^{0}||_{0}^{2}$ for all $n$.

Observe that for a noncentral chi-square, $\chi _{p}^{2}(r)$, $\Pr (\chi
_{p}^{2}(r)\leq t)$ is decreasing in the noncentrality parameter $r$ for
each $t$; thus $\Pr (\chi _{p}^{2}(r_{1})>t)>\Pr (\chi _{p}^{2}(r_{2})>t)$
for $r_{1}>r_{2}$. Therefore, the previous results imply that, for any $t$,
\begin{align*}
\lim_{n\rightarrow \infty }P_{n,Z^{\infty }}\left( ||u_{n}^{\ast
}||^{2}\times \widehat{QLR}_{n}(\phi _{0})\geq t\right) & =\Pr (\chi
_{1}^{2}(\kappa ^{2})\geq t) \\
& \leq \liminf_{n\rightarrow \infty }\Pr \left( \chi _{1}^{2}\left( \frac{%
||v_{n}^{\ast }||_{sd}^{2}}{||v_{n}^{0}||_{0}^{2}}\kappa ^{2}\right) \geq
t\right) \\
& =\liminf_{n\rightarrow \infty }P_{n,Z^{\infty }}(\widehat{QLR}%
_{n}^{0}(\phi _{0})\geq t).
\end{align*}%
\textit{Q.E.D.}

\medskip

\noindent \textbf{Proof of Theorem \ref{thm:BSQLR-loc-alt}}: The proof of
\textbf{Result (1)}\ is similar to that of Theorem \ref{thm:bootstrap}, so
we only present a sketch here. By assumptions \ref{ass:LAQ}(i) and \ref%
{ass:LAQ_B}(i) under local alternative, it follows that
\begin{align*}
& 0.5\left( \widehat{Q}_{n}^{B}(\widehat{\alpha }_{n}^{R,B}(-\frac{\mathbb{Z}%
_{n}^{\omega -1}(\boldsymbol{\alpha }_{n})}{B_{n}^{\omega }}))-\widehat{Q}%
_{n}^{B}(\widehat{\alpha }_{n}^{R,B})\right) \\
& =-\frac{\mathbb{Z}_{n}^{\omega -1}(\boldsymbol{\alpha }_{n})}{%
B_{n}^{\omega }}\{\mathbb{Z}_{n}^{\omega }(\boldsymbol{\alpha }_{n})+\langle
u_{n}^{\ast },\widehat{\alpha }_{n}^{R,B}-\boldsymbol{\alpha }_{n}\rangle \}+%
\frac{(\mathbb{Z}_{n}^{\omega -1}(\boldsymbol{\alpha }_{n}))^{2}}{%
2B_{n}^{\omega }}+o_{P_{V^{\infty }|Z^{\infty
}}}(r_{n}^{-1}),~wpa1(P_{n,Z^{\infty }}),
\end{align*}%
where $r_{n}^{-1}=\max \left\{ \left( -\frac{\mathbb{Z}_{n}^{\omega -1}(%
\boldsymbol{\alpha }_{n})}{B_{n}^{\omega }}\right) ^{2},\left\vert -\frac{%
\mathbb{Z}_{n}^{\omega -1}(\boldsymbol{\alpha }_{n})}{B_{n}^{\omega }}%
\right\vert n^{-1/2},o(n^{-1})\right\} =O_{P_{V^{\infty }|Z^{\infty
}}}(n^{-1}),~wpa1(P_{n,Z^{\infty }})$ under assumption \ref{ass:LAQ_B}%
(i)(ii) with $\boldsymbol{\alpha }_{n}$ (instead of $\alpha _{0}$).

By similar calculations to those in the proof of Result (1) of Theorem \ref%
{thm:bootstrap} (equation (\ref{boot-correct-center})),%
\begin{equation*}
\sqrt{n}\langle u_{n}^{\ast },\widehat{\alpha }_{n}^{R,B}-\widehat{\alpha }%
_{n}\rangle =o_{P_{V^{\infty }|Z^{\infty }}}(1),~wpa1(P_{n,Z^{\infty }}),
\end{equation*}%
i.e., the restricted bootstrap estimator $\widehat{\alpha }_{n}^{R,B}$
centers at $\widehat{\alpha }_{n}$, regardless of the local alternative. Thus%
\begin{equation*}
\langle u_{n}^{\ast },\widehat{\alpha }_{n}^{R,B}-\boldsymbol{\alpha }%
_{n}\rangle =\langle u_{n}^{\ast },\widehat{\alpha }_{n}^{R,B}-\widehat{%
\alpha }_{n}\rangle +\langle u_{n}^{\ast },\widehat{\alpha }_{n}-\boldsymbol{%
\alpha }_{n}\rangle =\langle u_{n}^{\ast },\widehat{\alpha }_{n}-\boldsymbol{%
\alpha }_{n}\rangle +o_{P_{V^{\infty }|Z^{\infty
}}}(n^{-1/2}),~wpa1(P_{n,Z^{\infty }}).
\end{equation*}%
This result and equation (\ref{eqn:LA-LAR}) (i.e., $\mathbb{Z}_{n}(%
\boldsymbol{\alpha }_{n})+\langle u_{n}^{\ast },\widehat{\alpha }_{n}-%
\boldsymbol{\alpha }_{n}\rangle =o_{P_{n,Z^{\infty }}}(n^{-1/2})$) imply
that
\begin{eqnarray*}
& & 0.5\left( \widehat{Q}_{n}^{B}(\widehat{\alpha }_{n}^{R,B}(-\frac{\mathbb{Z}%
_{n}^{\omega -1}(\boldsymbol{\alpha }_{n})}{B_{n}^{\omega }}))-\widehat{Q}%
_{n}^{B}(\widehat{\alpha }_{n}^{R,B})\right) \\
& =&-\frac{\mathbb{Z}_{n}^{\omega -1}(\boldsymbol{\alpha }_{n})}{%
B_{n}^{\omega }}\{\mathbb{Z}_{n}^{\omega }(\boldsymbol{\alpha }_{n})+\langle
u_{n}^{\ast },\widehat{\alpha }_{n}-\boldsymbol{\alpha }_{n}\rangle \}+\frac{%
(\mathbb{Z}_{n}^{\omega -1}(\boldsymbol{\alpha }_{n}))^{2}}{2B_{n}^{\omega }}%
+o_{P_{V^{\infty }|Z^{\infty }}}(r_{n}^{-1}),~wpa1(P_{n,Z^{\infty }}) \\
& =& -\frac{\mathbb{Z}_{n}^{\omega -1}(\boldsymbol{\alpha }_{n})}{%
B_{n}^{\omega }}\{\mathbb{Z}_{n}^{\omega -1}(\boldsymbol{\alpha }%
_{n})+o_{P_{n,Z^{\infty }}}(n^{-1/2})\}+\frac{(\mathbb{Z}_{n}^{\omega -1}(%
\boldsymbol{\alpha }_{n}))^{2}}{2B_{n}^{\omega }}+o_{P_{V^{\infty
}|Z^{\infty }}}(r_{n}^{-1}),~wpa1(P_{n,Z^{\infty }}) \\
& =& -\frac{(\mathbb{Z}_{n}^{\omega -1}(\boldsymbol{\alpha }_{n}))^{2}}{%
2B_{n}^{\omega }}\times \left( 1+o_{P_{V^{\infty }|Z^{\infty }}}(1)\right)
~wpa1(P_{n,Z^{\infty }}).
\end{eqnarray*}

Following the proof of Result (1) of Theorem \ref{thm:bootstrap} step 3 with
$\mathbb{Z}_{n}^{\omega -1}(\boldsymbol{\alpha }_{n})$ replacing $\mathbb{Z}%
_{n}^{\omega -1}$, we obtain:
\begin{equation*}
\frac{\widehat{QLR}_{n}^{B}(\widehat{\phi }_{n})}{\sigma _{\omega }^{2}}%
=\left( \sqrt{n}\frac{\mathbb{Z}_{n}^{\omega -1}(\boldsymbol{\alpha }_{n})}{%
\sigma _{\omega }\sqrt{B_{n}^{\omega }}}\right) ^{2}\times \left(
1+o_{P_{V^{\infty }|Z^{\infty }}}(1)\right) =O_{P_{V^{\infty }|Z^{\infty
}}}(1),~wpa1(P_{n,Z^{\infty }}).
\end{equation*}%
This shows that, since for the bootstrap SQLR the \textquotedblleft null
hypothesis is $\phi (\alpha )=\widehat{\phi }_{n}\equiv \phi (\widehat{%
\alpha }_{n})$\textquotedblright , it always centers correctly.

By similar calculations to those in the proof of Result (2) of Theorem \ref%
{thm:bootstrap}, the law of $\left( \sqrt{n}\frac{\mathbb{Z}_{n}^{\omega -1}(%
\boldsymbol{\alpha }_{n})}{\sigma _{\omega }\sqrt{B_{n}^{\omega }}}\right)
^{2}$ is asymptotically (and $wpa1(P_{n,Z^{\infty }})$) equal to the law of $%
\left( \frac{\mathbb{Z}}{||u_{n}^{\ast }||}\right) ^{2}$ where $\mathbb{Z}%
\sim N(0,1)$. This implies that the $a$-th quantile of the distribution of $%
\frac{\widehat{QLR}_{n}^{B}(\widehat{\phi }_{n})}{\sigma _{\omega }^{2}}$, $%
\widehat{c}_{n}(a)$, is uniformly bounded $wpa1(P_{n,Z^{\infty }})$. Also,
following the proof of Result (2) of Theorem \ref{thm:bootstrap} we obtain:%
\begin{equation*}
\sup_{t\in \mathbb{R}}\left\vert P_{V^{\infty }|Z^{\infty }}\left( \frac{%
\widehat{QLR}_{n}^{B}(\widehat{\phi }_{n})}{\sigma _{\omega }^{2}}\leq t\mid
Z^{n}\right) -P_{Z^{\infty }}\left( \widehat{QLR}_{n}(\phi _{0})\leq t\mid
H_{0}\right) \right\vert =o_{P_{V^{\infty }|Z^{\infty
}}}(1)~wpa1(P_{n,Z^{\infty }}).
\end{equation*}

This and Theorem \ref{thm:chi2_localt} (and the fact that $||u_{n}^{\ast
}||\leq c<\infty $) immediately imply \textbf{Results (2)}. \textit{Q.E.D.}

\medskip

\noindent \textbf{Proof of Theorem \ref{thm:t-localt}}: The proof is
analogous to that of Theorems \ref{thm:VE} and \ref{thm:chi2_localt} so we
only present a sketch here.

Under our assumptions, Theorem \ref{thm:VE} still holds under the local
alternatives $\boldsymbol{\alpha }_{n}$. Observe that, with $\boldsymbol{%
\alpha }_{n}=\alpha _{0}+d_{n}\Delta _{n}\in \mathcal{N}_{osn}$ and $%
d_{n}=o(1)$,
\begin{align*}
\mathcal{T}_{n}& \equiv \sqrt{n}\frac{\phi (\widehat{\alpha }_{n})-\phi _{0}%
}{||\widehat{v}_{n}^{\ast }||_{n,sd}}=\sqrt{n}\frac{\phi (\widehat{\alpha }%
_{n})-\phi _{0}}{||v_{n}^{\ast }||_{sd}}\times (1+o_{P_{n,Z^{\infty }}}(1))
\\
& =\sqrt{n}\langle u_{n}^{\ast },\widehat{\alpha }_{n}-\alpha _{0}\rangle
\times (1+o_{P_{n,Z^{\infty }}}(1))+o_{P_{n,Z^{\infty }}}(1) \\
& =\left( -\sqrt{n}\mathbb{Z}_{n}(\boldsymbol{\alpha }_{n})+\sqrt{n}\frac{%
d_{n}\kappa \left( 1+o(1)\right) }{||v_{n}^{\ast }||_{sd}}\right) \times
(1+o_{P_{n,Z^{\infty }}}(1))+o_{P_{n,Z^{\infty }}}(1),
\end{align*}%
where the second line follows from assumption \ref{ass:phi}; the third line
follows from equation (\ref{eqn:LA-LAR}), and $\sqrt{n}\mathbb{Z}_{n}(%
\boldsymbol{\alpha }_{n})\Rightarrow N(0,1)$ under the local alternatives
(i.e., assumption \ref{ass:LAQ}(ii) under the alternatives).

For \textbf{Result (1)}, under local alternatives with $%
d_{n}=n^{-1/2}||v_{n}^{\ast }||_{sd}$ we have:
\begin{equation*}
\mathcal{T}_{n}=-\left( \sqrt{n}\mathbb{Z}_{n}(\boldsymbol{\alpha }%
_{n})-\kappa \left( 1+o(1)\right) \right) \times \left( 1+o_{P_{n,Z^{\infty
}}}(1)\right) +o_{P_{n,Z^{\infty }}}(1),\text{ and }\mathcal{W}_{n}\equiv
\left( \mathcal{T}_{n}\right) ^{2}\Rightarrow \chi _{1}^{2}(\kappa ^{2}).
\end{equation*}

For \textbf{Result (2), }under local alternatives with $\sqrt{n}\frac{d_{n}}{%
||v_{n}^{\ast }||_{sd}}\rightarrow \infty $ we have:
\begin{equation*}
\mathcal{W}_{n}\equiv \left( \mathcal{T}_{n}\right) ^{2}=\left(
O_{P_{n,Z^{\infty }}}(1)-\sqrt{n}\frac{d_{n}\kappa \left( 1+o(1)\right) }{%
||v_{n}^{\ast }||_{sd}}\right) ^{2}\times \left( 1+o_{P_{n,Z^{\infty
}}}(1)\right) +o_{P_{n,Z^{\infty }}}(1)\rightarrow \infty \text{ }%
~wpa1(P_{n,Z^{\infty }})\text{.}
\end{equation*}%
\textit{Q.E.D.}

\medskip

\noindent \textbf{Proof of Theorem \ref{thm:waldB_con}} For \textbf{Result
(1)}, following the proofs of Theorems \ref{thm:bootstrap_2}(1) and \ref%
{thm:BSQLR-loc-alt}, we have: under local alternatives $\boldsymbol{\alpha }%
_{n}$ defined in (\ref{eq:locaalt-1}), for $j=1,2,$%
\begin{equation*}
\widehat{W}_{j,n}^{B}=-\sqrt{n}\frac{\mathbb{Z}_{n}^{\omega -1}(\boldsymbol{%
\alpha }_{n})}{\sigma _{\omega }}+o_{P_{V^{\infty }|Z^{\infty
}}}(1)~wpa1(P_{n,Z^{\infty }}).
\end{equation*}%
By similar calculations to those in the proof of Theorem \ref%
{thm:bootstrap_2}(1), the law of $\sqrt{n}\frac{\mathbb{Z}_{n}^{\omega -1}(%
\boldsymbol{\alpha }_{n})}{\sigma _{\omega }}$ is asymptotically (and $%
wpa1(P_{n,Z^{\infty }})$) equal to the law of $\mathbb{Z}\sim N(0,1)$. Then
under the local alternatives $\boldsymbol{\alpha }_{n}$,
\begin{equation}
\sup_{t\in \mathbb{R}}\left\vert P_{V^{\infty }|Z^{\infty }}\left( \widehat{W%
}_{j,n}^{B}\leq t\mid Z^{n}\right) -P_{Z^{\infty }}\left( \widehat{W}%
_{n}\leq t\right) \right\vert =o_{P_{V^{\infty }|Z^{\infty
}}}(1)~wpa1(P_{n,Z^{\infty }}),  \label{boot-wald-consistent}
\end{equation}%
where $\lim_{n\rightarrow \infty }P_{Z^{\infty }}\left( \widehat{W}_{n}\leq
t\right) =\Phi (t)$ (i.e., the standard normal cdf). Thus the $a$-th
quantile of the distribution of $\left( \widehat{W}_{j,n}^{B}\right) ^{2}$, $%
\widehat{c}_{j,n}(a)$, is uniformly bounded $wpa1(P_{n,Z^{\infty }})$.

For \textbf{Result (2a), }by Theorem \ref{thm:t-localt}(2), Result (1)
(i.e., equation (\ref{boot-wald-consistent})) and the continuous mapping
theorem, we have:%
\begin{align*}
& P_{n,Z^{\infty }}\left( \mathcal{W}_{n}\geq \widehat{c}_{j,n}(1-\tau
)\right) -P_{V^{\infty }|Z^{\infty }}\left( \left( \widehat{W}%
_{j,n}^{B}\right) ^{2}\geq \widehat{c}_{j,n}(1-\tau )\mid Z^{n}\right) \\
& =\Pr \left( \chi _{1}^{2}(\kappa ^{2})\geq \widehat{c}_{j,n}(1-\tau
)\right) -\Pr \left( \chi _{1}^{2}\geq \widehat{c}_{j,n}(1-\tau )\right)
+o_{P_{V^{\infty }|Z^{\infty }}}(1)~wpa1(P_{n,Z^{\infty }}).
\end{align*}%
Thus by the definition of $\widehat{c}_{j,n}(1-\tau )$ we obtain:
\begin{eqnarray*}
P_{n,Z^{\infty }}\left( \mathcal{W}_{n}\geq \widehat{c}_{j,n}(1-\tau
)\right) &=& \tau +\Pr \left( \chi _{1}^{2}(\kappa ^{2})\geq \widehat{c}%
_{j,n}(1-\tau )\right) -\Pr \left( \chi _{1}^{2}\geq \widehat{c}%
_{j,n}(1-\tau )\right) \\
&& +o_{P_{V^{\infty }|Z^{\infty
}}}(1)~wpa1(P_{n,Z^{\infty }}).
\end{eqnarray*}

\textbf{Result (2b) }directly follows from Theorem \ref{thm:t-localt}(2),
equation (\ref{boot-wald-consistent}) and the continuous mapping theorem.
\textit{Q.E.D.}

\subsection{Proofs for Section \protect\ref{subsec-A6} on asymptotic theory
under increasing dimension of $\protect\phi $}

\begin{lemma}
\label{lem:D_BDD} Let Assumption \ref{ass:sieve}(iv) hold. Then: there exist
positive finite constants $c,C$ such that
\begin{equation*}
c^{2}I_{d(n)}\leq \mathbb{D}_{n}^{2}\leq C^{2}I_{d(n)},
\end{equation*}%
where $I_{d(n)}$ is the $d(n)\times d(n)$ identity and for matrices $A\leq B$
means that $B-A$ is positive semi-definite.
\end{lemma}

\noindent \textbf{Proof of Lemma \ref{lem:D_BDD}}. By Assumption \ref%
{ass:sieve}(iv), the eigenvalues of $\Sigma _{0}(x)$ and $\Sigma (x)$ are
bounded away from zero and infinity uniformly in $x$. Therefore, for any
matrix $A$,
\begin{equation*}
A^{\prime }\Sigma ^{-1}(x)\Sigma _{0}(x)\Sigma ^{-1}(x)A\geq dA^{\prime
}\Sigma ^{-1}(x)A
\end{equation*}%
and
\begin{equation*}
A^{\prime }\Sigma ^{-1}(x)\Sigma _{0}(x)\Sigma ^{-1}(x)A\leq DA^{\prime
}\Sigma ^{-1}(x)A
\end{equation*}%
for some finite constant $0<d\leq D<\infty $, and for all $x$. Taking expectations at
both sides and choosing $A^{\prime }\equiv \frac{dm(x,\alpha _{0})}{d\alpha }%
[\mathbf{v}_{n}^{\ast }]^{\prime }$, these displays imply that
\begin{equation*}
\Omega _{sd,n}\geq d\Omega _{n}~and~\Omega _{sd,n}\leq D\Omega _{n}
\end{equation*}%
Thus
\begin{equation*}
\mathbb{D}_{n}^{2}=\Omega _{sd,n}^{1/2}\Omega _{n}^{-1}\Omega _{sd,n}\Omega
_{n}^{-1}\Omega _{sd,n}^{1/2}\geq d\{\Omega _{sd,n}^{1/2}\Omega
_{n}^{-1}\Omega _{sd,n}^{1/2}\}\geq d^{2}\Omega _{sd,n}^{1/2}\Omega
_{sd,n}^{-1}\Omega _{sd,n}^{1/2}=d^{2}I_{d(n)}.
\end{equation*}%
Similarly, $\mathbb{D}_{n}^{2}\leq D^{2}I_{d(n)}$. \textit{Q.E.D.}

\medskip

\begin{lemma}
\label{lem:Z_bdd} Let $\mathcal{T}_{n}^{M}\equiv \{t\in \mathbb{R}%
^{d(n)}:||t||_{e}\leq M_{n}n^{-1/2}\sqrt{d(n)}\}$. Then:
\begin{equation*}
||\Omega _{sd,n}^{-1/2}\mathbf{Z}_{n}||_{e}=O_{P}\left( n^{-1/2}\sqrt{d(n)}%
\right) \quad \text{and}\quad \Omega _{sd,n}^{-1/2}\mathbf{Z}_{n}\in
\mathcal{T}_{n}^{M}\text{ wpa1.}
\end{equation*}
\end{lemma}

\noindent \textbf{Proof of Lemma \ref{lem:Z_bdd}. } Let $\Omega
_{sd,n}^{-1/2}\mathbf{Z}_{n}\equiv n^{-1}\sum_{i=1}^{n}\zeta _{in}$ where $%
\zeta _{in}\in \mathbb{R}^{d(n)}$. Observe that $E[\zeta _{in}\zeta
_{in}^{\prime }]=I_{d(n)}$. It follows that
\begin{equation*}
E_{P}[(\Omega _{sd,n}^{-1/2}\mathbf{Z}_{n})^{\prime }(\Omega _{sd,n}^{-1/2}%
\mathbf{Z}_{n})]=tr\left\{ E_{P}[\Omega _{sd,n}^{-1/2}\mathbf{Z}_{n}\mathbf{Z%
}_{n}^{\prime }\Omega _{sd,n}^{-1/2}]\right\} =n^{-2}\sum_{i=1}^{n}tr\left\{
E_{P}[\zeta _{in}\zeta _{in}^{\prime }]\right\} =n^{-1}d(n),
\end{equation*}%
and thus the desired result follows by the Markov inequality. \textit{Q.E.D.}

\medskip

\begin{lemma}
\label{lem:1order_vv} Let Conditions for Lemma \ref{thm:ThmCONVRATEGRAL} and
Assumption \ref{ass:LAQ_vv} hold. Denote $\tilde{\gamma}_{n}\equiv \sqrt{%
s_{n}}(1+b_{n})+a_{n}$. Then:

(1) $\left\Vert \Omega _{sd,n}^{-1/2}\{\mathbf{Z}_{n}+\langle \mathbf{v}%
_{n}^{\ast \prime },\widehat{\alpha }_{n}-\alpha _{0}\rangle \}\right\Vert
_{e}=O_{P}(\sqrt{d(n)}\tilde{\gamma}_{n})=o_{P}(n^{-1/2})$;

(2) further let Assumption \ref{ass:phivv1} hold. Then
\begin{equation*}
\left\Vert \Omega _{sd,n}^{-1/2}\{\mathbf{Z}_{n}+\phi (\widehat{\alpha }%
_{n})-\phi (\alpha _{0})\}\right\Vert _{e}=o_{P}(n^{-1/2}).
\end{equation*}
\end{lemma}

\noindent \textbf{Proof of Lemma \ref{lem:1order_vv}}: For \textbf{Result
(1), }note that $||t||_{e}^{2}=\sum_{l=1}^{d(n)}|t_{l}|^{2}$ and if we
obtain $|t_{l}|=O_{P}(\tilde{\gamma}_{n})$ for $\tilde{\gamma}_{n}$
uniformly over $l$, then $||t||_{e}^{2}=O_{P}(d(n)\tilde{\gamma}_{n}^{2})$.

The rest of the proof follows closely the proof of Theorem \ref%
{thm:theta_anorm} so we only present the main steps. By definition of the
approximate PSMD estimator $\widehat{\alpha }_{n}$, and Assumption \ref%
{ass:LAQ_vv}(i),
\begin{equation*}
0\leq t^{\prime }\Omega _{sd,n}^{-1/2}\left( \mathbf{Z}_{n}+\langle \mathbf{v%
}_{n}^{\ast \prime },\widehat{\alpha }_{n}-\alpha _{0}\rangle \right) +\frac{%
1}{2}t^{\prime }\mathbb{B}_{n}t+O_{P}(r^{-1}(t)).
\end{equation*}%
We now choose $t=\sqrt{s_{n}}e$ where $e\in
\{(1,0,...,0),(0,1,0,...,0),...,(0,...,1)\}$, it is easy to see that this $%
t\in \mathcal{T}_{n}^{M}$, and thus the display above implies
\begin{equation*}
0\leq e^{\prime }\Omega _{sd,n}^{-1/2}\left( \mathbf{Z}_{n}+\langle \mathbf{v%
}_{n}^{\ast \prime },\widehat{\alpha }_{n}-\alpha _{0}\rangle \right) +O_{P}(%
\tilde{\gamma}_{n}).
\end{equation*}%
By changing the sign of $t$, it follows that
\begin{equation*}
\left\vert e^{\prime }\Omega _{sd,n}^{-1/2}\left( \mathbf{Z}_{n}+\langle
\mathbf{v}_{n}^{\ast \prime },\widehat{\alpha }_{n}-\alpha _{0}\rangle
\right) \right\vert =O_{P}(\tilde{\gamma}_{n}).
\end{equation*}%
Observe that the RHS holds uniformly over $e$, thus, since $e\in
\{(1,0,...,0),(0,1,0,...,0),...,(0,...,1)\}$, it follows that
\begin{equation*}
\left\Vert \Omega _{sd,n}^{-1/2}\left( \mathbf{Z}_{n}+\langle \mathbf{v}%
_{n}^{\ast \prime },\widehat{\alpha }_{n}-\alpha _{0}\rangle \right)
\right\Vert _{e}=O_{P}(\sqrt{d(n)}\tilde{\gamma}_{n})=o_{P}(n^{-1/2}),
\end{equation*}%
where the second equal sign is due to Assumption \ref{ass:LAQ_vv}(ii).

For \textbf{Result (2).} In view of Result (1), it suffices to show that
\begin{equation*}
\left\Vert \Omega _{sd,n}^{-1/2}\{\phi (\widehat{\alpha }_{n})-\phi (\alpha
_{0})-\langle \mathbf{v}_{n}^{\ast \prime },\widehat{\alpha }_{n}-\alpha
_{0}\rangle \}\right\Vert _{e}=o_{P}(n^{-1/2}).
\end{equation*}%
Following the proof of Theorem \ref{thm:theta_anorm} we have:%
\begin{equation*}
\langle \mathbf{v}_{n}^{\ast \prime },\widehat{\alpha }_{n}-\alpha
_{0}\rangle =\frac{d\phi (\alpha _{0})}{d\alpha }[\widehat{\alpha }%
_{n}-\alpha _{0,n}]=\frac{d\phi (\alpha _{0})}{d\alpha }[\widehat{\alpha }%
_{n}-\alpha _{0}]-\frac{d\phi (\alpha _{0})}{d\alpha }[\alpha _{0,n}-\alpha
_{0}].
\end{equation*}%
Since Assumption \ref{ass:phivv1}(ii)(iii) (with $t=0$) implies that
\begin{equation*}
\left\Vert \Omega _{sd,n}^{-1/2}\{\phi (\widehat{\alpha }_{n})-\phi (\alpha
_{0})-\frac{d\phi (\alpha _{0})}{d\alpha }[\widehat{\alpha }_{n}-\alpha
_{0}]+\frac{d\phi (\alpha _{0})}{d\alpha }[\alpha _{0,n}-\alpha
_{0}]\}\right\Vert _{e}=O_{P}(c_{n}),
\end{equation*}%
the desired result now follows from Assumption \ref{ass:phivv1}(iv) of $%
c_{n}=o(n^{-1/2})$. \textit{Q.E.D.}

\medskip

\noindent \textbf{Proof of Theorem \ref{thm:wald_vv}.} Throughout the proof
let $\hat{W}_{n}\equiv n(\phi (\widehat{\alpha }_{n})-\phi (\alpha
_{0}))^{\prime }\Omega _{sd,n}^{-1}(\phi (\widehat{\alpha }_{n})-\phi
(\alpha _{0}))$. By Lemma \ref{lem:1order_vv}(2),
\begin{equation*}
T_{n}\equiv (\phi (\widehat{\alpha }_{n})-\phi (\alpha _{0})+\mathbf{Z}%
_{n})^{\prime }\Omega _{sd,n}^{-1}(\phi (\widehat{\alpha }_{n})-\phi (\alpha
_{0})+\mathbf{Z}_{n})=o_{P}(n^{-1}).
\end{equation*}

Observe that
\begin{align*}
& |(\phi (\widehat{\alpha }_{n})-\phi (\alpha _{0}))^{\prime }\Omega
_{sd,n}^{-1}(\phi (\widehat{\alpha }_{n})-\phi (\alpha _{0}))-(\mathbf{Z}%
_{n})^{\prime }\Omega _{sd,n}^{-1}(\mathbf{Z}_{n})| \\
& \leq T_{n}+2||(\phi (\widehat{\alpha }_{n})-\phi (\alpha _{0})+\mathbf{Z}%
_{n})^{\prime }\Omega _{sd,n}^{-1/2}||_{e}\times ||\Omega _{sd,n}^{-1/2}%
\mathbf{Z}_{n}||_{e} \\
& =o_{P}(n^{-1})+2||(\phi (\widehat{\alpha }_{n})-\phi (\alpha _{0})+\mathbf{%
Z}_{n})^{\prime }\Omega _{sd,n}^{-1/2}||_{e}\times ||\Omega _{sd,n}^{-1/2}%
\mathbf{Z}_{n}||_{e} =o_{P}(n^{-1})+o_{P}(n^{-1}\sqrt{d(n)})
\end{align*}%
where the last equality is due to Lemmas \ref{lem:Z_bdd} and \ref%
{lem:1order_vv}(2). Therefore we obtain \textbf{Result (1)}:
\begin{equation*}
\hat{W}_{n}=(\sqrt{n}\mathbf{Z}_{n})^{\prime }\Omega _{sd,n}^{-1}(\sqrt{n}%
\mathbf{Z}_{n})+o_{P}(\sqrt{d(n)})\equiv \mathbf{W}_{n}+o_{P}(\sqrt{d(n)}).
\end{equation*}

\textbf{Result (2)} follows directly from Result (1) when $d(n)=d$ is fixed
and finite.

\textbf{Result (3)} follows from Result (1) and the following property:
\begin{equation*}
\Xi _{n}\equiv (2d(n))^{-1/2}(\mathbf{W}_{n}-d(n))\Rightarrow N(0,1)
\end{equation*}%
where $\mathbf{W}_{n}\equiv (\sqrt{n}\mathbf{Z}_{n})^{\prime }\Omega
_{sd,n}^{-1}(\sqrt{n}\mathbf{Z}_{n})$, or formally,
\begin{equation*}
\sup_{f\in BL_{1}(\mathbb{R})}|E[f(\Xi _{n})]-E[f(\mathbb{Z})]|=o(1)
\end{equation*}%
where $\mathbb{Z}\sim N(0,1)$ and $BL_{1}(\mathbb{R})$ is the space of
bounded (by 1) Lipschitz functions from $\mathbb{R}$ to $\mathbb{R}$.

By triangle inequality it suffices to show that
\begin{equation}
\sup_{f\in BL_{1}(\mathbb{R})}|E[f(\Xi _{n})]-E[f(\xi _{n})]|=o(1)
\label{eqn:SQLR_vv_1}
\end{equation}%
and
\begin{equation}
\sup_{f\in BL_{1}(\mathbb{R})}|E[f(\xi _{n})]-E[f(\mathbb{Z})]|=o(1)
\label{eqn:SQLR_vv_2}
\end{equation}%
where $\xi _{n}\equiv (2d(n))^{-1/2}(\sum_{j=1}^{d(n)}\mathbb{Z}%
_{j}^{2}-d(n))$ with $\mathbb{Z}_{j}\sim N(0,1)$ and independent across $%
j=1,...,d(n)$. We now show that both equations hold.\newline

\textbf{Equation}  \ref{eqn:SQLR_vv_1}. Let $t\mapsto \nu _{M}(t)\equiv \min
\{t^{\prime }t,M\}$ for some $M>0$. Observe
\begin{eqnarray*}
& &\left\vert E[f(\Xi _{n})]-E\left[ f\left( (2d(n))^{-1/2}(\nu _{M}(\Omega
_{sd,n}^{-1/2}\sqrt{n}\mathbf{Z}_{n})-d(n))\right) \right] \right\vert \\
&=& \left\vert E\left[ f\left( (2d(n))^{-1/2}(\nu _{\infty }(\Omega
_{sd,n}^{-1/2}\sqrt{n}\mathbf{Z}_{n})-d(n))\right) -f\left(
(2d(n))^{-1/2}(\nu _{M}(\Omega _{sd,n}^{-1/2}\sqrt{n}\mathbf{Z}%
_{n})-d(n))\right) \right] \right\vert \\
&=& \left\vert \int_{\{z:nz^{\prime }\Omega _{sd,n}^{-1}z>M\}}\left[ f\left(
\frac{\nu _{\infty }(\Omega _{sd,n}^{-1/2}\sqrt{n}z_{n})-d(n)}{\sqrt{2d(n)}}%
\right) -f\left( \frac{M-d(n)}{\sqrt{2d(n)}}\right) \right] P_{Z^{\infty
}}(dz)\right\vert \\
&\leq & 2P_{Z^{\infty }}\left( (\sqrt{n}\mathbf{Z}_{n})^{\prime }\Omega
_{sd,n}^{-1}(\sqrt{n}\mathbf{Z}_{n})>M\right)
\end{eqnarray*}%
where the last line follows from the fact that $f$ is bounded by 1.
Therefore, by the Markov inequality, for any $\epsilon $, there exists a $M$
such that
\begin{eqnarray*}
	\left\vert E[f(\Xi _{n})]-E\left[ f\left( (2d(n))^{-1/2}(\nu
	_{M}(\Omega _{sd,n}^{-1/2}\sqrt{n}\mathbf{Z}_{n})-d(n))\right) \right]
	\right\vert <\epsilon
\end{eqnarray*}
for sufficiently large $n$. A similar result holds
if we replace $\Omega _{sd,n}^{-1/2}\sqrt{n}\mathbf{Z}_{n}$ by $\mathcal{Z}%
_{n}=(\mathbb{Z}_{1},...,\mathbb{Z}_{d(n)})^{\prime }$ with $\mathbb{Z}%
_{j}\sim N(0,1)$ and independent across $j=1,...,d(n)$. Therefore, in order
to show equation \ref{eqn:SQLR_vv_1}, it suffices to show
\begin{equation*}
\sup_{f\in BL_{1}(\mathbb{R})}|E\left[ f\left( \Xi _{M,n}\right) \right]
-E[f(\xi _{M,n})]|=o(1)
\end{equation*}%
where $\Xi _{M,n}\equiv (2d(n))^{-1/2}(\nu _{M}(\Omega _{sd,n}^{-1/2}\sqrt{n}%
\mathbf{Z}_{n})-d(n))$ and $\xi _{M,n}\equiv (2d(n))^{-1/2}(\nu _{M}(%
\mathcal{Z}_{n})-d(n))$.

Since $f$ is uniformly bounded and continuous, it is clear that in order to
show the previous display, it suffices to show that
\begin{align}  \label{eqn:SQLR_vv_3}
(2d(n))^{-1/2}|\nu_{M}(\Omega^{-1/2}_{sd,n}\sqrt{n}\mathbf{Z}_{n}) - \nu_{M}(%
\mathcal{Z}_{n})| = o_{P}(1).
\end{align}

It turns out that $|\nu _{M}(t)-\nu _{M}(r)|\leq 2\sqrt{M}||t-r||_{e}$, so $%
t\mapsto \nu _{M}(t)$ is Lipschitz (and uniformly bounded). So in order to
show equation \ref{eqn:SQLR_vv_3} it is sufficient to show that for any $%
\delta >0$, there exists a $N(\delta )$ such that
\begin{equation*}
\Pr \left( (2d(n))^{-1/2}||\Omega _{sd,n}^{-1/2}\sqrt{n}\mathbf{Z}_{n}-%
\mathcal{Z}_{n}||_{e}>\delta \right) <\delta
\end{equation*}%
for all $n\geq N(\delta )$. Note that $\Omega _{sd,n}^{-1/2}\sqrt{n}\mathbf{Z%
}_{n}=\frac{1}{\sqrt{n}}\sum_{i=1}^{n}\Psi _{n}(Z_{i})$, with $\Psi
_{n}(z)\equiv \left( \frac{dm(x,\alpha _{0})}{d\alpha }[\mathbf{v}_{n}^{\ast
}]\Omega _{sd,n}^{-1/2}\right) ^{\prime }\rho (z,\alpha _{0})$, and that $%
\mathcal{Z}_{n}$ can be cast as $\frac{1}{\sqrt{n}}\sum_{i=1}^{n}\mathcal{Z}%
_{n,i}$ with $\mathcal{Z}_{n,i}\sim N(0,I_{d(n)})$, iid across $i=1,...,n$.
Following the arguments in Section 10.4 of \cite{Pollard_90}, we obtain: for
any $\delta >0$,
\begin{equation*}
\Pr \left( \left\Vert \sqrt{n}\mathbf{Z}_{n}^{\prime }\Omega _{sd,n}^{-1/2}-%
\frac{1}{\sqrt{n}}\sum_{i=1}^{n}\mathcal{Z}_{n,i}\right\Vert _{e}>3\delta
\right) \leq \mathcal{Y}_{d(n)}\left( \frac{\mu _{3,n}nd(n)^{5/2}}{(\delta
\sqrt{n})^{3}}\right) ,
\end{equation*}%
for any $n$, where $x\mapsto \mathcal{Y}_{d(n)}(x)\equiv Cx\times (1+|\log
(1/x)|/d(n))$ and $\mu _{3,n}\equiv E\left[ \left\Vert \left( \frac{%
dm(X,\alpha _{0})}{d\alpha }[\mathbf{v}_{n}^{\ast }]\Omega
_{sd,n}^{-1/2}\right) ^{\prime }\rho (Z,\alpha _{0})\right\Vert _{e}^{3}%
\right] $. Therefore,%
\begin{equation*}
\Pr \left( (2d(n))^{-1/2}||\Omega _{sd,n}^{-1/2}\sqrt{n}\mathbf{Z}_{n}-%
\mathcal{Z}_{n}||_{e}>\delta \right) \leq \mathcal{Y}_{d(n)}\left( \frac{\mu
_{3,n}nd(n)^{5/2}}{(\delta /3)^{3}d(n)^{3/2}n^{3/2}}\right) =\mathcal{Y}%
_{d(n)}\left( n^{-1/2}d(n)\frac{\mu _{3,n}}{(\delta /3)^{3}8}\right)
\rightarrow 0
\end{equation*}%
provided that $d(n)=o(\sqrt{n}\mu _{3,n}^{-1})$ which is assumed in the
Theorem Result (3).\newline

\textbf{Equation \ref{eqn:SQLR_vv_2}.} Observe that $\xi _{n}\equiv
(2d(n))^{-1/2}(\sum_{j=1}^{d(n)}\mathbb{Z}_{j}^{2}-d(n))$ with $\mathbb{Z}%
_{j}\sim N(0,1)$ i.i.d. across $j=1,...,d(n)$, $E[(\mathbb{Z}_{l}^{2}-1)]=0$
and $E[(\mathbb{Z}_{l}^{2}-1)^{2}]=2$. Thus, $\xi _{n}\Rightarrow N(0,1)$ by
a standard CLT. \textit{Q.E.D.}

\medskip

In the following we recall that $\alpha (t)\equiv \alpha +\mathbf{v}%
_{n}^{\ast }(\Omega _{sd,n})^{-1/2}t$ for $t\in \mathbb{R}^{d(n)}$.

\begin{lemma}
\label{lem:t_ex} Let all conditions for Theorem \ref{thm:chi2_vv}(1) hold.
Then there exists a $t_{n}$ (possibly random) such that: (1) $t_{n}\in
\mathcal{T}_{n}^{M}$ wpa1, (2) $\widehat{\alpha }_{n}(t_{n})\in \mathcal{A}%
_{k(n)}^{R}=\{\alpha \in \mathcal{A}_{k(n)}:\phi (\alpha )=\phi _{0}\}$
wpa1, and (3)
\begin{equation*}
0\leq n\{\widehat{Q}_{n}(\widehat{\alpha }_{n}(t_{n}))-\widehat{Q}_{n}(%
\widehat{\alpha }_{n})\}\leq (\sqrt{n}\Omega _{sd,n}^{-1/2}\mathbf{Z}%
_{n})^{\prime }\mathbb{D}_{n}(\sqrt{n}\Omega _{sd,n}^{-1/2}\mathbf{Z}%
_{n})+o_{P}(\sqrt{d(n)}).
\end{equation*}
\end{lemma}

\noindent \textbf{Proof of Lemma \ref{lem:t_ex}}\textsc{:} To show \textbf{%
Parts (1) and (2),} we define the following mappings:
\begin{equation*}
t\in \mathbb{R}^{d(n)}\mapsto \varphi _{n}(t)\equiv \Omega
_{sd,n}^{-1/2}\left\{ \phi \left( \widehat{\alpha }_{n}(t)\right) -\phi
(\alpha _{0})-\frac{d\phi (\alpha _{0})}{d\alpha }[\widehat{\alpha }%
_{n}(t)-\alpha _{0}]\right\}
\end{equation*}%
and $t\in \mathbb{R}^{d(n)}\mapsto \tau _{n}(t)\equiv -\mathbb{D}_{n}\Omega
_{sd,n}^{-1/2}\left\{ \langle \mathbf{v}_{n}^{\ast \prime },\widehat{\alpha }%
_{n}-\alpha _{0}\rangle +\frac{d\phi (\alpha _{0})}{d\alpha }[\alpha
_{0n}-\alpha _{0}]+\Omega _{sd,n}^{1/2}t\right\} $. Under our assumptions,
both mappings are continuous in $t$ (a.s.) and thus $\Phi _{n}\equiv \varphi
_{n}\circ \tau _{n}$ is also continuous in $t$ (a.s.). Given $%
c_{n}=o(n^{-1/2})$ satisfying assumption \ref{ass:phivv1}(iv), we define $%
\mathbb{T}_{n}\equiv \{t\in \mathbb{R}^{d(n)}:||t||_{e}\leq L_{n}c_{n}\}$
where $(L_{n})_{n}$ is a positive real valued sequence diverging to infinity
slowly such that $L_{n}c_{n}=o(n^{-1/2})$ (such a sequence exists by
assumption \ref{ass:phivv1}(iv)).

Let $S_{n}\equiv \{Z^{n}:\sup_{t\in \mathbb{T}_{n}}||\Phi _{n}(t)||_{e}\leq
L_{n}c_{n}\}$. By Lemmas \ref{lem:D_BDD}, \ref{lem:Z_bdd} and \ref%
{lem:1order_vv}(2), and assumption \ref{ass:phivv1}(iii), we have that for
any $t\in \mathbb{T}_{n}$,
\begin{equation*}
||\tau _{n}(t)||_{e}\leq O_{P}(\sqrt{d(n)}\{\tilde{\gamma}%
_{n}+n^{-1/2}\})+O(c_{n})+||\mathbb{D}_{n}t||_{e}=O_{P}(n^{-1/2}\sqrt{d(n)}%
)+O(L_{n}c_{n})
\end{equation*}%
where $\tilde{\gamma}_{n}\equiv \sqrt{|s_{n}|}(1+b_{n})+a_{n}=o(n^{-1/2})$
(by assumption \ref{ass:LAQ_vv}(ii)). Hence $\tau _{n}(t)\in \mathcal{T}%
_{n}^{M}$ for all $t\in \mathbb{T}_{n}$. This implies, by assumption \ref%
{ass:phivv1}(i)(ii), $P(S_{n})\rightarrow 1$.

Moreover, these results imply that $||\Phi _{n}(t)||_{e}\leq L_{n}c_{n}$ for
all $t\in \mathbb{T}_{n}$ and $Z^{n}\in S_{n}$. This implies that $\{\Phi
_{n}(t):t\in \mathbb{T}_{n}\}\subseteq \mathbb{T}_{n}$ wpa1.

For any given $n$, $\mathbb{T}_{n}$ is compact and convex in $\mathbb{R}%
^{d(n)}$ and since $\Phi _{n}$ is continuous and maps $\mathbb{T}_{n}$ into
itself (wpa1), by Brouwer fixed point theorem, wpa1 there exists a $\hat{t}%
_{n}\in \mathbb{T}_{n}$ such that $\Phi _{n}(\hat{t}_{n})=\hat{t}_{n}$.
Therefore,
\begin{equation*}
\hat{t}_{n}=\varphi _{n}\circ \tau _{n}(\hat{t}_{n})=\Omega
_{sd,n}^{-1/2}\left\{ \phi \left( \widehat{\alpha }_{n}(\tau _{n}(\hat{t}%
_{n}))\right) -\phi (\alpha _{0})-\frac{d\phi (\alpha _{0})}{d\alpha }[%
\widehat{\alpha }_{n}(\tau _{n}(\hat{t}_{n}))-\alpha _{0}]\right\} .
\end{equation*}%
Since\begin{eqnarray*}
	\widehat{\alpha }_{n}(\tau _{n}(\hat{t}_{n})) &=& \widehat{\alpha }_{n}+%
	\mathbf{v}_{n}^{\ast }\Omega _{sd,n}^{-1/2}\tau _{n}(\hat{t}_{n})\\
	&= &\widehat{%
		\alpha }_{n}-\mathbf{v}_{n}^{\ast }\Omega _{sd,n}^{-1/2}\mathbb{D}_{n}\Omega
	_{sd,n}^{-1/2}\left( \langle \mathbf{v}_{n}^{\ast \prime },\widehat{\alpha }%
	_{n}-\alpha _{0}\rangle +\frac{d\phi (\alpha _{0})}{d\alpha }[\alpha
	_{0n}-\alpha _{0}]+\Omega _{sd,n}^{1/2}\hat{t}_{n}\right)
\end{eqnarray*}
and $\Omega
_{sd,n}^{-1/2}\mathbb{D}_{n}\Omega _{sd,n}^{-1/2}=\Omega _{n}^{-1}$, we
obtain%
\begin{eqnarray*}
\frac{d\phi (\alpha _{0})}{d\alpha }[\widehat{\alpha }_{n}(\tau _{n}(\hat{t}%
_{n}))-\alpha _{0}] &= &\frac{d\phi (\alpha _{0})}{d\alpha }[\widehat{\alpha }%
_{n}-\alpha _{0}] \\
&& -\langle \mathbf{v}_{n}^{\ast \prime },\mathbf{v}_{n}^{\ast
}\rangle \Omega _{n}^{-1}\left( \langle \mathbf{v}_{n}^{\ast \prime },%
\widehat{\alpha }_{n}-\alpha _{0}\rangle +\frac{d\phi (\alpha _{0})}{d\alpha
}[\alpha _{0n}-\alpha _{0}]+\Omega _{sd,n}^{1/2}\hat{t}_{n}\right) .
\end{eqnarray*}%
Since $\langle \mathbf{v}_{n}^{\ast \prime },\mathbf{v}_{n}^{\ast }\rangle
=\Omega _{n}$ we obtain
\begin{eqnarray*}
\hat{t}_{n} &=&\Omega _{sd,n}^{-1/2}\left\{ \phi \left( \widehat{\alpha }%
_{n}(\tau _{n}(\hat{t}_{n}))\right) -\phi (\alpha _{0})-\frac{d\phi (\alpha
_{0})}{d\alpha }[\widehat{\alpha }_{n}-\alpha _{0}]+\langle \mathbf{v}%
_{n}^{\ast \prime },\widehat{\alpha }_{n}-\alpha _{0}\rangle +\frac{d\phi
(\alpha _{0})}{d\alpha }[\alpha _{0n}-\alpha _{0}]+\Omega _{sd,n}^{1/2}\hat{t%
}_{n}\right\} \\
&=&\Omega _{sd,n}^{-1/2}\left\{ \phi \left( \widehat{\alpha }_{n}(\tau _{n}(%
\hat{t}_{n}))\right) -\phi (\alpha _{0})\right\} +\hat{t}_{n}\text{.}
\end{eqnarray*}%
Thus $\Omega _{sd,n}^{-1/2}\left\{ \phi \left( \widehat{\alpha }_{n}(\tau
_{n}(\hat{t}_{n}))\right) -\phi (\alpha _{0})\right\} =0$ wpa1 iff $\phi
\left( \widehat{\alpha }_{n}(\tau _{n}(\hat{t}_{n}))\right) -\phi (\alpha
_{0})=0$ wpa1. Also, since $\tau _{n}(\hat{t}_{n})\in \mathcal{T}_{n}^{M}$
wpa1, Parts (1) and (2) hold with $t_{n}\equiv \tau _{n}(\hat{t}_{n})$.%
\newline

To show \textbf{Part (3),} recall that $\widehat{\alpha }_{n}\in \mathcal{N}%
_{osn}$ wpa1 and $\widehat{\alpha }_{n}(t_{n})\in \mathcal{A}_{k(n)}^{R}$ by
Parts (1) and (2) with $t_{n}\equiv \tau _{n}(\hat{t}_{n})$. We can rewrite $%
t_{n}$ as
\begin{equation*}
t_{n}\equiv -\mathbb{D}_{n}\Omega _{sd,n}^{-1/2}\left\{ \langle \mathbf{v}%
_{n}^{\ast \prime },\widehat{\alpha }_{n}-\alpha _{0}\rangle +A_{n}(\hat{t}%
_{n})\right\} \text{\quad with }A_{n}(t)\equiv \frac{d\phi (\alpha _{0})}{%
d\alpha }[\alpha _{0n}-\alpha _{0}]+\Omega _{sd,n}^{1/2}t\text{.}
\end{equation*}

Observe that $||t_{n}||_{e}=O_{P}(\sqrt{d(n)}n^{-1/2})$, so by Assumption %
\ref{ass:LAQ_vv}(i) and the definition of $\widehat{\alpha }_{n}$,
\begin{align*}
& 0\leq n\left[ \widehat{Q}_{n}(\widehat{\alpha }_{n}(t_{n}))-\widehat{Q}%
_{n}(\widehat{\alpha }_{n})\right] \\
=& n(t_{n})^{\prime }\Omega _{sd,n}^{-1/2}\{\langle \mathbf{v}_{n}^{\ast
\prime },\widehat{\alpha }_{n}-\alpha _{0}\rangle +\mathbf{Z}%
_{n}\}+0.5n\{t_{n}^{\prime }\mathbb{B}_{n}t_{n}\}+n\times
O_{P}(s_{n}+||t_{n}||_{e}a_{n}+||t_{n}||_{e}^{2}b_{n}) \\
\leq & n(t_{n})^{\prime }\Omega _{sd,n}^{-1/2}\{\langle \mathbf{v}_{n}^{\ast
\prime },\widehat{\alpha }_{n}-\alpha _{0}\rangle +\mathbf{Z}%
_{n}\}+0.5n\{t_{n}^{\prime }\mathbb{D}_{n}^{-1}t_{n}\}+n\times
O_{P}(s_{n}+||t_{n}||_{e}a_{n}+||t_{n}||_{e}^{2}b_{n})
\end{align*}%
where the third line follows from the fact that $\sup_{t:||t||_{e}=1}|t^{%
\prime }\{\mathbb{B}_{n}-\mathbb{D}_{n}^{-1}\}t|=O_{P}(b_{n})$ by
assumption, and thus we have: $t^{\prime }\mathbb{B}_{n}t\leq |t^{\prime }\{%
\mathbb{B}_{n}-\mathbb{D}_{n}^{-1}\}t|+t^{\prime }\mathbb{D}_{n}^{-1}t\leq
||t||_{e}^{2}O_{P}(b_{n})+t^{\prime }\mathbb{D}_{n}^{-1}t$ uniformly over $%
t\in \mathbb{R}^{d(n)}$ with $||t||_{e}=1$.

By the fact that $\Omega _{sd,n}^{-1/2}\mathbb{D}_{n}\mathbb{D}_{n}^{-1}%
\mathbb{D}_{n}\Omega _{sd,n}^{-1/2}=\Omega _{sd,n}^{-1/2}\mathbb{D}%
_{n}\Omega _{sd,n}^{-1/2}=\Omega _{n}^{-1}$, the definition of $t_{n}$ and
straightforward algebra, the previous display implies that
\begin{eqnarray*}
&  &n\left[ \widehat{Q}_{n}(\widehat{\alpha }_{n}(t_{n}))-\widehat{Q}_{n}(%
\widehat{\alpha }_{n})\right] \\
&\leq & -0.5n(\langle \mathbf{v}_{n}^{\ast \prime },\widehat{\alpha }%
_{n}-\alpha _{0}\rangle )^{\prime }\Omega _{n}^{-1}(\langle \mathbf{v}%
_{n}^{\ast \prime },\widehat{\alpha }_{n}-\alpha _{0}\rangle )-n(\langle
\mathbf{v}_{n}^{\ast \prime },\widehat{\alpha }_{n}-\alpha _{0}\rangle
)^{\prime }\Omega _{n}^{-1}(\mathbf{Z}_{n})-n(A_{n}(\hat{t}_{n}))^{\prime
}\Omega _{n}^{-1}\mathbf{Z}_{n} \\
&& +0.5n(A_{n}(\hat{t}_{n}))^{\prime }\Omega _{n}^{-1}(A_{n}(\hat{t}%
_{n}))+n\times O_{P}(s_{n}+||t_{n}||_{e}a_{n}+||t_{n}||_{e}^{2}b_{n}) \\
&\leq & n(\mathbf{Z}_{n})^{\prime }\Omega _{n}^{-1}(\mathbf{Z}_{n})-n(\langle
\mathbf{v}_{n}^{\ast \prime },\widehat{\alpha }_{n}-\alpha _{0}\rangle +%
\mathbf{Z}_{n})^{\prime }\Omega _{n}^{-1}(\mathbf{Z}_{n})-n(A_{n}(\hat{t}%
_{n}))^{\prime }\Omega _{n}^{-1}\mathbf{Z}_{n}+0.5n(A_{n}(\hat{t}%
_{n}))^{\prime }\Omega _{n}^{-1}(A_{n}(\hat{t}_{n})) \\
&& +n\times O_{P}(s_{n}+||t_{n}||_{e}a_{n}+||t_{n}||_{e}^{2}b_{n})
\end{eqnarray*}%
where second line follows because $(\langle \mathbf{v}_{n}^{\ast \prime },%
\widehat{\alpha }_{n}-\alpha _{0}\rangle )^{\prime }\Omega _{n}^{-1}(\langle
\mathbf{v}_{n}^{\ast \prime },\widehat{\alpha }_{n}-\alpha _{0}\rangle )\geq
0$ and straightforward algebra. Observe that%
\begin{eqnarray*}
&&\left[ (A_{n}(\hat{t}_{n}))^{\prime }\Omega _{n}^{-1}(A_{n}(\hat{t}_{n}))%
\right] ^{1/2} \\
&=&\left[ (A_{n}(\hat{t}_{n}))^{\prime }\Omega _{sd,n}^{-1/2}\mathbb{D}%
_{n}\Omega _{sd,n}^{-1/2}(A_{n}(\hat{t}_{n}))\right] ^{1/2} \\
&=&O_{P}\left( ||\mathbb{D}_{n}^{1/2}\Omega _{sd,n}^{-1/2}\frac{d\phi
(\alpha _{0})}{d\alpha }[\alpha _{0n}-\alpha _{0}]||_{e}+||\mathbb{D}%
_{n}^{1/2}\hat{t}_{n}||_{e}\right) \\
&=&O_{P}\left( ||\Omega _{sd,n}^{-1/2}\frac{d\phi (\alpha _{0})}{d\alpha }%
[\alpha _{0n}-\alpha _{0}]||_{e}+||\hat{t}_{n}||_{e}\right) \\
&=&O_{P}\left( c_{n}(1+L_{n})\right) =o_{P}\left( n^{-1/2}\right)
\end{eqnarray*}%
where the first equation follows from Lemma \ref{lem:D_BDD}, and the last
equality follows from assumption \ref{ass:phivv1}(iii) and the results from
Parts (1) and (2). Also
\begin{equation*}
\left[ (\mathbf{Z}_{n})^{\prime }\Omega _{n}^{-1}(\mathbf{Z}_{n})\right]
^{1/2}=\left[ (\Omega _{sd,n}^{-1/2}\mathbf{Z}_{n})^{\prime }\mathbb{D}%
_{n}(\Omega _{sd,n}^{-1/2}\mathbf{Z}_{n})\right] ^{1/2}=O_{P}(n^{-1/2}\sqrt{%
d(n)})
\end{equation*}%
by Lemmas \ref{lem:D_BDD} and \ref{lem:Z_bdd}. Thus $(A_{n}(\hat{t}%
_{n}))^{\prime }\Omega _{n}^{-1}\mathbf{Z}_{n}=o_{P}(n^{-1/2})O_{P}(n^{-1/2}%
\sqrt{d(n)})=o_{P}(n^{-1}\sqrt{d(n)})$. Similarly, $(\langle \mathbf{v}%
_{n}^{\ast \prime },\widehat{\alpha }_{n}-\alpha _{0}\rangle +\mathbf{Z}%
_{n})^{\prime }\Omega _{n}^{-1}(\mathbf{Z}_{n})=o_{P}(n^{-1}\sqrt{d(n)})$ by
Lemmas \ref{lem:D_BDD} and \ref{lem:1order_vv}(1). Since $%
||t_{n}||_{e}=O_{P}(\sqrt{d(n)}n^{-1/2})$, we have $%
n(s_{n}+(||t_{n}||_{e}a_{n}+(||t_{n}||_{e}^{2})b_{n})=O_{P}(ns_{n}+\sqrt{d(n)%
}n^{1/2}a_{n}+\sqrt{d(n)}\sqrt{d(n)}b_{n})=o_{P}(n^{-1}\sqrt{d(n)})$ under
assumption \ref{ass:LAQ_vv}(ii). Therefore,
\begin{eqnarray*}
n\left[ \widehat{Q}_{n}(\widehat{\alpha }_{n}(t_{n}))-\widehat{Q}_{n}(%
\widehat{\alpha }_{n})\right] &\leq & n(\mathbf{Z}_{n})^{\prime }\Omega
_{n}^{-1}(\mathbf{Z}_{n})+o_{P}(\sqrt{d(n)})\\
&= &(\sqrt{n}\Omega _{sd,n}^{-1/2}%
\mathbf{Z}_{n})^{\prime }\mathbb{D}_{n}(\sqrt{n}\Omega _{sd,n}^{-1/2}\mathbf{%
Z}_{n})+o_{P}(\sqrt{d(n)})
\end{eqnarray*}%
\textit{Q.E.D.}

\medskip

\noindent \textbf{Proof of Theorem \ref{thm:chi2_vv}}\textsc{:} The proof is
very similar to that of Theorem \ref{thm:chi2} and we only provide main
steps here.

\textsc{Step 1.} Similar to Steps 1 and 2 in the proof of Theorem \ref%
{thm:chi2}, by the definitions of $\widehat{\alpha }_{n}^{R}$ and $\widehat{%
\alpha }_{n}$ and Assumption \ref{ass:LAQ_vv}(i), it follows that for any
(possibly random) $t\in \mathcal{T}_{n}^{M}$,
\begin{eqnarray*}
0.5\widehat{QLR}_{n}(\phi _{0})& \geq & 0.5n\left( \widehat{Q}_{n}(\widehat{%
\alpha }_{n}^{R})-\widehat{Q}_{n}(\widehat{\alpha }_{n}^{R}(t))\right)
-o_{P}(1) \\
& = &-n\left( t^{\prime }\Omega _{sd,n}^{-1/2}\{\mathbf{Z}_{n}+\langle \mathbf{%
v}_{n}^{\ast \prime },\widehat{\alpha }_{n}^{R}-\alpha _{0}\rangle
\}+0.5t^{\prime }\mathbb{B}_{n}t\right)\\
& &+ O_{P}(s_{n}n+n||t||_{e}a_{n}+n||t||_{e}^{2}b_{n}).
\end{eqnarray*}

By Assumption \ref{ass:phivv1}(i)(ii),
\begin{equation*}
\left\Vert \Omega _{sd,n}^{-1/2}\left( \underbrace{\phi (\widehat{\alpha }%
_{n}^{R})-\phi (\alpha _{0})}_{=0}-\frac{d\phi (\alpha _{0})}{d\alpha }[%
\widehat{\alpha }_{n}^{R}-\alpha _{0}]\right) \right\Vert _{e}=O_{P}(c_{n}).
\end{equation*}%
Hence, by Assumption \ref{ass:phivv1}(iii),
\begin{equation}
\left\Vert \Omega _{sd,n}^{-1/2}\langle \mathbf{v}_{n}^{\ast \prime },%
\widehat{\alpha }_{n}^{R}-\alpha _{0}\rangle \right\Vert _{e}=O_{P}(c_{n}).
\label{A6(iii)}
\end{equation}%
Since $\sup_{t:||t||_{e}=1}|t^{\prime }\{\mathbb{B}_{n}-\mathbb{D}%
_{n}^{-1}\}t|=O_{P}(b_{n})$ by assumption, we have: $t^{\prime }\mathbb{B}%
_{n}t\leq |t^{\prime }\{\mathbb{B}_{n}-\mathbb{D}_{n}^{-1}\}t|+t^{\prime }%
\mathbb{D}_{n}^{-1}t\leq ||t||_{e}^{2}O_{P}(b_{n})+t^{\prime }\mathbb{D}%
_{n}^{-1}t$ uniformly over $t\in \mathbb{R}^{d(n)}$ with $||t||_{e}=1$.
This, Assumption \ref{ass:LAQ_vv}(i) and equation (\ref{A6(iii)}) together
imply that
\begin{equation*}
0.5\widehat{QLR}_{n}(\phi _{0})\geq -n\left( t^{\prime }\Omega _{sd,n}^{-1/2}%
\mathbf{Z}_{n}+0.5t^{\prime }\mathbb{D}_{n}^{-1}t\right)
+O_{P}(s_{n}n+n||t||_{e}(a_{n}+c_{n})+n||t||_{e}^{2}b_{n}).
\end{equation*}%
In the above display we let $t^{\prime }=-\mathbf{Z}_{n}^{\prime }\Omega
_{sd,n}^{-1/2}\mathbb{D}_{n}$, which, by Lemmas \ref{lem:D_BDD} and \ref%
{lem:Z_bdd}, is an admissible choice and $||t||_{e}=O_{P}\left( n^{-1/2}%
\sqrt{d(n)}\right) $. Observe that $t_{n}^{\prime }\Omega _{sd,n}^{-1/2}%
\mathbf{Z}_{n}=-\mathbf{Z}_{n}^{\prime }\Omega _{sd,n}^{-1/2}\mathbb{D}%
_{n}\Omega _{sd,n}^{-1/2}\mathbf{Z}_{n}$ and $t_{n}^{\prime }\mathbb{D}%
_{n}^{-1}t_{n}=\mathbf{Z}_{n}^{\prime }\Omega _{sd,n}^{-1/2}\mathbb{D}%
_{n}\Omega _{sd,n}^{-1/2}\mathbf{Z}_{n}$, we obtain:
\begin{eqnarray*}
0.5\widehat{QLR}_{n}(\phi _{0}) &\geq &0.5(\sqrt{n}\Omega _{sd,n}^{-1/2}%
\mathbf{Z}_{n})^{\prime }\mathbb{D}_{n}(\sqrt{n}\Omega _{sd,n}^{-1/2}\mathbf{%
Z}_{n})+O_{P}\left( s_{n}n+n^{1/2}\sqrt{d(n)}(a_{n}+c_{n})+d(n)b_{n}\right)
\\
&=&0.5(\sqrt{n}\Omega _{sd,n}^{-1/2}\mathbf{Z}_{n})^{\prime }\mathbb{D}_{n}(%
\sqrt{n}\Omega _{sd,n}^{-1/2}\mathbf{Z}_{n})+o_{P}(\sqrt{d(n)}),
\end{eqnarray*}%
where the last equal sign is due to Assumptions \ref{ass:phivv1}(iv) and \ref%
{ass:LAQ_vv}(ii).

\textsc{Step 2.} Similar to Step 3 in the proof of Theorem \ref{thm:chi2},
by the definitions of $\widehat{\alpha }_{n}^{R}$ and $\widehat{\alpha }_{n}$
and the result that $\widehat{\alpha }_{n}(t_{n})\in \mathcal{A}_{k(n)}^{R}$
(Lemma \ref{lem:t_ex}), with $t_{n}$ given in Lemma \ref{lem:t_ex}, we
obtain:
\begin{equation*}
0.5\widehat{QLR}_{n}(\phi _{0})\leq 0.5n\left( \widehat{Q}_{n}(\widehat{%
\alpha }_{n}(t_{n}))-\widehat{Q}_{n}(\widehat{\alpha }_{n})\right) +o_{P}(1).
\end{equation*}%
By lemma \ref{lem:t_ex}(3), it follows that%
\begin{equation*}
0.5\widehat{QLR}_{n}(\phi _{0})\leq 0.5(\sqrt{n}\Omega _{sd,n}^{-1/2}\mathbf{%
Z}_{n})^{\prime }\mathbb{D}_{n}(\sqrt{n}\Omega _{sd,n}^{-1/2}\mathbf{Z}%
_{n})+o_{P}(\sqrt{d(n)}).
\end{equation*}

\textsc{Step 3.} The results in steps 1 and 2 together imply that
\begin{equation*}
\widehat{QLR}_{n}(\phi _{0})=(\sqrt{n}\Omega _{sd,n}^{-1/2}\mathbf{Z}%
_{n})^{\prime }\mathbb{D}_{n}(\sqrt{n}\Omega _{sd,n}^{-1/2}\mathbf{Z}%
_{n})+o_{P}\left( \sqrt{d(n)}\right) ,
\end{equation*}%
which establishes \textbf{Result (1)}.\newline

\textbf{Result (2)} directly follows from Result (1) and the fact that $%
\mathbb{D}_{n}=I_{d(n)}$, $\Omega _{sd,n}=\Omega _{0,n}$ when $\Sigma
=\Sigma _{0}$.

\textbf{Result (3)} follows from Result (2), $\Omega _{sd,n}=\Omega _{0,n}$
when $\Sigma =\Sigma _{0}$, and the following property of $\mathbf{W}%
_{n}\equiv n\mathbf{Z}_{n}^{\prime }\Omega _{sd,n}^{-1}\mathbf{Z}_{n}:$
\begin{equation*}
(2d(n))^{-1/2}\left( \mathbf{W}_{n}-d(n)\right) \Rightarrow N(0,1),
\end{equation*}%
which has been established in the proof of Theorem \ref{thm:wald_vv} Result
(3). \textit{Q.E.D.}

\subsection{Proofs for Section \protect\ref{subsec-A2} on series LS
estimator $\widehat{m}$ and its bootstrap version}

\noindent \textbf{Proof of Lemma \ref{lem:suff_mcond_boot}}: For \textbf{%
Result (1)}, since
\begin{eqnarray*}
M_{n}(Z^{n}) &\equiv& P_{V^{\infty }|Z^{\infty }}\left( \sup_{\mathcal{N}%
_{osn}}\frac{\overline{\tau }_{n}}{n}\sum_{i=1}^{n}\left\Vert \widehat{m}%
^{B}(X_{i},\alpha )-\widetilde{m}(X_{i},\alpha )-\widehat{m}%
^{B}(X_{i},\alpha _{0})\right\Vert _{e}^{2}\geq M\mid Z^{n}\right) \\
&\leq & P_{V^{\infty }|Z^{\infty }}\left( \sup_{\mathcal{N}_{osn}}\frac{%
\overline{\tau }_{n}}{n}\sum_{i=1}^{n}\left\Vert \widehat{m}%
^{B}(X_{i},\alpha )-\widehat{m}(X_{i},\alpha )-\{\widehat{m}%
^{B}(X_{i},\alpha _{0})-\widehat{m}(X_{i},\alpha _{0})\}\right\Vert
_{e}^{2}\geq \frac{M}{2}\mid Z^{n}\right) \\
&& +P_{V^{\infty }|Z^{\infty }}\left( \sup_{\mathcal{N}_{osn}}\frac{\overline{%
\tau }_{n}}{n}\sum_{i=1}^{n}\left\Vert \widehat{m}(X_{i},\alpha )-\widetilde{%
m}(X_{i},\alpha )-\widehat{m}(X_{i},\alpha _{0})\right\Vert _{e}^{2}\geq
\frac{M}{2}\mid Z^{n}\right) \\
&\equiv & M_{1,n}(Z^{n})+M_{2,n}(Z^{n}),
\end{eqnarray*}%
we have: for all $\delta >0$, there is a $M(\delta )>0$ such that for all $%
M\geq M(\delta )$,%
\begin{equation*}
P_{Z^{\infty }}\left( M_{n}(Z^{n})\geq 2\delta \right) \leq P_{Z^{\infty
}}\left( M_{1,n}(Z^{n})\geq \delta \right) +P_{Z^{\infty }}\left(
M_{2,n}(Z^{n})\geq \delta \right) .
\end{equation*}%
By following the proof of Lemma C.3(ii) of \cite{CP_WP07}, we have that $%
P_{Z^{\infty }}\left( M_{2,n}(Z^{n})\geq \delta \right) <\delta /2$
eventually. Thus, to establish Result (1), it suffices to bound%
\begin{equation*}
P_{Z^{\infty }}\left( \{M_{1,n}(Z^{n})\geq \delta \}\cap \{\lambda _{\min
}((P^{\prime }P)/n)>c\}\right) +P_{Z^{\infty }}(\lambda _{\min }((P^{\prime
}P)/n)\leq c).
\end{equation*}%
By Assumption \ref{ass:m_ls}(ii)(iii) and theorem 1 in \cite{NEWEY_97} $%
\lambda _{\min }((P^{\prime }P)/n)\geq c>0$ with probability $P_{Z^{\infty
}} $ approaching one, hence $P_{Z^{\infty }}(\lambda _{\min }((P^{\prime
}P)/n)\leq c)<\delta /4$ eventually. To bound the term corresponding to $%
M_{1,n}$, we note that\footnote{%
To ease the notational burden in the proof, we assume $d_{\rho }=1$; when $%
d_{\rho }>1$ the same proof steps hold, component by component.}
\begin{align*}
& \sum_{i=1}^{n}\left\Vert \widehat{m}^{B}\left( X_{i},\alpha \right) -%
\widehat{m}\left( X_{i},\alpha \right) -\{\widehat{m}^{B}(X_{i},\alpha _{0})-%
\widehat{m}\left( X_{i},\alpha _{0}\right) \}\right\Vert _{e}^{2} \\
& =\sum_{i=1}^{n}\Delta \zeta ^{B}(\alpha )^{\prime }P(P^{\prime
}P)^{-}p^{J_{n}}(X_{i})p^{J_{n}}(X_{i})^{\prime }(P^{\prime }P)^{-}P^{\prime
}\Delta \zeta ^{B}(\alpha ) \\
& =\Delta \zeta ^{B}(\alpha )^{\prime }P(P^{\prime }P)^{-}P^{\prime }\Delta
\zeta ^{B}(\alpha ) \\
& \leq \frac{1}{\lambda _{min}((P^{\prime }P)/n)}\{n^{-1}\Delta \zeta
^{B}(\alpha )^{\prime }PP^{\prime }\Delta \zeta ^{B}(\alpha )\};
\end{align*}%
where $\Delta \zeta ^{B}(\alpha )=((\omega _{1}-1)\Delta \rho (Z_{1},\alpha
),...,(\omega _{n}-1)\Delta \rho (Z_{n},\alpha ))^{\prime }$ with $\Delta
\rho (Z,\alpha )\equiv \rho (Z,\alpha )-\rho (Z,\alpha _{0})$. It is thus
sufficient to show that, for large enough $n$,
\begin{equation}
P_{Z^{\infty }}\left( P_{V^{\infty }|Z^{\infty }}\left( \sup_{\mathcal{N}%
_{osn}}\frac{\overline{\tau }_{n}}{n^{2}}\Delta \zeta ^{B}(\alpha )^{\prime
}PP^{\prime }\Delta \zeta ^{B}(\alpha )\geq M\mid Z^{n}\right) \geq \delta
\right) <\delta ,  \label{eqn:M_1n}
\end{equation}%
which is established in Lemma \ref{lem:M_1n}.

For \textbf{Result (2)}, recall that $\ell _{n}^{B}(x,\alpha )\equiv
\widetilde{m}(x,\alpha )+\widehat{m}^{B}(x,\alpha _{0})$. By similar
calculations to those in \cite{AC_Emetrica03} (p. 1824) it follows
\begin{align*}
& E_{P_{V^{\infty }}}\left[ n^{-1}\sum_{i=1}^{n}\left\Vert \widehat{m}%
^{B}(X_{i},\alpha _{0})\right\Vert _{e}^{2}\right] \\
=& E_{P_{V^{\infty }}}\left[ p^{J_{n}}(X_{i})^{\prime }(P^{\prime
}P)^{-}P^{\prime }E_{P_{V^{\infty }|X^{\infty }}}\left[ \rho ^{B}(\alpha
_{0})\rho ^{B}(\alpha _{0})^{\prime }|X^{n}\right] P(P^{\prime
}P)^{-}p^{J_{n}}(X_{i})\right]
\end{align*}%
where $\rho ^{B}(\alpha )\equiv (\rho ^{B}(V_{1},\alpha ),...,\rho
^{B}(V_{n},\alpha ))^{\prime }$ with $\rho ^{B}(V_{i.},\alpha )\equiv \omega
_{i}\rho (Z_{i},\alpha )$. Note that {\small {%
\begin{align*}
E_{P_{V|X^{\infty }}}[\rho ^{B}(V_{i},\alpha _{0})\rho ^{B}(V_{j},\alpha
_{0})^{\prime }|X^{n}]& =E_{P_{\Omega }}[\omega _{i}\omega
_{j}E_{P_{V|X}}[\rho (Z_{i},\alpha _{0})\rho (Z_{j},\alpha _{0})^{\prime
}|X_{i},X_{j}]] \\
& =0\quad \text{for all }i\neq j\text{,}
\end{align*}%
}}and%
\begin{equation*}
E_{P_{V|X^{\infty }}}[\rho ^{B}(V_{i},\alpha _{0})\rho ^{B}(V_{i},\alpha
_{0})^{\prime }|X^{n}]=\sigma _{\omega }^{2}\Sigma _{0}(X_{i})\text{.}
\end{equation*}%
So under Assumption \ref{ass:Wboot} or \ref{ass:Wboot_e}, Assumptions \ref%
{ass:sieve}(iv) and \ref{ass:m_ls}(ii), applying Markov inequality we
obtain: for all $\delta >0$, there is a $M(\delta )>0$ such that for all $%
M\geq M(\delta )$,
\begin{equation*}
P_{Z^{\infty }}\left( P_{V^{\infty }|Z^{\infty }}\left( \frac{J_{n}}{n}%
n^{-1}\sum_{i=1}^{n}\left\Vert \widehat{m}^{B}(X_{i},\alpha _{0})\right\Vert
_{e}^{2}\geq M\mid Z^{n}\right) \geq \delta \right) <\delta .
\end{equation*}%
To establish Result (2), with $(\tau _{n}^{\prime })^{-1}=\max \{\frac{J_{n}%
}{n},b_{m,J_{n}}^{2},(M_{n}\delta _{n})^{2}\}$, it remains to show that%
\begin{equation}
P_{Z^{\infty }}\left( \sup_{\mathcal{N}_{osn}}\tau _{n}^{\prime
}n^{-1}\sum_{i=1}^{n}\left\Vert \widetilde{m}(X_{i},\alpha )\right\Vert
_{e}^{2}\geq M\right) <\delta .  \label{mtildeb}
\end{equation}%
By Lemma SM.1 of \cite{CP_WPS}, under Assumptions \ref{ass:m_ls} and \ref%
{ass:rho_Donsker}(i), we have: there are finite constants $c,c^{\prime }>0$
such that, for all $\delta >0$, there is a $N(\delta )$, such that for all $%
n\geq N(\delta )$,
\begin{equation*}
P_{Z^{\infty }}\left( \forall \alpha \in \mathcal{N}_{osn}:~cE_{P_{X}}\left[
||\widetilde{m}\left( X,\alpha \right) ||_{e}^{2}\right] \leq \frac{1}{n}%
\sum_{i=1}^{n}||\widetilde{m}\left( X_{i},\alpha \right) ||_{e}^{2}\leq
c^{\prime }E_{P_{X}}\left[ ||\widetilde{m}\left( X,\alpha \right) ||_{e}^{2}%
\right] \right) >1-\delta .
\end{equation*}%
Thus to show (\ref{mtildeb}), it suffices to show that
\begin{equation*}
\sup_{\mathcal{N}_{osn}}\tau _{n}^{\prime }E_{P_{X}}\left[ \left\Vert
\widetilde{m}(X,\alpha )\right\Vert _{e}^{2}\right] =O(1).
\end{equation*}%
By Assumption \ref{ass:m_ls}(ii) it follows
\begin{align*}
\sup_{\alpha \in \mathcal{N}_{osn}}E_{P_{X}}\left[ \left\Vert \widetilde{m}%
(X,\alpha )\right\Vert _{e}^{2}\right] & \leq \sup_{\mathcal{N}%
_{osn}}\left\{ E_{P_{X}}\left[ \left\Vert \widetilde{m}(X,\alpha
)-m(X,\alpha )\right\Vert _{e}^{2}\right] +E_{P_{X}}\left[ \left\Vert
m(X,\alpha )\right\Vert _{e}^{2}\right] \right\} \\
& \leq const.\sup_{\alpha \in \mathcal{N}_{osn}}\max \left\{
b_{m,J_{n}}^{2},||\alpha -\alpha _{0}||^{2}\right\} =O((\tau _{n}^{\prime
})^{-1}),
\end{align*}%
where the last inequality follows from Assumptions \ref{ass:m_ls}%
(ii)(iii)(iv) and \ref{ass:weak_equiv}. We thus obtain Result (2).\newline

For \textbf{Result (3)}, we note that
\begin{equation*}
\frac{1}{n}\sum_{i=1}^{n}\left\Vert \widehat{m}^{B}\left( X_{i},\alpha
\right) \right\Vert _{\widehat{\Sigma }^{-1}}^{2}-\frac{1}{n}%
\sum_{i=1}^{n}\left\Vert \ell _{n}^{B}\left( X_{i},\alpha \right)
\right\Vert _{\widehat{\Sigma }^{-1}}^{2}=R_{1n}^{B}(\alpha
)+2R_{2n}^{B}(\alpha ),
\end{equation*}%
where
\begin{equation*}
R_{1n}^{B}(\alpha )\equiv \frac{1}{n}\sum_{i=1}^{n}\left\Vert \widehat{m}%
^{B}\left( X_{i},\alpha \right) -\ell _{n}^{B}\left( X_{i},\alpha \right)
\right\Vert _{\widehat{\Sigma }^{-1}}^{2},\text{ }R_{2n}^{B}(\alpha )\leq
\sqrt{R_{1n}^{B}}\sqrt{\frac{1}{n}\sum_{i=1}^{n}\left\Vert \ell
_{n}^{B}\left( X_{i},\alpha \right) \right\Vert _{\widehat{\Sigma }^{-1}}^{2}%
}.
\end{equation*}

By Result (1) and Assumption \ref{ass:VE}(iii), we have:
\begin{equation*}
P_{Z^{\infty }}\left( P_{V^{\infty }|Z^{\infty }}\left( \sup_{\mathcal{N}%
_{osn}}\overline{\tau }_{n}R_{1n}^{B}(\alpha )\geq M\mid Z^{n}\right) \geq
\delta \right) <\delta
\end{equation*}%
with $\overline{\tau }_{n}^{-1}=\delta _{n}^{2}(M_{n}\delta _{s,n})^{2\kappa
}C_{n}$. By Results (1) and (2), and Assumption \ref{ass:VE}(iii), we have:
\begin{equation*}
P_{Z^{\infty }}\left( P_{V^{\infty }|Z^{\infty }}\left( \sup_{\mathcal{N}%
_{osn}}\tilde{\tau}_{n}R_{2n}^{B}(\alpha )\geq M\mid Z^{n}\right) \geq
\delta \right) <\delta
\end{equation*}%
with $\tilde{\tau}_{n}^{-1}\equiv M_{n}\delta _{n}^{2}(M_{n}\delta
_{s,n})^{\kappa }\sqrt{C_{n}}$. By Assumption \ref{ass:rho_Donsker}(iii) and
the fact that $L_{n}$ diverges, we obtain the desired result. \textit{Q.E.D.}

In the following we state Lemma \ref{lem:M_1n} and its proof.

\begin{lemma}
\label{lem:M_1n} Let Assumptions \ref{ass:weak_equiv}(i)(ii), \ref{ass:m_ls}%
(iii), \ref{ass:rho_Donsker}(i)(ii) and either \ref{ass:Wboot} or \ref%
{ass:Wboot_e} hold. Then: for all $\delta >0$, there is a $M(\delta )>0$
such that for all $M\geq M(\delta )$,
\begin{equation*}
P_{Z^{\infty }}\left( P_{V^{\infty }|Z^{\infty }}\left( \sup_{\mathcal{N}%
_{osn}}\frac{\overline{\tau }_{n}}{n^{2}}\Delta \zeta ^{B}(\alpha )^{\prime
}PP^{\prime }\Delta \zeta ^{B}(\alpha )\geq M\mid Z^{n}\right) \geq \delta
\right) <0.5\delta
\end{equation*}%
eventually, with $\overline{\tau }_{n}^{-1}\equiv (\delta _{n})^{2}\left(
M_{n}\delta _{s,n}\right) ^{2\kappa }C_{n}$, where $\Delta \zeta ^{B}(\alpha
)=((\omega _{1}-1)\Delta \rho (Z_{1},\alpha ),...,(\omega _{n}-1)\Delta \rho
(Z_{n},\alpha ))^{\prime }$ and $\Delta \rho (Z,\alpha )\equiv \rho
(Z,\alpha )-\rho (Z,\alpha _{0})$.
\end{lemma}

\noindent \textbf{Proof of Lemma \ref{lem:M_1n}}: Denote
\begin{equation*}
M_{1n}^{\prime }(Z^{n})\equiv P_{V^{\infty }|Z^{\infty }}\left( \sup_{%
\mathcal{N}_{osn}}\frac{\overline{\tau }_{n}}{n^{2}}\Delta \zeta ^{B}(\alpha
)^{\prime }PP^{\prime }\Delta \zeta ^{B}(\alpha )\geq M\mid Z^{n}\right) .
\end{equation*}%
By the Markov inequality
\begin{equation*}
M_{1n}^{\prime }(Z^{n})\leq M^{-1}E_{P_{V^{\infty }|Z^{\infty }}}\left[
\sup_{\mathcal{N}_{osn}}\frac{\overline{\tau }_{n}}{n^{2}}\Delta \zeta
^{B}(\alpha )^{\prime }PP^{\prime }\Delta \zeta ^{B}(\alpha )\right] .
\end{equation*}%
Hence it is sufficient to bound
\begin{eqnarray*}
P_{Z^{\infty }}\left( M_{1n}^{\prime }(Z^{n})\geq \delta \right) &\leq &%
\frac{1}{M\delta }E_{P_{V^{\infty }}}\left[ \sup_{\mathcal{N}_{osn}}\frac{%
\overline{\tau }_{n}}{n^{2}}\Delta \zeta ^{B}(\alpha )^{\prime }PP^{\prime
}\Delta \zeta ^{B}(\alpha )\right] \\
&=&\frac{\overline{\tau }_{n}}{nM\delta }\sum_{j=1}^{J_{n}}E_{P_{V^{\infty
}}}\left[ \sup_{\mathcal{N}_{osn}}\left( \frac{1}{\sqrt{n}}%
\sum_{i=1}^{n}(\omega _{i}-1)f_{j}(Z_{i},\alpha )\right) ^{2}\right] ,
\end{eqnarray*}%
where the first inequality follows from the law of iterated expectations and
the Markov inequality, and the second equality is due to the notation $%
f_{j}(z,\alpha )\equiv p_{j}(x)\{\rho (z,\alpha )-\rho (z,\alpha _{0})\}$.

Under assumption \ref{ass:Wboot}, $\{(\omega _{i}-1)f_{j}(Z_{i},\alpha
)\}_{i=1}^{n}$ are independent, and thus, by proposition A.1.6 in \cite%
{VdV-W_book96} (VdV-W),
\begin{eqnarray*}
& & \frac{\overline{\tau }_{n}}{nM\delta }\sum_{j=1}^{J_{n}}E_{P_{V^{\infty }}}%
\left[ \sup_{\mathcal{N}_{osn}}\left( n^{-1/2}\sum_{i=1}^{n}(\omega
_{i}-1)f_{j}(Z_{i},\alpha )\right) ^{2}\right] \\
& \leq & \frac{\overline{\tau }_{n}}{nM\delta }\sum_{j=1}^{J_{n}}\left(
E_{P_{V^{\infty }}}\left[ \sup_{\mathcal{N}_{osn}}\left\vert
n^{-1/2}\sum_{i=1}^{n}(\omega _{i}-1)f_{j}(Z_{i},\alpha )\right\vert \right]
+\sqrt{E[\max_{i\leq n}\sup_{\mathcal{N}_{osn}}\left\vert n^{-1/2}(\omega
_{i}-1)f_{j}(Z_{i},\alpha )\right\vert ^{2}]}\right) ^{2}.
\end{eqnarray*}%
The second term in the RHS is bounded above by
\begin{equation*}
\sqrt{nn^{-1}E_{P_{V^{\infty }}}[(\omega _{i}-1)^{2}\sup_{\mathcal{N}%
_{osn}}\left\vert f_{j}(Z_{i},\alpha )\right\vert ^{2}]}\leq \sqrt{%
E_{P_{\omega }}[(\omega _{i}-1)^{2}]E_{P_{Z^{\infty }}}[\sup_{\mathcal{N}%
_{osn}}\left\vert f_{j}(Z_{i},\alpha )\right\vert ^{2}]}=O((M_{n}\delta
_{s,n})^{\kappa })
\end{equation*}%
by Assumptions \ref{ass:m_ls}(iii), \ref{ass:rho_Donsker}(ii) and \ref%
{ass:Wboot}. Hence, under assumption \ref{ass:Wboot} we need to control
\begin{equation}
\frac{\overline{\tau }_{n}}{nM\delta }\sum_{j=1}^{J_{n}}\left(
E_{P_{V^{\infty }}}\left[ \sup_{\mathcal{N}_{osn}}\left\vert
n^{-1/2}\sum_{i=1}^{n}(\omega _{i}-1)f_{j}(Z_{i},\alpha )\right\vert \right]
\right) ^{2}+O\left( \frac{\overline{\tau }_{n}J_{n}}{nM\delta }(M_{n}\delta
_{s,n})^{2\kappa }\right) .  \label{eqn:bdd_iid}
\end{equation}

Under Assumption \ref{ass:Wboot_e}, $((\omega _{i}-1)f_{j}(Z_{i},\alpha
))_{i}$ are \emph{not} independent. So we need to take some additional steps
to arrive to an equation of the form of (\ref{eqn:bdd_iid}). Under
Assumption \ref{ass:Wboot_e}, it follows
\begin{align*}
& \frac{\overline{\tau }_{n}}{M\delta }\sum_{j=1}^{J_{n}}E_{P_{V^{\infty }}}%
\left[ \sup_{\mathcal{N}_{osn}}\left( n^{-1}\sum_{i=1}^{n}(\omega
_{i}-1)f_{j}(Z_{i},\alpha )\right) ^{2}\right] \\
=& \frac{\overline{\tau }_{n}}{M\delta }\sum_{j=1}^{J_{n}}E_{P_{V^{\infty }}}%
\left[ \sup_{\mathcal{N}_{osn}}\left( n^{-1}\sum_{i=1}^{n}\omega
_{i}f_{j}(Z_{i},\alpha )-n^{-1}\sum_{i=1}^{n}f_{j}(Z_{i},\alpha )\right) ^{2}%
\right] \\
=& \frac{\overline{\tau }_{n}}{M\delta }\sum_{j=1}^{J_{n}}E_{P_{Z^{\infty
}}\times P_{\hat{Z}^{\infty }}}\left[ \sup_{\mathcal{N}_{osn}}\left(
n^{-1}\sum_{i=1}^{n}(\delta _{\hat{Z}_{i}}-\mathbb{P}_{n})[f_{j}(\cdot
,\alpha )]\right) ^{2}\right] ,
\end{align*}%
where the last line follows from the fact that $\omega _{i}$ are the number
of times the variable $Z_{i}$ appear on the bootstrap sample. Thus, the
distribution of $\omega _{i}\delta _{Z_{i}}$ is the same as that of $\delta
_{\hat{Z}_{i}}$ where $(\hat{Z}_{i})_{i}$ is the bootstrap sample, i.e., an
i.i.d. sample from $\mathbb{P}_{n}\equiv n^{-1}\sum_{i=1}^{n}\delta _{Z_{i}}$%
. By a slight adaptation of lemma 3.6.6 in VdV-W (allowing for square of the
norm), it follows
\begin{equation*}
E_{P_{Z^{\infty }}\times P_{\hat{Z}^{\infty }}}\left[ \sup_{\mathcal{N}%
_{osn}}\left( n^{-1}\sum_{i=1}^{n}(\delta _{\hat{Z}_{i}}-\mathbb{P}%
_{n})[f_{j}(\cdot ,\alpha )]\right) ^{2}\right] \leq E_{P_{Z^{\infty }}}%
\left[ E_{P_{\tilde{N}^{\infty }}}\left[ \sup_{\mathcal{N}_{osn}}\left(
n^{-1}\sum_{i=1}^{n}\tilde{N}_{i}\delta _{Z_{i}}[f_{j}(\cdot ,\alpha
)]\right) ^{2}\right] \right] ,
\end{equation*}%
where $\tilde{N}_{i}=N_{i}-N_{i}^{\prime }$ with $N_{i}$ and $N_{i}^{\prime
} $ being iid Poisson variables with parameter $0.5$ ($P_{\tilde{N}^{\infty
}}$ is the corresponding probability). Note that now, $\{\tilde{N}%
_{i}f_{j}(Z_{i},\alpha )\}_{i=1}^{n}$ are independent. So by proposition
A.1.6 in VdV-W,
\begin{align*}
& \frac{\overline{\tau }_{n}}{nM\delta }\sum_{j=1}^{J_{n}}E_{Q}\left[ \sup_{%
\mathcal{N}_{osn}}\left( n^{-1/2}\sum_{i=1}^{n}\tilde{N}_{i}f_{j}(Z_{i},%
\alpha )\right) ^{2}\right] \\
& \leq \frac{\overline{\tau }_{n}}{nM\delta }\sum_{j=1}^{J_{n}}\left( E_{Q}%
\left[ \sup_{\mathcal{N}_{osn}}\left\vert n^{-1/2}\sum_{i=1}^{n}\tilde{N}%
_{i}f_{j}(Z_{i},\alpha )\right\vert \right] +\sqrt{E[\max_{i\leq n}\sup_{%
\mathcal{N}_{osn}}\left\vert n^{-1/2}\tilde{N}_{i}f_{j}(Z_{i},\alpha
)\right\vert ^{2}]}\right) ^{2},
\end{align*}%
where $Q\equiv P_{Z^{\infty }}\times P_{\tilde{N}^{\infty }}$. By
Cauchy-Schwarz inequality, the second term in the RHS is bounded above by
\begin{equation*}
\sqrt{nn^{-1}E_{Q}[|\tilde{N}|^{2}\sup_{\mathcal{N}_{osn}}|f_{j}(Z,\alpha
)|^{2}]}\leq \sqrt{E_{P_{\tilde{N}}}[|\tilde{N}|^{2}]E_{P_{Z}}[\sup_{%
\mathcal{N}_{osn}}|f_{j}(Z,\alpha )|^{2}]}=O((M_{n}\delta _{s,n})^{\kappa })
\end{equation*}%
by Assumptions \ref{ass:m_ls}(iii) and \ref{ass:rho_Donsker}(ii) and $E[|%
\tilde{N}|^{2}]<\infty $. Therefore, under Assumption \ref{ass:Wboot_e} we
need to control
\begin{equation}
\frac{\overline{\tau }_{n}}{nM\delta }\sum_{j=1}^{J_{n}}\left( E_{Q}\left[
\sup_{\mathcal{N}_{osn}}\left\vert n^{-1/2}\sum_{i=1}^{n}\tilde{N}%
_{i}f_{j}(Z_{i},\alpha )\right\vert \right] \right) ^{2}+O\left( \frac{%
\overline{\tau }_{n}J_{n}}{nM\delta }(M_{n}\delta _{s,n})^{2\kappa }\right) .
\label{eqn:bdd_e}
\end{equation}

Applying lemma 2.9.1 of VdV-W we can bound the leading terms in equations (%
\ref{eqn:bdd_iid}) and (\ref{eqn:bdd_e}) respectively as follows,
\begin{align}
& \frac{\overline{\tau }_{n}}{nM\delta }\sum_{j=1}^{J_{n}}E_{P_{V^{\infty }}}%
\left[ \sup_{\mathcal{N}_{osn}}\left\vert n^{-1/2}\sum_{i=1}^{n}(\omega
_{i}-1)\delta _{Z_{i}}[f_{j}(\cdot ,\alpha )]\right\vert \right]  \notag
\label{eqn:VDVW_iid_1} \\
\leq & \frac{\overline{\tau }_{n}}{nM\delta }\sum_{j=1}^{J_{n}}\left\{
\int_{0}^{\infty }\sqrt{P(|\omega -1|\geq t)}dt\right\} \max_{1\leq l\leq
n}E_{P_{Z^{\infty }}\times P_{\epsilon ^{\infty }}}\left[ \sup_{\mathcal{N}%
_{osn}}\left\vert l^{-1/2}\sum_{i=1}^{l}\epsilon _{i}\delta
_{Z_{i}}[f_{j}(\cdot ,\alpha )]\right\vert \right] ,
\end{align}%
and
\begin{align}
& \frac{\overline{\tau }_{n}}{nM\delta }\sum_{j=1}^{J_{n}}E_{P_{Z^{\infty }}}%
\left[ E_{P_{\tilde{N}}}\left[ \sup_{\mathcal{N}_{osn}}\left\vert
n^{-1/2}\sum_{i=1}^{n}\tilde{N}_{i}\delta _{Z_{i}}[f_{j}(\cdot ,\alpha
)]\right\vert \right] \right]  \notag  \label{eqn:VDVW_e_1} \\
\leq & \frac{\overline{\tau }_{n}}{nM\delta }\sum_{j=1}^{J_{n}}\left\{
\int_{0}^{\infty }\sqrt{P(|\tilde{N}|\geq t)}dt\right\} \max_{1\leq l\leq
n}E_{P_{Z^{\infty }}\times P_{\epsilon ^{\infty }}}\left[ \sup_{\mathcal{N}%
_{osn}}\left\vert l^{-1/2}\sum_{i=1}^{l}\epsilon _{i}\delta
_{Z_{i}}[f_{j}(\cdot ,\alpha )]\right\vert \right] ,
\end{align}%
where $(\epsilon _{i})_{i=1}^{n}$ is a sequence of Rademacher random
variables.
Note that $\left\{ \int_{0}^{\infty }\sqrt{P(|\omega -1|\geq t)}dt\right\}
<\infty $ (under Assumption \ref{ass:Wboot}), and also $\left\{
\int_{0}^{\infty }\sqrt{P(|\tilde{N}|\geq t)}dt\right\} \leq 2\sqrt{2}$ (see
VdV-W p. 351). Hence in both cases we need to bound
\begin{align}
& \frac{\overline{\tau }_{n}}{nM\delta }\sum_{j=1}^{J_{n}}\left( \max_{1\leq
l\leq n}E_{P_{Z^{\infty }}\times P_{\epsilon ^{\infty }}}\left[ \sup_{%
\mathcal{N}_{osn}}\left\vert l^{-1/2}\sum_{i=1}^{l}\epsilon _{i}\delta
_{Z_{i}}[f_{j}(\cdot ,\alpha )]\right\vert \right] \right) ^{2}  \notag \\
\leq & \frac{\overline{\tau }_{n}}{nM\delta }\sum_{j=1}^{J_{n}}\left(
\max_{1\leq l\leq n}E_{P_{Z^{\infty }}\times P_{\epsilon ^{\infty }}}\left[
\sup_{\mathcal{N}_{osn}}\left\vert l^{-1/2}\sum_{i=1}^{l}\epsilon _{i}\delta
_{Z_{i}}[\bar{f}_{j}(\cdot ,\alpha )]\right\vert \right] \right.  \notag \\
& \left. +\max_{1\leq l\leq n}E_{P_{Z^{\infty }}\times P_{\epsilon ^{\infty
}}}\left[ \sup_{\mathcal{N}_{osn}}\left\vert l^{-1/2}\sum_{i=1}^{l}\epsilon
_{i}E_{P_{Z}}[f_{j}(Z,\alpha )]\right\vert \right] \right) ^{2}  \notag \\
\leq & 2T_{1,n}+2T_{2,n},  \label{eqn:M_1n_0}
\end{align}%
where $\bar{f}_{j}(\cdot ,\alpha )=f_{j}(\cdot ,\alpha ) - E[f_{j}(Z ,\alpha )]$,
\begin{equation*}
T_{1,n}=\frac{\overline{\tau }_{n}}{nM\delta }\sum_{j=1}^{J_{n}}\left(
\max_{1\leq l\leq n}E_{P_{Z^{\infty }}\times P_{\epsilon ^{\infty }}}\left[
\sup_{\mathcal{N}_{osn}}\left\vert l^{-1/2}\sum_{i=1}^{l}\epsilon _{i}\delta
_{Z_{i}}[\bar{f}_{j}(\cdot ,\alpha )]\right\vert \right] \right) ^{2}
\end{equation*}%
and%
\begin{equation*}
T_{2,n}=\frac{\overline{\tau }_{n}}{nM\delta }\sum_{j=1}^{J_{n}}\left(
\max_{1\leq l\leq n}E_{P_{Z^{\infty }}\times P_{\epsilon ^{\infty }}}\left[
\sup_{\mathcal{N}_{osn}}\left\vert l^{-1/2}\sum_{i=1}^{l}\epsilon
_{i}E_{P_{Z}}[f_{j}(Z,\alpha )]\right\vert \right] \right) ^{2}.
\end{equation*}%
To bound the term $T_{2,n}$, we note that
\begin{align*}
& \max_{1\leq j\leq J_{n}}\max_{1\leq l\leq n}E_{P_{Z^{\infty }}\times
P_{\epsilon ^{\infty }}}\left[ \sup_{\mathcal{N}_{osn}}\left\vert
l^{-1/2}\sum_{i=1}^{l}\epsilon _{i}E_{P_{Z}}[f_{j}(Z,\alpha )]\right\vert %
\right] \\
& =\max_{1\leq j\leq J_{n}}\max_{1\leq l\leq n}\sup_{\mathcal{N}%
_{osn}}|E_{P_{Z}}[f_{j}(Z,\alpha )]|E_{P_{\epsilon ^{\infty }}}\left[
\left\vert l^{-1/2}\sum_{i=1}^{l}\epsilon _{i}\right\vert \right] \\
& \leq \max_{1\leq j\leq J_{n}}\max_{1\leq l\leq n}\sup_{\mathcal{N}%
_{osn}}|E_{P_{X}}[p_{j}(X)\Delta m(X,\alpha )]|\sqrt{E_{P_{\epsilon ^{\infty
}}}\left[ \left( l^{-1/2}\sum_{i=1}^{l}\epsilon _{i}\right) ^{2}\right] } \\
& \leq \max_{1\leq j\leq J_{n}}\max_{1\leq l\leq n}\left( \sqrt{%
E_{P_{Z}}[|p_{j}(X)|^{2}]}\sup_{\mathcal{N}_{osn}}\sqrt{E_{P_{X}}[|\Delta
m(X,\alpha )|^{2}]}\sqrt{E_{P_{\epsilon ^{\infty }}}\left[
l^{-1}\sum_{i=1}^{l}\left( \epsilon _{i}\right) ^{2}\right] }\right) \\
& =O(M_{n}\delta _{n}),
\end{align*}%
where $\Delta m(X,\alpha )\equiv m(X,\alpha )-m(X,\alpha _{0})$ and the
inequality follows from Cauchy-Schwarz and the fact that $\epsilon _{i}$ are
independent, and the last two equal signs are due to Assumptions \ref%
{ass:weak_equiv}(i)(ii) and \ref{ass:m_ls}(iii). Thus $T_{2,n}\leq
const.\times (M_{n}\delta _{n})^{2}\frac{\overline{\tau }_{n}J_{n}}{nM\delta
}$.

To bound the term $T_{1,n}$, we note that by the \textquotedblleft
desymmetrization lemma\textquotedblright\ 2.3.6 in VdV-W (note that $\bar{f}%
_{j}(Z_{i},\alpha )$ are centered),
\begin{equation*}
T_{1,n}\leq const.\times \frac{\overline{\tau }_{n}}{nM\delta }%
\sum_{j=1}^{J_{n}}\max_{1\leq l\leq n}\left( E_{P_{Z^{\infty }}}\left[ \sup_{%
\mathcal{N}_{osn}}\left\vert l^{-1/2}\sum_{i=1}^{l}\bar{f}_{j}(Z_{i},\alpha
)\right\vert \right] \right) ^{2}.
\end{equation*}%
By \cite{VdV-W_book96} theorem 2.14.2, we have (up to some omitted
constant), for all $j$,
\begin{align*}
& E_{P_{Z^{\infty }}}\left[ \sup_{\alpha \in \mathcal{N}_{osn}}\left\vert
l^{-1/2}\sum_{i=1}^{l}\bar{f}_{j}(Z_{i},\alpha )\right\vert \right] \\
\leq & \left\{ (M_{n}\delta _{s,n})^{\kappa }\int_{0}^{1}\sqrt{1+\log
N_{[]}(w(M_{n}\delta _{s,n})^{\kappa },\mathcal{E}_{ojn},||\cdot
||_{L^{2}(f_{Z})})}dw\right\}
\end{align*}%
where $\mathcal{E}_{ojn}=\{p_{j}(\cdot )(\rho (\cdot ,\alpha )-\rho (\cdot
,\alpha _{0}))-E[p_{j}(\cdot )(\rho (\cdot ,\alpha )-\rho (\cdot ,\alpha
_{0})]):\alpha \in \mathcal{N}_{osn}\}$.

Given any $w>0$, let $(\{g_{l}^{m},g_{u}^{m}\})_{m=1,...,N(w)}$ be the $%
||.||_{L^{2}(f_{Z})}$-norm brackets of $\mathcal{O}_{on}$. If $\{\rho (\cdot
,\alpha )-\rho (\cdot ,\alpha _{0})\}\in \mathcal{O}_{on}$ belongs to a
bracket $\{g_{l}^{m},g_{u}^{m}\}$, then, since $|p_{j}(x)|<const<\infty $ by
Assumption \ref{ass:m_ls}(iii),
\begin{equation*}
g_{l}^{m}(Z)\leq p_{j}(X)\{\Delta \rho (Z,\alpha )\}\leq g_{u}^{m}(Z)
\end{equation*}%
(where $\{g_{l}^{m},g_{u}^{m}\}$ are transformations of the original ones,
given by $(1\{p_{j}>0\}g_{l}^{m}+1\{p_{j}\leq 0\}g_{u}^{m})p_{j}$ and $%
(1\{p_{j}>0\}g_{u}^{m}+1\{p_{j}\leq 0\}g_{l}^{m})p_{j}$ and since $%
|p_{j}(x)|<const<\infty $ the $||.||_{L^{2}(f_{Z})}$-norm of the new
brackets is given by $\delta \times 2const$. We keep the same notation and
omit the constant "$2const$" to ease the notational burden), and from the
previous calculations it is easy to see that
\begin{equation*}
\{g_{l}^{m}(Z)-E[g_{u}^{m}(Z)]\}\leq p_{j}(X)\Delta \rho (Z,\alpha
)-E[p_{j}(X)\Delta \rho (Z,\alpha )]\leq \{g_{u}^{m}(Z)-E[g_{l}^{m}(Z)]\}.
\end{equation*}%
So functions of the form $%
(\{(g_{l}^{m}(Z)-E[g_{u}^{m}(Z)]),(g_{u}^{m}(Z)-E[g_{l}^{m}(Z)])%
\})_{m=1,...,N(w)}$ form $||.||_{L^{2}(f_{V})}$-norm brackets on $\mathcal{E}%
_{ojn}$. By construction, $N_{[]}(w,\mathcal{E}_{ojn},||.||_{L^{2}(f_{Z})})%
\leq N(w)$. Hence (up to some omitted constants)
\begin{align*}
& E_{P_{Z^{\infty }}}\left[ \sup_{\alpha \in \mathcal{N}_{osn}}\left\vert
l^{-1/2}\sum_{i=1}^{l}\bar{f}_{j}(Z_{i},\alpha )\right\vert \right] \\
\leq & (M_{n}\delta _{s,n})^{\kappa }\max_{j=1,...,J_{n}}\left\{ \int_{0}^{1}%
\sqrt{1+\log N_{[]}(w(M_{n}\delta _{s,n})^{\kappa },\mathcal{O}_{on},||\cdot
||_{L^{2}(f_{Z})})}dw\right\} \\
\leq & (M_{n}\delta _{s,n})^{\kappa }\sqrt{C_{n}},
\end{align*}%
where the last inequality follows from assumption \ref{ass:rho_Donsker}(ii).
Notice that the above RHS does not depend on $l$ nor on $j$, so we obtain
\begin{equation}
\max_{1\leq j\leq J_{n}}\max_{1\leq l\leq n}E_{P_{Z^{\infty }}}\left[
\sup_{\alpha \in \mathcal{N}_{osn}}\left( l^{-1/2}\sum_{i=1}^{l}\bar{f}%
_{j}(Z_{i},\alpha )\right) ^{2}\right] \leq const.\times (M_{n}\delta
_{s,n})^{2\kappa }C_{n}  \label{non-b-bound}
\end{equation}%
and hence $T_{1,n}\leq const.\times (M_{n}\delta _{s,n})^{2\kappa }C_{n}%
\frac{\overline{\tau }_{n}J_{n}}{nM\delta }$.

Note that $\max \left\{ (M_{n}\delta _{n})^{2},(M_{n}\delta _{s,n})^{2\kappa
}\right\} =(M_{n}\delta _{s,n})^{2\kappa }$ (by assumption) and that $%
\overline{\tau }_{n}^{-1}\equiv \frac{J_{n}}{n}\left( M_{n}\delta
_{s,n}\right) ^{2\kappa }C_{n}$, the desired result follows. \textit{Q.E.D.}

\medskip

\noindent \textbf{Proof of Lemma \ref{lem:T1n}}: Denote
\begin{equation*}
T_{nI}^{B}\equiv \sup_{\mathcal{N}_{osn}}\left\vert \frac{1}{n}%
\sum_{i=1}^{n}\left( \frac{d\widetilde{m}(X_{i},\alpha )}{d\alpha }%
[u_{n}^{\ast }]\right) ^{\prime }\widehat{\Sigma }(X_{i})^{-1}\ell
_{n}^{B}(X_{i},\alpha )-\frac{1}{n}\sum_{i=1}^{n}\left( \frac{%
dm(X_{i},\alpha _{0})}{d\alpha }[u_{n}^{\ast }]\right) ^{\prime }\Sigma
(X_{i})^{-1}\ell _{n}^{B}(X_{i},\alpha )\right\vert ,
\end{equation*}%
and
\begin{equation*}
T_{nII}^{B}\equiv \sup_{\mathcal{N}_{osn}}\left\vert \frac{1}{n}%
\sum_{i=1}^{n}\left( \frac{dm(X_{i},\alpha _{0})}{d\alpha }[u_{n}^{\ast
}]\right) ^{\prime }\Sigma (X_{i})^{-1}\ell _{n}^{B}(X_{i},\alpha )-\left\{
\mathbb{Z}_{n}^{\omega }+\langle u_{n}^{\ast },\alpha -\alpha _{0}\rangle
\right\} \right\vert .
\end{equation*}%
It suffices to show that for all $\delta >0$, there is $N(\delta )$ such
that for all $n\geq N(\delta )$,
\begin{equation}
P_{Z^{\infty }}\left( P_{V^{\infty }|Z^{\infty }}\left( \sqrt{n}%
T_{nI}^{B}\geq \delta \mid Z^{n}\right) \geq \delta \right) <\delta
\label{eqn:proof_T1n_1}
\end{equation}%
and
\begin{equation}
P_{Z^{\infty }}\left( P_{V^{\infty }|Z^{\infty }}\left( \sqrt{n}%
T_{nII}^{B}\geq \delta \mid Z^{n}\right) \geq \delta \right) <\delta .
\label{eqn:proof_T1n_2}
\end{equation}

We first verify equation (\ref{eqn:proof_T1n_1}). Note that
\begin{align*}
T_{nI}^{B}\leq & \sup_{\mathcal{N}_{osn}}\left\vert \frac{1}{n}%
\sum_{i=1}^{n}\left( \frac{d\widetilde{m}(X_{i},\alpha )}{d\alpha }%
[u_{n}^{\ast }]-\frac{dm(X_{i},\alpha _{0})}{d\alpha }[u_{n}^{\ast }]\right)
^{\prime }\widehat{\Sigma }(X_{i})^{-1}\ell _{n}^{B}(X_{i},\alpha
)\right\vert \\
& +\sup_{\mathcal{N}_{osn}}\left\vert \frac{1}{n}\sum_{i=1}^{n}\left( \frac{%
dm(X_{i},\alpha _{0})}{d\alpha }[u_{n}^{\ast }]\right) ^{\prime }\{\widehat{%
\Sigma }(X_{i})^{-1}-\Sigma (X_{i})^{-1}\}\ell _{n}^{B}(X_{i},\alpha
)\right\vert \\
\equiv & T_{nIa}^{B}+T_{nIb}^{B}.
\end{align*}%
By Assumption \ref{ass:VE}(iii) and the Cauchy-Schwarz inequality, it
follows that, for some $C\in (0,\infty )$,
\begin{eqnarray*}
& & P_{Z^{\infty }}\left( P_{V^{\infty }|Z^{\infty }}\left( \sqrt{n}%
T_{nIa}^{B}\geq \delta \mid Z^{n}\right) \geq \delta \right) \\
& \leq & P_{Z^{\infty }}\left( P_{V^{\infty }|Z^{\infty }}\left( \sup_{\mathcal{N}%
_{osn}}\sqrt{\frac{\sum_{i=1}^{n}\left\Vert \frac{dm(X_{i},\alpha _{0})}{%
d\alpha }[u_{n}^{\ast }]-\frac{d\widetilde{m}(X_{i},\alpha )}{d\alpha }%
[u_{n}^{\ast }]\right\Vert _{e}^{2}}{n}}\sqrt{\frac{\sum_{i=1}^{n}\left\Vert
\ell _{n}^{B}(X_{i},\alpha )\right\Vert _{e}^{2}}{n}}\geq \frac{C\delta }{%
\sqrt{n}}\mid Z^{n}\right) \geq \delta \right) \\
& & +P_{Z^{\infty }}\left( \lambda _{min}(\widehat{\Sigma }(X))<c\right) .
\end{eqnarray*}%
 The second term in the RHS vanishes eventually, so we focus on the first
term. It follows
\begin{eqnarray*}
& & P_{Z^{\infty }}\left( P_{V^{\infty }|Z^{\infty }}\left( \sup_{\mathcal{N}%
_{osn}}\sqrt{\frac{\sum_{i=1}^{n}\left\Vert \frac{dm(X_{i},\alpha _{0})}{%
d\alpha }[u_{n}^{\ast }]-\frac{d\widetilde{m}(X_{i},\alpha )}{d\alpha }%
[u_{n}^{\ast }]\right\Vert _{e}^{2}}{n}}\sqrt{\frac{1}{n}\sum_{i=1}^{n}\left%
\Vert \ell _{n}^{B}(X_{i},\alpha )\right\Vert _{e}^{2}}\geq \frac{C\delta }{%
\sqrt{n}}\mid Z^{n}\right) \geq \delta \right) \\
& \leq & P_{Z^{\infty }}\left( P_{V^{\infty }|Z^{\infty }}\left( \sup_{%
\mathcal{N}_{osn}}\sqrt{\frac{1}{n}\sum_{i=1}^{n}\left\Vert \frac{%
dm(X_{i},\alpha _{0})}{d\alpha }[u_{n}^{\ast }]-\frac{d\widetilde{m}%
(X_{i},\alpha )}{d\alpha }[u_{n}^{\ast }]\right\Vert _{e}^{2}}\sqrt{\frac{Mn%
}{\tau _{n}^{\prime }}}\geq C\delta \mid Z^{n}\right) \geq 0.5\delta \right)
\\
& + & P_{Z^{\infty }}\left( P_{V^{\infty }|Z^{\infty }}\left( \sup_{\mathcal{N}%
_{osn}}\sqrt{\frac{\tau _{n}^{\prime }}{n}\sum_{i=1}^{n}\left\Vert \ell
_{n}^{B}(X_{i},\alpha )\right\Vert _{e}^{2}}\geq \sqrt{M}\mid Z^{n}\right)
\geq 0.5\delta \right) .
\end{eqnarray*}%
By Lemma \ref{lem:suff_mcond_boot}(2) the second term on the RHS is less
than $0.5\delta $ eventually (with $\left( \tau _{n}^{\prime }\right)
^{-1}=const.(M_{n}\delta _{n})^{2}$). Regarding the first term, note that
\begin{align*}
& \sup_{\mathcal{N}_{osn}}\sqrt{\frac{1}{n}\sum_{i=1}^{n}\left\Vert \frac{%
dm(X_{i},\alpha _{0})}{d\alpha }[u_{n}^{\ast }]-\frac{d\widetilde{m}%
(X_{i},\alpha )}{d\alpha }[u_{n}^{\ast }]\right\Vert _{e}^{2}}\sqrt{\frac{n}{%
\tau _{n}^{\prime }}} \\
& \leq \sup_{\mathcal{N}_{osn}}\sqrt{\frac{1}{n}\sum_{i=1}^{n}\left\Vert
\frac{d\widetilde{m}(X_{i},\alpha _{0})}{d\alpha }[u_{n}^{\ast }]-\frac{d%
\widetilde{m}(X_{i},\alpha )}{d\alpha }[u_{n}^{\ast }]\right\Vert
_{e}^{2}\times \frac{n}{\tau _{n}^{\prime }}} \\
& +\sup_{\mathcal{N}_{osn}}\sqrt{\frac{1}{n}\sum_{i=1}^{n}\left\Vert \frac{%
dm(X_{i},\alpha _{0})}{d\alpha }[u_{n}^{\ast }]-\frac{d\widetilde{m}%
(X_{i},\alpha _{0})}{d\alpha }[u_{n}^{\ast }]\right\Vert _{e}^{2}\times
\frac{n}{\tau _{n}^{\prime }}} \\
& \leq \sup_{\mathcal{N}_{osn}}\sqrt{\frac{1}{n}\sum_{i=1}^{n}\left\Vert
\frac{dm(X_{i},\alpha _{0})}{d\alpha }[u_{n}^{\ast }]-\frac{dm(X_{i},\alpha )%
}{d\alpha }[u_{n}^{\ast }]\right\Vert _{e}^{2}\times \frac{n}{\tau
_{n}^{\prime }}}+o_{P_{Z^{\infty }}}(1),
\end{align*}%
by the LS projection property and the definition of $\widetilde{m}$, as well
as by the Markov inequality and Assumption \ref{ass:anor-mtilde}(i). Next,
by the Markov inequality and Assumption \ref{ass:cont_diffm}(ii), we have:%
\begin{eqnarray*}
&&P_{Z^{\infty }}\left( \sup_{\mathcal{N}_{osn}}\sqrt{\frac{1}{n}%
\sum_{i=1}^{n}\left\Vert \frac{dm(X_{i},\alpha _{0})}{d\alpha }[u_{n}^{\ast
}]-\frac{dm(X_{i},\alpha )}{d\alpha }[u_{n}^{\ast }]\right\Vert _{e}^{2}}%
\sqrt{\frac{n}{\tau _{n}^{\prime }}}\geq 0.5\delta \right) \\
&\leq &\frac{2}{\delta }\sqrt{E_{P_{Z^{\infty }}}\left[ \sup_{\mathcal{N}%
_{osn}}\left\Vert \frac{dm(X,\alpha _{0})}{d\alpha }[u_{n}^{\ast }]-\frac{%
dm(X,\alpha )}{d\alpha }[u_{n}^{\ast }]\right\Vert _{e}^{2}\right] \times
\frac{n}{\tau _{n}^{\prime }}}\rightarrow 0.
\end{eqnarray*}%
Thus, we established that
\begin{equation*}
P_{Z^{\infty }}\left( P_{V^{\infty }|Z^{\infty }}\left( \sqrt{n}%
T_{nIa}^{B}\geq \delta \mid Z^{n}\right) \geq \delta \right) <\delta \text{%
\quad eventually.}
\end{equation*}%
By similar arguments, Assumptions \ref{ass:VE}(iii) and \ref{ass:rho_Donsker}%
(iv), Lemma \ref{lem:suff_mcond_boot}(2), and that $\frac{1}{n}%
\sum_{i=1}^{n}\left\Vert \frac{dm(X_{i},\alpha _{0})}{d\alpha }[u_{n}^{\ast
}]\right\Vert _{e}^{2}$ is bounded in probability, it can be shown that
\begin{equation*}
P_{Z^{\infty }}\left( P_{V^{\infty }|Z^{\infty }}\left( \sqrt{n}%
T_{nIb}^{B}\geq \delta \mid Z^{n}\right) \geq \delta \right)
<\delta,~eventually.
\end{equation*}%
Therefore, we establish equation (\ref{eqn:proof_T1n_1}).

For equation (\ref{eqn:proof_T1n_2}), let $g(X,u_{n}^{\ast })\equiv \left(
\frac{dm(X,\alpha _{0})}{d\alpha }[u_{n}^{\ast }]\right) ^{\prime }\Sigma
^{-1}(X)$. Then
\begin{align*}
T_{nII}^{B}& \leq \sup_{\mathcal{N}_{osn}}\left\vert \frac{1}{n}%
\sum_{i=1}^{n}g(X_{i},u_{n}^{\ast })\widetilde{m}(X_{i},\alpha )-\langle
u_{n}^{\ast },\alpha -\alpha _{0}\rangle \right\vert +\left\vert \frac{1}{n}%
\sum_{i=1}^{n}g(X_{i},u_{n}^{\ast })\widehat{m}^{B}(X_{i},\alpha _{0})-%
\mathbb{Z}_{n}^{\omega }\right\vert \\
& \equiv T_{nIIa}+T_{nIIb}^{B}.
\end{align*}%
Thus to show equation (\ref{eqn:proof_T1n_2}) it suffices to show that $%
\sqrt{n}T_{nIIa}=o_{P_{Z^{\infty }}}(1)$ and that
\begin{equation}
P_{Z^{\infty }}\left( P_{V^{\infty }|Z^{\infty }}\left( \sqrt{n}%
T_{nIIb}^{B}\geq \delta \mid Z^{n}\right) \geq \delta \right) <\delta \quad
\text{eventually}.  \label{IIbB}
\end{equation}

First we consider the term $T_{nIIa}$. This part of proof is similar to
those in \cite{AC_Emetrica03}, \cite{AC_JOE07} and \cite{CP_WP07a} for their
regular functional $\lambda ^{\prime }\theta $ case, and hence we shall be
brief. By the orthogonality properties of the LS projection and the
definition of $\widetilde{m}(X_{i},\alpha )$ and $\widetilde{g}%
(X_{i},u_{n}^{\ast })$, we have:
\begin{equation*}
n^{-1}\sum_{i=1}^{n}g(X_{i},u_{n}^{\ast })\widetilde{m}(X_{i},\alpha
)=n^{-1}\sum_{i=1}^{n}\widetilde{g}(X_{i},u_{n}^{\ast })m(X_{i},\alpha ).
\end{equation*}%
By Cauchy-Schwarz inequality,
\begin{eqnarray*}
&&\sup_{\mathcal{N}_{osn}}\left\vert \frac{1}{n}\sum_{i=1}^{n}\{\widetilde{g}%
(X_{i},u_{n}^{\ast })-g(X_{i},u_{n}^{\ast })\}\{m(X_{i},\alpha
)-m(X_{i},\alpha _{0})\}\right\vert \\
&\leq &\sqrt{\frac{1}{n}\sum_{i=1}^{n}||\widetilde{g}(X_{i},u_{n}^{\ast
})-g(X_{i},u_{n}^{\ast })||_{e}^{2}}\sup_{\mathcal{N}_{osn}}\sqrt{\frac{1}{n}%
\sum_{i=1}^{n}||m(X_{i},\alpha )-m(X_{i},\alpha _{0})||_{e}^{2}}.
\end{eqnarray*}%
By assumption \ref{ass:anor-mtilde}(iii),
\begin{eqnarray*}
\sqrt{n}\sup_{\mathcal{N}_{osn}}%
\frac{1}{n}\sum_{i=1}^{n}\left\{ ||m(X_{i},\alpha )-m(X_{i},\alpha
_{0})||_{e}^{2}-E_{P_{X}}[||m(X_{1},\alpha )-m(X_{1},\alpha
_{0})||_{e}^{2}]\right\} =o_{P}(1).
\end{eqnarray*}
Thus, since $\sup_{\mathcal{N}%
_{osn}}E_{P_{X}}[||m(X_{1},\alpha )-m(X_{1},\alpha
_{0})||_{e}^{2}]=O(M_{n}^{2}\delta _{n}^{2})$, it follows
\begin{equation*}
\sup_{\mathcal{N}_{osn}}\frac{1}{n}\sum_{i=1}^{n}||m(X_{i},\alpha
)-m(X_{i},\alpha _{0})||_{e}^{2}=O_{P_{Z^{\infty }}}\left( (M_{n}\delta
_{n})^{2}+o_{P_{Z^{\infty }}}(n^{-1/2})\right) .
\end{equation*}%
This, Assumption \ref{ass:anor-mtilde}(ii) and $\delta _{n}=o(n^{-1/4})$ (by
assumption \ref{ass:rho_Donsker}(iv)) imply that%
\begin{eqnarray*}
&&\sup_{\mathcal{N}_{osn}}\left\vert \frac{1}{n}\sum_{i=1}^{n}\{\widetilde{g}%
(X_{i},u_{n}^{\ast })-g(X_{i},u_{n}^{\ast })\}\{m(X_{i},\alpha
)-m(X_{i},\alpha _{0})\}\right\vert \\
&\leq &o_{P_{Z^{\infty }}}(\frac{1}{\sqrt{n}M_{n}\delta _{n}})\times
O_{P_{Z^{\infty }}}\left( \sqrt{(M_{n}\delta _{n})^{2}+o(n^{-1/2})}\right)
=o_{P_{Z^{\infty }}}(n^{-1/2})
\end{eqnarray*}%
Therefore,
\begin{equation*}
\sqrt{n}T_{nIIa}=\sqrt{n}\sup_{\mathcal{N}_{osn}}\left\vert \frac{1}{n}%
\sum_{i=1}^{n}g(X_{i},u_{n}^{\ast })m(X_{i},\alpha )-\langle u_{n}^{\ast
},\alpha -\alpha _{0}\rangle \right\vert +o_{P_{Z^{\infty }}}(n^{-1/2}).
\end{equation*}%
By assumption \ref{ass:anor-mtilde}(iv), $\sqrt{n}\sup_{\mathcal{N}%
_{osn}}\left\vert \frac{1}{n}\sum_{i=1}^{n}g(X_{i},u_{n}^{\ast
})m(X_{i},\alpha )-E_{P_{X}}\left[ g(X_{1},u_{n}^{\ast })\{m(X_{1},\alpha
)-m(X_{1},\alpha _{0})\right] \right\vert =o_{P_{Z^{\infty }}}(1)$. Thus, by
Assumption \ref{ass:cont_diffm}(iv), we conclude that $\sqrt{n}%
T_{nIIa}=o_{P_{Z^{\infty }}}(1)$.

Next we consider the term $T_{nIIb}^{B}$. By the orthogonality properties of
the LS projection,
\begin{equation*}
n^{-1}\sum_{i=1}^{n}g(X_{i},u_{n}^{\ast })\widehat{m}^{B}(X_{i},\alpha
_{0})=n^{-1}\sum_{i=1}^{n}\widetilde{g}(X_{i},u_{n}^{\ast })\rho
^{B}(V_{i},\alpha _{0}),
\end{equation*}%
where $\rho ^{B}(V_{i},\alpha _{0})\equiv \omega _{i,n}\rho (Z_{i},\alpha
_{0})$ and $\left\{ \omega _{i,n}\right\} _{i=1}^{n}$ is independent of $%
\left\{ Z_{i}\right\} _{i=1}^{n}$.

Hence, by applying the Markov inequality twice, it follows that
\begin{eqnarray*}
&&P_{Z^{\infty }}\left( P_{V^{\infty }|Z^{\infty }}\left( \sqrt{n}%
T_{nIIb}^{B}\geq \delta \mid Z^{n}\right) \geq \delta \right) \\
&\leq &\delta ^{-4}E_{P_{V^{\infty }}}\left[ n^{-1}\left(
\sum_{i=1}^{n}\{g(X_{i},u_{n}^{\ast })-\widetilde{g}(X_{i},u_{n}^{\ast
})\}\rho ^{B}(V_{i},\alpha _{0})\right) ^{2}\right] .
\end{eqnarray*}%
Regarding the cross-products terms where $i\neq j$, note that
\begin{eqnarray*}
& & E_{P_{V^{\infty }}}\left[ \{g(X_{j},u_{n}^{\ast })-\widetilde{g}%
(X_{j},u_{n}^{\ast })\}\{g(X_{i},u_{n}^{\ast })-\widetilde{g}%
(X_{i},u_{n}^{\ast })\}\rho ^{B}(V_{i},\alpha _{0})\rho ^{B}(V_{j},\alpha
_{0})\right] \\
&=& E_{P_{V^{\infty }}}\left[ \{g(X_{j},u_{n}^{\ast })-\widetilde{g}%
(X_{j},u_{n}^{\ast })\}\{g(X_{i},u_{n}^{\ast })-\widetilde{g}%
(X_{i},u_{n}^{\ast })\}E_{P_{V^{\infty }|X^{\infty }}}\left[ \rho
^{B}(V_{i},\alpha _{0})\rho ^{B}(V_{j},\alpha _{0})\mid X^{n}\right] \right]
\\
&=& E_{P_{V^{\infty }}}\left[ \{g(X_{j},u_{n}^{\ast })-\widetilde{g}%
(X_{j},u_{n}^{\ast })\}\{g(X_{i},u_{n}^{\ast })-\widetilde{g}%
(X_{i},u_{n}^{\ast })\}E_{P_{V^{\infty }|X^{\infty }}}\left[ \omega
_{i}\omega _{j}|X^{n}\right] E_{P_{Z^{\infty }|X^{\infty }}}\left[ \rho
(Z_{i},\alpha _{0})\rho (Z_{j},\alpha _{0})\mid X^{n}\right] \right] \\
&=& 0,
\end{eqnarray*}%
since $E_{P_{Z^{\infty }|X^{\infty }}}\left[ \rho (Z_{i},\alpha _{0})\rho
(Z_{j},\alpha _{0})\mid X^{n}\right] =E_{P_{Z|X}}\left[ \rho (Z_{i},\alpha
_{0}) \vert X_{i}\right] E_{P_{Z|X}}\left[ \rho (Z_{j},\alpha _{0})\vert X_{j}\right]
=0$ for $i\neq j$. Thus, it suffices to study
\begin{align*}
& \delta ^{-4}E_{P_{V^{\infty }}}\left[ n^{-1}\sum_{i=1}^{n}\left(
g(X_{i},u_{n}^{\ast })-\widetilde{g}(X_{i},u_{n}^{\ast })\right) ^{2}\left(
\rho ^{B}(V_{i},\alpha _{0})\right) ^{2}\right] \\
& =\delta ^{-4}n^{-1}\sum_{i=1}^{n}E_{P_{V^{\infty }}}\left[ \left(
g(X_{i},u_{n}^{\ast })-\widetilde{g}(X_{i},u_{n}^{\ast })\right)
^{2}E_{P_{V^{\infty }|X^{\infty }}}\left[ \left( \omega _{i}\rho
(Z_{i},\alpha _{0})\right) ^{2}\mid X^{n}\right] \right] .
\end{align*}%
By the original-sample $\left\{ Z_{i}\right\} _{i=1}^{n}$ being i.i.d., $%
\left\{ \omega _{i,n}\right\} _{i=1}^{n}$ being independent of $\left\{
Z_{i}\right\} _{i=1}^{n}$, Assumption \ref{ass:sieve}(iv) and the fact that $%
\sigma _{\omega }^{2}<\infty $, we can majorize the previous expression (up
to an omitted constant) by
\begin{equation*}
\delta ^{-4}E_{P_{V^{\infty }}}\left[ \left( g(X_{i},u_{n}^{\ast })-%
\widetilde{g}(X_{i},u_{n}^{\ast })\right) ^{2}\right] =o(1),
\end{equation*}%
where the last equality is due to Assumption \ref{ass:anor-mtilde}(ii).
Hence we established equation (\ref{IIbB}). The desired result now follows.
\textit{Q.E.D.}

\medskip

\noindent \textbf{Proof of Lemma \ref{lem:T3n}}: By the Cauchy-Schwarz
inequality and Assumption \ref{ass:VE}(iii), it suffices to show that
\begin{equation*}
P_{Z^{\infty }}\left( P_{V^{\infty }|Z^{\infty }}\left( \sup_{\mathcal{N}%
_{osn}}\sqrt{n^{-1}\sum_{i=1}^{n}\left\Vert \frac{d^{2}\widetilde{m}%
(X_{i},\alpha )}{d\alpha ^{2}}[u_{n}^{\ast },u_{n}^{\ast }]\right\Vert
_{e}^{2}}\sup_{\mathcal{N}_{osn}}\sqrt{n^{-1}\sum_{i=1}^{n}\left\Vert \ell
_{n}^{B}(X_{i},\alpha )\right\Vert _{e}^{2}}\geq \delta \mid Z^{n}\right)
\geq \delta \right) <\delta .
\end{equation*}%
By Lemma \ref{lem:suff_mcond_boot}(2), it suffices to show that
\begin{equation*}
P_{Z^{\infty }}\left( \sup_{\mathcal{N}_{osn}}\sqrt{n^{-1}\sum_{i=1}^{n}%
\left\Vert \frac{d^{2}\widetilde{m}(X_{i},\alpha )}{d\alpha ^{2}}%
[u_{n}^{\ast },u_{n}^{\ast }]\right\Vert _{e}^{2}}\geq \frac{\delta }{%
M_{n}\delta _{n}}\right) <\delta .
\end{equation*}%
By Markov inequality and the LS projection properties, the LHS of the
previous equation can be bounded above by
\begin{equation*}
\frac{M_{n}^{2}\delta _{n}^{2}}{\delta ^{2}}E_{P_{X}}\left[ \sup_{\mathcal{N}%
_{osn}}\left\Vert \frac{d^{2}\widetilde{m}(X,\alpha )}{d\alpha ^{2}}%
[u_{n}^{\ast },u_{n}^{\ast }]\right\Vert _{e}^{2}\right] \leq \frac{%
M_{n}^{2}\delta _{n}^{2}}{\delta ^{2}}E_{P_{X}}\left[ \sup_{\mathcal{N}%
_{osn}}\left\Vert \frac{d^{2}m(X,\alpha )}{d\alpha ^{2}}[u_{n}^{\ast
},u_{n}^{\ast }]\right\Vert _{e}^{2}\right] <\delta
\end{equation*}%
eventually, which is satisfied given Assumption \ref{ass:cont_diffm}(iii).
The desired result follows. \textit{Q.E.D.}

\medskip

\noindent \textbf{Proof of Lemma \ref{lem:T2n}}: For \textbf{Result (1)}, we
first want to show that
\begin{equation}
\sup_{\mathcal{N}_{osn}}\left\vert \frac{1}{n}\sum_{i=1}^{n}\left\{
\left\Vert \frac{d\widetilde{m}(X_{i},\alpha )}{d\alpha }[u_{n}^{\ast
}]\right\Vert _{\widehat{\Sigma }^{-1}}^{2}-\left\Vert \frac{dm(X_{i},\alpha
_{0})}{d\alpha }[u_{n}^{\ast }]\right\Vert _{\Sigma ^{-1}}^{2}\right\}
\right\vert \leq T_{n,I}+T_{n,II}+T_{n,III}=o_{P_{Z^{\infty }}}(1)
\label{eqn:proof_T2n_1}
\end{equation}%
where
\begin{equation*}
T_{n,I}=\sup_{\mathcal{N}_{osn}}\left\vert \frac{1}{n}\sum_{i=1}^{n}\left\{
\left\Vert \frac{d\widetilde{m}(X_{i},\alpha )}{d\alpha }[u_{n}^{\ast
}]\right\Vert _{\widehat{\Sigma }^{-1}}^{2}-\left\Vert \frac{d\widetilde{m}%
(X_{i},\alpha _{0})}{d\alpha }[u_{n}^{\ast }]\right\Vert _{\widehat{\Sigma }%
^{-1}}^{2}\right\} \right\vert ,
\end{equation*}%

\begin{equation*}
T_{n,II}=\left\vert \frac{1}{n}\sum_{i=1}^{n}\left\{ \left\Vert \frac{d%
\widetilde{m}(X_{i},\alpha _{0})}{d\alpha }[u_{n}^{\ast }]\right\Vert _{%
\widehat{\Sigma }^{-1}}^{2}-\left\Vert \frac{dm(X_{i},\alpha _{0})}{d\alpha }%
[u_{n}^{\ast }]\right\Vert _{\widehat{\Sigma }^{-1}}^{2}\right\} \right\vert
,
\end{equation*}%
\begin{equation*}
T_{n,III}=\left\vert \frac{1}{n}\sum_{i=1}^{n}\left\{ \left\Vert \frac{%
dm(X_{i},\alpha _{0})}{d\alpha }[u_{n}^{\ast }]\right\Vert _{\widehat{\Sigma
}^{-1}}^{2}-\left\Vert \frac{dm(X_{i},\alpha _{0})}{d\alpha }[u_{n}^{\ast
}]\right\Vert _{\Sigma ^{-1}}^{2}\right\} \right\vert .
\end{equation*}%
Therefore, to prove equation (\ref{eqn:proof_T2n_1}), it suffices to show
that
\begin{equation*}
T_{n,j}=o_{P_{Z^{\infty }}}(1)\text{ for }j\in \{I,II,III\}.
\end{equation*}

Note that for $||\cdot ||_{L^{2}(P_{n})}$ with $P_{n}$ being the empirical
measure, $|||a||_{L^{2}(P_{n})}^{2}-||b||_{L^{2}(P_{n})}^{2}|\leq
||a-b||_{L^{2}(P_{n})}^{2}+2|\langle b,a-b\rangle _{L^{2}(P_{n})}|$. Now,
let $a\equiv \frac{d\widetilde{m}(X_{i},\alpha )}{d\alpha }[u_{n}^{\ast }]$
and $b\equiv \frac{d\widetilde{m}(X_{i},\alpha _{0})}{d\alpha }[u_{n}^{\ast
}]$. In order to show $T_{n,I}=o_{P_{Z^{\infty }}}(1)$, under Assumption \ref%
{ass:VE}(iii), it suffices to show
\begin{equation*}
\sqrt{n^{-1}\sum_{i=1}^{n}\left\Vert \frac{d\widetilde{m}(X_{i},\alpha _{0})%
}{d\alpha }[u_{n}^{\ast }]\right\Vert _{e}^{2}}\sup_{\mathcal{N}_{osn}}\sqrt{%
n^{-1}\sum_{i=1}^{n}\left\Vert \frac{d\widetilde{m}(X_{i},\alpha )}{d\alpha }%
[u_{n}^{\ast }]-\frac{d\widetilde{m}(X_{i},\alpha _{0})}{d\alpha }%
[u_{n}^{\ast }]\right\Vert _{e}^{2}}=o_{P_{Z^{\infty }}}(1).
\end{equation*}%
By the property of LS projection, we have:%
\begin{equation*}
n^{-1}\sum_{i=1}^{n}\left\Vert \frac{d\widetilde{m}(X_{i},\alpha _{0})}{%
d\alpha }[u_{n}^{\ast }]\right\Vert _{e}^{2}\leq
n^{-1}\sum_{i=1}^{n}\left\Vert \frac{dm(X_{i},\alpha _{0})}{d\alpha }%
[u_{n}^{\ast }]\right\Vert _{e}^{2}=O_{P_{Z^{\infty }}}(1)
\end{equation*}%
due to iid data, Markov inequality, the definition of $E_{P_{Z^{\infty }}}%
\left[ \left\Vert \frac{dm(X_{i},\alpha _{0})}{d\alpha }[u_{n}^{\ast
}]\right\Vert _{\Sigma ^{-1}}^{2}\right] $ and Assumption \ref{ass:sieve}%
(iv). Next, by the property of LS projection, we have:%
\begin{eqnarray*}
&&\sup_{\mathcal{N}_{osn}}n^{-1}\sum_{i=1}^{n}\left\Vert \frac{d\widetilde{m}%
(X_{i},\alpha )}{d\alpha }[u_{n}^{\ast }]-\frac{d\widetilde{m}(X_{i},\alpha
_{0})}{d\alpha }[u_{n}^{\ast }]\right\Vert _{e}^{2} \\
&\leq &\sup_{\mathcal{N}_{osn}}n^{-1}\sum_{i=1}^{n}\left\Vert \frac{%
dm(X_{i},\alpha )}{d\alpha }[u_{n}^{\ast }]-\frac{dm(X_{i},\alpha _{0})}{%
d\alpha }[u_{n}^{\ast }]\right\Vert _{e}^{2}=o_{P_{Z^{\infty }}}(1)
\end{eqnarray*}%
due to iid data, Markov inequality and Assumption \ref{ass:cont_diffm}(ii).
Thus we established $T_{n,I}=o_{P_{Z^{\infty }}}(1)$.

By similar algebra as before, in order to show $T_{n,II}=o_{P_{Z^{\infty
}}}(1)$, given Assumption \ref{ass:VE}(iii), it suffices to show
\begin{equation*}
\sqrt{n^{-1}\sum_{i=1}^{n}\left\Vert \frac{dm(X_{i},\alpha _{0})}{d\alpha }%
[u_{n}^{\ast }]\right\Vert _{e}^{2}}\sqrt{n^{-1}\sum_{i=1}^{n}\left\Vert
\frac{d\widetilde{m}(X_{i},\alpha _{0})}{d\alpha }[u_{n}^{\ast }]-\frac{%
dm(X_{i},\alpha _{0})}{d\alpha }[u_{n}^{\ast }]\right\Vert _{e}^{2}}%
=o_{P_{Z^{\infty }}}(1).
\end{equation*}%
The term $n^{-1}\sum_{i=1}^{n}\left\Vert \frac{dm(X_{i},\alpha _{0})}{%
d\alpha }[u_{n}^{\ast }]\right\Vert _{e}^{2}=O_{P_{Z^{\infty }}}(1)$ is due
to iid data, Markov inequality, the definition of $E_{P_{Z^{\infty }}}\left[
\left\Vert \frac{dm(X_{i},\alpha _{0})}{d\alpha }[u_{n}^{\ast }]\right\Vert
_{\Sigma ^{-1}}^{2}\right] $ and Assumption \ref{ass:sieve}(iv). The term $%
n^{-1}\sum_{i=1}^{n}\left\Vert \frac{d\widetilde{m}(X_{i},\alpha _{0})}{%
d\alpha }[u_{n}^{\ast }]-\frac{dm(X_{i},\alpha _{0})}{d\alpha }[u_{n}^{\ast
}]\right\Vert _{e}^{2}=o_{P_{Z^{\infty }}}(1)$ is due to iid data, Markov
inequality and Assumption \ref{ass:anor-mtilde}(i). Thus $%
T_{n,II}=o_{P_{Z^{\infty }}}(1)$.

Finally, $T_{n,III}=o_{P_{Z^{\infty }}}(1)$ follows from the fact that $%
n^{-1}\sum_{i=1}^{n}\left\Vert \frac{dm(X_{i},\alpha _{0})}{d\alpha }%
[u_{n}^{\ast }]\right\Vert _{e}^{2}=O_{P_{Z^{\infty }}}(1)$ and Assumption %
\ref{ass:VE}(iii). We thus established equation (\ref{eqn:proof_T2n_1}).
Since
\begin{equation*}
E_{P_{Z^{\infty }}}\left[ n^{-1}\sum_{i=1}^{n}\left\Vert \frac{%
dm(X_{i},\alpha _{0})}{d\alpha }[u_{n}^{\ast }]\right\Vert _{\Sigma
^{-1}}^{2}\right] =E_{P_{X}}\left[ g(X,u_{n}^{\ast })\Sigma
(X)g(X,u_{n}^{\ast })^{\prime }\right] \leq C<\infty ,
\end{equation*}%
we obtain Result (1).

\textbf{Result (2)} immediately follows from equation (\ref{eqn:proof_T2n_1}%
) and Assumption \ref{ass:LLN_triangular}. \textit{Q.E.D.}

\pagebreak

\section{Sieve Score Statistic and Score Bootstrap}

\label{app:appD}

In the main text we present the sieve Wald, SQLR statistics and their
bootstrap versions. Here we consider sieve score (or LM) statistic and its
bootstrap version. Both the sieve score test and score bootstrap only
require to compute the original-sample restricted PSMD estimator of $\alpha
_{0}$, and hence are computationally attractive.

Recall that $\widehat{\alpha }_{n}^{R}$ is the original-sample restricted
PSMD estimator (\ref{Rpsmd}). Let $\widehat{v}_{n}^{\ast R}$ be computed in
the same way as $\widehat{v}_{n}^{\ast }$ in Subsection \ref{sec:est_avar},
except that we use $\widehat{\alpha }_{n}^{R}$ instead of $\widehat{\alpha }%
_{n}$. And%
\begin{equation*}
||\widehat{v}_{n}^{\ast R}||_{n,sd}^{2}=n^{-1}\sum_{i=1}^{n}\left( \frac{d%
\widehat{m}(X_{i},\widehat{\alpha }_{n}^{R})}{d\alpha }[\widehat{v}%
_{n}^{\ast R}]\right) ^{\prime }\widehat{\Sigma }_{i}^{-1}\rho (Z_{i},%
\widehat{\alpha }_{n}^{R})\rho (Z_{i},\widehat{\alpha }_{n}^{R})^{\prime }%
\widehat{\Sigma }_{i}^{-1}\left( \frac{d\widehat{m}(X_{i},\widehat{\alpha }%
_{n}^{R})}{d\alpha }[\widehat{v}_{n}^{\ast R}]\right)
\end{equation*}%
Denote%
\begin{eqnarray*}
\widehat{S}_{n} &\equiv &\frac{1}{\sqrt{n}}\sum_{i=1}^{n}\left( \frac{d%
\widehat{m}(X_{i},\widehat{\alpha }_{n}^{R})}{d\alpha }[\widehat{v}%
_{n}^{\ast R}/||\widehat{v}_{n}^{\ast R}||_{n,sd}]\right) ^{\prime }\widehat{%
\Sigma }_{i}^{-1}\widehat{m}(X_{i},\widehat{\alpha }_{n}^{R}) \\
\widehat{S}_{1,n} &\equiv &\frac{1}{\sqrt{n}}\sum_{i=1}^{n}\left( \frac{d%
\widehat{m}(X_{i},\widehat{\alpha }_{n}^{R})}{d\alpha }[\widehat{v}%
_{n}^{\ast R}/||\widehat{v}_{n}^{\ast R}||_{n,sd}]\right) ^{\prime }\widehat{%
\Sigma }_{i}^{-1}\rho (Z_{i},\widehat{\alpha }_{n}^{R})
\end{eqnarray*}%
and
\begin{eqnarray*}
\widehat{S}_{n}^{B} &\equiv &\frac{1}{\sqrt{n}}\sum_{i=1}^{n}\left( \frac{d%
\widehat{m}(X_{i},\widehat{\alpha }_{n}^{R})}{d\alpha }[\widehat{v}%
_{n}^{\ast R}/||\widehat{v}_{n}^{\ast R}||_{n,sd}]\right) ^{\prime }\widehat{%
\Sigma }_{i}^{-1}\{\widehat{m}^{B}(X_{i},\widehat{\alpha }_{n}^{R})-\widehat{%
m}(X_{i},\widehat{\alpha }_{n}^{R})\} \\
\widehat{S}_{1,n}^{B} &\equiv &\frac{1}{\sqrt{n}}\sum_{i=1}^{n}\left( \frac{d%
\widehat{m}(X_{i},\widehat{\alpha }_{n}^{R})}{d\alpha }[\widehat{v}%
_{n}^{\ast R}/||\widehat{v}_{n}^{\ast R}||_{n,sd}]\right) ^{\prime }\widehat{%
\Sigma }_{i}^{-1}\{(\omega _{i,n}-1)\rho (Z_{i},\widehat{\alpha }_{n}^{R})\}.
\end{eqnarray*}%
Then
\begin{eqnarray*}
Var\left( \widehat{S}_{1,n}^{B}\mid Z^{n}\right) &=&\frac{\sigma _{\omega
}^{2}\sum_{i=1}^{n}\left( \frac{d\widehat{m}(X_{i},\widehat{\alpha }_{n}^{R})%
}{d\alpha }[\widehat{v}_{n}^{\ast R}]\right) ^{\prime }\widehat{\Sigma }%
_{i}^{-1}\rho (Z_{i},\widehat{\alpha }_{n}^{R})\rho (Z_{i},\widehat{\alpha }%
_{n}^{R})^{\prime }\widehat{\Sigma }_{i}^{-1}\left( \frac{d\widehat{m}(X_{i},%
\widehat{\alpha }_{n}^{R})}{d\alpha }[\widehat{v}_{n}^{\ast R}]\right) }{n||%
\widehat{v}_{n}^{\ast R}||_{n,sd}^{2}} \\
&=&\sigma _{\omega }^{2},
\end{eqnarray*}%
which coincides with that of $\widehat{S}_{1,n}$ (once adjusted by $\sigma
_{\omega }^{2}$).

Following the results in Subsection \ref{sec:est_avar} one can compute $%
\widehat{v}_{n}^{\ast R}$ in closed form, $\widehat{v}_{n}^{\ast R}=\bar{\psi%
}^{k(n)}(\cdot )^{\prime }\widetilde{D}_{n}^{-}\widetilde{\digamma }_{n}$
where
\begin{equation*}
\widetilde{\digamma }_{n}=\frac{d\phi (\widehat{\alpha }_{n}^{R})}{d\alpha }[%
\bar{\psi}^{k(n)}(\cdot )],\quad \widetilde{D}_{n}=n^{-1}\sum_{i=1}^{n}%
\left( \frac{d\widehat{m}(X_{i},\widehat{\alpha }_{n}^{R})}{d\alpha }[\bar{%
\psi}^{k(n)}(\cdot )^{\prime }]\right) ^{\prime }\widehat{\Sigma }%
_{i}^{-1}\left( \frac{d\widehat{m}(X_{i},\widehat{\alpha }_{n}^{R})}{d\alpha
}[\bar{\psi}^{k(n)}(\cdot )^{\prime }]\right) .
\end{equation*}%
And $||\widehat{v}_{n}^{\ast R}||_{n,sd}^{2}=\widetilde{\digamma }%
_{n}^{\prime }\widetilde{D}_{n}^{-}\widetilde{\mho }_{n}\widetilde{D}_{n}^{-}%
\widetilde{\digamma }_{n}$ with
\begin{equation*}
\widetilde{\mho }_{n}=\frac{1}{n}\sum_{i=1}^{n}\left( \frac{d\widehat{m}%
(X_{i},\widehat{\alpha }_{n}^{R})}{d\alpha }[\bar{\psi}^{k(n)}(\cdot
)^{\prime }]\right) ^{\prime }\widehat{\Sigma }_{i}^{-1}\rho (Z_{i},\hat{%
\alpha}_{n}^{R})\rho (Z_{i},\hat{\alpha}_{n}^{R})^{\prime }\widehat{\Sigma }%
_{i}^{-1}\left( \frac{d\widehat{m}(X_{i},\widehat{\alpha }_{n}^{R})}{d\alpha
}[\bar{\psi}^{k(n)}(\cdot )^{\prime }]\right) .
\end{equation*}

Therefore, the bootstrap sieve score statistic $\widehat{S}_{1,n}^{B}$ can
be expressed as
\begin{eqnarray*}
\widehat{S}_{1,n}^{B} &=&\frac{1}{\sqrt{n}}\sum_{i=1}^{n}\left( \frac{d%
\widehat{m}(X_{i},\widehat{\alpha }_{n}^{R})}{d\alpha }[\widehat{v}%
_{n}^{\ast R}/||\widehat{v}_{n}^{\ast R}||_{n,sd}]\right) ^{\prime }\widehat{%
\Sigma }_{i}^{-1}(\omega _{i,n}-1)\rho (Z_{i},\widehat{\alpha }_{n}^{R}) \\
&=&\left( \widetilde{\digamma }_{n}^{\prime }\widetilde{D}_{n}^{-}\widetilde{%
\mho }_{n}\widetilde{D}_{n}^{-}\widetilde{\digamma }_{n}\right) ^{-1/2}%
\widetilde{\digamma }_{n}^{\prime }\widetilde{D}_{n}^{-}\frac{1}{\sqrt{n}}%
\sum_{i=1}^{n}\left( \frac{d\widehat{m}(X_{i},\widehat{\alpha }_{n}^{R})}{%
d\alpha }[\bar{\psi}^{k(n)}(\cdot )^{\prime }]\right) ^{\prime }\widehat{%
\Sigma }_{i}^{-1}(\omega _{i,n}-1)\rho (Z_{i},\widehat{\alpha }_{n}^{R}).
\end{eqnarray*}%
For the case of IID weights, this expression is similar to that proposed in
\cite{KS2012} for parametric models, which suggests the potential higher
order refinements of the bootstrap sieve score test $\left( \widehat{S}%
_{1,n}^{B}\right) ^{2}$. We leave it to future research for bootstrap
refinement.

In the rest of this section, to simplify presentation, we assume that $%
\widehat{m}(x,\alpha )$ is a series LS estimator (\ref{mhat}) of $m(x,\alpha
)$. Then we have:
\begin{equation*}
\widehat{m}^{B}(x,\widehat{\alpha }_{n}^{R})-\widehat{m}(x,\widehat{\alpha }%
_{n}^{R})=\left( \sum_{j=1}^{n}(\omega _{j,n}-1)\rho (Z_{j},\widehat{\alpha }%
_{n}^{R})p^{J_{n}}(X_{j})^{\prime }\right) (P^{\prime }P)^{-}p^{J_{n}}(x).
\end{equation*}%
When $\widehat{\Sigma }=I$ then we have:%
\begin{equation*}
\widehat{S}_{n}=\frac{1}{\sqrt{n}}\sum_{i=1}^{n}\left( \frac{d\widehat{m}%
(X_{i},\widehat{\alpha }_{n}^{R})}{d\alpha }[\widehat{v}_{n}^{\ast R}/||%
\widehat{v}_{n}^{\ast R}||_{n,sd}]\right) ^{\prime }\rho (Z_{i},\widehat{%
\alpha }_{n}^{R})=\widehat{S}_{1,n}
\end{equation*}%
\begin{equation*}
\widehat{S}_{n}^{B}=\frac{1}{\sqrt{n}}\sum_{i=1}^{n}\left( \frac{d\widehat{m}%
(X_{i},\widehat{\alpha }_{n}^{R})}{d\alpha }[\widehat{v}_{n}^{\ast R}/||%
\widehat{v}_{n}^{\ast R}||_{n,sd}]\right) ^{\prime }(\omega _{i,n}-1)\rho
(Z_{i},\widehat{\alpha }_{n}^{R})=\widehat{S}_{1,n}^{B}.
\end{equation*}

Let $\{\epsilon _{n}\}_{n=1}^{\infty }$ and $\{\zeta _{n}\}_{n=1}^{\infty }$
be real valued positive sequences such that $\epsilon _{n}=o(1)$ and $\zeta
_{n}=o(1)$.

\begin{assumption}
\label{ass:SCORE_B} (i) $\max \{\epsilon _{n},n^{-1/4}\}M_{n}\delta
_{n}=o(n^{-1/2})$
\begin{equation*}
\sup_{\mathcal{N}_{osn}}\sup_{u\in \overline{\mathbf{V}}_{n}\colon
||u||=1}n^{-1}\sum_{i=1}^{n}\left\Vert \frac{d\widehat{m}(X_{i},\alpha )}{%
d\alpha }[u]-\frac{dm(X_{i},\alpha )}{d\alpha }[u]\right\Vert
_{e}^{2}=O_{P_{Z^{\infty }}}(\max \{n^{-1/2},\epsilon _{n}^{2}\});
\end{equation*}%
(ii) there is a continuous mapping $\Upsilon :\mathbb{R}_{+}\rightarrow
\mathbb{R}_{+}$ such that $\max \{ \Upsilon (\zeta
_{n}),n^{-1/4}\}M_{n}\delta _{n}=o(n^{-1/2})$ and%
\begin{equation*}
\sup_{\mathcal{N}_{osn}}\sup_{\overline{\mathbf{V}}_{n}\colon ||u_{n}^{\ast
}-u||\leq \zeta _{n}}n^{-1}\sum_{i=1}^{n}\left\Vert \frac{dm(X_{i},\alpha )}{%
d\alpha }[u_{n}^{\ast }]-\frac{dm(X_{i},\alpha )}{d\alpha }[u]\right\Vert
_{e}^{2}=O_{P_{Z^{\infty }}}(\max \{n^{-1/2},\left( \Upsilon (\zeta
_{n})\right) ^{2}\});
\end{equation*}%
(iii) $||\widehat{u}_{n}^{\ast R}-u_{n}^{\ast }||=O_{P_{Z^{\infty }}}(\zeta
_{n})$ where $\widehat{u}_{n}^{\ast R}\equiv \widehat{v}_{n}^{\ast R}/||%
\widehat{v}_{n}^{\ast R}||_{sd}$.
\end{assumption}

Assumption \ref{ass:SCORE_B}(i) can be obtained by similar conditions to
those imposed in \cite{AC_Emetrica03}. Assumption \ref{ass:SCORE_B}(ii) can
be established by controlling the entropy, as in VdV-W Chapter 2.11 and $E%
\left[ \left\Vert \frac{dm(X,\alpha )}{d\alpha }[u_{n}^{\ast }]-\frac{%
dm(X,\alpha )}{d\alpha }[u]\right\Vert _{e}^{2}\right] =o(1)$ for all $%
||u_{n}^{\ast }-u||<\zeta _{n}$; this result is akin to that in lemma 1 of
\cite{CLvK_Emetrica03}. However, Assumption \ref{ass:SCORE_B}(ii) can also
be obtained by weaker conditions, yielding a $\left( \Upsilon (\zeta
_{n})\right) ^{2}$ that is slower than $O(n^{-1/2})$ provided that $\Upsilon
(\zeta _{n})M_{n}\delta _{n}=o(n^{-1/2})$. In the proof we show that $||%
\widehat{u}_{n}^{\ast R}-u_{n}^{\ast }||=o_{P_{Z^{\infty }}}(1)$; faster
rates of convergence will relax the conditions needed to show part (ii).

\begin{theorem}
\label{thm:bootstrap_3} Let $\widehat{\alpha }_{n}^{R}$ be the restricted
PSMD estimator (\ref{Rpsmd}), and conditions for Lemma \ref%
{thm:ThmCONVRATEGRAL} and Proposition \ref{pro:conv-rate-RPSMDE} hold. Let
Assumptions \ref{ass:phi}, \ref{ass:m_ls} - \ref{ass:cont_diffm}, \ref%
{ass:LAQ}(ii), \ref{ass:VE}, \ref{ass:sigma-smooth} and \ref{ass:SCORE_B}
hold and that $n\delta _{n}^{2}\left( M_{n}\delta _{s,n}\right) ^{2\kappa
}C_{n}=o(1)$. Then, under the null hypothesis of $\phi (\alpha _{0})=\phi
_{0}$,

(1) $\widehat{S}_{n}=\sqrt{n}\mathbb{Z}_{n}+o_{P_{Z^{\infty
}}}(1)\Rightarrow N(0,1)$.

(2) Further, if conditions for Lemma \ref{lem:cons_boot} and Assumptions \ref%
{ass:LAQ_B}(ii), \ref{ass:Wboot} or \ref{ass:Wboot_e} hold, then:
\begin{eqnarray*}
\left\vert \mathcal{L}_{V^{\infty }|Z^{\infty }}(\sigma _{\omega }^{-1}%
\widehat{S}_{n}^{B}\mid Z^{n})-\mathcal{L}(\widehat{S}_{n})\right\vert
&=&o_{P_{Z^{\infty }}}(1),\quad \text{and} \\
\sup_{t\in \mathbb{R}}\left\vert P_{V^{\infty }|Z^{\infty }}(\sigma _{\omega
}^{-1}\widehat{S}_{n}^{B}\leq t|Z^{n})-P_{Z^{\infty }}(\widehat{S}_{n}\leq
t)\right\vert &=&o_{P_{V^{\infty }|Z^{\infty }}}(1)~wpa1(P_{Z^{\infty }}).
\end{eqnarray*}
\end{theorem}

\medskip

\noindent \textbf{Proof of Theorem \ref{thm:bootstrap_3}}\textsc{:} We first
note that by Lemma \ref{lem:Qdiff_B}, Assumptions \ref{ass:LAQ}(i) and \ref%
{ass:LAQ_B}(i) hold. Also, by Proposition \ref{pro:conv-rate-RPSMDE} we have
$\widehat{\alpha }_{n}^{R}\in \mathcal{N}_{osn}$ wpa1 under the null
hypothesis of $\phi (\alpha _{0})=\phi _{0}$. Under the null hypothesis, and
Assumption \ref{ass:phi}, we also have (see Step 1 in the proof of Theorem %
\ref{thm:chi2}):%
\begin{equation*}
\sqrt{n}\langle u_{n}^{\ast },\widehat{\alpha }_{n}^{R}-\alpha _{0}\rangle
=o_{P_{Z^{\infty }}}(1).
\end{equation*}

For \textbf{Result (1)}, we show that $\widehat{S}_{n}$ is asymptotically
standard normal under the null hypothesis in two steps.

\textsc{Step 1.} We first show that $\left\vert \frac{||\widehat{v}%
_{n}^{\ast R}||_{sd}}{||\widehat{v}_{n}^{\ast R}||_{n,sd}}-1\right\vert
=o_{P_{Z^{\infty }}}(1)$ and $||\widehat{u}_{n}^{\ast R}-u_{n}^{\ast
}||=o_{P_{Z^{\infty }}}(1)$, where $\widehat{u}_{n}^{\ast R}\equiv \widehat{v%
}_{n}^{\ast R}/||\widehat{v}_{n}^{\ast R}||_{sd}$ and $\widehat{v}_{n}^{\ast
R}$ is computed in the same way as that in Subsection \ref{sec:est_avar},
except that we use $\widehat{\alpha }_{n}^{R}$ instead of $\widehat{\alpha }%
_{n}$.

$\left\vert \frac{||\widehat{v}_{n}^{\ast R}||_{sd}}{||\widehat{v}_{n}^{\ast
R}||_{n,sd}}-1\right\vert =o_{P_{Z^{\infty }}}(1)$ can be established in the
same way as that of Theorem \ref{thm:VE}(1). Also, following the proof of
Theorem \ref{thm:VE}(1), we obtain:%
\begin{equation*}
\left\Vert \frac{\widehat{v}_{n}^{\ast R}-v_{n}^{\ast }}{||v_{n}^{\ast }||}%
\right\Vert =o_{P_{Z^{\infty }}}(1),\quad \frac{||\widehat{v}_{n}^{\ast R}||%
}{||v_{n}^{\ast }||_{sd}}=O_{P_{Z^{\infty }}}\left( 1\right) ,\quad
\sup_{v\in \overline{\mathbf{V}}_{n}}\left\vert \frac{\langle v_{n}^{\ast }-%
\widehat{v}_{n}^{\ast R},v\rangle }{||v||\times ||\widehat{v}_{n}^{\ast R}||}%
\right\vert =o_{P_{Z^{\infty }}}(1).
\end{equation*}%
This and Assumption \ref{ass:sieve}(iv) imply that $\left\vert \frac{\langle
\widehat{v}_{n}^{\ast R},\widehat{v}_{n}^{\ast R}-v_{n}^{\ast }\rangle }{||%
\widehat{v}_{n}^{\ast R}||_{sd}^{2}}\right\vert =o_{P_{Z^{\infty }}}(1)$ and
$\left\vert \frac{\langle v_{n}^{\ast },\widehat{v}_{n}^{\ast R}-v_{n}^{\ast
}\rangle }{||\widehat{v}_{n}^{\ast R}||_{sd}^{2}}\right\vert =\frac{%
||v_{n}^{\ast }||_{sd}}{||\widehat{v}_{n}^{\ast R}||_{sd}}\times
o_{P_{Z^{\infty }}}(1)$. Therefore,
\begin{equation*}
\left\vert \frac{||v_{n}^{\ast }||_{sd}^{2}}{||\widehat{v}_{n}^{\ast
R}||_{sd}^{2}}-1\right\vert \leq \left\vert \frac{\langle \widehat{v}%
_{n}^{\ast R},\widehat{v}_{n}^{\ast R}-v_{n}^{\ast }\rangle }{||\widehat{v}%
_{n}^{\ast R}||_{sd}^{2}}\right\vert +\left\vert \frac{\langle v_{n}^{\ast },%
\widehat{v}_{n}^{\ast R}-v_{n}^{\ast }\rangle }{||\widehat{v}_{n}^{\ast
R}||_{sd}^{2}}\right\vert =o_{P_{Z^{\infty }}}(1).
\end{equation*}%
and
\begin{equation*}
\left\vert \frac{||v_{n}^{\ast }||_{sd}}{||\widehat{v}_{n}^{\ast R}||_{sd}}%
-1\right\vert =o_{P_{Z^{\infty }}}(1).
\end{equation*}%
Thus%
\begin{eqnarray*}
||\widehat{u}_{n}^{\ast R}-u_{n}^{\ast }|| &=&\left\Vert \frac{\widehat{v}%
_{n}^{\ast R}}{||\widehat{v}_{n}^{\ast R}||_{sd}}-\frac{v_{n}^{\ast }}{%
||v_{n}^{\ast }||_{sd}}\right\Vert =\left\Vert \frac{\widehat{v}_{n}^{\ast R}%
}{||v_{n}^{\ast }||_{sd}}(1+o_{P_{Z^{\infty }}}(1))-\frac{v_{n}^{\ast }}{%
||v_{n}^{\ast }||_{sd}}\right\Vert \\
&=&\left\Vert \frac{\widehat{v}_{n}^{\ast R}-v_{n}^{\ast }}{||v_{n}^{\ast
}||_{sd}}\right\Vert +o_{P_{Z^{\infty }}}(\frac{||\widehat{v}_{n}^{\ast R}||%
}{||v_{n}^{\ast }||_{sd}})=o_{P_{Z^{\infty }}}(1).
\end{eqnarray*}

\textsc{Step 2.} We show that under the null hypothesis,
\begin{equation}
\widehat{S}_{n}=\sqrt{n}\mathbb{Z}_{n}+o_{P_{Z^{\infty }}}(1)\equiv \frac{1}{%
\sqrt{n}}\sum_{i=1}^{n}\left( \frac{dm(X_{i},\alpha _{0})}{d\alpha }%
[u_{n}^{\ast }]\right) ^{\prime }\Sigma (X_{i})^{-1}\rho (Z_{i},\alpha
_{0})+o_{P_{Z^{\infty }}}(1).  \label{score1}
\end{equation}%
By Step 1, it suffices to show that under the null hypothesis,%
\begin{equation*}
\overline{S}_{n}\equiv \frac{1}{\sqrt{n}}\sum_{i=1}^{n}\left( \frac{d%
\widehat{m}(X_{i},\widehat{\alpha }_{n}^{R})}{d\alpha }[\widehat{u}%
_{n}^{\ast R}]\right) ^{\prime }\widehat{\Sigma }^{-1}(X_{i})\widehat{m}%
(X_{i},\widehat{\alpha }_{n}^{R})=\sqrt{n}\mathbb{Z}_{n}+o_{P_{Z^{\infty
}}}(1).
\end{equation*}%
Recall that $\ell _{n}(x,\alpha )\equiv \widehat{m}(x,\alpha _{0})+%
\widetilde{m}(x,\alpha )$. We have:
\begin{align*}
& \left\vert \overline{S}_{n}-\frac{1}{\sqrt{n}}\sum_{i=1}^{n}\left( \frac{d%
\widehat{m}(X_{i},\widehat{\alpha }_{n}^{R})}{d\alpha }[\widehat{u}%
_{n}^{\ast R}]\right) ^{\prime }\widehat{\Sigma }(X_{i})^{-1}\ell _{n}(X_{i},%
\widehat{\alpha }_{n}^{R})\right\vert \\
\leq & \sqrt{n}\sqrt{n^{-1}\sum_{i=1}^{n}\left\Vert \widehat{\Sigma }^{-1/2}(X_{i})
\frac{d\widehat{m}(X_{i},\widehat{\alpha }_{n}^{R})}{d\alpha }[\widehat{u}%
_{n}^{\ast R}]\right\Vert _{e}^{2}}\sqrt{n^{-1}\sum_{i=1}^{n}\left\Vert
\widehat{m}(X_{i},\widehat{\alpha }_{n}^{R})-\ell _{n}(X_{i},\widehat{\alpha
}_{n}^{R})\right\Vert _{e}^{2}},
\end{align*}%
By Lemma \ref{lem:suff_mcond_boot}(1) and the assumption that $n\delta
_{n}^{2}(M_{n}\delta _{s,n})^{2\kappa }C_{n}=o(1)$, we have:%
\begin{equation*}
\sqrt{n^{-1}\sum_{i=1}^{n}\left\Vert \widehat{m}(X_{i},\widehat{\alpha }%
_{n}^{R})-\ell _{n}(X_{i},\widehat{\alpha }_{n}^{R})\right\Vert _{e}^{2}}%
=o_{P_{Z^{\infty }}}(n^{-1/2}).
\end{equation*}%
Also $n^{-1}\sum_{i=1}^{n}\left\Vert \widehat{\Sigma }^{-1/2}(X_{i})\frac{d\widehat{m}(X_{i},\widehat{%
\alpha }_{n}^{R})}{d\alpha }[\widehat{u}_{n}^{\ast R}]\right\Vert
_{e}^{2}\asymp 1$ by Step 1, assumptions \ref{ass:cont_diffm} and \ref{ass:SCORE_B}. Therefore
\begin{equation*}
\overline{S}_{n}=\frac{1}{\sqrt{n}}\sum_{i=1}^{n}\left( \frac{d\widehat{m}%
(X_{i},\widehat{\alpha }_{n}^{R})}{d\alpha }[\widehat{u}_{n}^{\ast
R}]\right) ^{\prime }\widehat{\Sigma }(X_{i})^{-1}\ell _{n}(X_{i},\widehat{%
\alpha }_{n}^{R})+o_{P_{Z^{\infty }}}(1).
\end{equation*}%
Assumption \ref{ass:SCORE_B}(i) implies that%
\begin{equation*}
n^{-1}\sum_{i=1}^{n}\left\Vert \frac{d\widehat{m}(X_{i},\widehat{\alpha }%
_{n}^{R})}{d\alpha }[\widehat{u}_{n}^{\ast R}]-\frac{dm(X_{i},\widehat{%
\alpha }_{n}^{R})}{d\alpha }[\widehat{u}_{n}^{\ast R}]\right\Vert
_{e}^{2}=O_{P_{Z^{\infty }}}(\max \left\{ n^{-1/2},\epsilon _{n}^{2}\right\}
).
\end{equation*}%
And $n^{-1}\sum_{i=1}^{n}\left\Vert \ell _{n}(X_{i},\widehat{\alpha }%
_{n}^{R})\right\Vert _{e}^{2}=O_{P_{Z^{\infty }}}((M_{n}\delta _{n})^{2})$
by Lemma \ref{lem:suff_mcond_boot}(2). These results, Assumption \ref%
{ass:SCORE_B}(i) and Assumption \ref{ass:VE}(iii) together lead to
\begin{eqnarray*}
&&\frac{1}{\sqrt{n}}\sum_{i=1}^{n}\left( \frac{d\widehat{m}(X_{i},\widehat{%
\alpha }_{n}^{R})}{d\alpha }[\widehat{u}_{n}^{\ast R}]\right) ^{\prime }%
\widehat{\Sigma }(X_{i})^{-1}\ell _{n}(X_{i},\widehat{\alpha }_{n}^{R}) \\
&=&\frac{1}{\sqrt{n}}\sum_{i=1}^{n}\left( \frac{dm(X_{i},\widehat{\alpha }%
_{n}^{R})}{d\alpha }[\widehat{u}_{n}^{\ast R}]\right) ^{\prime }\Sigma
(X_{i})^{-1}\ell _{n}(X_{i},\widehat{\alpha }_{n}^{R})+o_{P_{Z^{\infty }}}(1)
\\
&=&\frac{1}{\sqrt{n}}\sum_{i=1}^{n}\left( \frac{dm(X_{i},\widehat{\alpha }%
_{n}^{R})}{d\alpha }[u_{n}^{\ast }]\right) ^{\prime }\Sigma (X_{i})^{-1}\ell
_{n}(X_{i},\widehat{\alpha }_{n}^{R})+o_{P_{Z^{\infty }}}(1),
\end{eqnarray*}%
where the second equality is due to $||\widehat{u}_{n}^{\ast R}-u_{n}^{\ast
}||=O_{P_{Z^{\infty }}}(\zeta _{n})$ (Assumption \ref{ass:SCORE_B}(iii)) and
Assumption \ref{ass:SCORE_B}(ii).

Since $\widehat{\alpha }_{n}^{R}\in \mathcal{N}_{osn}$ wpa1 under the null
hypothesis, $\sqrt{n}\langle u_{n}^{\ast },\widehat{\alpha }_{n}^{R}-\alpha
_{0}\rangle =o_{P_{Z^{\infty }}}(1)$, and by analogous calculations to those
in the proof of Lemma \ref{lem:T1n}, we obtain:%
\begin{equation*}
\frac{1}{\sqrt{n}}\sum_{i=1}^{n}\left( \frac{dm(X_{i},\widehat{\alpha }%
_{n}^{R})}{d\alpha }[u_{n}^{\ast }]\right) ^{\prime }\Sigma (X_{i})^{-1}\ell
_{n}(X_{i},\widehat{\alpha }_{n}^{R})=\sqrt{n}\mathbb{Z}_{n}+o_{P_{Z^{\infty
}}}(1),
\end{equation*}%
and hence equation (\ref{score1}) holds. By Assumption \ref{ass:LAQ}(ii) we
have: $\widehat{S}_{n}\Rightarrow N(0,1)$ under the null hypothesis.

For \textbf{Result (2)}, we now show that $\widehat{S}_{n}^{B}$ also
converges weakly (in the sense of Bootstrap Section \ref{sec:bootstrap}) to
a standard normal under the null hypothesis. It suffices to show that
\begin{equation}
\widehat{S}_{n}^{B}=\frac{1}{\sqrt{n}}\sum_{i=1}^{n}(\omega _{i}-1)\left(
\frac{dm(X_{i},\alpha _{0})}{d\alpha }[u_{n}^{\ast }]\right) ^{\prime
}\Sigma (X_{i})^{-1}\rho (Z_{i},\alpha _{0})+o_{P_{V^{\infty }|Z^{\infty
}}}(1)\text{ }wpa1(P_{Z^{\infty }}).  \label{score2}
\end{equation}

Note that $\ell _{n}^{B}(X_{i},\widehat{\alpha }_{n}^{R})-\ell _{n}(X_{i},%
\widehat{\alpha }_{n}^{R}))=\widehat{m}^{B}(X_{i},\alpha _{0})-\widehat{m}%
(X_{i},\alpha _{0})$, and that $n^{-1}\sum_{i=1}^{n}||\widehat{m}%
^{B}(X_{i},\alpha _{0})-\widehat{m}(X_{i},\alpha
_{0})||_{e}^{2}=O_{P_{V^{\infty }|Z^{\infty }}}(J_{n}/n)$ wpa1($P_{Z^{\infty
}}$) (see the proof of Lemma \ref{lem:suff_mcond_boot}). We have, by
calculations similar to Step 2,
\begin{equation*}
\left\vert \widehat{S}_{n}^{B}-\frac{1}{\sqrt{n}}\sum_{i=1}^{n}\left( \frac{%
dm(X_{i},\widehat{\alpha }_{n}^{R})}{d\alpha }[u_{n}^{\ast }]\right)
^{\prime }\Sigma (X_{i})^{-1}\{\ell _{n}^{B}(X_{i},\widehat{\alpha }%
_{n}^{R})-\ell _{n}(X_{i},\widehat{\alpha }_{n}^{R})\}\right\vert
=o_{P_{V^{\infty }|Z^{\infty }}}(1)\text{ }wpa1(P_{Z^{\infty }}).
\end{equation*}

By analogous calculations to those in the proof of Lemma \ref{lem:T1n}, we
obtain equation (\ref{score2}). This and Result (1) and Assumption \ref%
{ass:LAQ_B}(ii) now imply that under the null and conditional on the data, $%
\sigma _{\omega }^{-1}\widehat{S}_{n}^{B}$ is also asymptotically standard
normally distributed. The last part of Result (2) can be established in the
same way as that of Theorem \ref{thm:bootstrap_2}(1), and is omitted.
\textit{Q.E.D.}

\end{document}